\documentclass{elsart}

\usepackage{amssymb}
\usepackage{amsfonts}
\usepackage{amsmath}
\usepackage{rotating}
\usepackage{graphicx,color}
\usepackage{subfigure}
\usepackage{bm}
\usepackage{color}
\usepackage{multirow}
\usepackage{enumitem}
\usepackage{tikz}
\usepackage{breakcites}
\usepackage{epstopdf}
\usepackage{morefloats}
\usepackage[pdfpagelabels, hypertexnames=true, plainpages=false, naturalnames=false, breaklinks=true]{hyperref}

\newcommand{\descrcell}[2]{%
  \scriptsize
  \begin{tabular}[t]{@{}c@{}}\normalsize#1\\\normalsize#2\end{tabular}%
}

\begin{document}
\setcounter{page}{1}
\allowdisplaybreaks[4]
\begin{frontmatter}

\title{Algebraic multiscale method for flow in heterogeneous porous media with embedded discrete fractures (F-AMS)}

\author[TuDelft]{Matei \c{T}ene\thanksref{label1}\corauthref{cor1}}
\author[PI]{Mohammed Saad Al Kobaisi\thanksref{label2}}
\author[TuDelft]{Hadi Hajibeygi\thanksref{label3}}

\address[TuDelft]{Department of Geoscience and Engineering, Delft University of Technology,\\P.O. Box 5048, 2600 GA Delft, The Netherlands.}
\address[PI]{The Petroleum Institute, P.O. Box: 2533, Abu Dhabi, United Arab Emirates.}

\thanks[label1]{m.tene@tudelft.nl}
\thanks[label2]{malkobaisi@pi.ac.ae}
\thanks[label3]{h.hajibeygi@tudelft.nl}

\corauth[cor1]{Corresponding author.}

\begin{abstract}
This paper introduces an Algebraic MultiScale method for simulation of flow in heterogeneous porous media with embedded discrete Fractures (F-AMS). First, multiscale coarse grids are independently constructed for both porous matrix and fracture networks. Then, a map between coarse- and fine-scale is obtained by algebraically computing basis functions with local support. In order to extend the localization assumption to the fractured media, four types of basis functions are investigated: (1) Decoupled-AMS, in which the two media are completely decoupled, (2) Frac-AMS and (3) Rock-AMS, which take into account only one-way transmissibilities, and (4) Coupled-AMS, in which the matrix and fracture interpolators are fully coupled.
In order to ensure scalability, the F-AMS framework permits full flexibility in terms of the resolution of the fracture coarse grids. Numerical results are presented for two- and three-dimensional heterogeneous test cases.
During these experiments, the performance of F-AMS, paired with ILU(0) as second-stage smoother in a convergent iterative procedure, is studied by monitoring CPU times and convergence rates.
Finally, in order to investigate the scalability of the method, an extensive benchmark study is conducted, where a commercial algebraic multigrid solver is used as reference. The results show that, given an appropriate coarsening strategy, F-AMS is insensitive to severe fracture and matrix conductivity contrasts, as well as the length of the fracture networks. Its unique feature is that a fine-scale mass conservative flux field can be reconstructed after any iteration, providing efficient approximate solutions in time-dependent simulations.
\end{abstract}


\end{frontmatter}
\endNoHyper

\section{Introduction}
\label{sec:intro}

In many geoscience applications, including hydrocarbon production and geothermal energy exploitation, the target formations are naturally fractured. These formations are often highly heterogeneous matrix rock is crossed by several networks of lower-dimensional highly-conductive fractures at multiple length scales \cite{Berkowitz1}. This raises important challenges for flow simulation, motivating the development of advanced modelling and numerical solution techniques. 

Among the proposed methods, the hierarchical fracture modelling approach allows for avoiding complexities associated with the discretization and dynamic nature of fracture geometries \cite{SeongFrac2}. In this approach, small-scale fractures (below the matrix grid resolution) are homogenized within the matrix rock, altering its effective permeability \cite{Sander}. The remaining fractures are then represented as explicit control volumes \cite{SeongFrac1,hadi-frac-jcp}. If the fracture and matrix grids are generated independently, then the formation is said to be discretized according to the Embedded Discrete Fracture Model (EDFM) \cite{Moinfar1}. Alternatively, the fracture cells can be constrained to lie at the interfaces of matrix cells, i.e. by employing Discrete Fracture Modelling methods (DFM), which require unstructured grids \cite{karimi-fard}. Both DFM and EDFM have been applied to reservoirs with complex fracture geometries \cite{stephanM1} and fluid physics \cite{rainer1,Moinfar2}. Recent developments include higher-order approximation schemes within finite-volume \cite{dfm-Ahmed2015} and finite-element \cite{Geiger09} methods.

Note that the total number of degrees of freedom (DOF), even after homogenizing small-scale fractures, is beyond the scope of classical simulation methods. This motivates the development of an efficient multiscale method for heterogeneous fractured porous media.

Multiscale finite element (MSFE) \cite{Hou97} and finite volume (MSFV) methods \cite{Jenny03} have been introduced and evolved mainly for heterogeneous, but non-fractured, porous media (see  \cite{Kippe06} for a comparison). Recent developments include efficient solution of the pressure equation for multi-component displacements within sequentially- \cite{SEON-BO,hadi-compositional-spej} and fully-implicit \cite{cusini,ADM} frameworks. They have also been extended to capture complex wells \cite{Patrick-well2,Wolfsteiner-well} and to the transport equations \cite{zhou-trans}. In addition, enriched multiscale methods are targeted at media with high conductivity contrasts \cite{Yalchin-enriched1,Yalchin-enriched2,Davide14}, with modifications to maintain their monotonicity \cite{Hesse08a,yixuan-monotonemsfv}. More recently, a multiscale formulation has been devised to support unstructured grids \cite{MsRSB-olav}. 

The first application of MSFV methods to fractured reservoirs was developed on the basis of coupling the matrix pressure to the average pressure in each fracture network, at coarse scale \cite{hadi-frac-jcp}. Based on their results for two-dimensional (2D) problems, this method was efficient for media with highly conductive fracture networks which span short spatial length scales (relative to that of the domain). However, convergence was observed to degrade for test cases with significant variations in the pressure distribution along the fracture network. In combination with streamline-based mixed formulations, multiscale methods have also been employed to 2D fractured reservoirs \cite{Natvig09}. More recently, a multiscale approach was developed for 2D reservoirs which assigned one coarse node at each fracture intersection only, with no coarse nodes in the matrix \cite{Sandve}. Note that none of these methods include 3D heterogeneous reservoirs nor has their performance been benchmarked against a commercial linear solver. More importantly, the literature is lacking a multiscale method which allows for flexible coarse grids inside the matrix as well as its embedded fractures and, thus,  able to accommodate heterogeneous cases with fracture networks of different length scales. 

This paper presents the development of an Algebraic MultiScale method for heterogeneous Fractured porous media (F-AMS) using EDFM. Given a partition of the fine-scale cells into primal and dual coarse blocks for both the matrix and fracture networks, the algorithm algebraically constructs the multiscale prolongation (mapping coarse- to fine-scale) and restriction (mapping fine- to coarse-scale) operators. The prolongation operator columns are the local basis functions, solved on dual-coarse cells, for both matrix and fractures.

F-AMS supports four different matrix-fracture coupling strategies, at the coarse-scale. First, the Decoupled-AMS basis functions are defined by neglecting the contribution of a medium's coarse solution (e.g., fractures) in the interpolated solution in the other (e.g., matrix), thus preserving sparsity in the resulting coarse-scale system. Then, two semi-coupled (one-way) strategies, Rock-AMS and Frac-AMS, are considered. The Rock-AMS approach constructs a prolongation operator in which the matrix coarse solutions also contribute in computing the interpolated fine-scale solution in neighbouring fractures. Similarly, Frac-AMS considers the influence of the fracture coarse solution when interpolating the pressure inside the surrounding porous rock. Finally, the fully coupled strategy, Coupled-AMS, is devised, where coarse-scale solutions from both media play a role in finding the fine-scale solution of each other. This last approach, although allowing for full fracture-matrix coupling, leads to a dense coarse-scale system with additional overhead during the associated algebraic (matrix-vector, matrix-matrix) operations. As such, for practical applications, the Coupled-AMS prolongation operator may require tuning via truncation, where values below a specified threshold are algebraically deleted, followed by a rescaling step, to maintain partition of unity. This option is also investigated in the paper.

To summarize, F-AMS allows for arbitrary coarse grid resolutions in both fractures and matrix, as well as all possible coarse-scale coupling between them. Furthermore, once these coarse grids are defined, the F-AMS procedure is formulated and implemented in algebraic form, in line with the previously published formulations of incompressible (AMS) \cite{yixuan-ams} and compressible (C-AMS) \cite{Tene-cams} flows. In the limit, if the Frac-AMS coupling strategy and only one coarse node per fracture network is used, F-AMS automatically reduces to the method described in \cite{hadi-frac-jcp}. However, this setup proves inefficient for many of the test cases in this paper. From a bottom-up perspective, F-AMS extends the AMS prolongation operator, as previously described in \cite{yixuan-ams}, with additional columns. Some of these columns correspond to the enriched fracture coarse domain, as explained above. The remainder represent local well basis functions for Peaceman wells \cite{Patrick-well2}. 

In order to test F-AMS method, a proof-of-concept implementation is developed with a focus on reservoirs defined on 3D structured grids with embedded vertical fracture plates (for the challenges associated with unstructured multiscale simulation see \cite{MsRSB-olav}). For the presented experiments, a finite-element (FE) restriction operator is employed to obtain a symmetric-positive-definite coarse system. If approximate (non-converged) F-AMS solutions are used, in the last iteration step, a finite-volume (FV) restriction operator is employed followed by mass-conservative reconstruction of fine-scale flux field. The performance of this in-house object oriented serial-processing simulator was measured based on CPU times, as well as convergence rates. Numerical test cases are considered in order to study the effect of the different components of the algorithm, namely, the coarsening ratios and basis function coupling strategies. F-AMS is developed as an accurate multiscale (approximate) pressure solver. In order to assess its scalability, however, its performance (based on CPU times) is compared against the commercial Algebraic MultiGrid (AMG) solver, SAMG \cite{SAMG}, as a preconditioner. Results of these systematic studies show that only few DOF per fracture network are necessary to obtain a good trade-off between convergence rate and computational expense. In conclusion, F-AMS is found efficient and scalable for solving flow in heterogeneous and naturally fractured porous media. Its development marks an important step forward towards the integration multiscale methods as ``black-box'' pressure solvers within existing reservoir simulators, with the possibility of extension to more complex physics and scenarios. 

The paper is structured as follows. First, the EDFM fine-scale discrete system is described in Section~\ref{sec:equations}. Then, the components of the F-AMS algorithm are detailed in Section~\ref{sec:fams}. Section~\ref{sec:results} consists of numerical results for both 2D and 3D test cases. Finally, conclusions and remarks make the subject of Section~\ref{sec:conclusions}.


\section{Fine-scale discretized system}
\label{sec:equations}

The mass-conservation equations for single-phase flow in fractured media, using Darcy's law, can be written as
\begin{align}
\bigg[\frac{\partial (\phi \rho)}{\partial t}-\nabla \cdot (\rho \lambda \cdot \nabla  p)\bigg]^m &= [\rho q]^{mw} + [\rho Q]^{m} + [\rho q]^{mf} \hspace{5mm} \text{on} \hspace{2mm} \Omega^m \subset R^n \label{pm}
\end{align}
for the matrix (superscript $^m$) and 
\begin{align}
\bigg[\frac{\partial (\phi \rho)}{\partial t}-\nabla \cdot (\rho \lambda \cdot  \nabla p)\bigg]^f &= [\rho q]^{fw}+ [\rho Q]^{f} + [\rho q]^{fm} \hspace{5mm} \text{on} \hspace{2mm} \Omega^f \subset R^{n-1}\label{pf}
\end{align}
for the fracture (superscript $^f$) spatial domains. Here, the phase mobility $\lambda = k^{*}/\mu$ consists of fluid viscosity ($\mu$) and effective fracture ($k^{*} = k^f$) or rock ($k^{*} = k^m$) permeability. Note that the latter can also account for the homogenized small-scale fractures, as described in the hierarchical fracture model \cite{SeongFrac1,Sander}. Also, $\rho$ and $\phi$ are the fluid density and rock porosity, respectively. Here, for simplicity, fracture porosity is always considered to be $1$. 
The $q^{mw}$ and $q^{fw}$ denote the matrix and fracture external source terms, respectively, i.e. from injection/production wells. For a perforated matrix control volume $V$, it reads
\begin{align}
q^{mw}_V~ =~ PI~ \lambda^m~ (p^w - p^m) / V~ \equiv~ \beta^m~ (p^w - p^m),
\end{align}
where $\beta^m = PI~ \lambda^m / V$, $q^{mw}_{tot} = \int_V q^{mw}_V dV$ is the total injection rate, and $PI$ is the productivity index \cite{Peaceman_well}. Similarly, the fracture-matrix coupling terms are modeled such that for a matrix volume $V$ intersecting with a fracture surface $A$ one obtains
\begin{align}
q^{mf}_V~ =~ CI~  \lambda^{f-m}~ (p^f - p^m) / V~ \equiv~ \eta^m~ (p^f - p^m)
\end{align}
and
\begin{align}
q^{fm}_A~ =~ CI~  \lambda^{f-m}~ (p^m - p^f) / A~ \equiv~ \eta^f~ (p^m - p^f),
\end{align}
where $\eta^m = CI~  \lambda^{f-m} / V$ and $\eta^f = CI~  \lambda^{f-m} / A$. This ensures the total flux between a fracture element of area $A$ and a matrix element of volume $V$ is equal \cite{hadi-frac-jcp}, i.e.,
\begin{align}
\int_{V} {q^{mf}_V} dV = - \int_{A} q^{fm}_A dA.
\end{align}
The $\lambda^{f-m}$ is the effective mobility at the interface between the fractures and their surrounding matrix. The $CI$ is the connectivity index, defined on a discrete system as
\begin{align}
CI_{ij}  = \frac{A_{ij}}{\langle d \rangle_{ij}},
\end{align}
where $A_{ij}$ is the area fraction of fracture element $i$ overlapping with the matrix element $j$, and $\langle d \rangle_{ij}$ is the average distance of the two elements \cite{hadi-frac-jcp}. Finally, the $Q^m$ and $Q^f$ terms describe other external source terms for matrix and fractures (e.g., gravity terms).

\begin{figure}[htb!]%
\subfigure[Fine-scale grid]{\includegraphics[width=0.4\textwidth,height=0.4\textwidth]{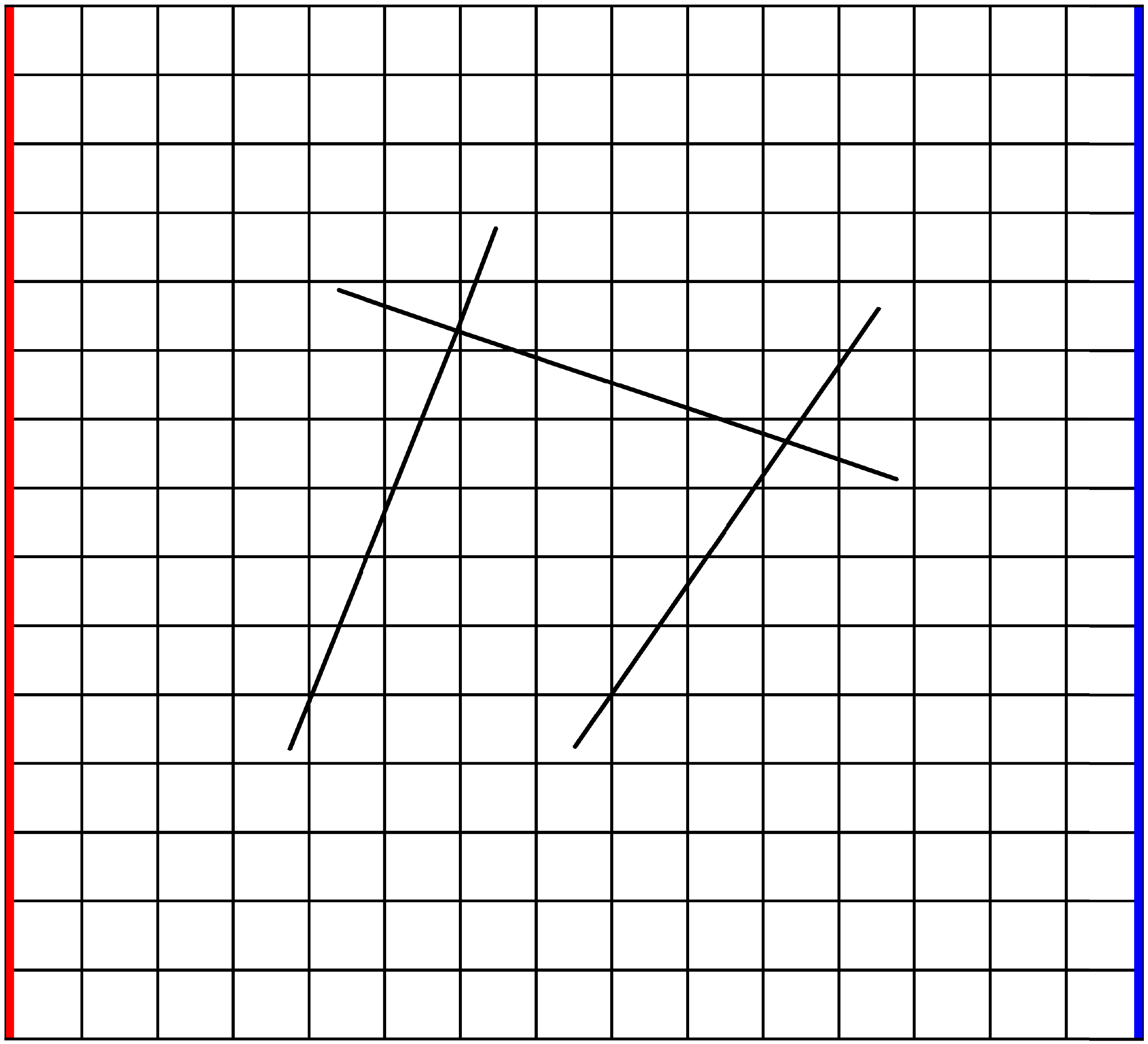}\label{fig:fine2d}}\hspace{1.7cm}%
\subfigure[Pressure solution]{\includegraphics[width=0.4\textwidth,height=0.4\textwidth]{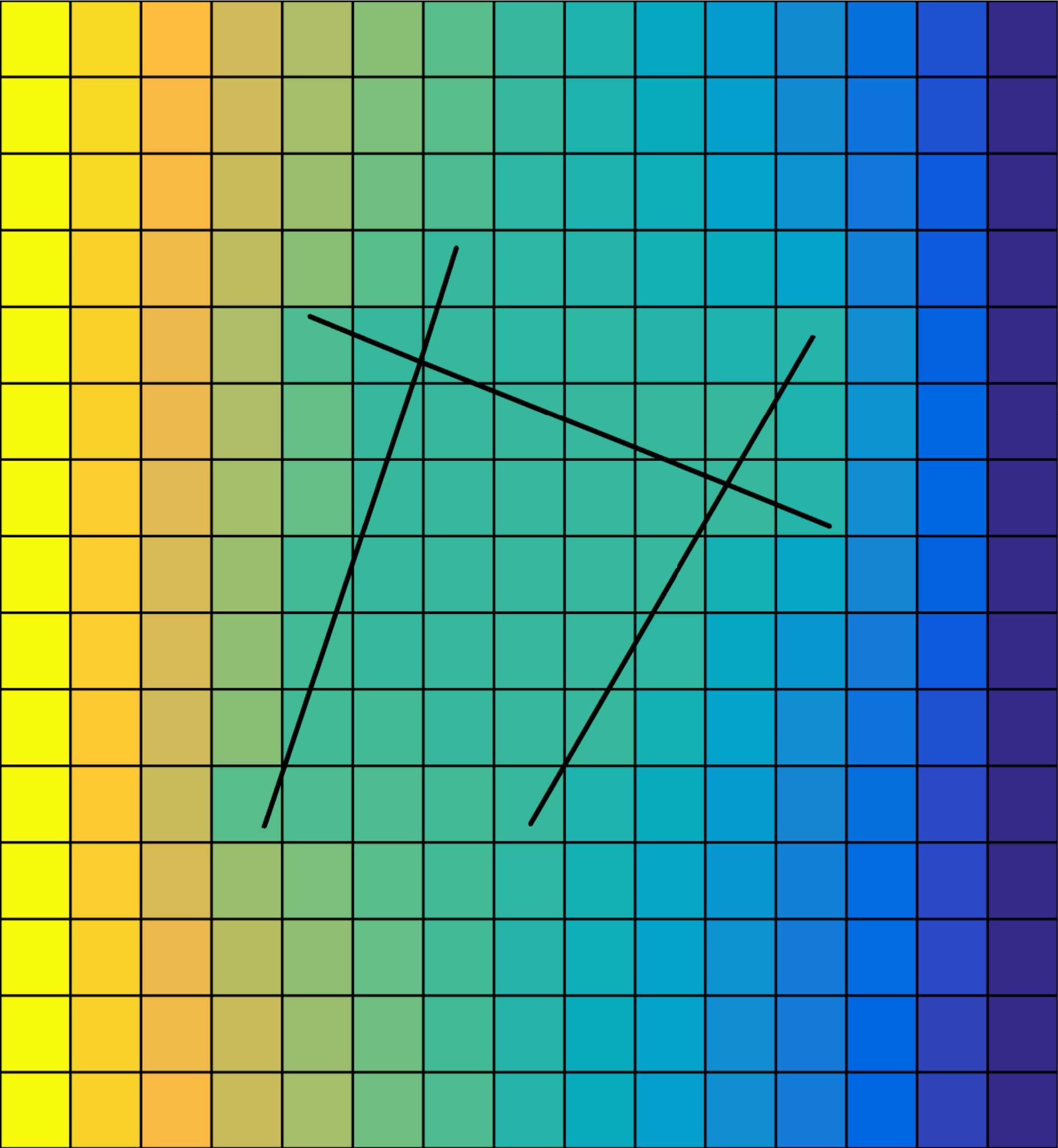}\hspace{0.15cm}\includegraphics[height=0.4\textwidth]{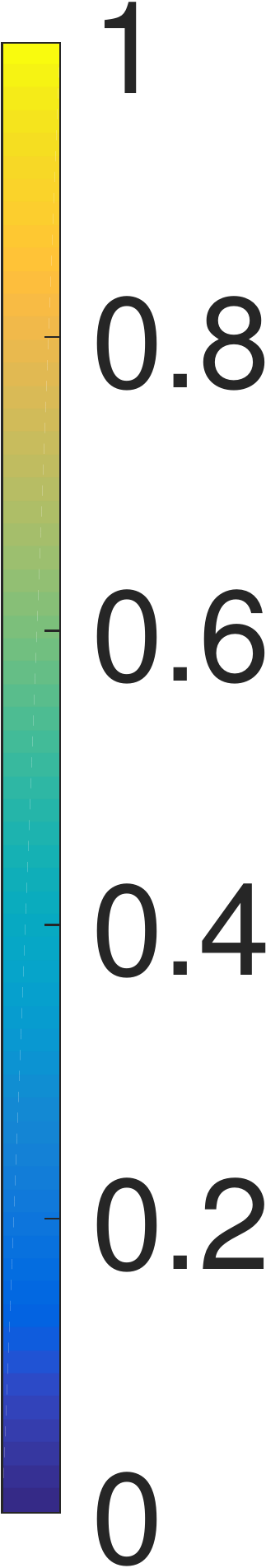}\label{fig:pres2d_homo}\hspace{-1cm}}%
\caption{Illustration of a 2D fine-scale computational grid (a) which contains $15 \times 15$ homogeneous matrix, $21$ fracture cells and two pressure-constrained wells at the West and East boundaries with values of $1$ and $0$, respectively. The pressure solution is plotted in (b), where fractures are 100 times more conductive than the matrix.}%
\label{fig:fineGrid}%
\end{figure}

These equations are to be solved for matrix and fracture pressures, $p^m$ and $p^f$, on the matrix $\Omega^m$ and fracture $\Omega^f$ domains, as depicted in Fig.~\ref{fig:fineGrid}. Note that a fracture network can consist of several fractures, which are represented in a lower dimensional space, i.e. $\Omega^f \subset R^{n-1}$, than the matrix (reservoir rock) $\Omega^m \subset R^{n}$. The main advantage of this type of formulation is that the matrix and fracture grids are independent and, thus, can be freely adapted to accommodate the appropriate physics for each medium. This is especially important in highly fractured reservoirs or when fractures are generated (and closed) during simulation, e.g., in geothermal formations \cite{KarvounisGeothermal}.

The incompressible single-phase pressure solution obtained using the EDFM approach for a 2D fractured reservoir model, shown in Fig.~\ref{fig:fine2d}, is provided in Fig.~\ref{fig:pres2d_homo}. Two pressure-constrained wells are placed on the East and West boundaries, and the reservoir rock is homogeneous. 

When non-linearities are present (e.g., compressible flows), the flow equations need to first be linearized, i.e.,
\begin{align}\label{dis-mf}
\bm A^{\nu} p^{\nu+1} \equiv {\left[\begin{array}{ccc}\bm A^{mm} & \bm A^{mf} & \bm A^{mw} \\ \bm A^{fm} & \bm A^{ff} & \bm A^{fw} \\ \bm A^{wm} & \bm A^{wf} & \bm A^{ww} \end{array}\right]^{\nu}}
\left[\begin{array}{c}p^m \\p^f\\p^w\end{array}\right]^{\nu+1}
=
\left[\begin{array}{c}q^m \\q^f \\ q^{w} \end{array}\right]^{\nu}
\equiv q^{\nu},
\end{align}
and then iteratively solved, in a Newton-Raphson loop, until the converged solution is achieved. Note that this system \eqref{dis-mf} shows an implicit treatment of the coupling between fracture and matrix through the $\bm A^{fm}$ entries, and that $\bm A$ can be non-symmetric, due to the compressibility effects \cite{Tene-cams,hadi-comp-jcp}.

Developing an efficient solution strategy for the linearized system \eqref{dis-mf} is quite challenging for several reasons. On the one hand, the size of this system can exceed several millions of unknowns for realistic test cases. On the other hand, the value of the condition numbers for the system matrix is worsened by high contrasts between reservoir properties (matrix permeability is highly heterogeneous over large scales, fractures are typically much more conductive than the matrix, etc.). 

Clearly a classical upscaling method cannot be employed here due to the highly resolved fractures, which play an important role in mass transport. This creates a niche for conservative multiscale methods, which have the important advantage of solving coarse-scale problems while honouring fine-scale data \cite{imsfv-jcp,Zhou} in an iterative error reduction strategy \cite{giuseppe-ib08,hadi-aimsfv-jcp,SeongTrJcp} which allows for mass-conservative flux reconstruction at any stage \cite{Jenny06}. Next, the development of the F-AMS method is presented. 

\section {Algebraic Multiscale Method for Heterogeneous Porous Media with Embedded Discrete Fractures (F-AMS)}
\label{sec:fams}

\begin{figure}[htb!]%
\centering%
\subfigure[2D coarse grid]{\includegraphics[width=0.42\textwidth,height=0.4\textwidth]{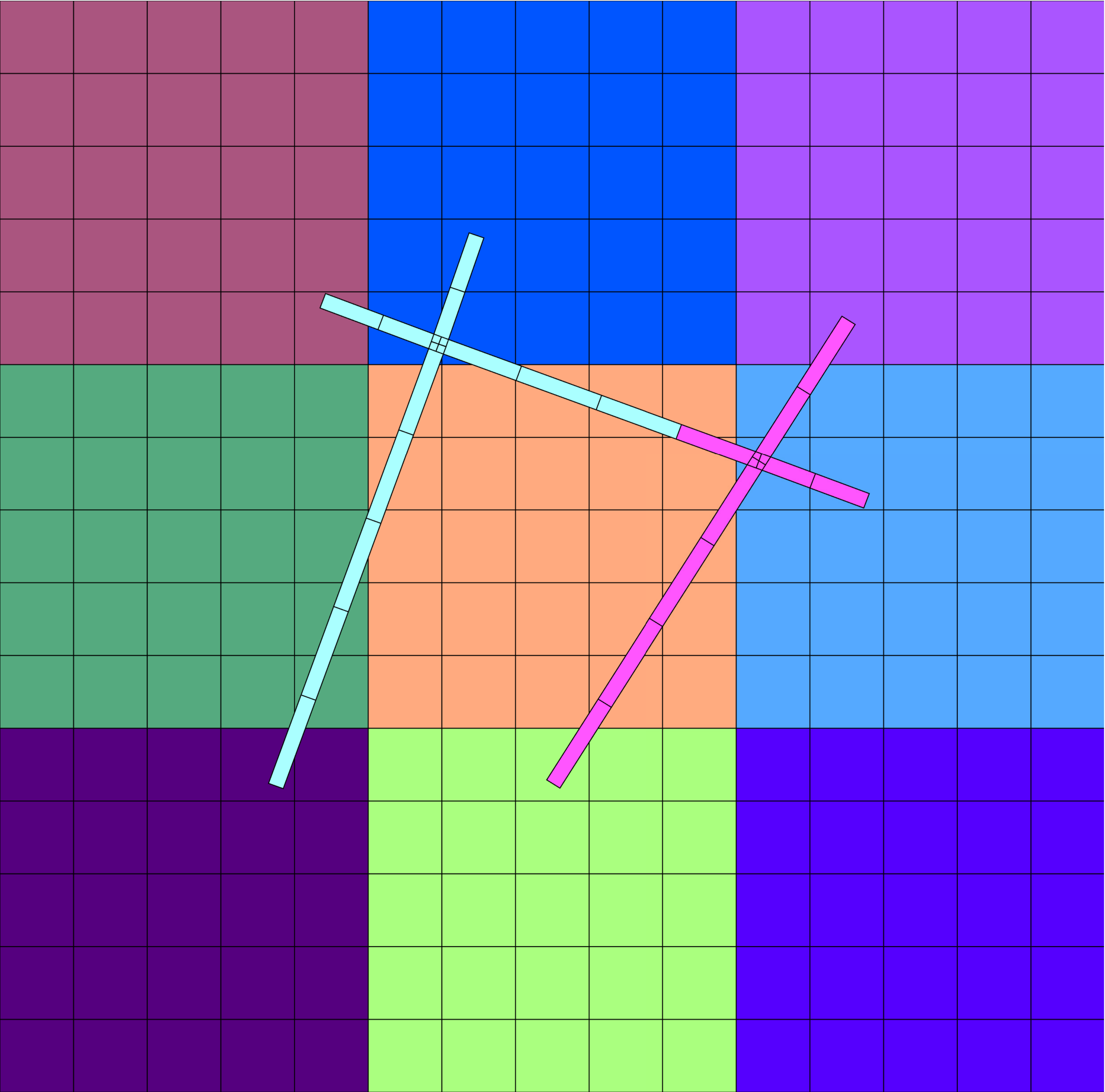}\label{fig:primal2d}}\hspace{1.5cm}%
\subfigure[2D dual-coarse grid]{\includegraphics[width=0.42\textwidth,height=0.4\textwidth]{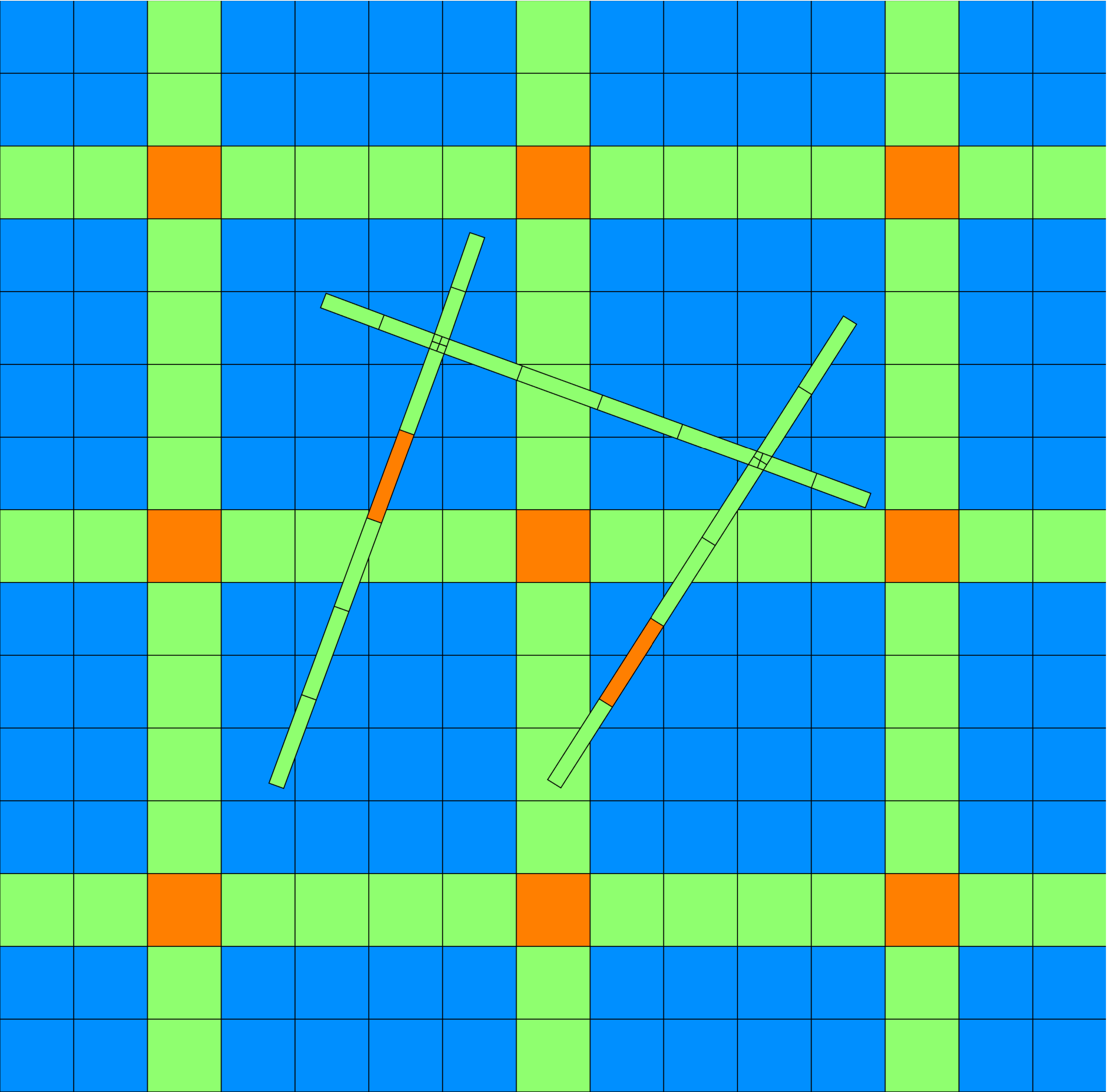}\label{fig:dual2d}}\\%
\subfigure[3D matrix coarse grid]{\includegraphics[width=0.42\textwidth,height=0.4\textwidth]{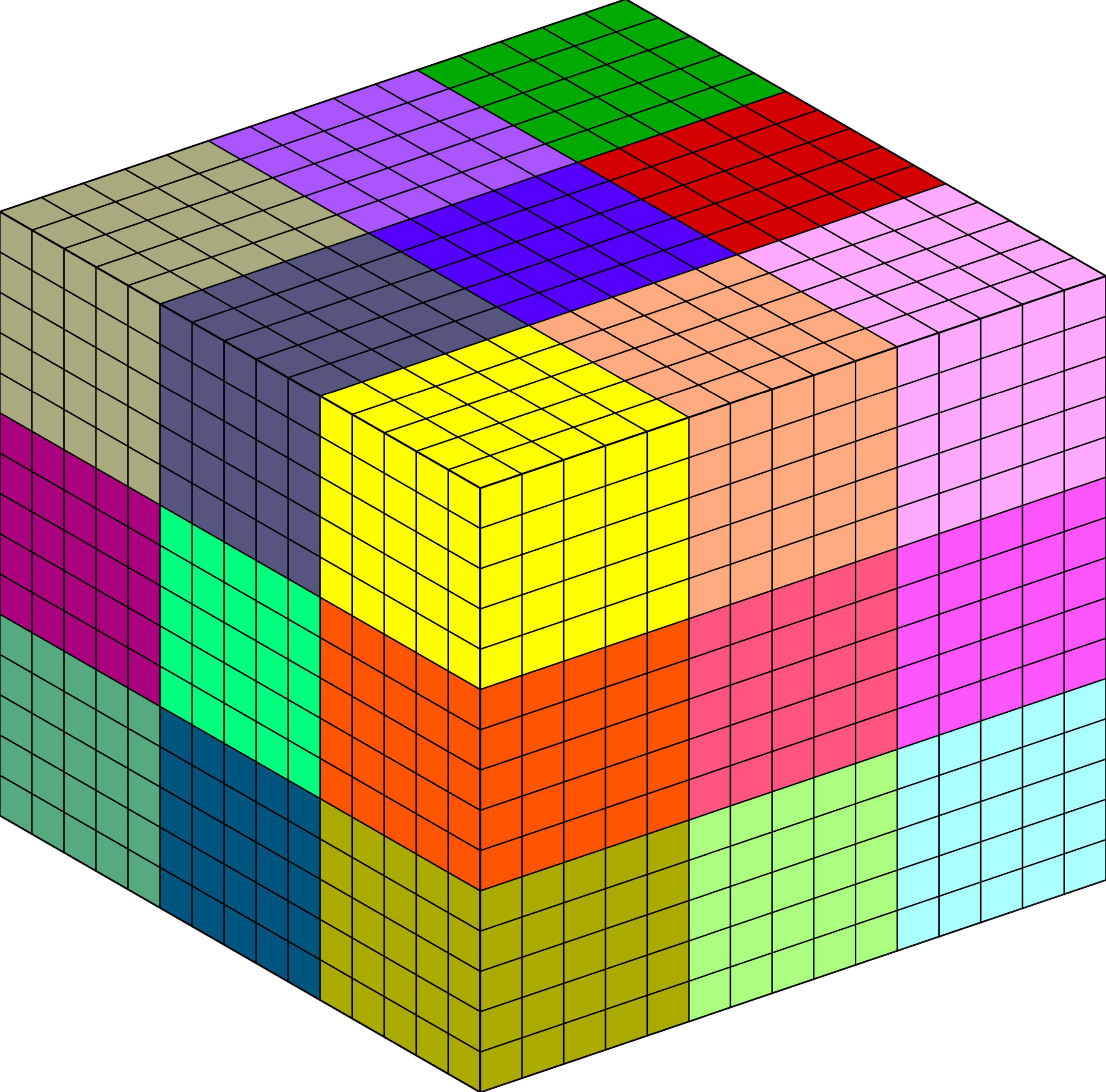}\label{fig:primal3d}}\hspace{1.5cm}%
\subfigure[3D matrix dual-coarse grid]{\includegraphics[width=0.42\textwidth,height=0.4\textwidth]{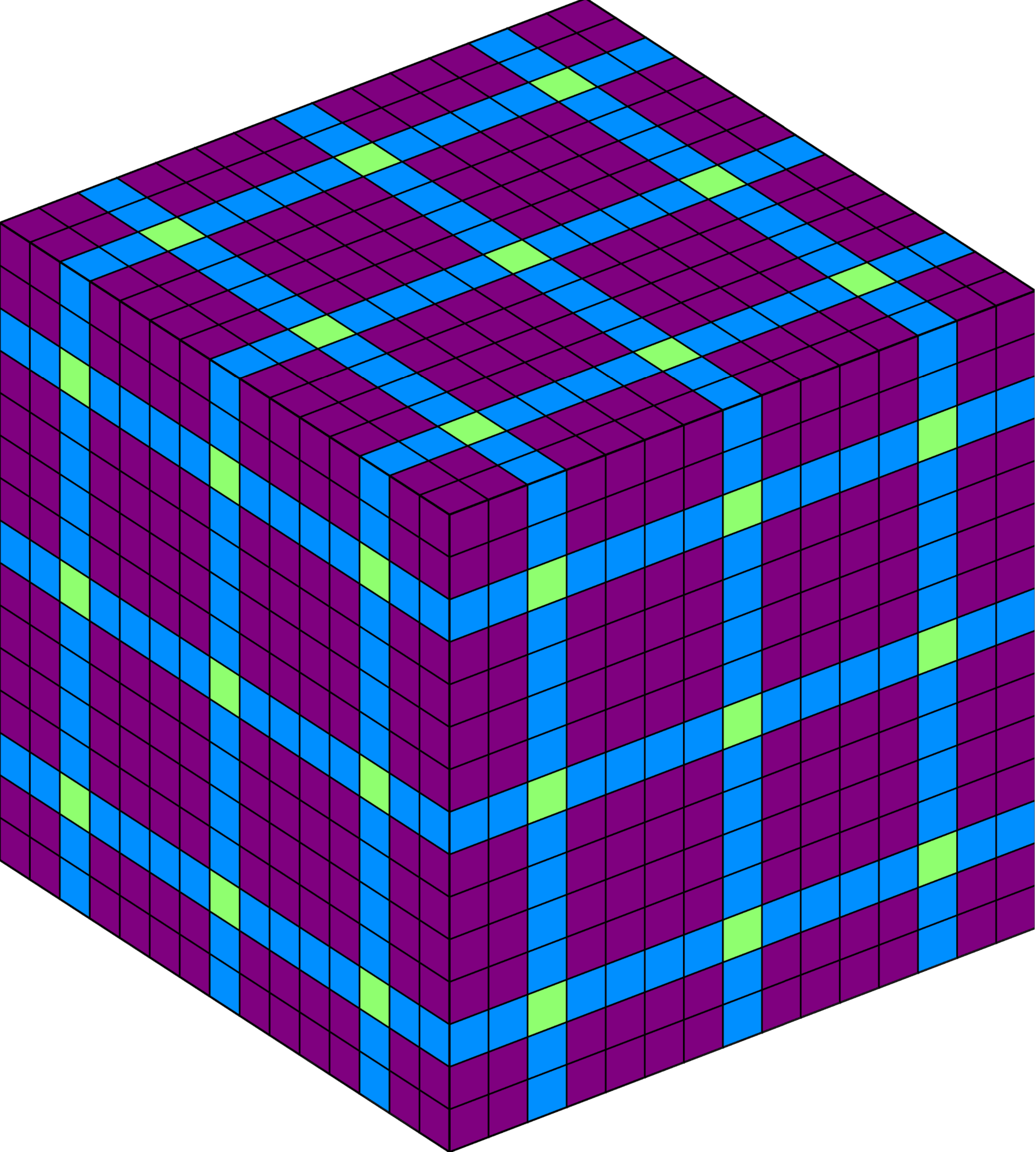}\label{fig:dual3d}}\\%
\subfigure[3D fracture coarse grid]{\includegraphics[width=0.42\textwidth,height=0.4\textwidth]{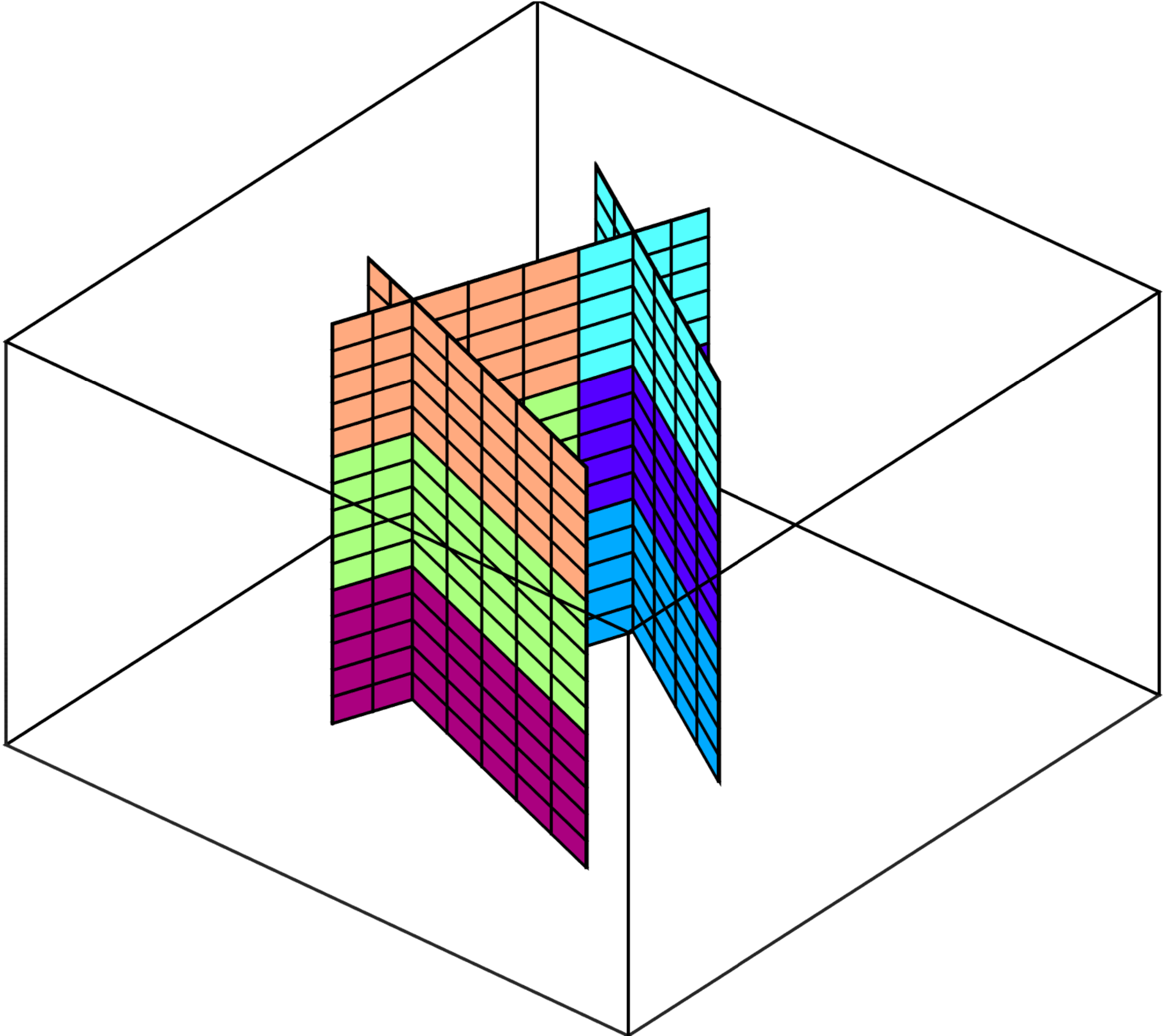}\label{fig:primal3dfrac}}\hspace{1.5cm}%
\subfigure[3D fracture dual-coarse grid]{\includegraphics[width=0.42\textwidth,height=0.4\textwidth]{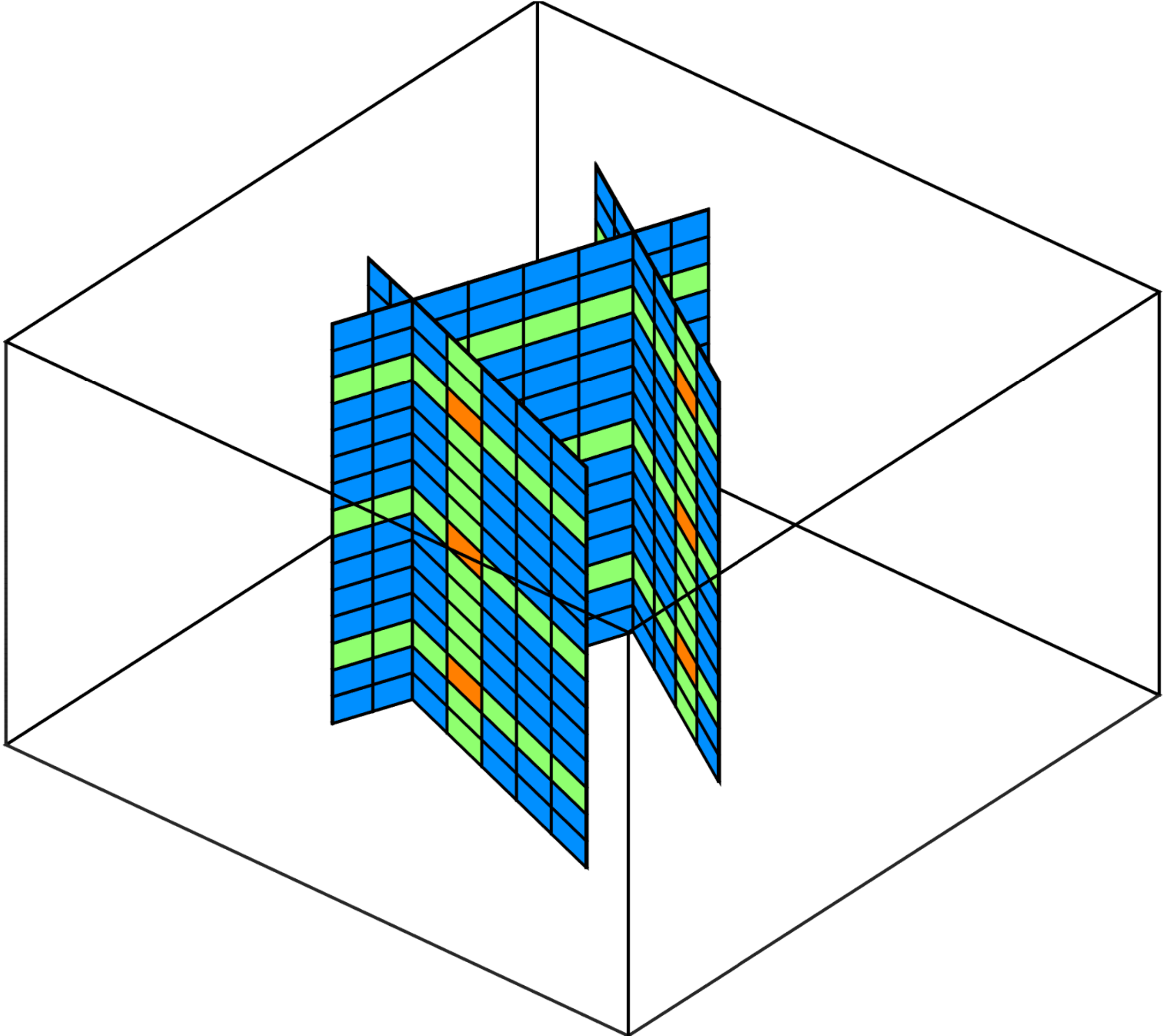}\label{fig:dual3dfrac}}%
\caption{F-AMS coarse grids defined on 2D (top) and 3D domains (middle and bottom). The primal grid (left) consists of non-overlapping coarse blocks, each shown in a different colour. The dual grid (right) is split into nodes, shown in orange, 1D blocks (edges) in green, 2D blocks (faces) in blue and 3D blocks (interiors) in purple. Note that the fracture aperture in (a) and (b) is magnified for clarity.}%
\label{fig:coarseGrids}%
\end{figure}

This section describes the F-AMS procedure, an efficient multiscale solution strategy for Eq.~\eqref{dis-mf}. Given a computational domain with $N_f$ fracture networks and $N_w$ wells, F-AMS first superimposes two coarse grids on top of both the matrix and fracture domains. The primal-coarse grid is a non-overlapping decomposition of the domain, inside which a fine-scale cell is selected as coarse node (Fig.~\ref{fig:primal2d} for 2D, and \ref{fig:primal3d} and \ref{fig:primal3dfrac} for 3D cases). By connecting the coarse nodes, a secondary overlapping coarse grid is obtained, which is called the dual coarse grid (Figs.~\ref{fig:dual2d}, \ref{fig:dual3d} and \ref{fig:dual3dfrac}). There exist $N_{cm}$ and $N_{dm}$ matrix primal-coarse and dual-coarse blocks and, similarly, each fracture network $f_i$ contains $N_{cf_i}$ and $N_{df_i}$ fracture primal-coarse and dual-coarse blocks. Note that $N_w$ (injection or production) wells exist in the domain, as driving forces for the flow.

F-AMS approximates the solution to Eq.~\eqref{dis-mf}, $p$, as a superposition of coarse-scale solutions ($\breve p$) using locally computed basis functions ($\Phi$), i.e.
\begin{align}\label{pm-ms}
p^m \approx {p'^m} = \sum_{i=1}^{N_{cm}} \Phi_i^{mm} \breve p_i^m + \sum_{i=1}^{N_{f}} \sum_{j=1}^{N_{cf_i}} \Phi_j^{mf_i} \breve p_j^{f_i} +  \sum_{k=1}^{N_{w}} \Phi_j^{mw} \breve p_k^{w},
\end{align}
for the matrix and 
\begin{align}\label{ms-f}
p^f \approx {p'^f} = \sum_{i=1}^{N_{cm}} \Phi_i^{fm} \breve p_i^m + \sum_{i=1}^{N_{f}} \sum_{j=1}^{N_{c{f_i}}} \Phi_j^{ff_i} \breve p_j^{f_i} + \sum_{k=1}^{N_{w}} \Phi_j^{fw} \breve p_k^{w}
\end{align}
for the fractures, respectively. The basis functions associated with matrix coarse cells (i.e., $\Phi^{m*}$) are $\Phi_i^{mm}$ for matrix-matrix effects, $\Phi_j^{mf}$ for the matrix-fracture coupling, and $\Phi_j^{mw}$ matrix-well interactions. These basis functions (interpolators) are employed in order to capture the effects of all the important factors (matrix, fractures, and wells) in the construction of a good approximation for the matrix pressure field $p^m$. Similarly for fractures, $\Phi^{f*}$ consists of the contributions from the matrix $\Phi^{fm}$, fractures $\Phi^{ff}$, and wells $\Phi^{fw}$, if present. 

One of the novel aspects of this work is that the pressure field inside fractures, $p^f$, is included explicitly in the multiscale formulation (Eq.~\eqref{ms-f}). This means that the fracture grid cells are also decomposed into primal and dual coarse blocks (Fig.~\ref{fig:coarseGrids}), similar to the matrix. Their solutions are also mapped to the coarse scale and back to the original resolution, again, similar to the matrix. More specifically, each fracture network $f_i$ is decomposed into $N_{cf_i}$ primal-coarse grid blocks, for which sets of basis functions are calculated. One could employ the same formulation for wells, i.e., discretize them into several fine-scale cells which can then be coarsened on the superimposed multiscale grids. However, for the sake of simplicity, in the experiments presented in this work, each well is assigned a single fine-scale DOF, which is mapped to the coarse-scale using the identity restriction operator, i.e.
\begin{align}
p^{w_i} = p'^{w_i} = \breve p^{w_i} \hspace{1cm} \forall i \in \{1, \cdots, N_w\}.
\end{align}

In algebraic notation, the superpositions \eqref{pm-ms} and \eqref{ms-f} can be expressed as
\begin{align}\label{pm-ams}
p^m \approx {p'}^m = \bm{\mathcal{P}}^m \breve p \equiv [\bm{\mathcal{P}}^{mm} ~ ~ \bm{\mathcal{P}}^{mf} ~ ~ \bm{\mathcal{P}}^{mw}] ~ [\breve p^m  ~ ~ \breve p^f  ~ ~ \breve p^w]^T
\end{align}
and 
\begin{align}
p^f \approx {p'}^f =  \bm{\mathcal{P}}^f \breve p \equiv [\bm{\mathcal{P}}^{fm} ~ ~  \bm{\mathcal{P}}^{ff}  ~ ~ \bm{\mathcal{P}}^{fw}] ~ [\breve p^m  ~ ~ \breve p^f  ~ ~ \breve p^w]^T,
\end{align}
respectively.

The basis functions are assembled in the columns of the multiscale \emph{prolongation operator}, $\bm{\mathcal{P}}$, with the dimension of $N_{fine} \times N_{coarse}$, where $N_{fine}$ and $N_{coarse}$ are the total number of fine- and coarse-scale control volumes, respectively. The part of $\bm{\mathcal{P}}$ corresponding to the matrix fine-cells is defined as
\begin{align}\label{pm-ams-pro}
\bm{\mathcal{P}}^m = 
\left[\begin{array}{ccc|ccccccc|ccc}
\vdots & \cdots & \vdots & \vdots & \cdots & \vdots& \cdots & \vdots& \cdots&\vdots&\vdots & \cdots & \vdots\\
\Phi^{mm}_1 & \cdots & \Phi^{mm}_{N_{cm}} & \Phi^{mf_1}_1 & \cdots & \Phi^{mf_1}_{N_{cf_1}}& \cdots & \Phi^{mf_{N_f}}_1& \cdots& \Phi^{mf_{N_f}}_{N_{cf_{N_f}}} &\Phi^{mw}_1 & \cdots & \Phi^{mw}_{N_w}\\
\vdots & \cdots & \vdots & \vdots & \cdots & \vdots & \cdots & \vdots& \cdots&\vdots&\vdots & \cdots & \vdots
\end{array}\right].
\end{align}
Notice the three sub-blocks which represent matrix-matrix, matrix-fracture, and matrix-well coupling. Similarly, the prolongation operator for fractures can be stated as
\begin{align}
\bm{\mathcal{P}}^f  =
\left[\begin{array}{ccc|ccccccc|ccc}
\vdots & \cdots & \vdots & \vdots & \cdots & \vdots& \cdots & \vdots& \cdots&\vdots&\vdots & \cdots & \vdots\\
\Phi^{fm}_1 & \cdots & \Phi^{fm}_{N_{cm}} & \Phi^{ff_1}_1 & \cdots & \Phi^{ff_1}_{N_{cf_1}}& \cdots & \Phi^{ff_{N_f}}_1& \cdots& \Phi^{ff_{N_f}}_{N_{cf_{N_f}}}&\Phi^{fw}_1 & \cdots & \Phi^{fw}_{N_w}\\
\vdots & \cdots & \vdots & \vdots & \cdots & \vdots & \cdots & \vdots& \cdots&\vdots& \vdots& \cdots&\vdots
\end{array}\right].
\end{align}

Algebraically, the complete F-AMS prolongation operator reads
\begin{align}\label{prolongation_single}
\bm{\mathcal{P}} = 
\left[\begin{array}{c}
\bm{\mathcal{P}}^{m} \\ \bm{\mathcal{P}}^{f} \\ \bm{\mathcal{P}}^{w} 
\end{array} \right]
=
\left[\begin{array}{ccc}
\bm{\mathcal{P}}^{mm} & \bm{\mathcal{P}}^{mf} & \bm{\mathcal{P}}^{mw} \\
\bm{\mathcal{P}}^{fm}  & \bm{\mathcal{P}}^{ff} & \bm{\mathcal{P}}^{fw}\\
\bm{\mathcal{P}}^{wm} & \bm{\mathcal{P}}^{wf} & \bm{\mathcal{P}}^{ww}\\   
\end{array} \right],
\end{align}
where $\bm{\mathcal{P}}^{wm}$ and $\bm{\mathcal{P}}^{wf}$ are set to zero, while $\bm{\mathcal{P}}^{ww}$ is the identity matrix. 

Note that the prolongation operator, as described in Eq.~\eqref{prolongation_single}, allows full flexibility in consideration of the fracture-matrix coupling in the interpolated solution, i.e. via the values in $\bm{\mathcal{P}}^{mf}$ and $\bm{\mathcal{P}}^{fm}$. This leads to the definition of four operators, differentiated by the coupling strategy they employ:
\begin{enumerate}
\item \textbf{Decoupled-AMS:} $ \bm{\mathcal{P}}^{mf} =  0$ and $\bm{\mathcal{P}}^{fm} = 0$
\item \textbf{Frac-AMS:} only $\bm{\mathcal{P}}^{fm} = 0$.
\item \textbf{Rock-AMS:} only $\bm{\mathcal{P}}^{mf} = 0$.
\item \textbf{Coupled-AMS:} $\bm{\mathcal{P}}^{mf}$ and $\bm{\mathcal{P}}^{fm}$ both non-zero.
\end{enumerate}
The first option, i.e., Decoupled-AMS, constructs the most sparse $\bm{\mathcal{P}}$ and thus has an efficient setup phase. The fourth option, i.e., Coupled-AMS, can lead to more accurate multiscale simulations, however, it can severely increase the density of the multiscale operators. In such cases, one may be able to obtain a trade-off between the quality of the prolongation operator and its sparsity via truncation, followed by a rescaling of the rows to ensure partition of unity. A CPU-based study considering the overhead introduced by the density of $\bm{\mathcal{P}}$ is presented in Section~\ref{sec:results}. 

In order to construct the coarse-scale system, F-AMS also needs the specification of a restriction operator, which is a map from fine- to coarse-scale (dimension $N_{coarse} \times N_{fine}$). Due to its algebraic formulation, F-AMS can accommodate multiscale finite volumes (MSFV), multiscale finite elements (MSFE) or even a hybrid multiscale finite elements and volumes restriction (MSMIX). More specifically, the MSMIX employs a FV-based restriction for part of the domain (e.g., fractures or wells), and FE for the rest (e.g., matrix rock). It is important to note that, after any MSFV stage, it is possible to construct a mass-conservative flux field for both matrix and fractures. As such, in multiphase simulations, if iterations for pressure Eq.~\eqref{dis-mf} are stopped before full convergence is achieved, MSFV needs to be employed before solving transport equations. On the other hand, MSFE leads to a symmetric-positive-definite (SPD) coarse system if the fine-scale system matrix is also SPD \cite{yixuan-ams,Tene-cams}, and is the option used during all numerical experiments presented in this work. Note that MSMIX can be tweaked to achieve the desired compromise between MSFV and MSFE. 

In the following sub-sections, first, the formulation of the local basis functions is explained. Then, the F-AMS system, and finally the simulation strategy is described in detail.

\subsection{Basis function formulations}
\label{sec:basis}

As stated before, F-AMS constructs a non-overlapping partition on the given fine-scale computational domain for both matrix rock and fracture cells, i.e. the primal-coarse grid. Then, by connecting the coarse nodes, the overlapping decomposition of the domain, i.e., dual-coarse grid, is obtained. Following the original description of the MSFV basis functions \cite{Jenny03,Jenny06} and its algebraic description \cite{Zhou-tams,yixuan-ams,Tene-cams}, local basis functions are calculated for each coarse node $i$, corresponding to each dual block, by respecting the \emph{wirebasket} hierarchy \cite{wirebasket}. First the pressure in the coarse nodes (also called vertices, shown in orange in Fig.~\ref{fig:dual2d}) is set to
\begin{align}
\delta_{ij} =%
\begin{cases}
1, & \text{ if } i=j \\
0, & \text{ if } i \neq j
\end{cases}.
\label{kronecker-basis}
\end{align}
Then the dual blocks in the neighbourhood of node $i$ are resolved, in sequence, as follows: first all the 1D dual blocks (or edges, shown in green in Fig.~\ref{fig:dual2d}), followed by the 2D (or faces, shown in blue in Fig.~\ref{fig:dual2d}) and, finally, if applicable, the 3D dual blocks (or interiors, shown in purple in Fig.~\ref{fig:dual3d}). The fact that each dual block (e.g. edge) neglects the transmissibilities to neighbouring cells belonging to blocks of inferior rank in the wirebasket hierarchy (i.e. faces and interiors), constitutes the \emph{localization assumption} \cite{yixuan-ams}, which ensures that each basis function has a limited support.

The pressure values obtained in the manner described above, for each coarse node $i$, are assembled in column $i$ of the prolongation operator to form basis function $\Phi_i$.

It is important to note that by having independent fine-scale grids for each media, a matrix cell (say from a face block) can be directly connected to fracture cells belonging to dual blocks of any type (vertex, edge or face). This is an important difference from non-fractured media, where the two-point flux approximation (TPFA) stencil on structured grids ensured that any dual block would have connections only to blocks of directly superior or inferior rank in the wirebasket hierarchy (e.g., the external neighbours of a face cell are either edge or interior cells). As such, the multscale localization assumption needs to be extended to account for the connection between the two media. In the F-AMS framework, this leads to definition of basis functions which account for different degrees of coupling between the matrix and its perforating fractures. 

In order to provide a compact definition of the various basis functions (matrix, fracture, well), paired with one of the four different coupling strategies considered (Decoupled-AMS, Rock-AMS, Frac-AMS, Coupled-AMS), the following ``skeleton'' is introduced
\begin{align}
- \nabla \cdot (\lambda^* \cdot \nabla \Phi^{*\bullet})%
+ \sum_{j \in \text{conn}_{mf}^*} \eta^*_j~ \xi(\Phi^{*\bullet})%
+ \sum_{j \in \text{perf}_w^*} \beta^*_j~ (\Phi^{*\bullet} - \Phi^{w\bullet}) %
= 0, \label{skeleton-basis}
\end{align}
which is solved for all basis functions $\forall \Phi^{*\bullet} \in \{\Phi^{mm},\Phi^{mf},\Phi^{mw},\Phi^{fm},\Phi^{ff},\Phi^{fw}\}$, subject to the localization assumption within each domain. Recall that $\bm{\mathcal{P}}^{wm}$ and $\bm{\mathcal{P}}^{wf}$ are zero, while $\bm{\mathcal{P}}^{ww}$ is the identity matrix. The skeleton expression \eqref{skeleton-basis} is based on the incompressible pressure equation, since it was found computationally efficient, even when compressibility is involved \cite{Tene-cams}. In its definition, $\text{perf}_w^*$ represents the set of (matrix or fracture) cells perforated by wells. Moreover, $\text{conn}_{mf}^*$ is the set of all cells with cross-media connectivities from the corresponding (matrix or fracture) domain. Finally, the $\xi(\Phi^{*\bullet})$ function gives the type of matrix-fracture coupling captured by the basis function, and will be specified separately for each strategy, as follows:

\begin{enumerate}[itemsep=0.5cm]

\item \textbf{Decoupled-AMS:} all basis functions have no-flow boundary conditions between the matrix and fracture domains, i.e.,
\begin{align}
\xi(\Phi^{*\bullet}) = 0 \hspace{5mm} \forall \Phi^{*\bullet} \in \{\Phi^{mm},\Phi^{mf},\Phi^{mw},\Phi^{fm},\Phi^{ff},\Phi^{fw}\},
\end{align}
which means that the fracture-matrix coupling term is completely omitted in Eq~\eqref{skeleton-basis}. The prolongation is solved algebraically, as described in Appendix~\ref{sec:app-decoupledAMS}. Figure~\ref{fig:decoupledAMS} illustrates the step-by-step procedure for a fracture and matrix basis function belonging to the 2D reservoir from Fig.~\ref{fig:fineGrid}. Note that the support of each of the interpolators is restricted to their containing medium. Finally, the Decoupled-AMS approach can be seen as applying the original AMS to separate sub-domains (i.e., matrix and fractures), having them coupled only at the coarse-scale system.

\begin{figure}[htb!]%
\centering%
\setlength{\fboxsep}{0pt}%
\setlength{\fboxrule}{1pt}%
\setlength{\unitlength}{1cm}%
\begin{picture}(15,0.5)%
  \put(5.4,0.2){\large Decoupled-AMS}%
\end{picture}\\%
\subfigure[$\Phi^{mm}$ vertices]{\fbox{\includegraphics[width=0.3\textwidth,height=0.3\textwidth]{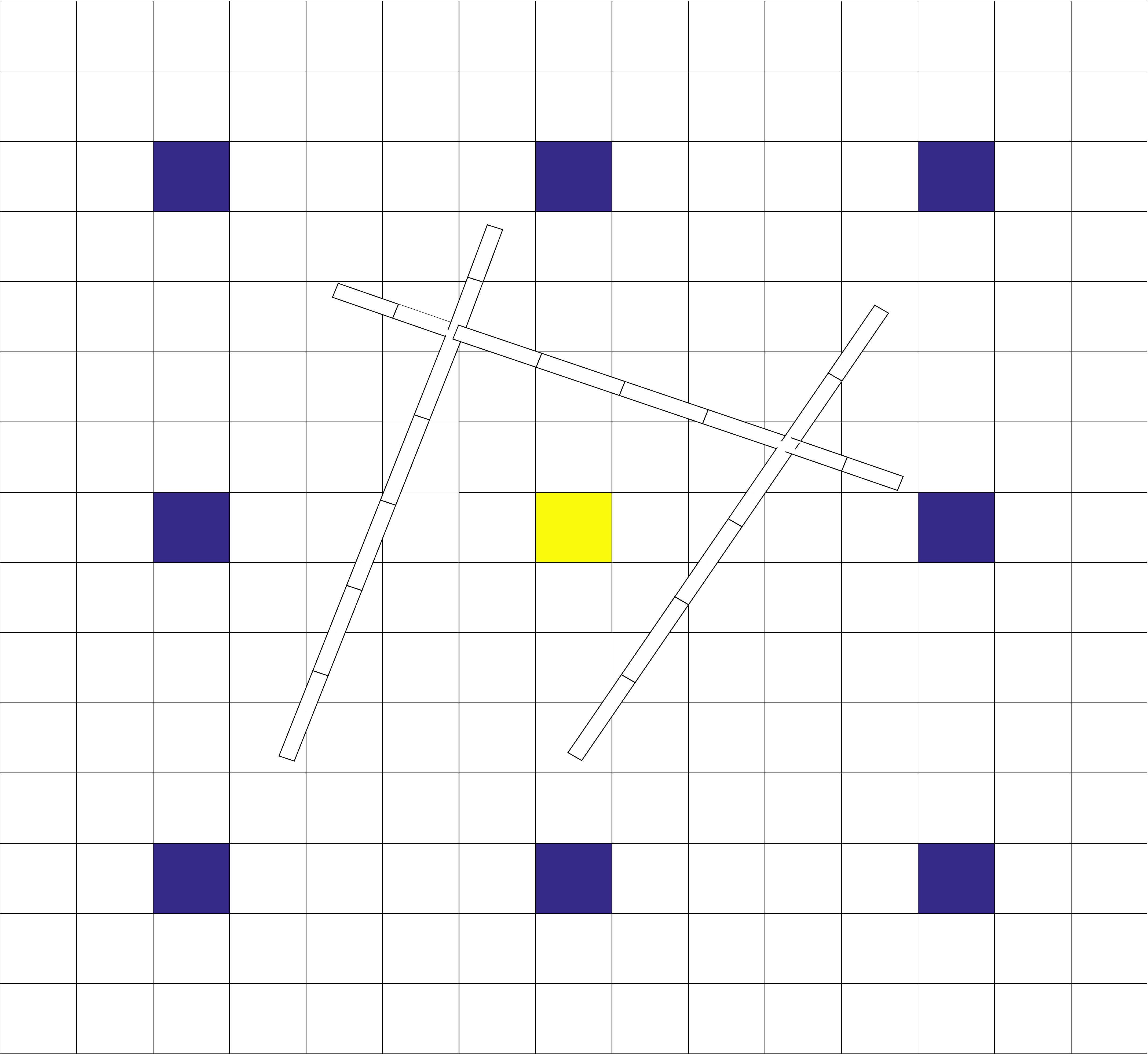}\label{fig:decoupledAMS_rock_vertices}}}\hspace{0.3cm}%
\subfigure[$\Phi^{mm}$ edges]{\fbox{\includegraphics[width=0.3\textwidth,height=0.3\textwidth]{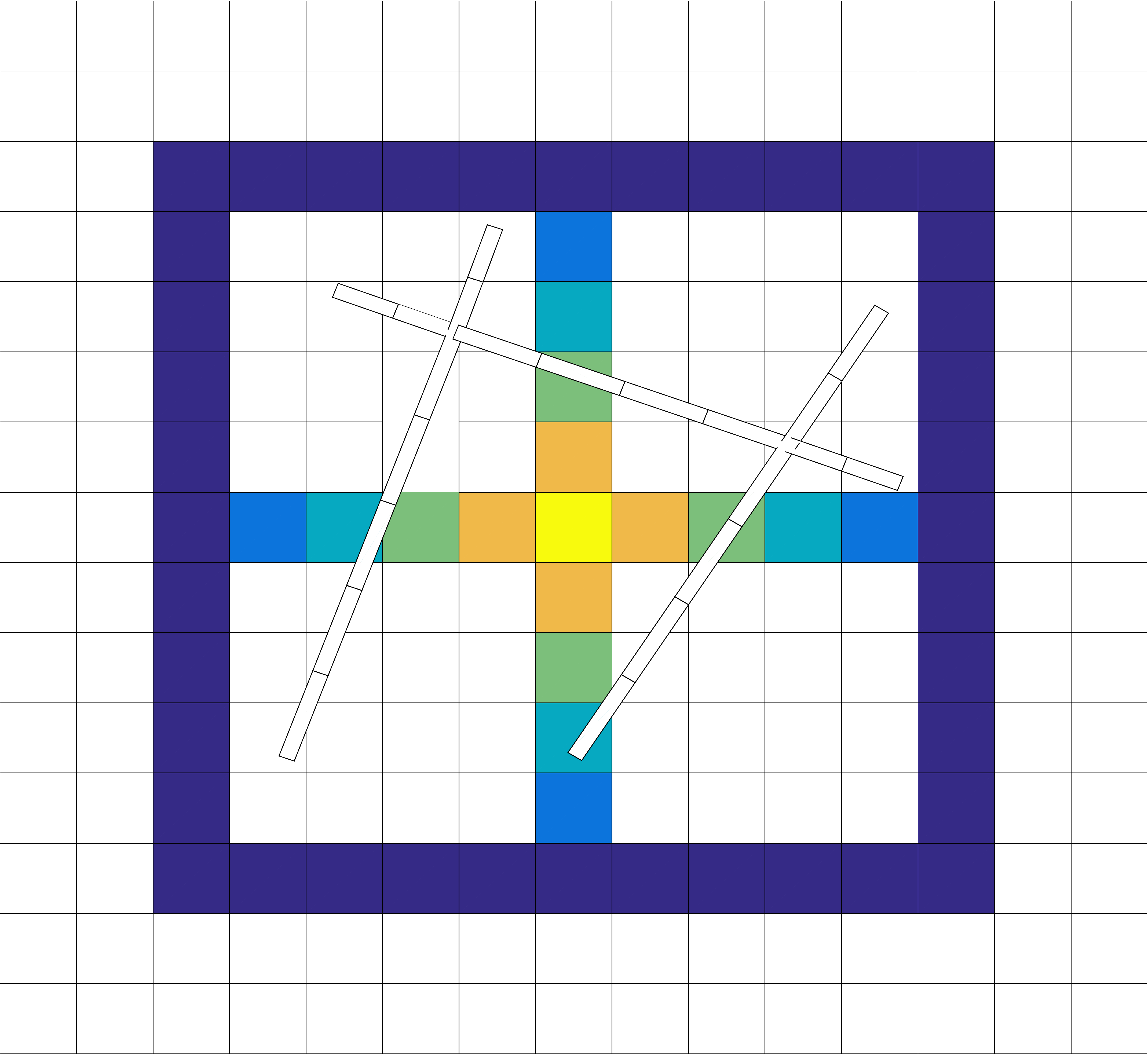}\label{fig:decoupledAMS_rock_edges}}}\hspace{0.3cm}%
\subfigure[$\Phi^{mm}$ faces]{\fbox{\includegraphics[width=0.3\textwidth,height=0.3\textwidth]{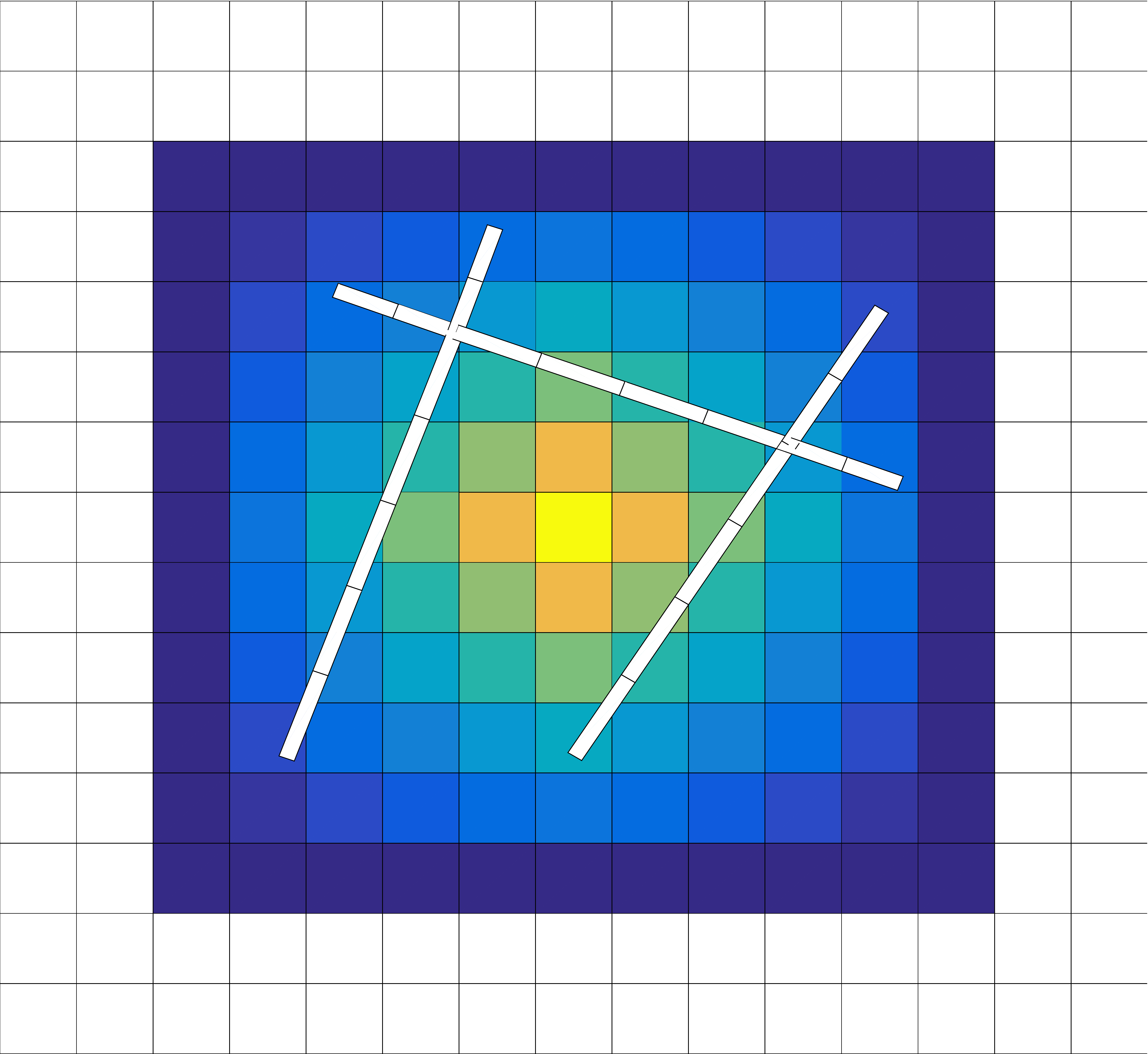}}\hspace{0.3cm}\includegraphics[height=0.3\textwidth]{Fig/method/colorbar.pdf}\label{fig:decoupledAMS_rock}\hspace{-0.8cm}}\\%
\subfigure[$\Phi^{ff}$ vertices]{\fbox{\includegraphics[width=0.3\textwidth,height=0.3\textwidth]{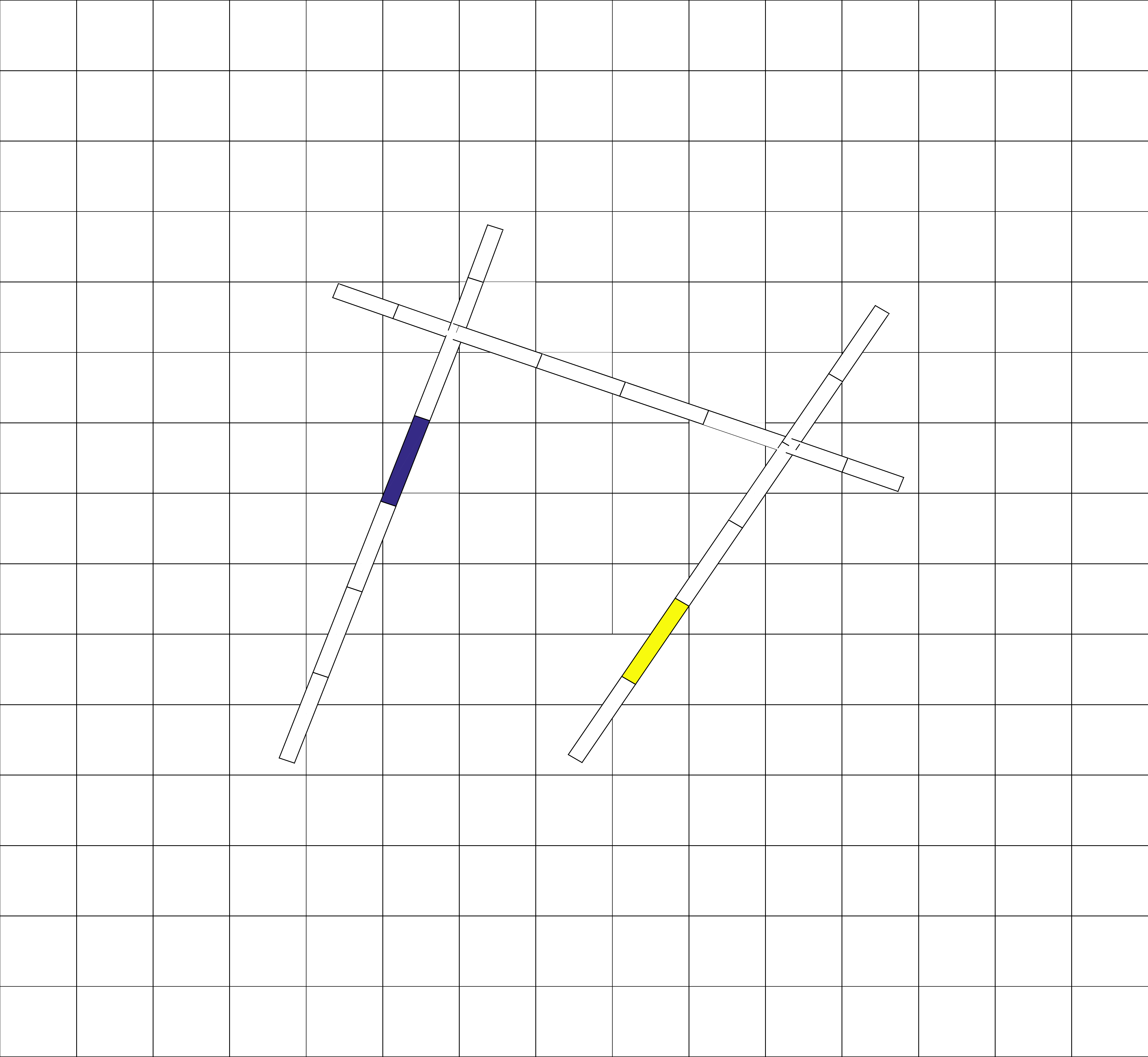}\label{fig:decoupledAMS_frac2_vertices}}}\hspace{1cm}%
\subfigure[$\Phi^{ff}$ edges]{\fbox{\includegraphics[width=0.3\textwidth,height=0.3\textwidth]{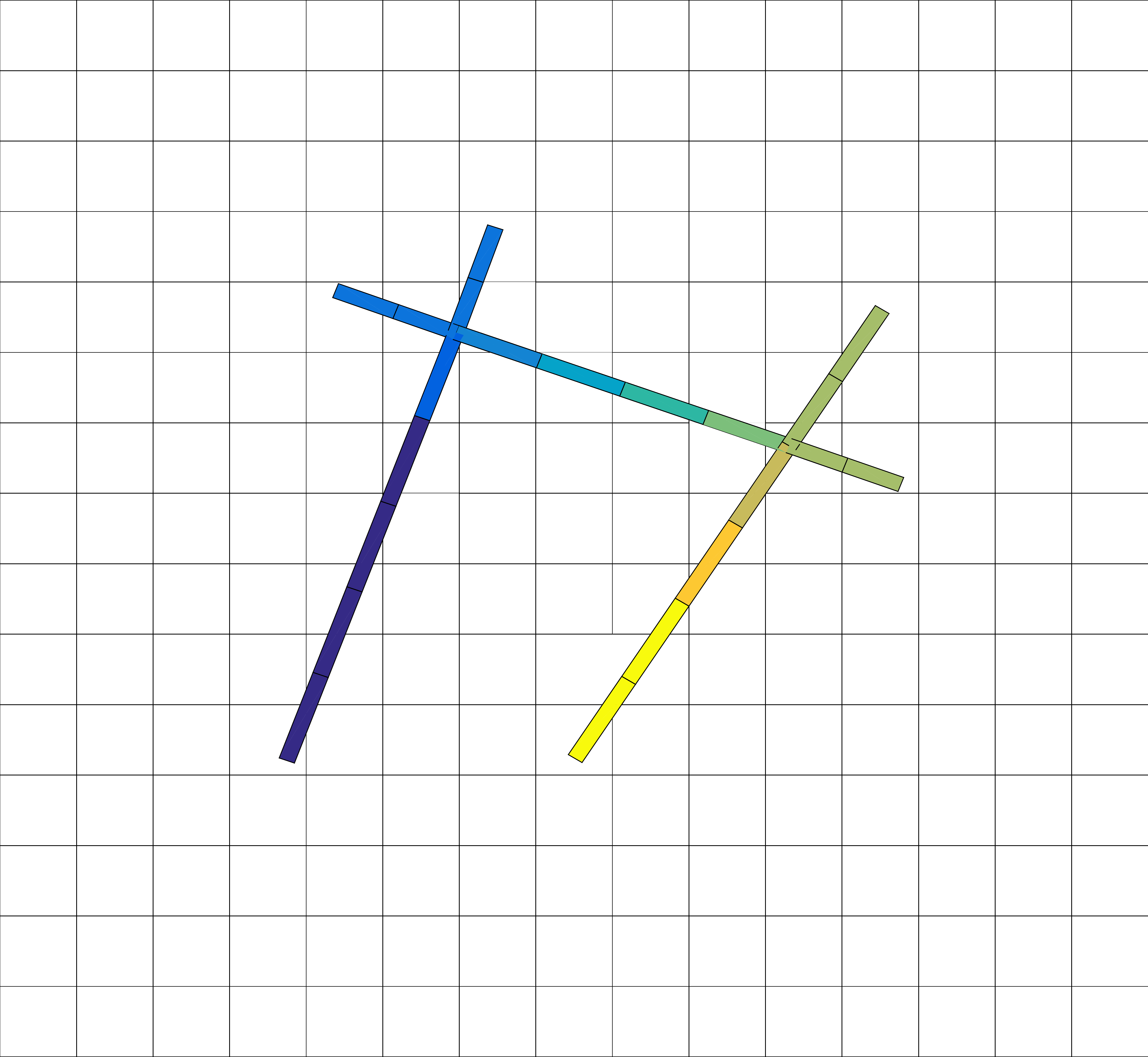}}\hspace{0.3cm}\includegraphics[height=0.3\textwidth]{Fig/method/colorbar.pdf}\label{fig:decoupledAMS_frac2}\hspace{-0.8cm}}%
\caption{Step-by-step construction of a Decoupled-AMS matrix basis function (top) and fracture function (bottom). In this strategy, $\Phi^{mf}$ and $\Phi^{fm}$ are explicitly set to $0$ and the two media have no-flow boundary conditions towards each other. Note that the fracture aperture is magnified for clarity.}%
\label{fig:decoupledAMS}%
\end{figure}%

\item  \textbf{Frac-AMS:} the fracture basis functions in the fracture domain, $\Phi^{ff}$, are first computed subject to no-flow conditions towards the matrix, the same as in the Decoupled-AMS (Fig.~\ref{fig:decoupledAMS_frac2_vertices} and \ref{fig:decoupledAMS_frac2}), i.e. by substituting
\begin{align}
\xi(\Phi^{ff}) = 0
\end{align}
in Eq.~\eqref{skeleton-basis}. Then, the obtained values are fixed and used as Dirichlet boundary conditions while computing $\Phi^{mf}$, for which the matrix-fracture transmissibility is taken into account, i.e.,
\begin{align}
\xi(\Phi^{mf}) = \Phi^{mf} - \Phi^{ff}.
\end{align}

\noindent On the other hand, the matrix basis functions, $\Phi^{mm}$, are solved by setting $\Phi^{fm} = 0$ as Dirichlet boundary condition, i.e., 
\begin{align}
\xi(\Phi^{mm}) = \Phi^{mm}.
\end{align}

\noindent This procedure is performed algebraically, as described in Appendix~\ref{sec:app-fracAMS}. Note from Fig.~\ref{fig:fracAMS} that, after this computation, the fracture functions have non-zero values in the matrix, while the support of the matrix basis functions is restricted to the rock domain.

\begin{figure}[htb!]%
\centering%
\setlength{\fboxsep}{0pt}%
\setlength{\fboxrule}{1pt}%
\setlength{\unitlength}{1cm}%
\begin{picture}(15,0.5)%
  \put(5.6,0.2){\large Frac-AMS}%
\end{picture}\\%
\subfigure[$\Phi^{mm}$ vertices]{\fbox{\includegraphics[width=0.3\textwidth,height=0.3\textwidth]{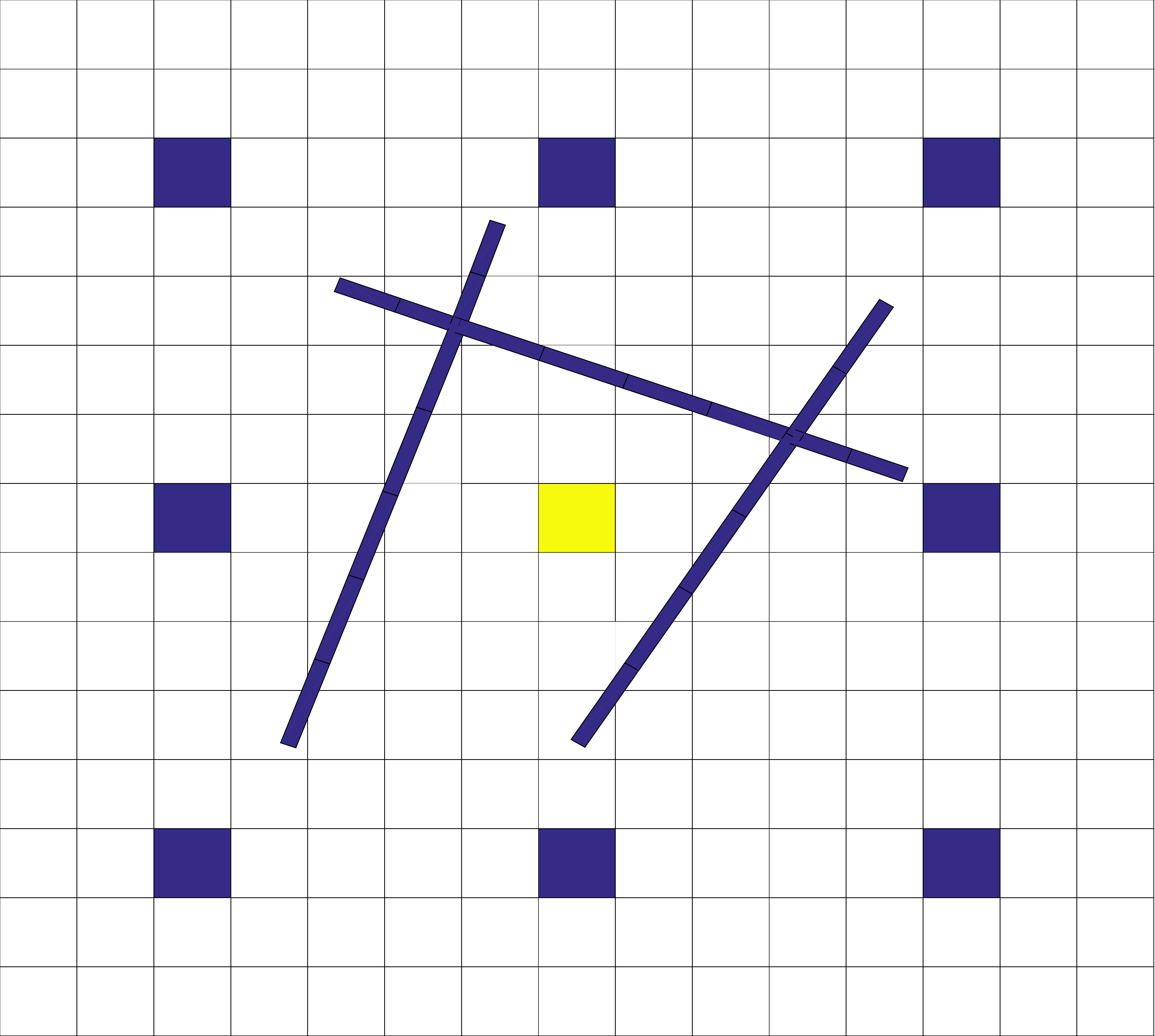}\label{fig:fracAMS_rock_vertices}}}\hspace{0.3cm}%
\subfigure[$\Phi^{mm}$ edges]{\fbox{\includegraphics[width=0.3\textwidth,height=0.3\textwidth]{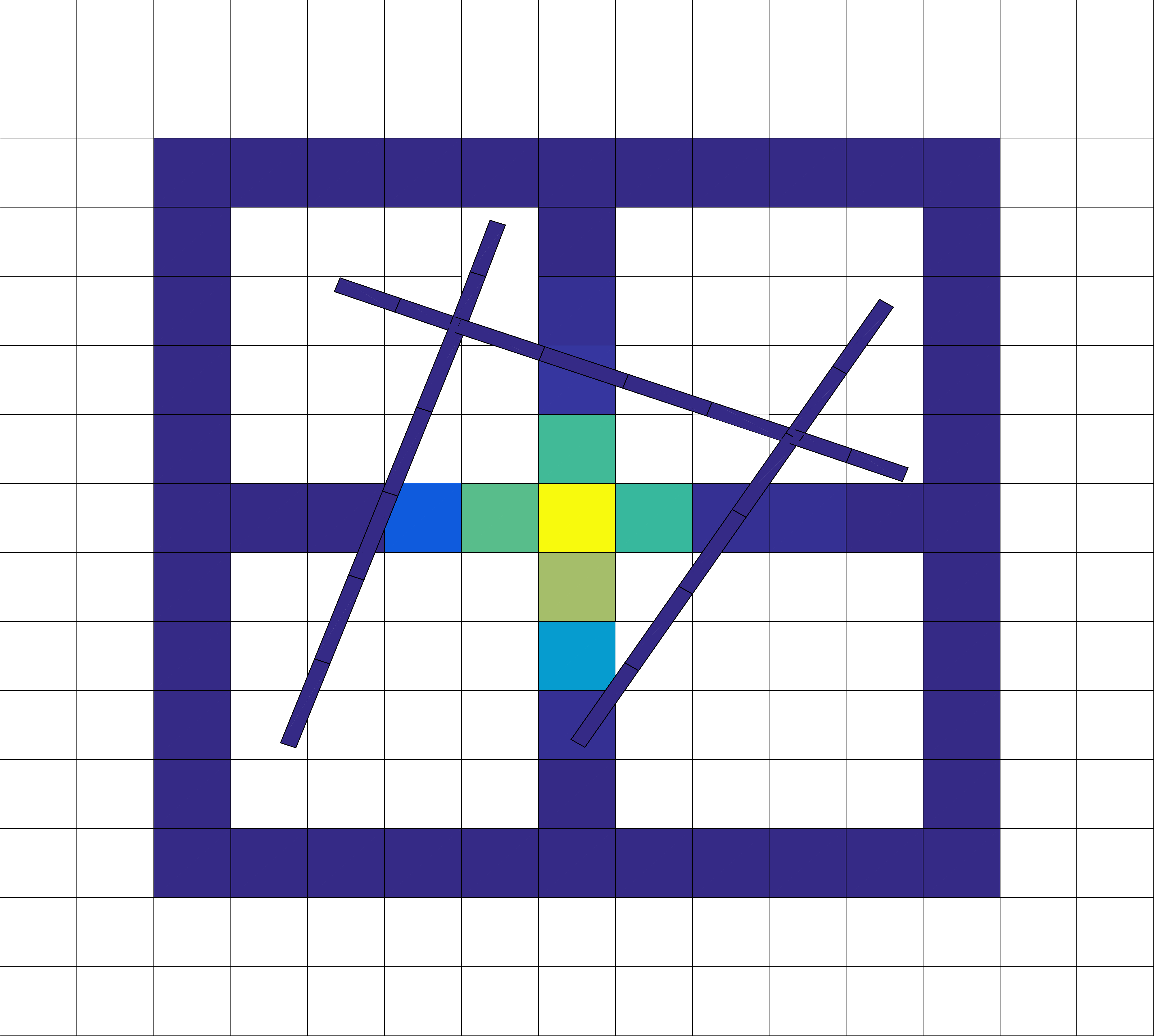}\label{fig:fracAMS_rock_edges}}}\hspace{0.3cm}%
\subfigure[$\Phi^{mm}$ faces]{\fbox{\includegraphics[width=0.3\textwidth,height=0.3\textwidth]{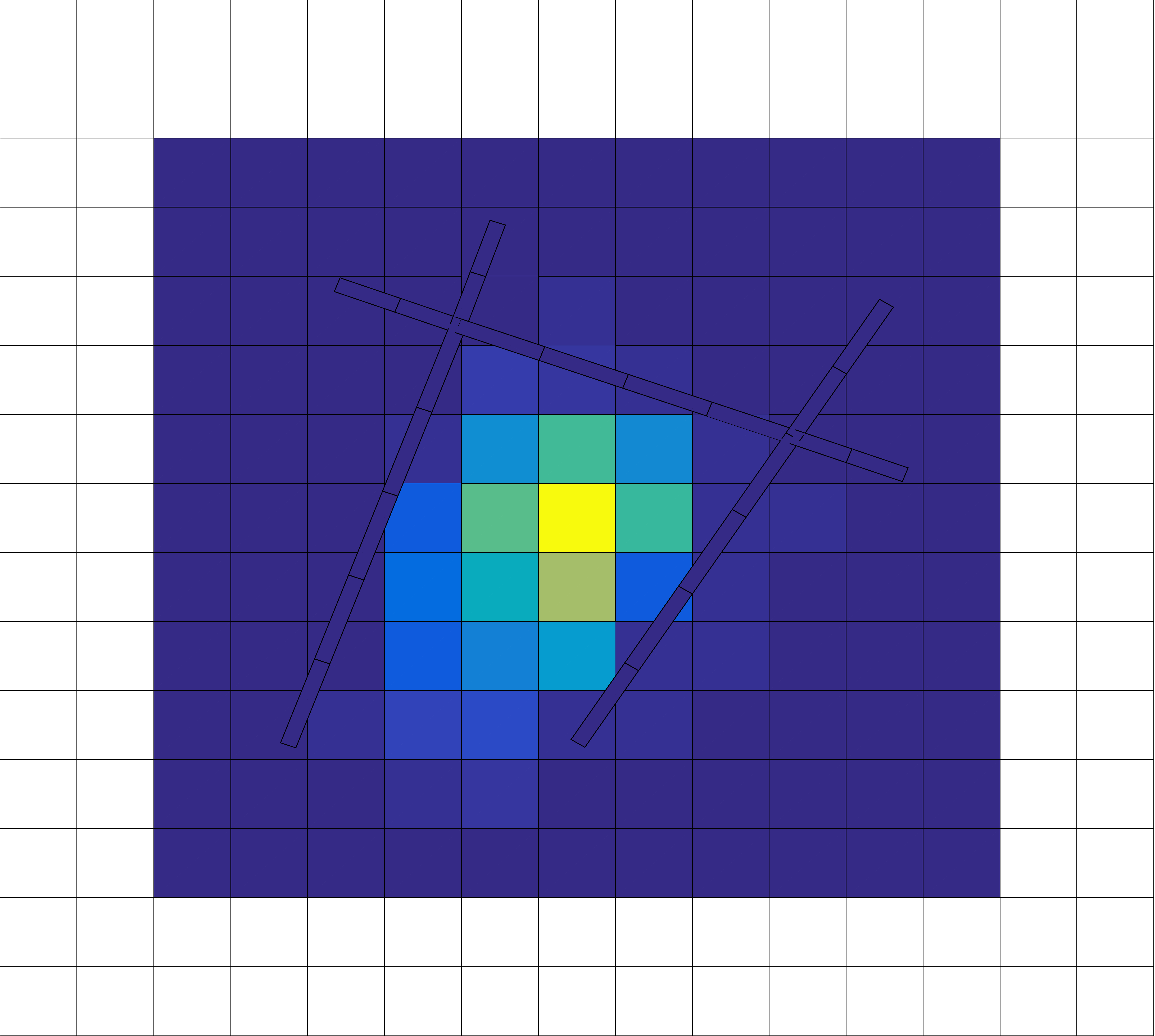}}\hspace{0.3cm}\includegraphics[height=0.3\textwidth]{Fig/method/colorbar.pdf}\label{fig:fracAMS_rock}\hspace{-0.8cm}}\\%
\subfigure[$\Phi^{mf}$ vertices]{\fbox{\includegraphics[width=0.3\textwidth,height=0.3\textwidth]{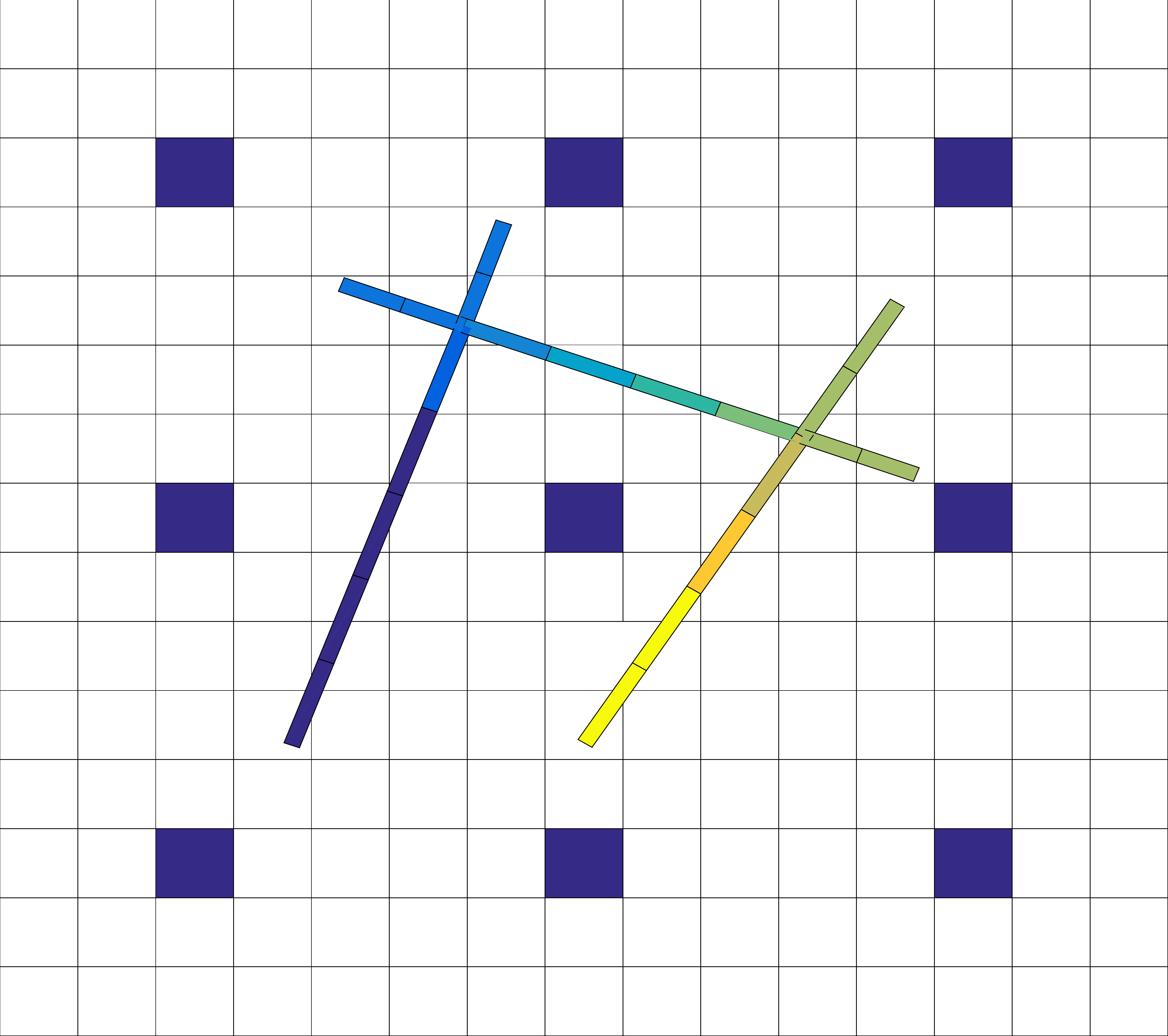}\label{fig:fracAMS_frac2_vertices}}}\hspace{0.3cm}%
\subfigure[$\Phi^{mf}$ edges]{\fbox{\includegraphics[width=0.3\textwidth,height=0.3\textwidth]{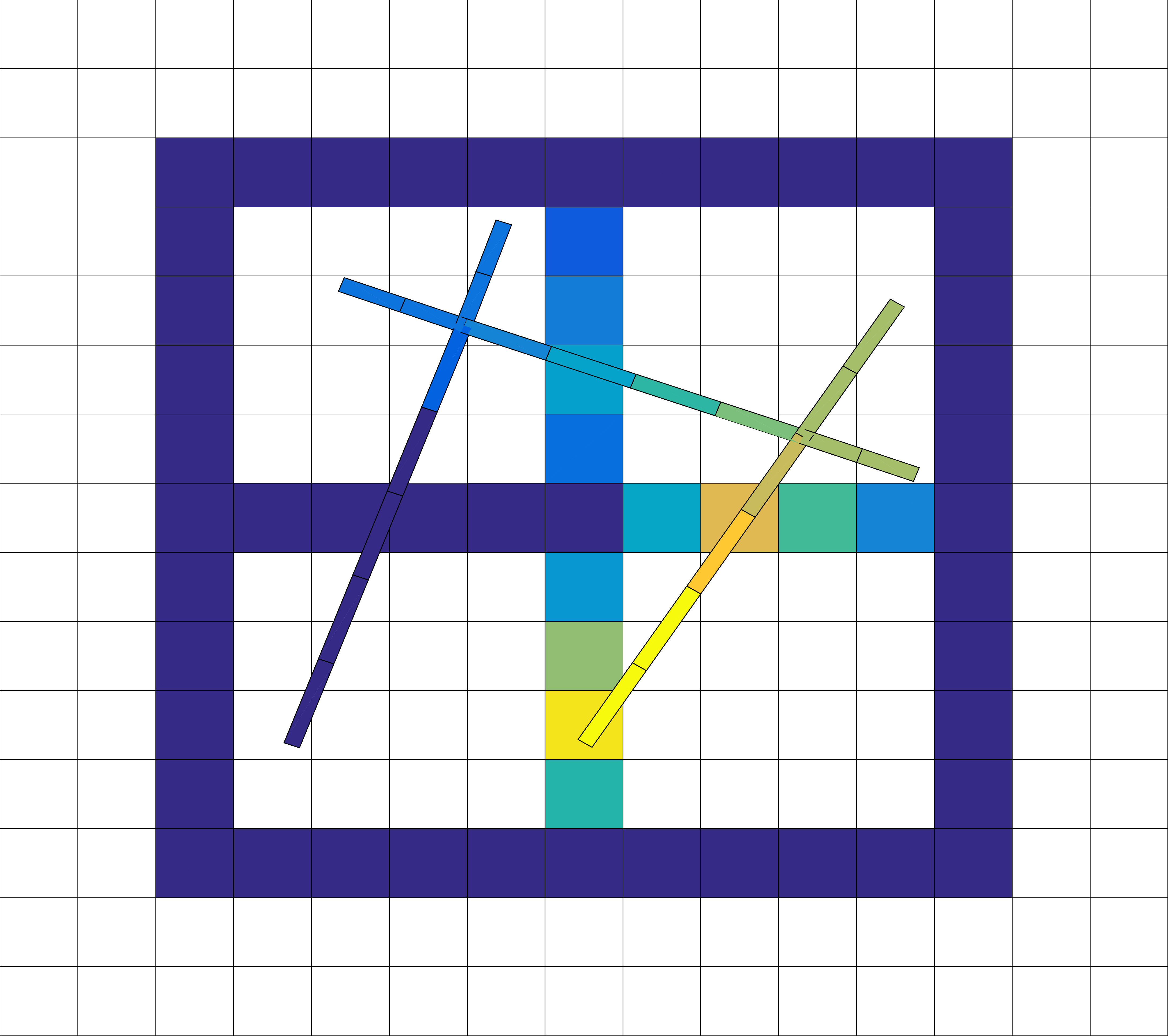}\label{fig:fracAMS_frac2_edges}}}\hspace{0.3cm}%
\subfigure[$\Phi^{mf}$ faces]{\fbox{\includegraphics[width=0.3\textwidth,height=0.3\textwidth]{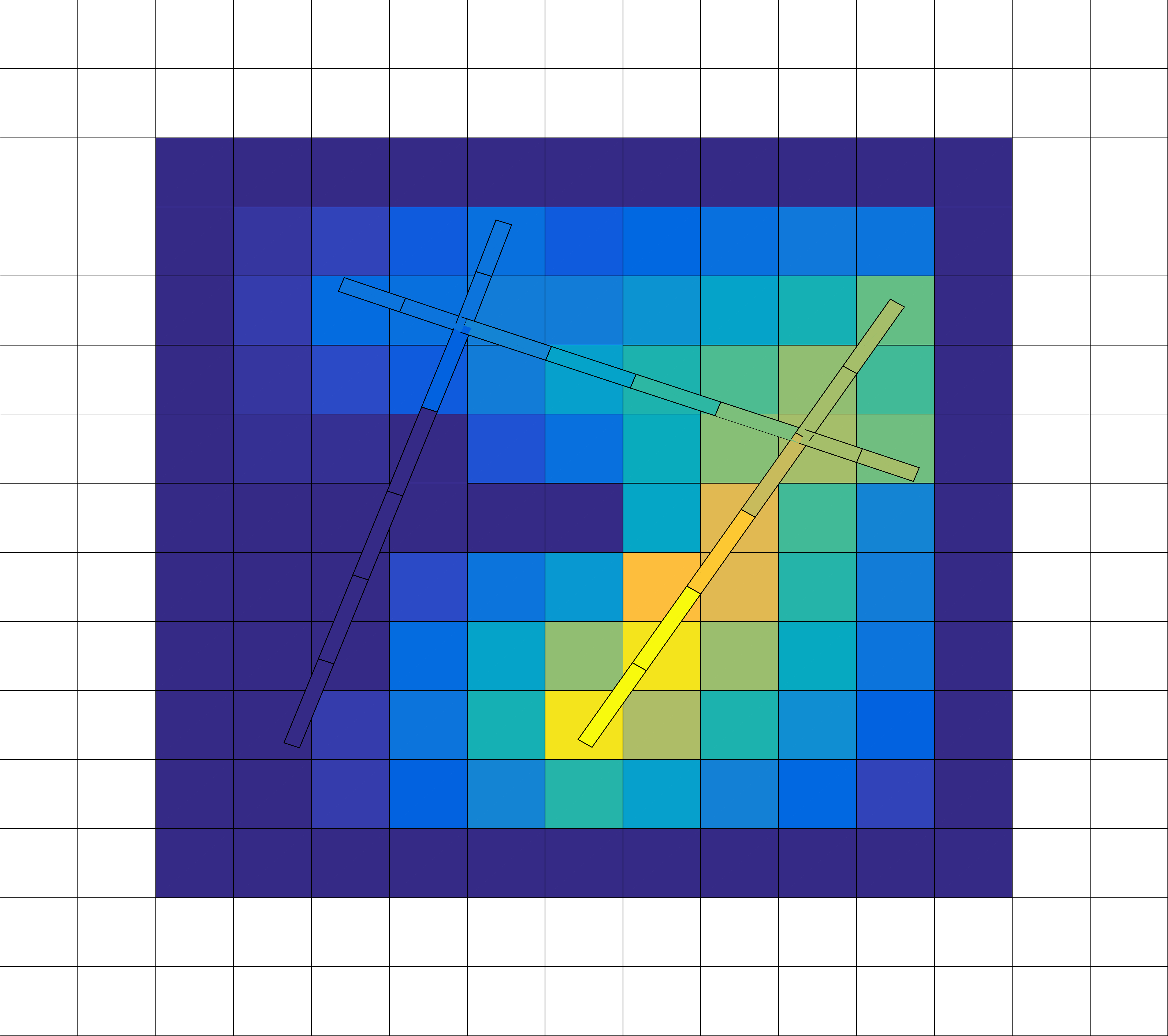}}\hspace{0.3cm}\includegraphics[height=0.3\textwidth]{Fig/method/colorbar.pdf}\label{fig:fracAMS_frac}\hspace{-0.8cm}}%
\caption{Step-by-step construction of a Frac-AMS matrix basis function (top) and fracture function (bottom). In this strategy, $\Phi^{fm}= 0$  and used as Dirichlet boundary condition while solving $\Phi^{mm}$. Then, $\Phi^{ff}$ is obtained similar as in Decoupled-AMS (Figs.~\ref{fig:decoupledAMS_frac2_vertices}-\ref{fig:decoupledAMS_frac2}) and used as Dirichlet boundary condition for $\Phi^{mf}$.}%
\label{fig:fracAMS}%
\end{figure}%

\item \textbf{Rock-AMS:} First, $\Phi^{mm}$ is computed with no-flow to the fractures, as with Decoupled-AMS (Figs.~\ref{fig:decoupledAMS_rock_vertices}-\ref{fig:decoupledAMS_rock}), i.e.,
\begin{align}
\xi(\Phi^{mm}) = 0.
\end{align} 

\noindent Then, the values are fixed and used as Dirichlet boundaries while solving $\Phi^{fm}$, for which the fracture-matrix connections are taken into account in Eq.~\ref{skeleton-basis}, i.e.,
\begin{align}
\xi(\Phi^{fm}) = \Phi^{fm} - \Phi^{mm}.
\end{align}

\noindent For the fracture functions, $\Phi^{mf} = 0$ which is used as Dirichlet condition to compute $\Phi^{ff}$, i.e.,
\begin{align}
\xi(\Phi^{ff}) = \Phi^{ff}.
\end{align}

\noindent Appendix~\ref{sec:app-rockAMS} presents the algebraic procedure corresponding to this coupling strategy. Note from Fig.~\ref{fig:rockAMS} that, in Rock-AMS, the matrix basis functions have non-zero values in the fractures, while the opposite does not hold.

\begin{figure}[htb!]%
\centering%
\setlength{\fboxsep}{0pt}%
\setlength{\fboxrule}{1pt}%
\setlength{\unitlength}{1cm}%
\begin{picture}(15,0.5)%
  \put(5.6,0.2){\large Rock-AMS}%
\end{picture}\\%
\subfigure[$\Phi^{fm}$ vertices]{\fbox{\includegraphics[width=0.3\textwidth,height=0.3\textwidth]{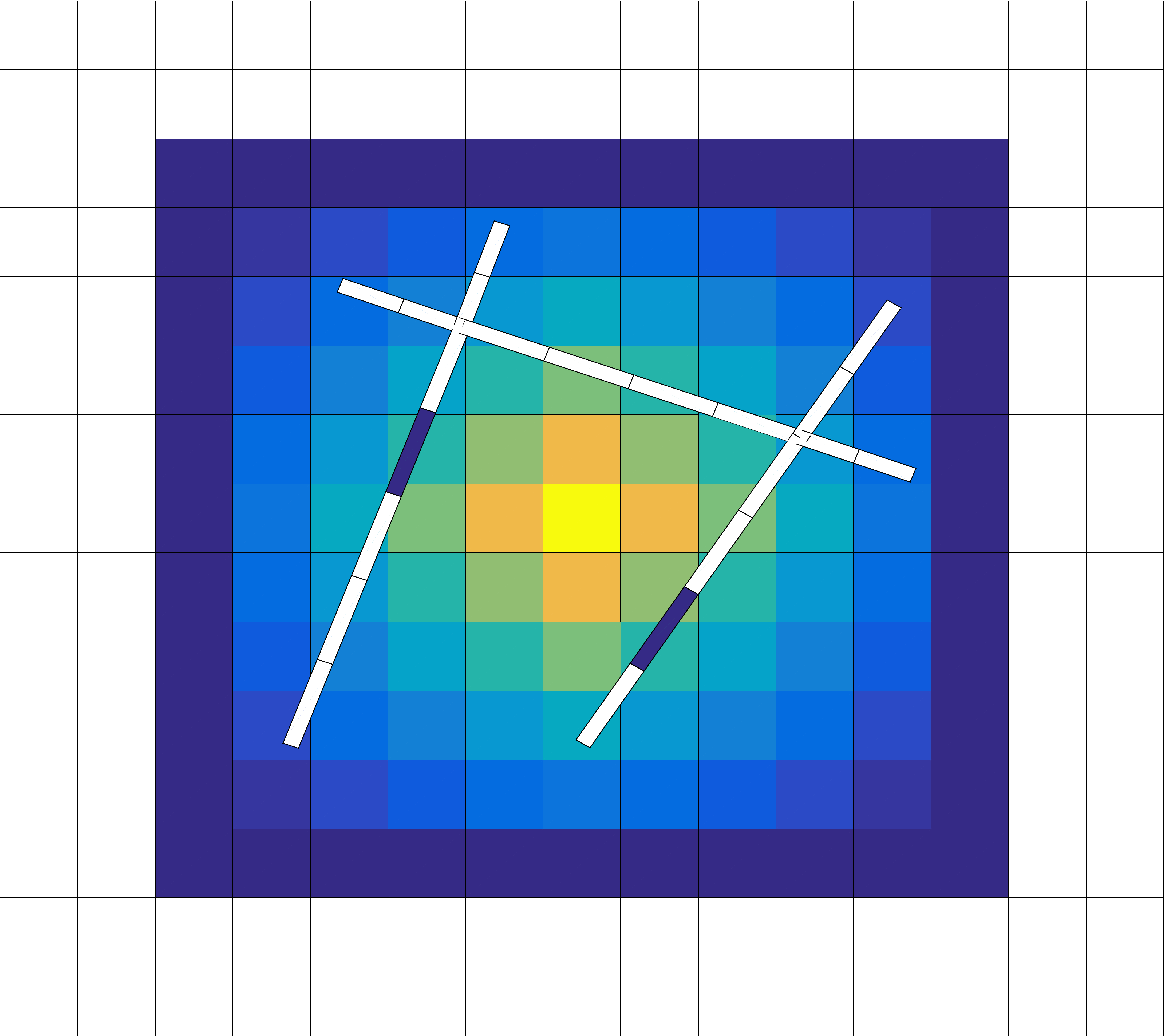}\label{fig:rockAMS_rock_faces}}}\hspace{1cm}%
\subfigure[$\Phi^{fm}$ edges]{\fbox{\includegraphics[width=0.3\textwidth,height=0.3\textwidth]{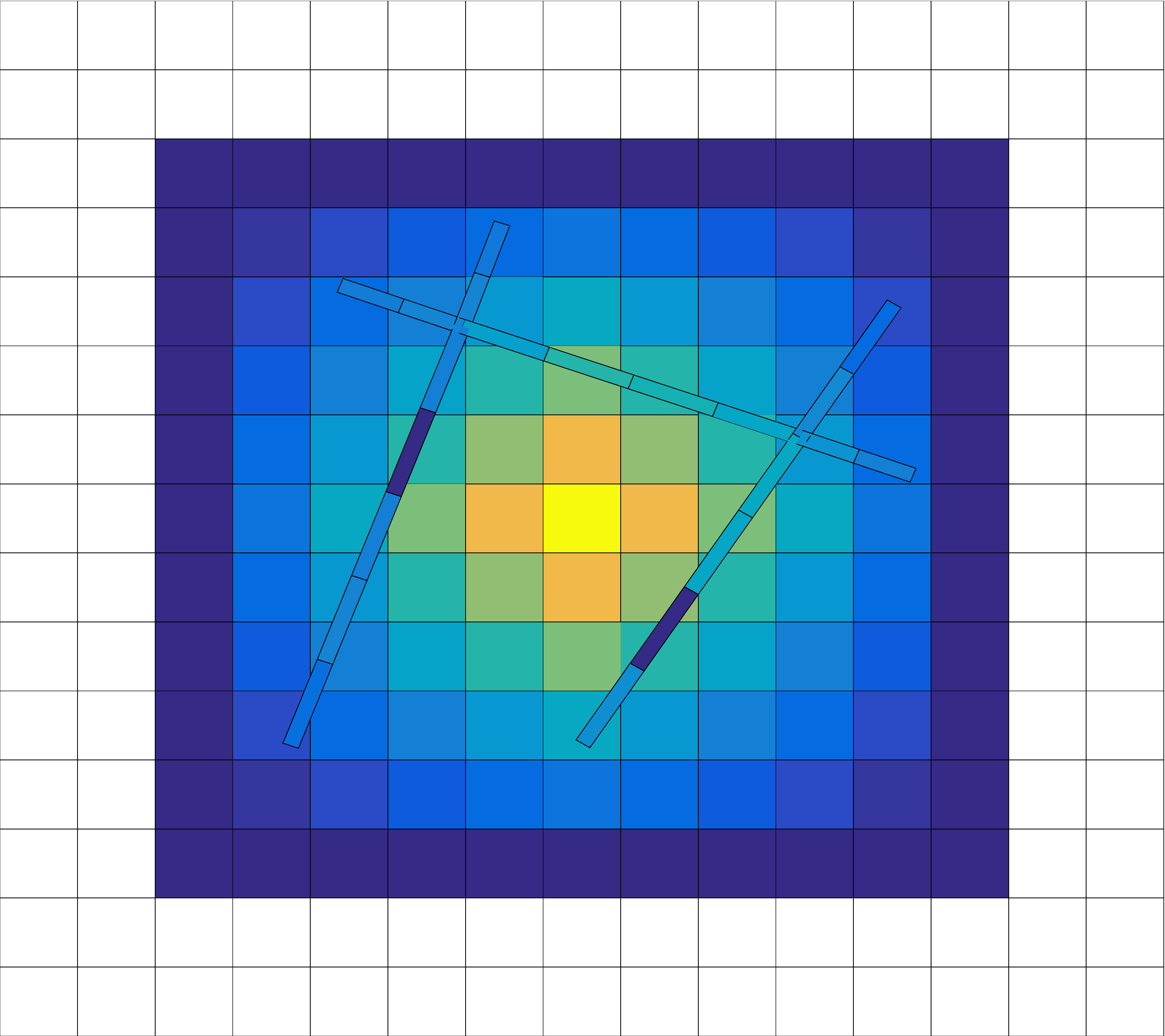}}\hspace{0.3cm}\includegraphics[height=0.3\textwidth]{Fig/method/colorbar.pdf}\label{fig:rockAMS_rock}\hspace{-0.8cm}}\\%
\subfigure[$\Phi^{ff}$ vertices]{\fbox{\includegraphics[width=0.3\textwidth,height=0.3\textwidth]{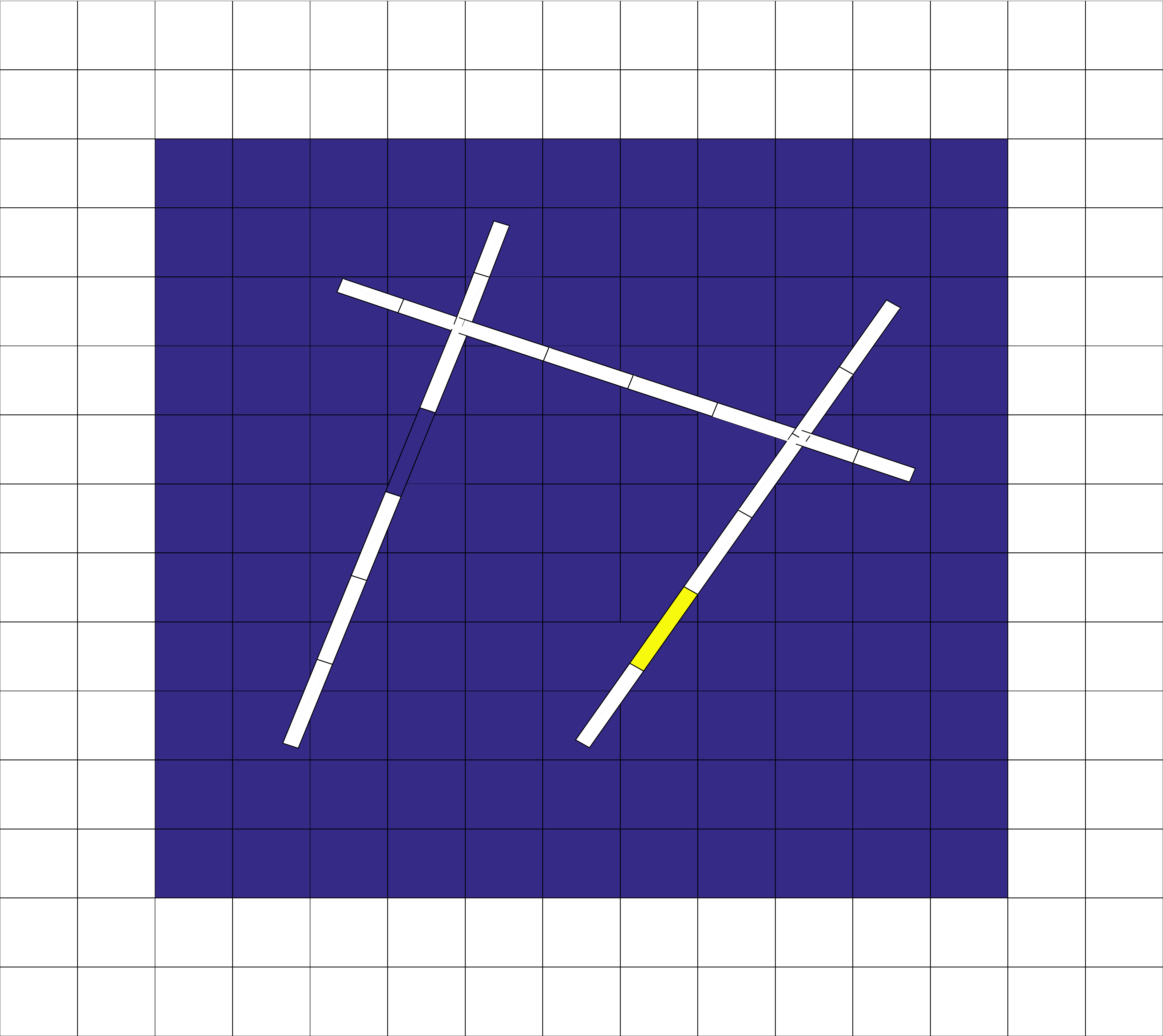}\label{fig:rockAMS_frac2_vertices}}}\hspace{1cm}%
\subfigure[$\Phi^{ff}$ edges]{\fbox{\includegraphics[width=0.3\textwidth,height=0.3\textwidth]{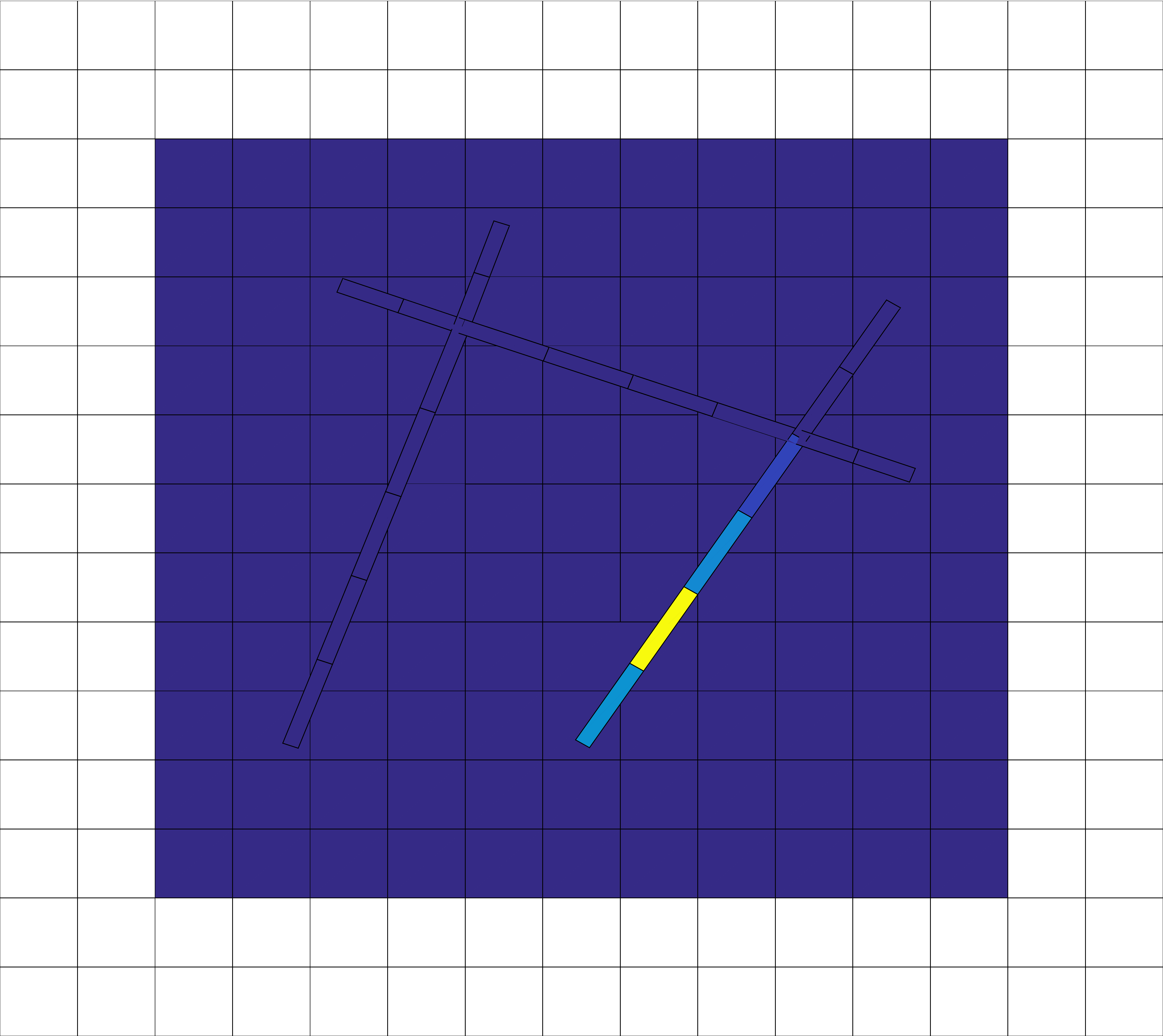}}\hspace{0.3cm}\includegraphics[height=0.3\textwidth]{Fig/method/colorbar.pdf}\label{fig:rockAMS_frac2}\hspace{-0.8cm}}%
\caption{Step-by-step construction of a Rock-AMS matrix basis function (top) and fracture function (bottom). In this strategy, the $\Phi^{mm}$ are obtained similar to Decoupled-AMS (Figs.~\ref{fig:decoupledAMS_rock_vertices}-\ref{fig:decoupledAMS_rock}) and used as Dirichlet boundary conditions while solving $\Phi^{fm}$. Then, $\Phi^{mf} = 0$ and used as Dirichlet boundary conditions for $\Phi^{ff}$.}%
\label{fig:rockAMS}%
\end{figure}%

\item \textbf{Coupled-AMS:} In order to preserve the two-way coupling between fractures and matrix, adjacent dual blocks of the same type are merged (e.g. fracture edges with the matrix edges they perforate), as shown in Fig.~\ref{fig:mergedDuals}. To clarify, two blocks are considered adjacent if there is a non-zero transmissibility between a cell from one of them and a cell from the other.%

\begin{figure}[htb!]%
\centerline{%
\subfigure[Merged edges in 2D]{\raisebox{0.7cm}{\includegraphics[width=0.3\textwidth,height=0.3\textwidth]{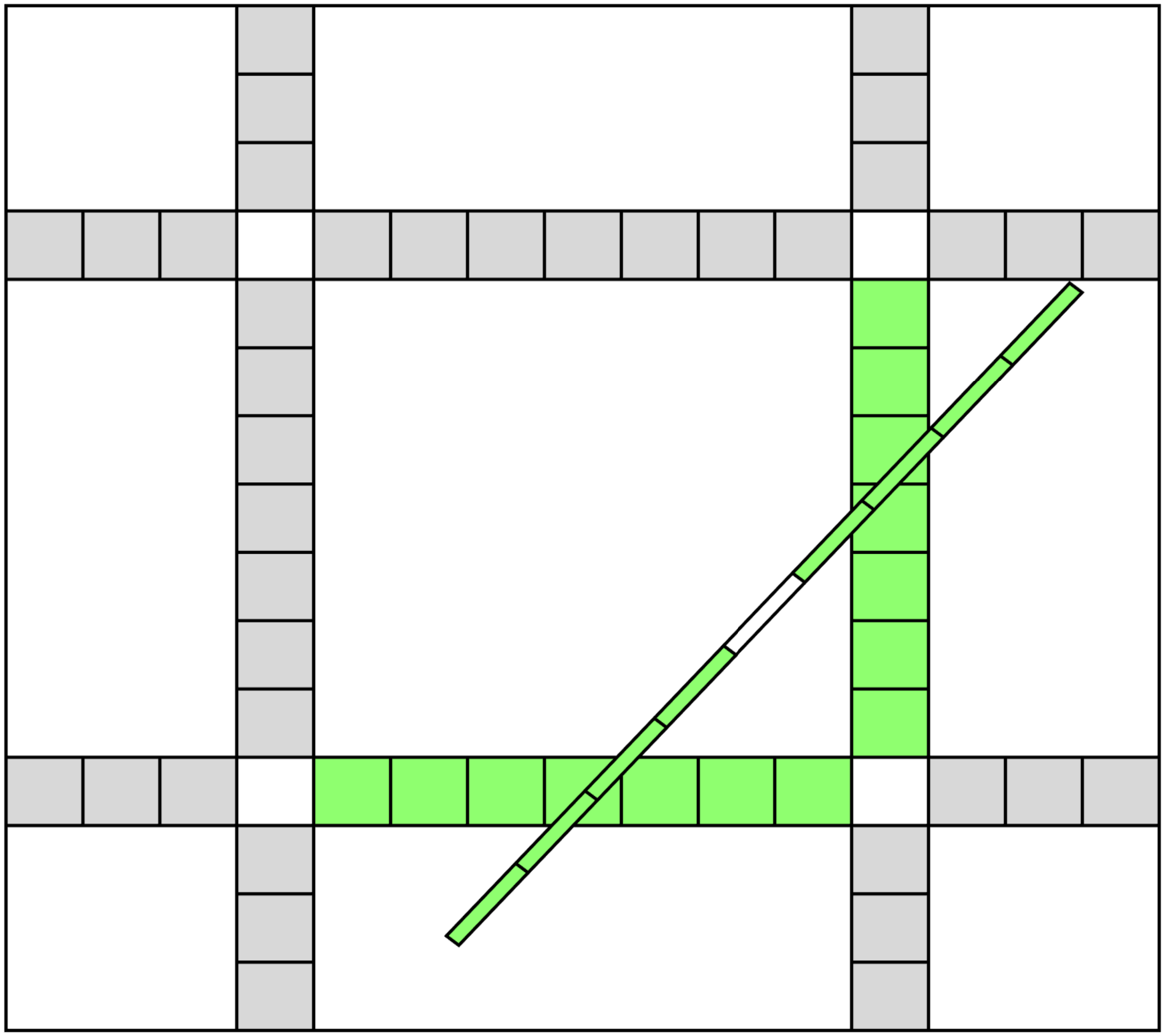}\label{fig:merged_edges_2d}}}\hspace{0.1cm}%
\subfigure[Merged edges in 3D]{\includegraphics[width=0.38\textwidth,height=0.4\textwidth]{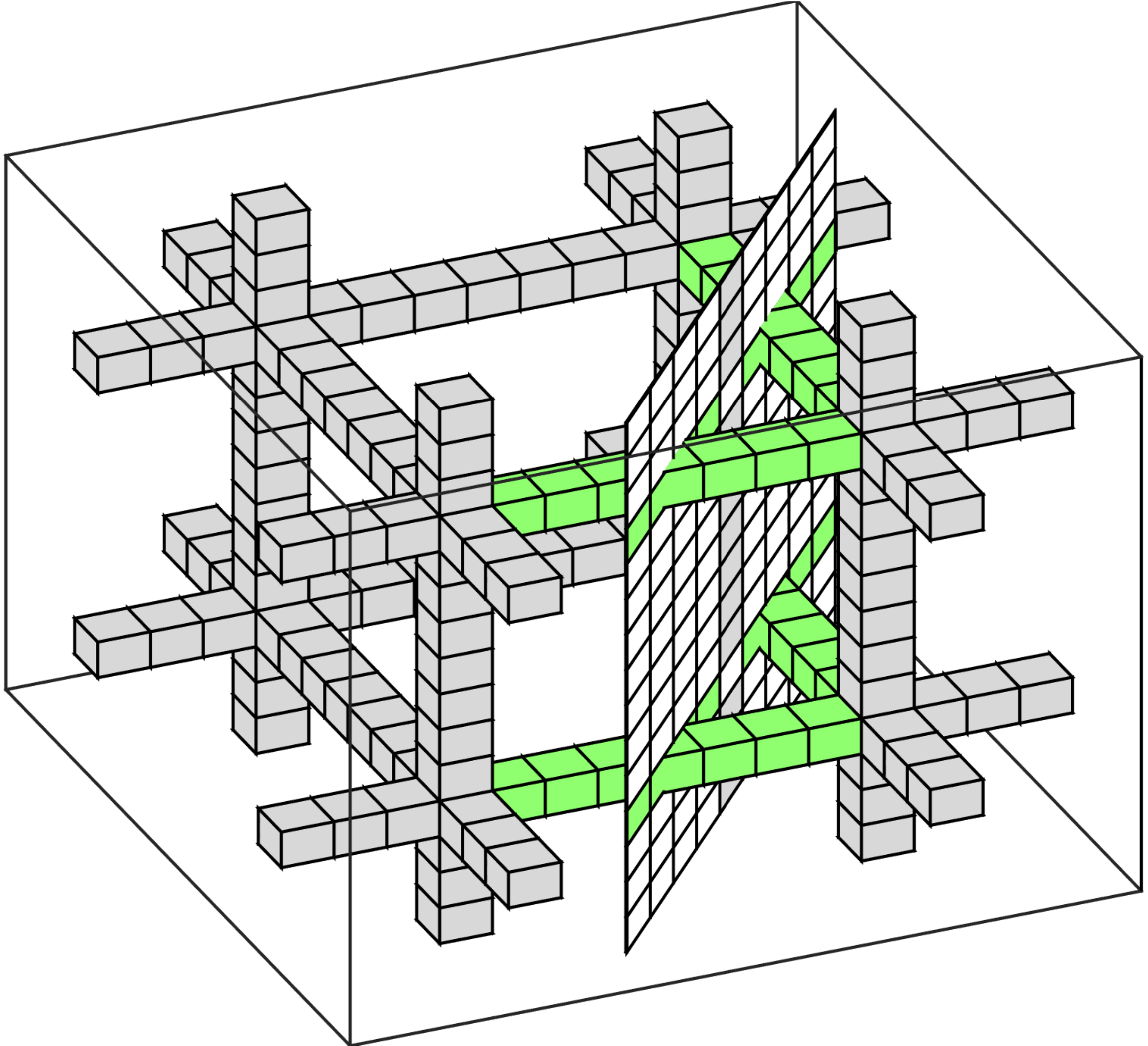}\label{fig:merged_edges_3d}}\hspace{0.1cm}%
\subfigure[Merged faces in 3D]{\includegraphics[width=0.38\textwidth,height=0.4\textwidth]{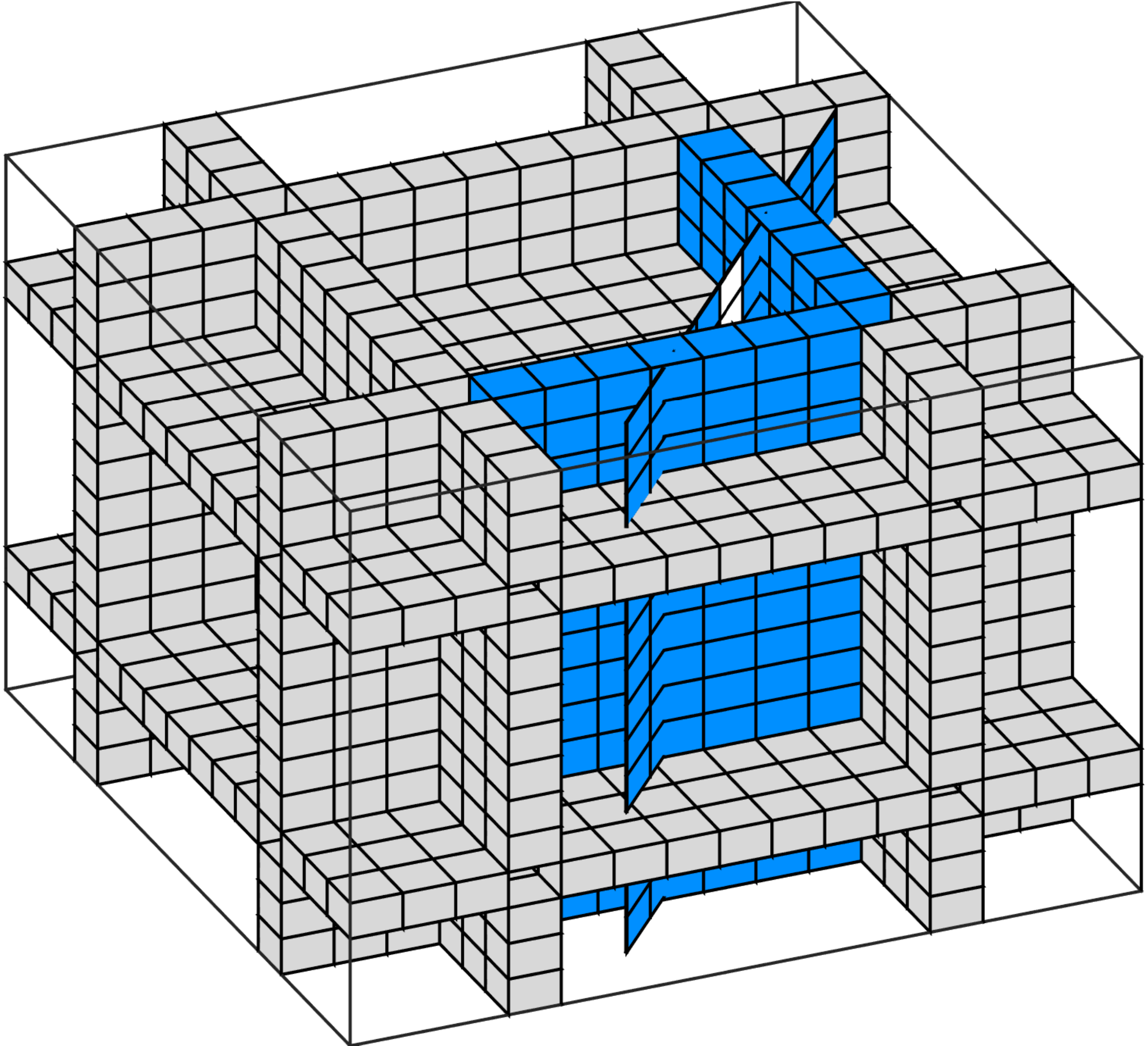}\label{fig:merged_faces_3d}}}%
\caption{The merging of adjacent dual blocks of the same type in 2D (a) and 3D (b and c), in order to preserve the full coupling between fractures and matrix in the Coupled-AMS basis functions.}%
\label{fig:mergedDuals}%
\end{figure}%

\begin{figure}[htb!]%
\centering%
\setlength{\fboxsep}{0pt}%
\setlength{\fboxrule}{1pt}%
\setlength{\unitlength}{1cm}%
\begin{picture}(15,0.7)%
  \put(5.4,0.2){\large Coupled-AMS}%
\end{picture}\\%
\subfigure[$\Phi^{mm}$ \& $\Phi^{fm}$ vertices]{\fbox{\includegraphics[width=0.3\textwidth,height=0.3\textwidth]{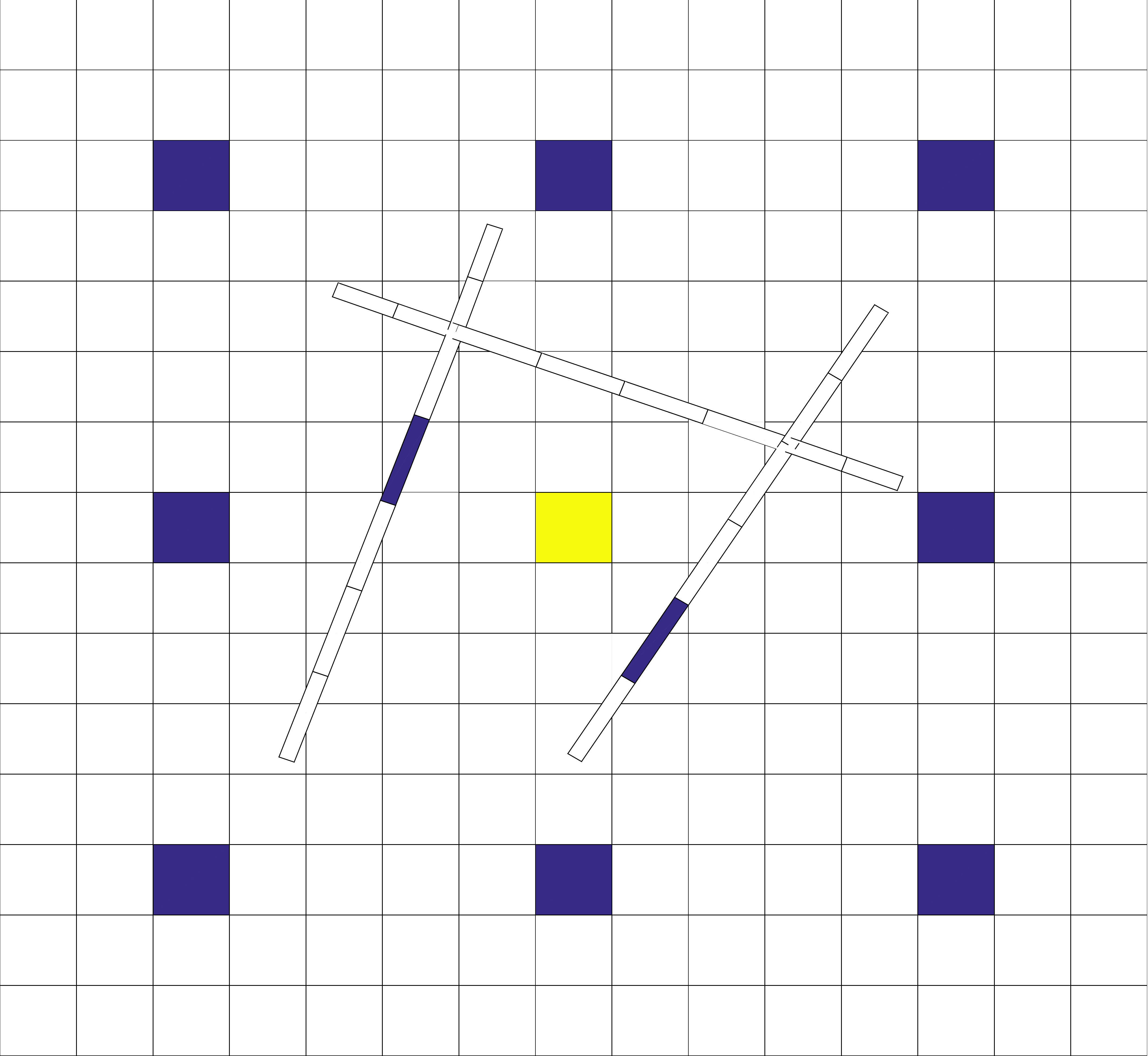}\label{fig:coupledAMS_rock_vertices}}}\hspace{0.3cm}%
\subfigure[$\Phi^{mm}$ \& $\Phi^{fm}$ edges]{\fbox{\includegraphics[width=0.3\textwidth,height=0.3\textwidth]{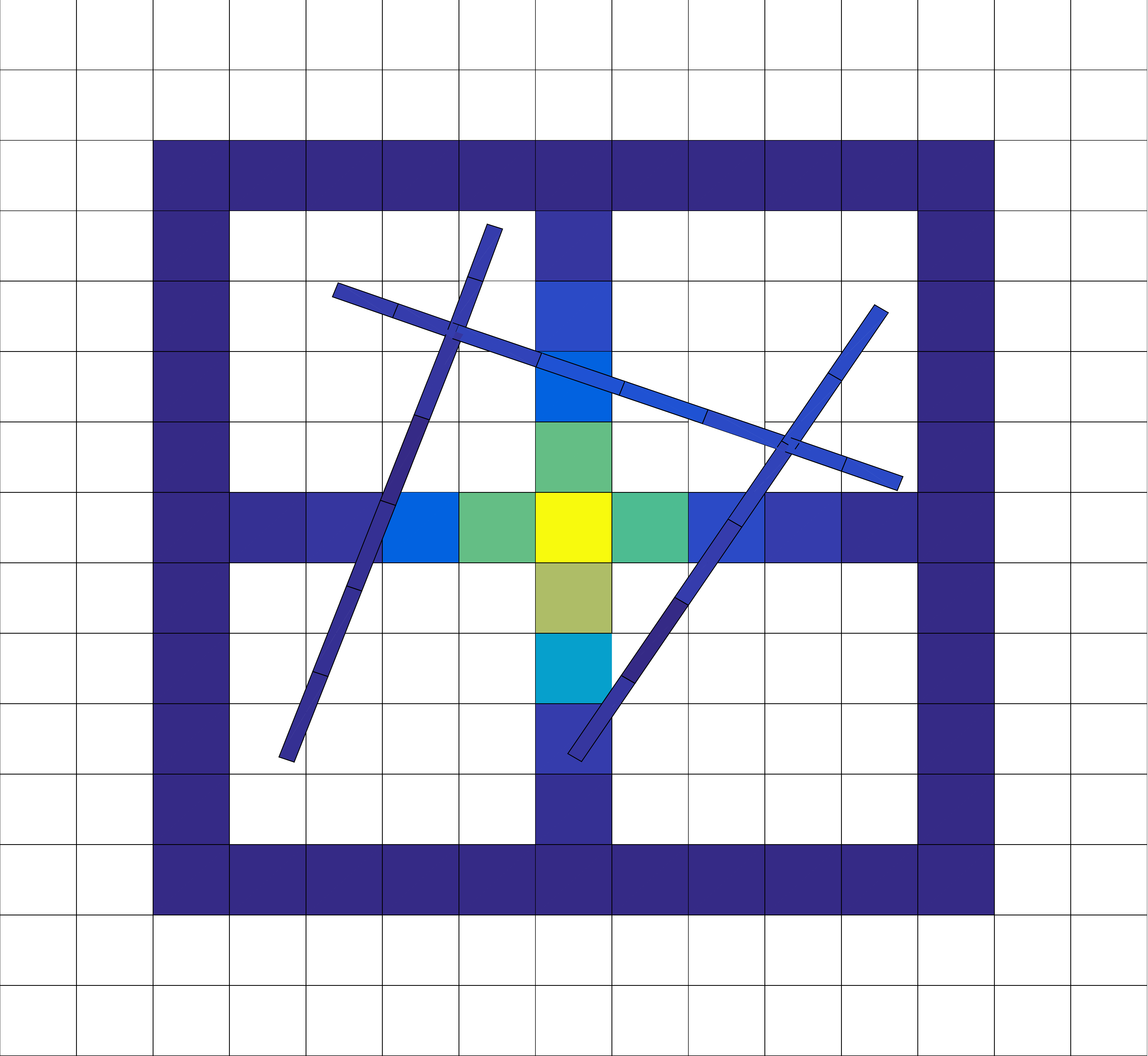}\label{fig:coupledAMS_rock_edges}}}\hspace{0.3cm}%
\subfigure[$\Phi^{mm}$ \& $\Phi^{fm}$ faces]{\fbox{\includegraphics[width=0.3\textwidth,height=0.3\textwidth]{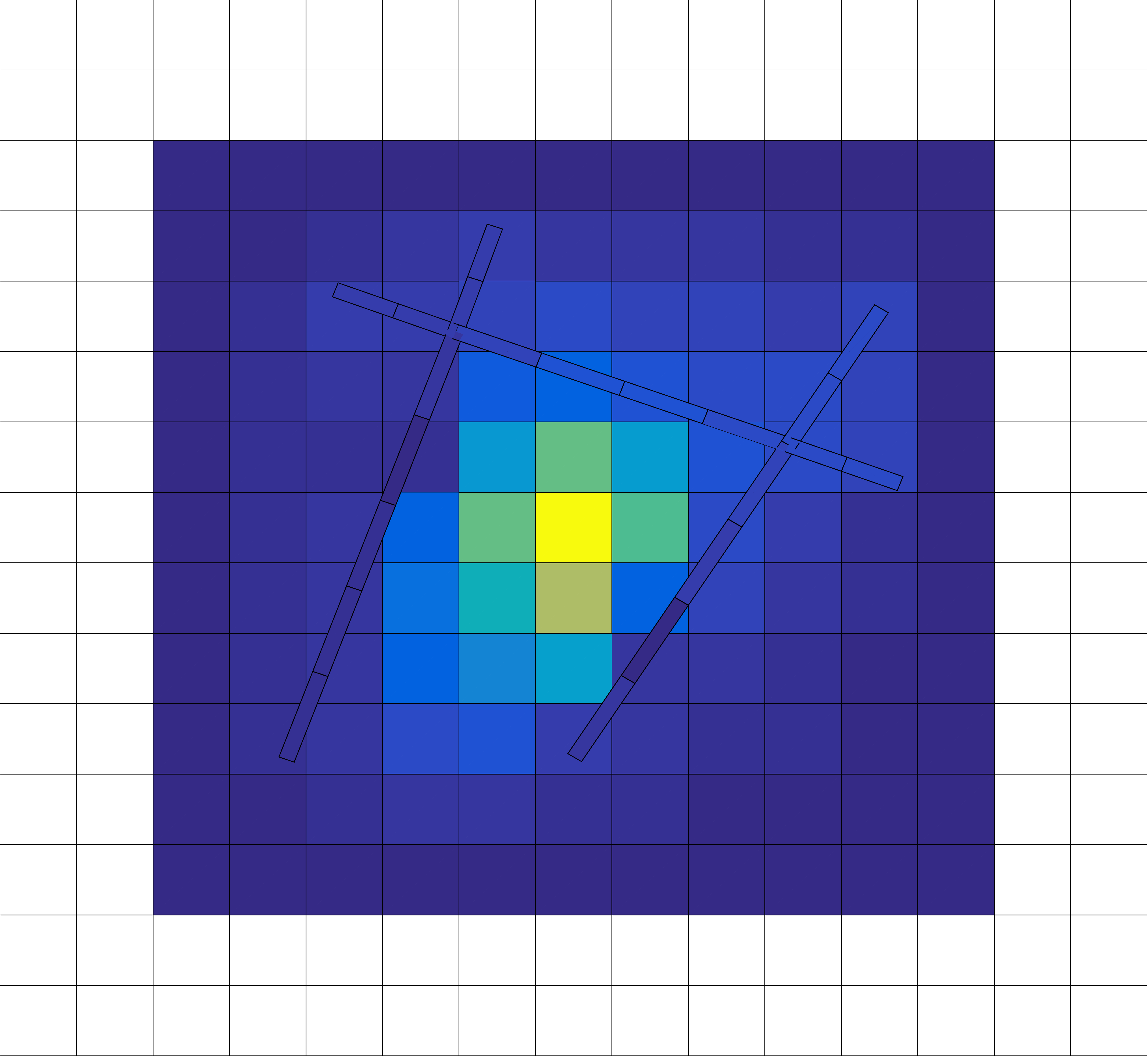}}\hspace{0.3cm}\includegraphics[height=0.3\textwidth]{Fig/method/colorbar.pdf}\label{fig:coupledAMS_rock}\hspace{-0.8cm}}\\%
\subfigure[$\Phi^{ff}$ \& $\Phi^{mf}$ vertices]{\fbox{\includegraphics[width=0.3\textwidth,height=0.3\textwidth]{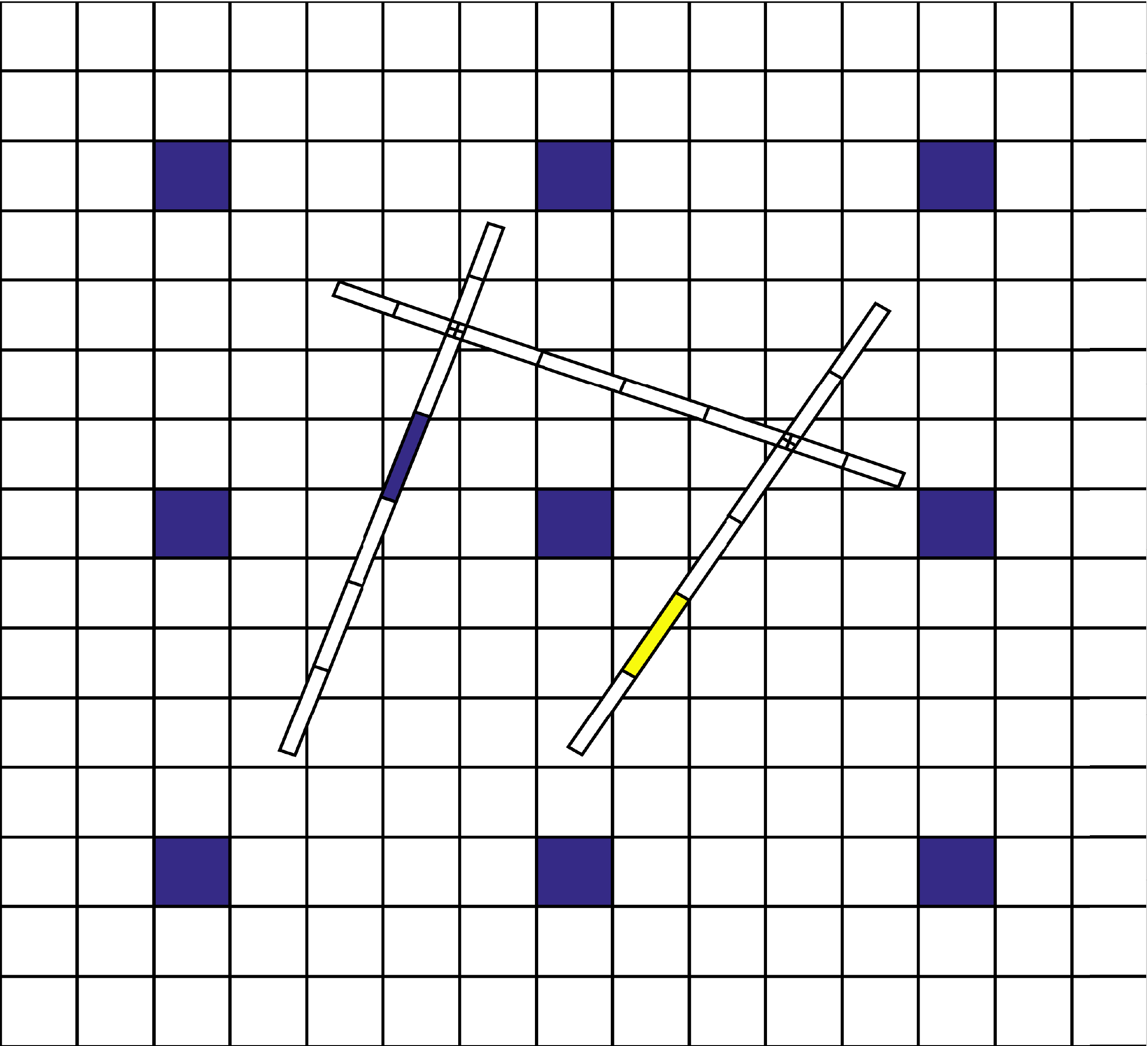}\label{fig:coupledAMS_frac2_vertices}}}\hspace{0.3cm}%
\subfigure[$\Phi^{ff}$ \& $\Phi^{mf}$ edges]{\fbox{\includegraphics[width=0.3\textwidth,height=0.3\textwidth]{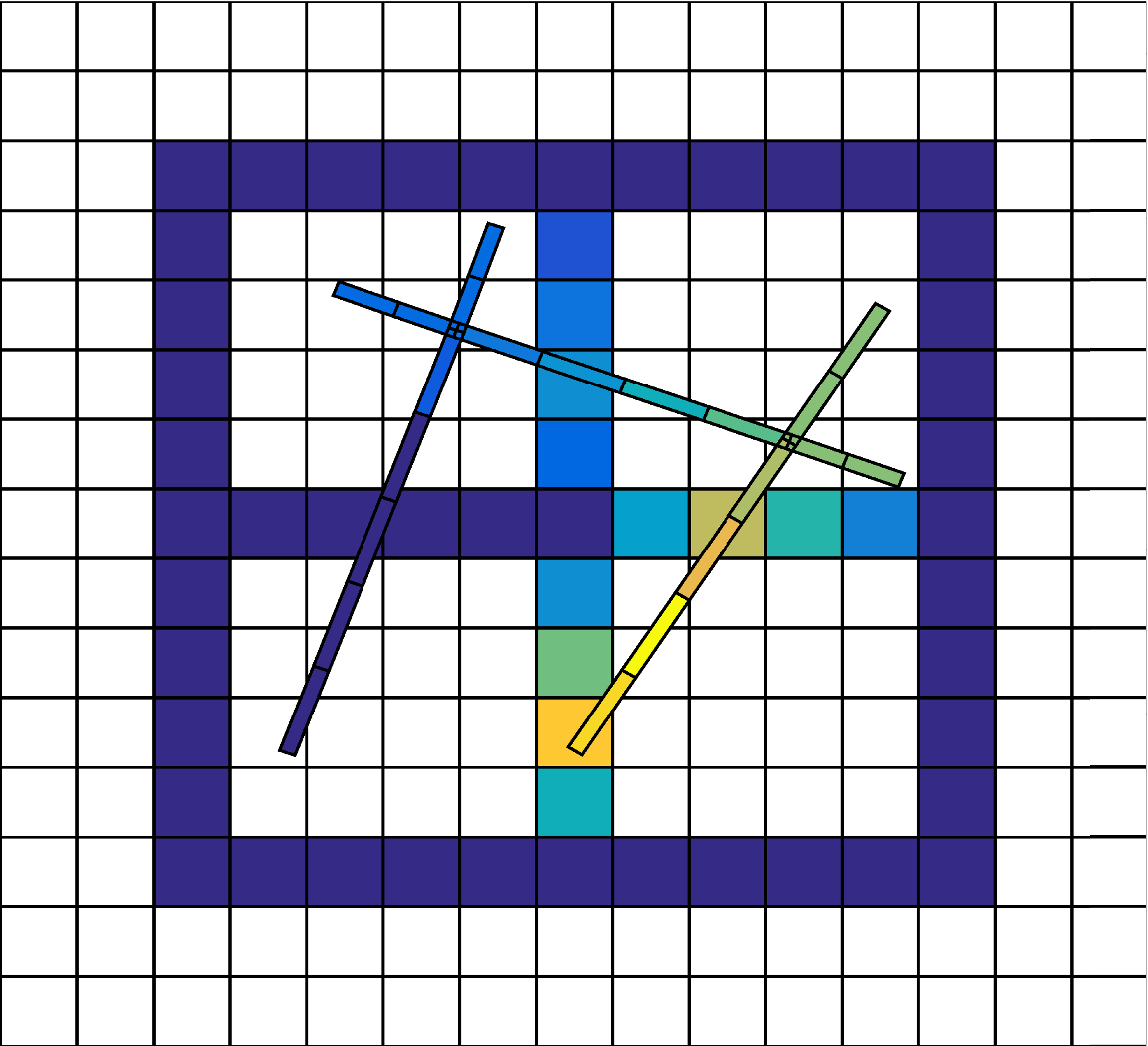}\label{fig:coupledAMS_frac2_edges}}}\hspace{0.3cm}%
\subfigure[$\Phi^{ff}$ \& $\Phi^{mf}$ faces]{\fbox{\includegraphics[width=0.3\textwidth,height=0.3\textwidth]{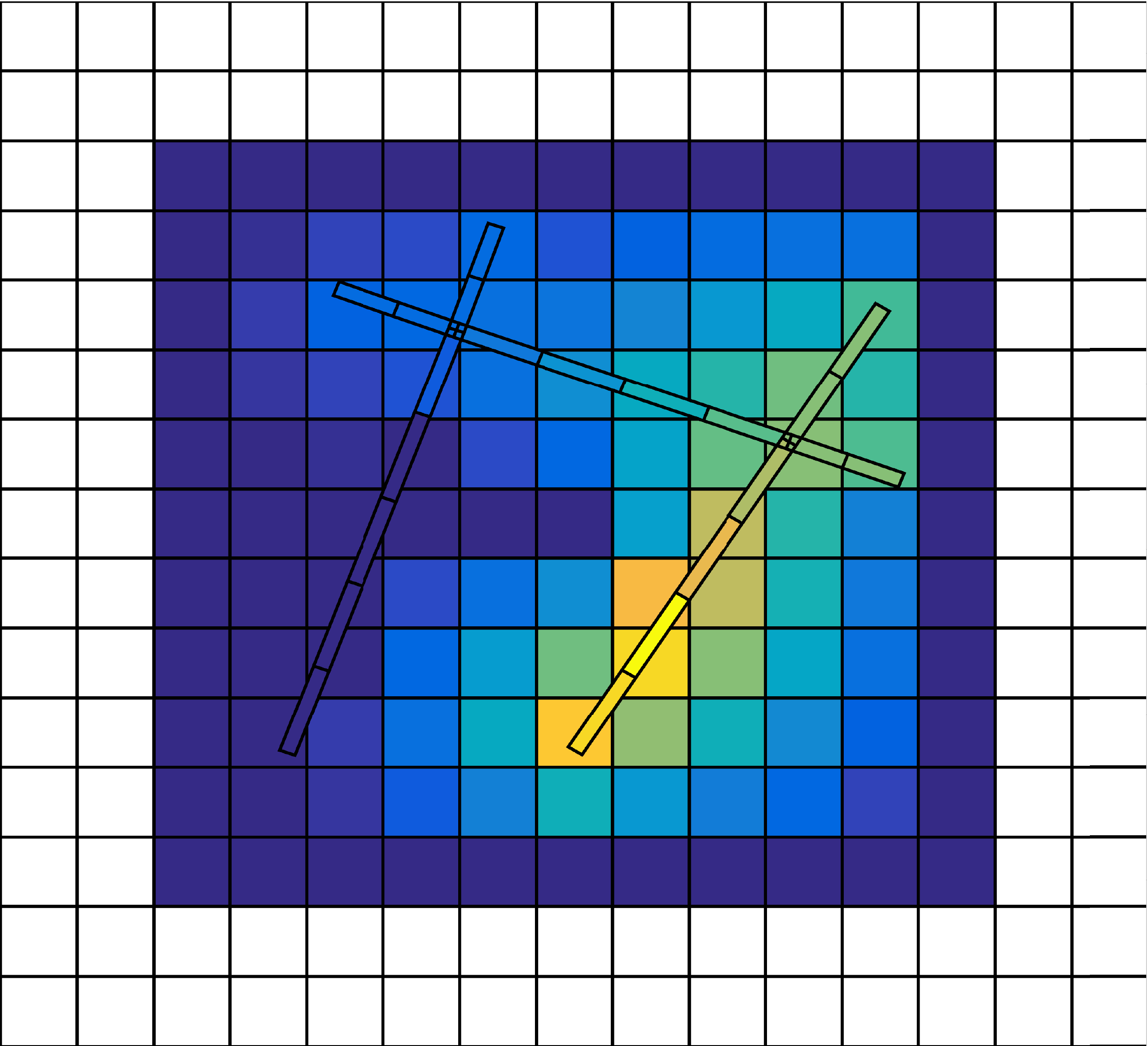}}\hspace{0.3cm}\includegraphics[height=0.3\textwidth]{Fig/method/colorbar.pdf}\label{fig:coupledAMS_frac2}\hspace{-0.8cm}}%
\caption{Step-by-step construction of a Coupled-AMS matrix basis function (top) and fracture function (bottom). In this strategy, both media preserve their connectivity towards each other.}%
\label{fig:coupledAMS}%
\end{figure}%

\noindent On the resulting dual-coarse grid, computation of the basis functions follow the usual wirebasket hierarchy, with full consideration of the coupling, i.e., using 
\begin{align}
\xi(\Phi^{mm}) = \Phi^{mm} - \Phi^{fm}, & & \xi(\Phi^{fm}) = \Phi^{fm} - \Phi^{mm}, \nonumber \\ 
\xi(\Phi^{ff}) = \Phi^{ff} - \Phi^{mf} & \text{\hspace{9mm}and} & \xi(\Phi^{mf}) = \Phi^{mf} - \Phi^{ff}, 
\end{align}
for matrix and fracture cells, respectively (see Appendix~\ref{sec:app-coupledAMS} for the detailed algebraic procedure). Figure~\ref{fig:coupledAMS} presents an illustration of the Coupled-AMS basis functions. Note that fracture functions have non-zero values in the matrix, and vice versa, the matrix basis functions also have non-zero values inside fractures.

\end{enumerate}

It is worth mentioning that, by construction, all four coupling strategies result in basis functions which form a partition of unity.

The consideration of wells is similar to what was described in the literature for 2D problems \cite{Patrick-well2}, but extended here to 3D problems and integrated within the F-AMS framework. Each well is represented as a single coarse node and a well function is computed locally. The resulting values are assembled in the designated column of the prolongation operator. Note that, even for well perforations in the matrix, the corresponding well functions can have non-zero values inside fractures as well if either Coupled-AMS or Rock-AMS are employed.

As mentioned before, the basis functions have local support. However, in the Coupled-AMS case, the merging of dual blocks (Fig.~\ref{fig:mergedDuals}) can increase this support substantially in cases with a high density of interconnected (or long) fractures. This, in turn, can lead to a dense prolongation operator, with a severe impact on computational performance. To overcome this, one can impose a limiting criterion on the merging of the dual blocks, for example a maximum number of fine grid cells. Alternatively, or in combination with the previous method, one can discard the non-zeros from $\bm{\mathcal{P}}$ which lie below a specified threshold, and rescale the rows accordingly to preserve the partition of unity. The latter choice is studied in detail in Section~\ref{sec:results}.

It is important to note that F-AMS basis functions are computed at the beginning of a time-dependent simulation, and adaptively updated only if the fine-scale properties change beyond a threshold value \cite{hadi-aimsfv-jcp}. Next, the F-AMS solution algorithm is described.

\subsection{The F-AMS solution algorithm}
\label{sec:algorithm}

In addition to the prolongation, in order to obtain the coarse-scale pressure system, the restriction operator $\bm{\mathcal{R}}$ (i.e., map from fine to coarse scale) is now defined. As previously described, F-AMS identifies three distinct types of features in the domain, i.e., matrix, fractures, and wells, therefore, $\bm{\mathcal{R}}$ can be defined in a much more general form than in the previous studies \cite{Zhou-tams,yixuan-ams,Tene-cams}.

The first option is to apply a FV-based restriction to all domains, i.e. the MSFV restriction operator, $\bm{\mathcal{R}}^{FV}$, where the entry at row $i$ and column $j$ is $1$ only if the fine-scale cell $j$ (either from the matrix, fractures, or wells) belongs to primal-coarse block $i$. MSFV ensures mass conservation, at the coarse-scale, thus allowing the reconstruction of a fine-scale conservative flux field. However, in previous works \cite{yixuan-ams,yixuan-monotonemsfv}, it has been found sensitive to the heterogeneity contrast in the domain. Alternatively, one can construct a Galerkin-FE-based restriction operator, as $\bm{\mathcal{R}}^{FE} = \bm{\mathcal{P}}^T$, traditionally called MSFE, which leads to a symmetric-positive-definite (SPD) coarse linear system, if the fine-scale system \eqref{dis-mf} is also SPD. Finally, one can consider a third option, where some of the features (e.g., matrix) are restricted according to FE, while, for the rest, FV is used. This will be referred to as the MSMIX restriction operator. Note that, for multiphase flow problems, if the pressure system is not solved to machine accuracy, a final iteration with $\bm{\mathcal{R}}^{FV}$ needs to be employed, followed by a conservative flux reconstruction stage, in order to facilitate the solution of the transport equations \cite{hadi-comp-jcp,zhou-trans}.

Although F-AMS can be used as a single-sweep multiscale solver, where the approximate solution, $p'$, is used with no iterations, previous studies have shown that an iterative procedure is needed for highly-heterogeneous reservoirs \cite{imsfv-jcp}. To this end, one needs to pair the F-AMS multiscale step with a fine-scale smoother, which ensures error reduction to any desired level. The F-AMS algorithm can now be compiled, as shown in Table~\ref{fams-algorithm}.

\begingroup
\renewcommand{\arraystretch}{1.2}
\begin{table}[htb!]
\caption{F-AMS solution algorithm}
\begin{tabular}{|l|}
\hline%
Repeat the following steps until convergence to the desired accuracy is reached:\\
  \hspace{1mm} 1. \textbf{Initialize:} $p^\nu \leftarrow p^{\nu+1}$ \\
  \hspace{1mm} 2. \textbf{Update linear system entries:} $\bm A^{\nu}$ and $q^\nu$ in \eqref{dis-mf}\\
  \hspace{1mm} 3. \textbf{Update residual:} $r^\nu = q^\nu - \bm A^\nu p^\nu$\\
  \hspace{1mm} 4. \textbf{Compute (or adaptively update) $\bm{\mathcal{P}}$}: follow either coupling strategy from\\ Subsection~\ref{sec:basis}. \\
  \hspace{1mm} 5. \textbf{Multiscale Stage:} $\delta p^ {\nu +1/2} = \bm{\mathcal{P}} (\bm{\mathcal{R}} \bm A^{\nu} \bm{\mathcal{P}})^{-1} \bm{\mathcal{R}} \ r^{\nu} $\\
  \hspace{1mm} 6. \textbf{Update residual} $r^{\nu+1/2} = r^{\nu} - \bm A^{\nu}\delta p^{\nu+1/2}$\\
  \hspace{1mm} 7. \textbf{Smoothing Stage:} $\delta p^ {\nu +2/2} = \bm M^{-1}_S \  r^{\nu+1/2}$\\
  \hspace{1mm} 8. \textbf{Update solution:} $p^{\nu+1} = (p^\nu + \delta p^{\nu +1/2} +  \delta p^{\nu +2/2})$\\
  \hline
\end{tabular}
\label{fams-algorithm}
\end{table}
\endgroup

The smoothing operator, $\bm M^{-1}_S$, approximates the inverse of the complete fine-scale linear operator, $\bm A^\nu$, via ILU(0) decomposition \cite{Saad}. Note that the contrast between the matrix and fracture transmissibility values is usually severe, leading to a high condition number in Eq.~\eqref{dis-mf}. In such cases, F-AMS can easily be extended to include another smoothing stage, which employs iterations on the sub-block systems corresponding to each media, i.e. $\bm A^{mm}$ and $\bm A^{ff}$ from \eqref{dis-mf}. A detailed study of the impact of such a smoothing stage is beyond the scope of this paper, and makes the object of future research.

Next, numerical results are presented in order to study the effect of each component on the performance of the F-AMS algorithm. Then, the scalability of F-AMS is studied, as linear solver, in a CPU benchmark, where the SAMG commercial solver \cite{SAMG} is used as reference.

\section{Numerical Results}
\label{sec:results}

The aim of this section is to investigate the performance of F-AMS while simulating flow through fractured porous media. First, a 2D reservoir with heterogeneous matrix rock and a relatively complex fracture network is considered. A distance-based graph algorithm is introduced, in order to automate the fracture coarsening process. The convergence behaviour of F-AMS is studied, considering each of the four different coupling strategies introduced in the previous section. Then, the same fracture network is extruded along the Z axis and embedded in a 3D heterogeneous domain, for which simulations are performed considering different coarsening strategies, fracture-matrix conductivity ratios, fracture densities and domain sizes. CPU times are measured in detail for both the setup and solution stages in all test cases, and compared to those obtained using the industrial-grade SAMG solver \cite{SAMG}. Finally, the same 3D reservoir is used to investigate the effect of heterogeneous fracture conductivities, spanning several orders of magnitude.

During the upcoming simulations, special attention is given to the conductivity contrast between the matrix and the fracture domains. The transmissibility ratio $T_{ratio}$ is introduced as
\begin{align}
T_{ratio} = \frac{\langle T \rangle_{frac}}{\langle T \rangle_{rock}},
\label{tratio}
\end{align}
i.e. the ratio between the average fracture $\langle T \rangle_{frac}$ and matrix $ \langle T \rangle_{rock}$ transmissibility values, respectively. 

It is important to note that, in all test cases, F-AMS employs a FE restriction operator. Furthermore, the coarse-scale linear system and the basis functions in each dual block are all solved using a direct solver, based on LU decomposition, from the PETSc package \cite{petsc}.

For some experiments, a detailed breakdown of the CPU time spent in each stage of the F-AMS algorithm is presented. In the legends of the corresponding bar plots, the ``Initialization'' refers to the time spent on allocation of memory for the various data structures, the ``Operators'' represents the computation of basis functions and construction of the restriction and prolongation matrices (Step 4 in Table~\ref{fams-algorithm}). Also, ``Fine linsys. constr.'' denotes computation of the transmissibility values and fine-scale linear system assembly. In addition, the matrix multiplications resulting in the coarse-scale system are labelled as ``Coarse linsys, constr.'', while ``Solution''  stands for the solution of the coarse system followed by the interpolation (Step 5 in Table~\ref{fams-algorithm}). Finally, ``Smoother'' accounts for Step 7.

\subsection{Distance-based fracture coarsening}

A fracture network can be represented as a graph, in which the fracture lines (plates) are the arcs, while the nodes are the locations at which fractures intersect. This leads to a quasi-unstructured grid, where the complexity mostly revolves around the representation of the intersections. In this work, each intersection is assigned a pressure value, which is explicitly represented in the fine-scale linear system \eqref{dis-mf} via an equation describing the conservation of flux coming from/going into the fracture control volumes it connects.

\begin{figure}[htb!]%
\setlength{\fboxsep}{0pt}%
\setlength{\fboxrule}{0.2pt}%
\setlength{\unitlength}{1cm}%
\centering%
\subfigure[$d_{min} = \infty$, 1 DOF]{\fbox{\includegraphics[width=0.45\textwidth,height=0.4\textwidth]{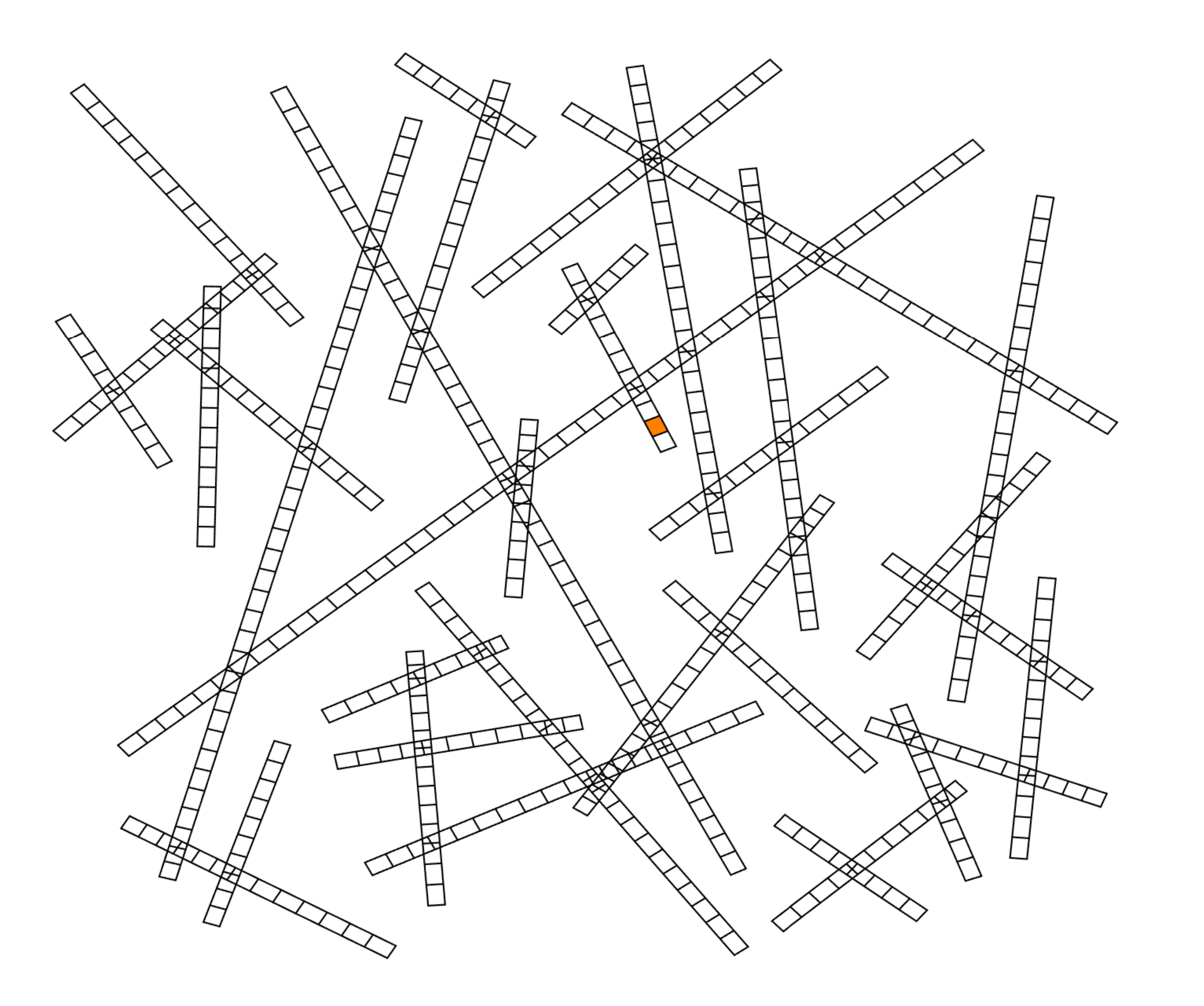}}\label{fig:1dof}}\hspace{1.3cm}%
\subfigure[$d_{min} = 20$, 15 DOF]{\fbox{\includegraphics[width=0.45\textwidth,height=0.4\textwidth]{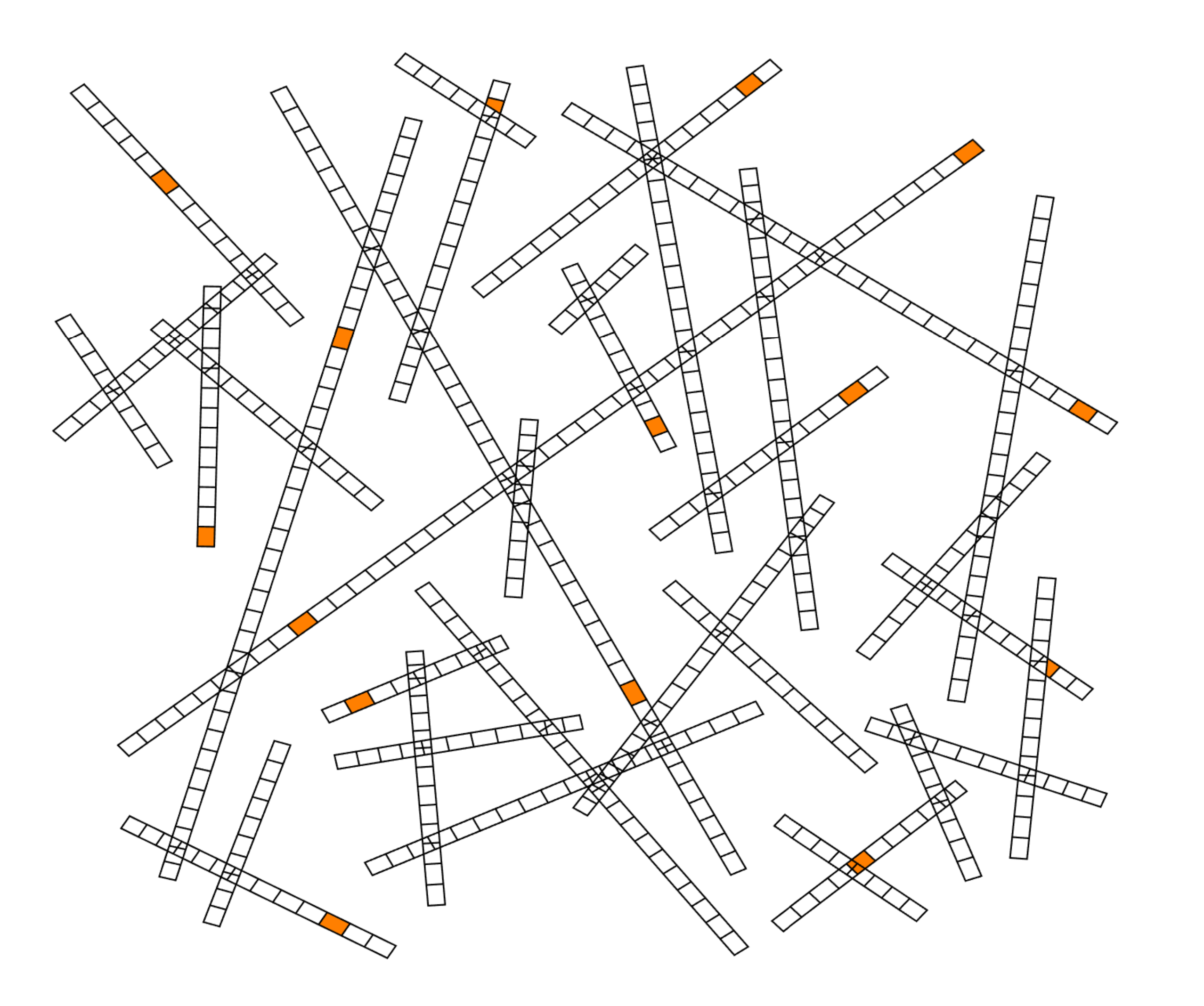}}\label{fig:15dof}}%
\caption{Coarse nodes obtained by the distance-based coarsening algorithm on a network with 575 fine-scale cells. Note that the aperture is magnified for clarity.}%
\label{fig:fracCoarsening}%
\end{figure}

As previously described, F-AMS requires primal- and dual-coarse grids in the (quasi-unstructured) fracture domain. In order to hide this complexity from the user, this paper introduces a distance-based algorithm for the automatic coarsening of fracture networks, as described in Table~\ref{fractureCoarsening}. The network's graph is traversed in a breadth-first order such that a distance of at least $d_{min}$ cells is guaranteed between each pair of resulting coarse nodes. As such, $d_{min}$ can be seen as a fracture coarsening factor. Note that choosing $d_{min} = \infty$ results in a single coarse node, as shown in Fig.~\ref{fig:1dof}. Moreover, Fig.~\ref{fig:15dof} depicts the result of the coarsening algorithm for $d_{min} = 20$ cells on a fairly complex fracture network.

\begingroup
\renewcommand{\arraystretch}{1.2}
\begin{table}[htb!]
\caption{Distance-based fracture coarsening algorithm.}
\begin{tabular}{|l|}
\hline
Repeat the following for each fracture network, $f_i$, which has $N_{cells}$ fine cells:\\
  \hspace{1mm} 1. Choose $d_{min}$, the minimum distance between two fracture coarse nodes.\\
  \hspace{1mm} 2. Initialize three empty queues: $Q_{vertex}, Q_1, Q_2$.\\
  \hspace{1mm} 3. Initialize two vectors of length $N_{cells}$: $level$ (set to $\infty$) and $primal$ (set to $0$).\\
  \hspace{1mm} 4. $N_{primal} \leftarrow 0$\\
  \hspace{1mm} 5. Choose a starting cell from $f_i$ and add it to $Q_{vertex}$.\\
  \hspace{1mm} Repeat until $Q_{vertex}$ is empty:\\
  \hspace{6mm} 6. $vertex \leftarrow$ extract top of $Q_{vertex}$\\
  \hspace{6mm} 7. $N_{primal} \leftarrow N_{primal} + 1$\\
  \hspace{6mm} 8. Create primal block number $N_{primal}$ with $vertex$ as its coarse node.\\
  \hspace{6mm} 9. $primal[vertex] \leftarrow N_{primal}$ and $level[vertex] \leftarrow 0$.\\
  \hspace{6mm} 10. Add $vertex$ to $Q_1$.\\
  \hspace{6mm} For $dist$ from $1$ up to and including $d_{min}$:\\
  \hspace{11mm} Repeat until $Q_1$ is empty:\\
  \hspace{16mm} 11. $cell \leftarrow$ extract top of $Q_1$\\
  \hspace{16mm} For each neighbour of $cell$, $neigh_j$, with $level[neigh_j] > dist$:\\
  \hspace{21mm} 12. remove $neigh_j$ from $Q_{vertex}$, if it is present (i.e. $level[neigh_j]$ = $d_{min}$)\\
  \hspace{21mm} 13. $level[neigh_j] \leftarrow dist$\\
  \hspace{21mm} 14. $primal[neigh_j] \leftarrow N_{primal}$\\
  \hspace{21mm} 15. add $neigh_j$ to $Q_2$.\\
  \hspace{11mm} 16. swap $Q_1$ and $Q_2$.\\
  \hspace{6mm} 17. empty $Q_1$ into $Q_{vertex}$.\\~\\
  \hspace{1mm} At this point, $primal[i]$ gives the index of fine cell $i$'s primal block, while $level[i]$ is \\
  \hspace{1mm} the distance from cell $i$ to the nearest vertex. The fine cells which were not marked\\
  \hspace{1mm} as vertices will form edges on the dual-coarse grid.\\
  \hline
\end{tabular}
\label{fractureCoarsening}
\end{table}
\endgroup

\subsection{F-AMS convergence}

\begin{figure}[htb!]%
\subfigure[$log_{10}(k^m)$ and fracture map]{\raisebox{0.5cm}{\includegraphics[width=0.45\textwidth,height=0.39\textwidth]{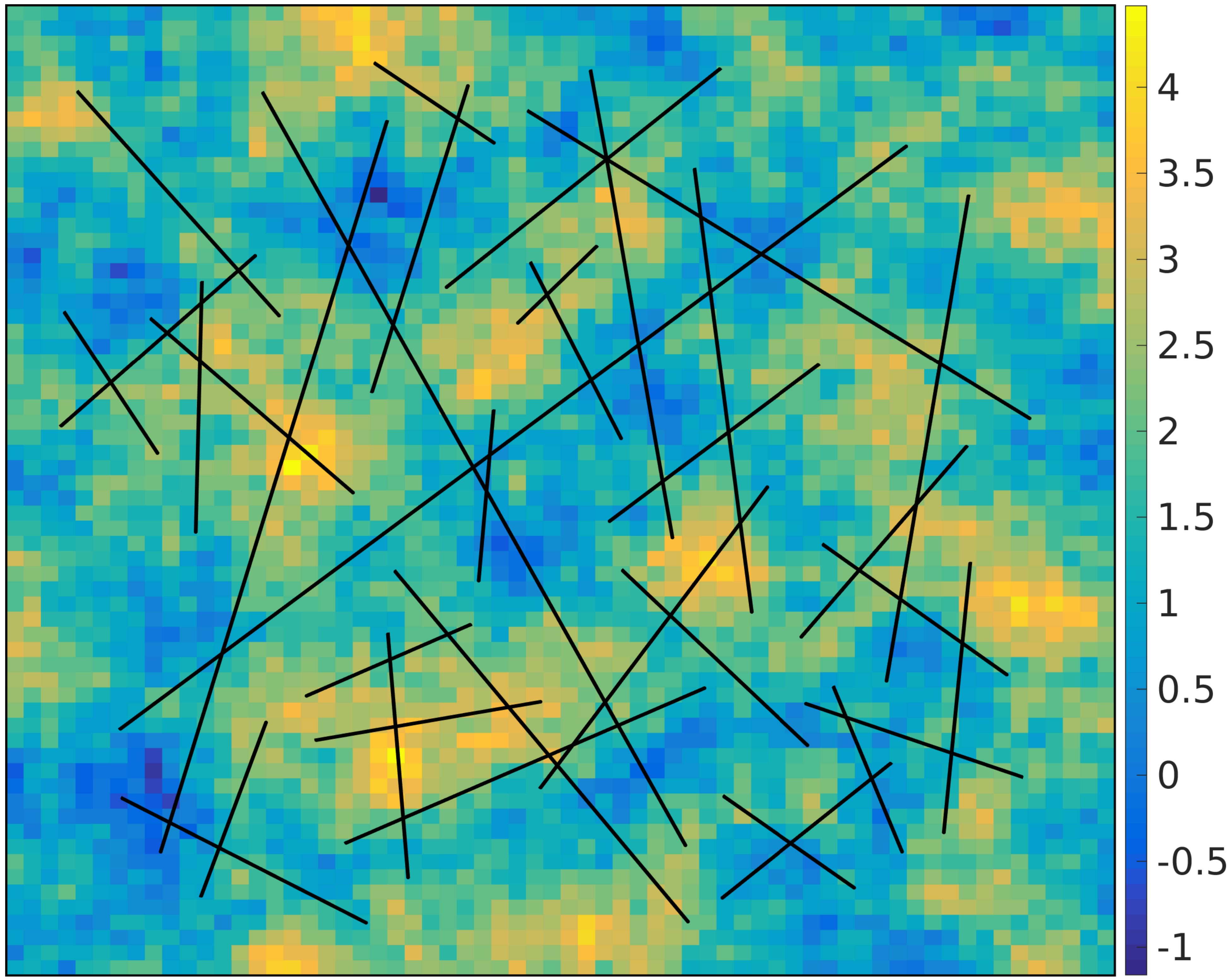}}\label{fig:perm2d}}\hspace{1cm}%
\subfigure[Pressure solution]{\includegraphics[width=0.5\textwidth,height=0.45\textwidth]{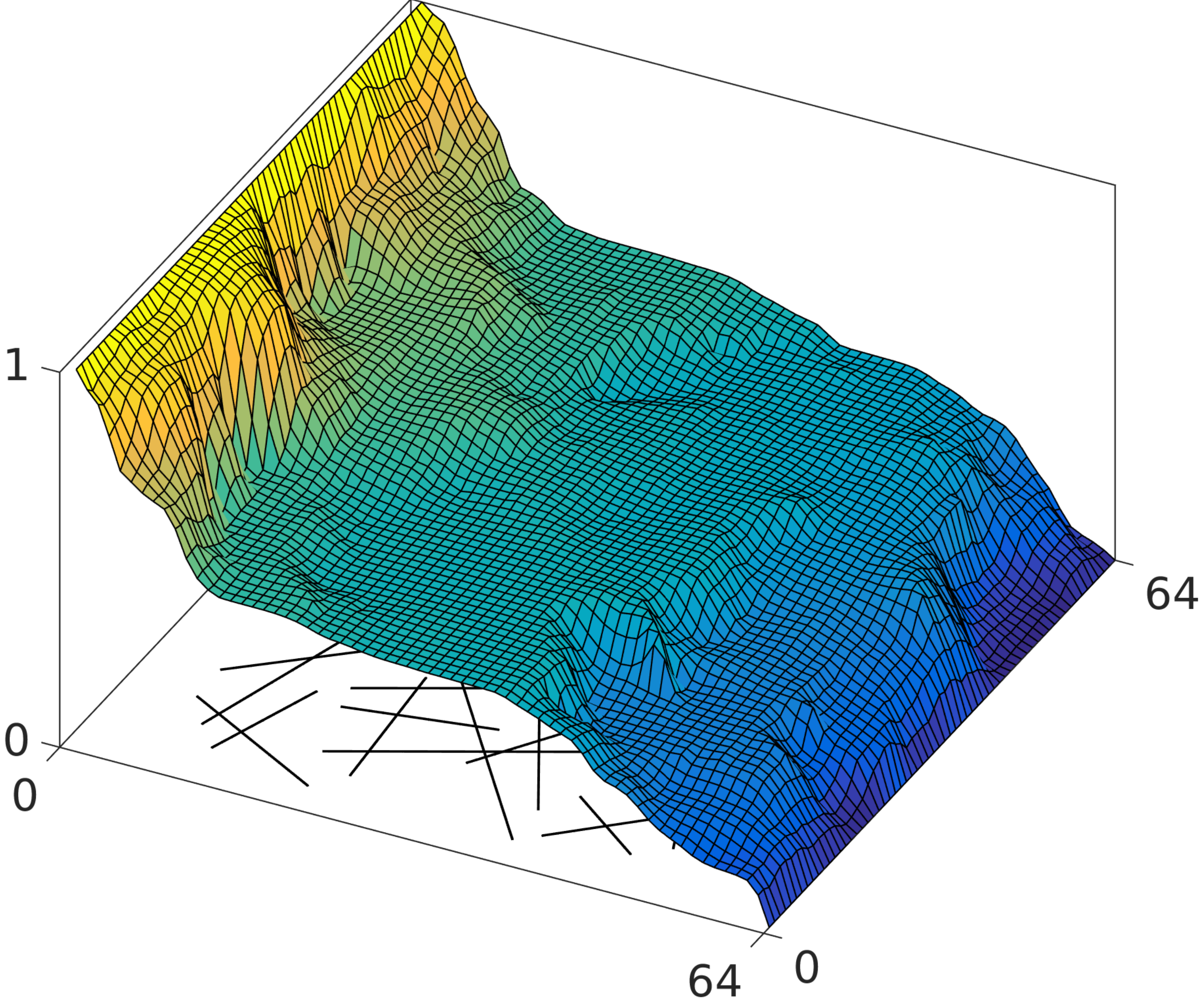}\label{fig:pres2d}}%
\caption{2D Test Case: heterogeneous matrix rock with $64 \times 64$ matrix and $575$ fracture grid cells (a). Two pressure-constrained wells at the West and East boundaries are placed, resulting in the pressure solution shown in (b) for $T_{ratio} = 100$.}%
\label{fig:2dcase}%
\end{figure}

The fracture network from Fig.~\ref{fig:fracCoarsening} was embedded into a heterogeneous (patchy) matrix rock with two pressure-constrained wells added on the West and East boundaries. This 2D test case, depicted in Fig.~\ref{fig:2dcase}, was used to study the convergence properties of F-AMS, with the four coupling strategies presented before.

\begin{figure}[htb!]%
\subfigure[Decoupled-AMS 1 fracture DOF,\newline $\|\epsilon\|_2 = 8.6557$]{\includegraphics[width=0.5\textwidth,height=0.4\textwidth]{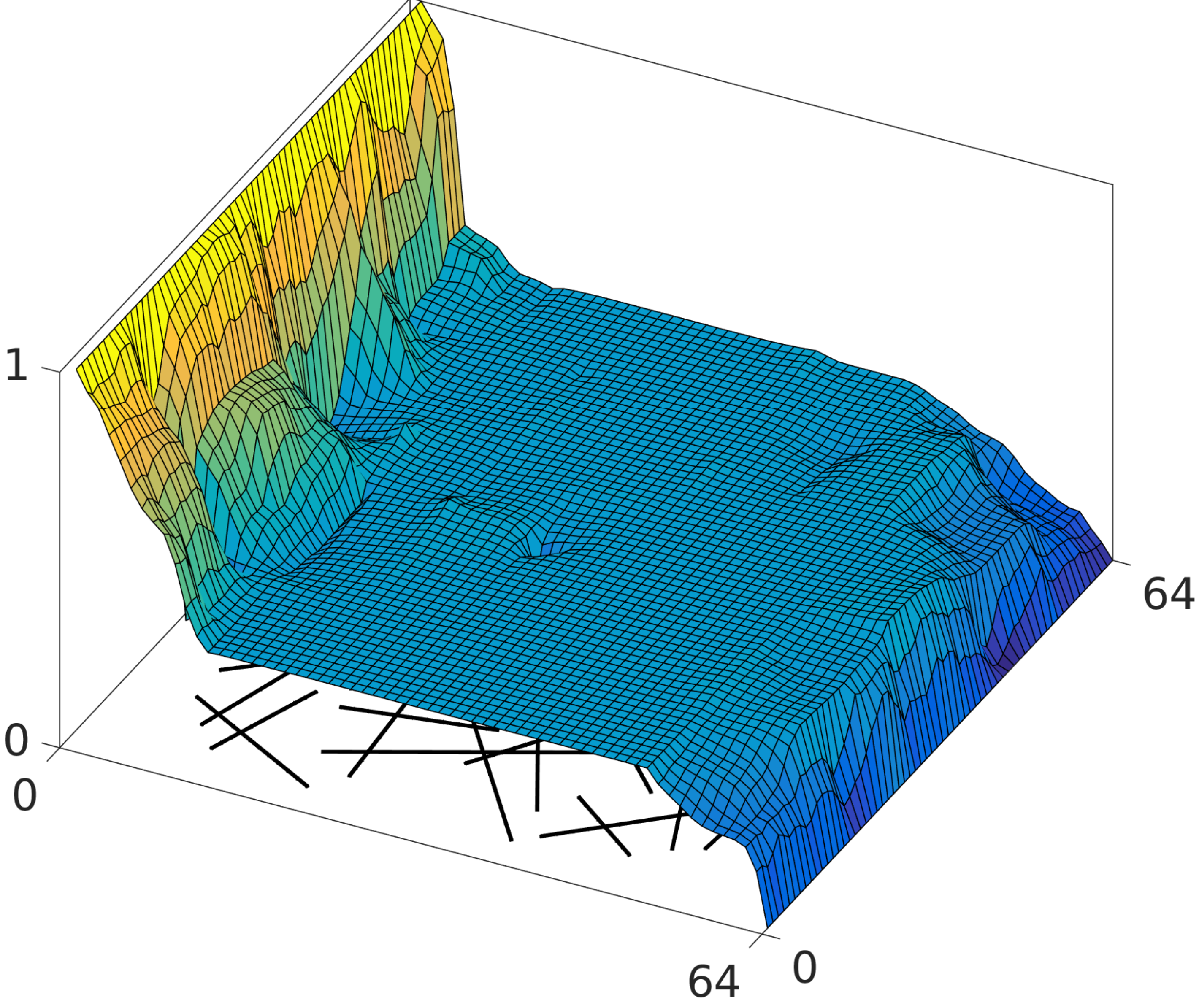}\label{fig:decoupledAMS_1dof}}\hspace{0.3cm}%
\subfigure[Frac-AMS 1 fracture DOF,\newline $\|\epsilon\|_2 = 6.1719$]{\includegraphics[width=0.5\textwidth,height=0.4\textwidth]{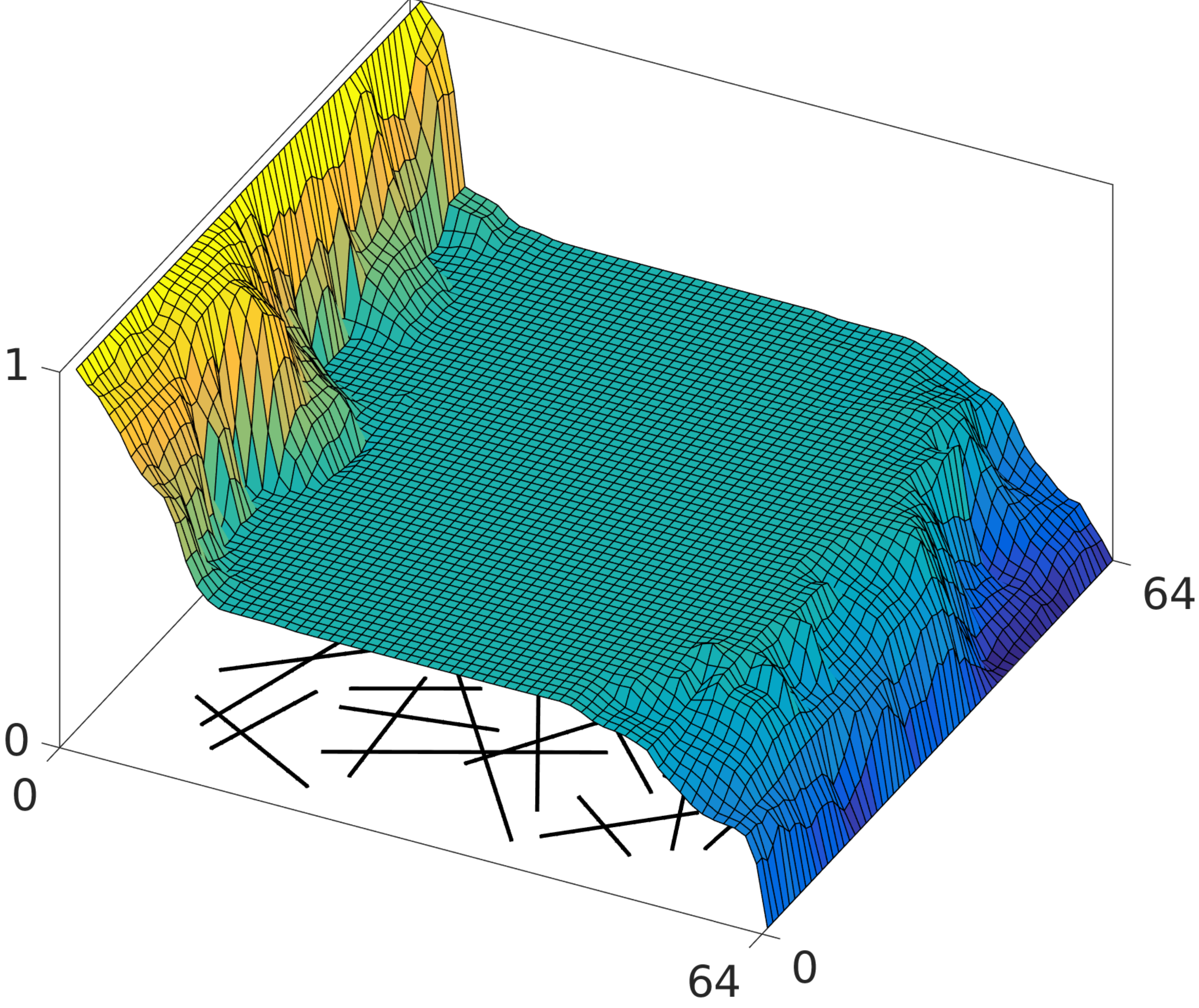}\label{fig:fracAMS_1dof}}\\%
\subfigure[Rock-AMS 1 fracture DOF,\newline $\|\epsilon\|_2 = 2.6708$]{\includegraphics[width=0.5\textwidth,height=0.4\textwidth]{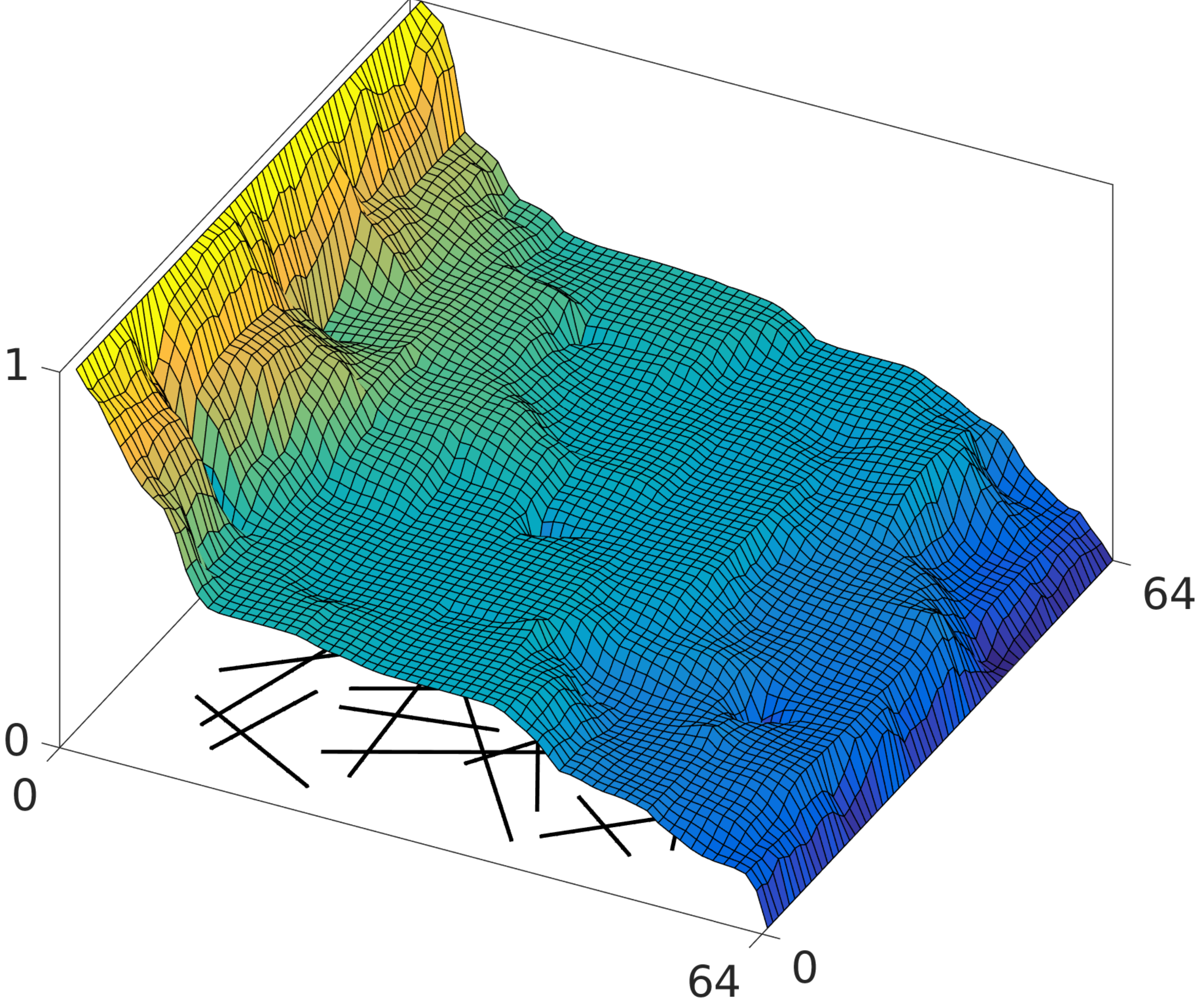}\label{fig:rockAMS_1dof}}\hspace{0.3cm}%
\subfigure[Coupled-AMS 1 fracture DOF,\newline $\|\epsilon\|_2 = 4.4906$]{\includegraphics[width=0.5\textwidth,height=0.4\textwidth]{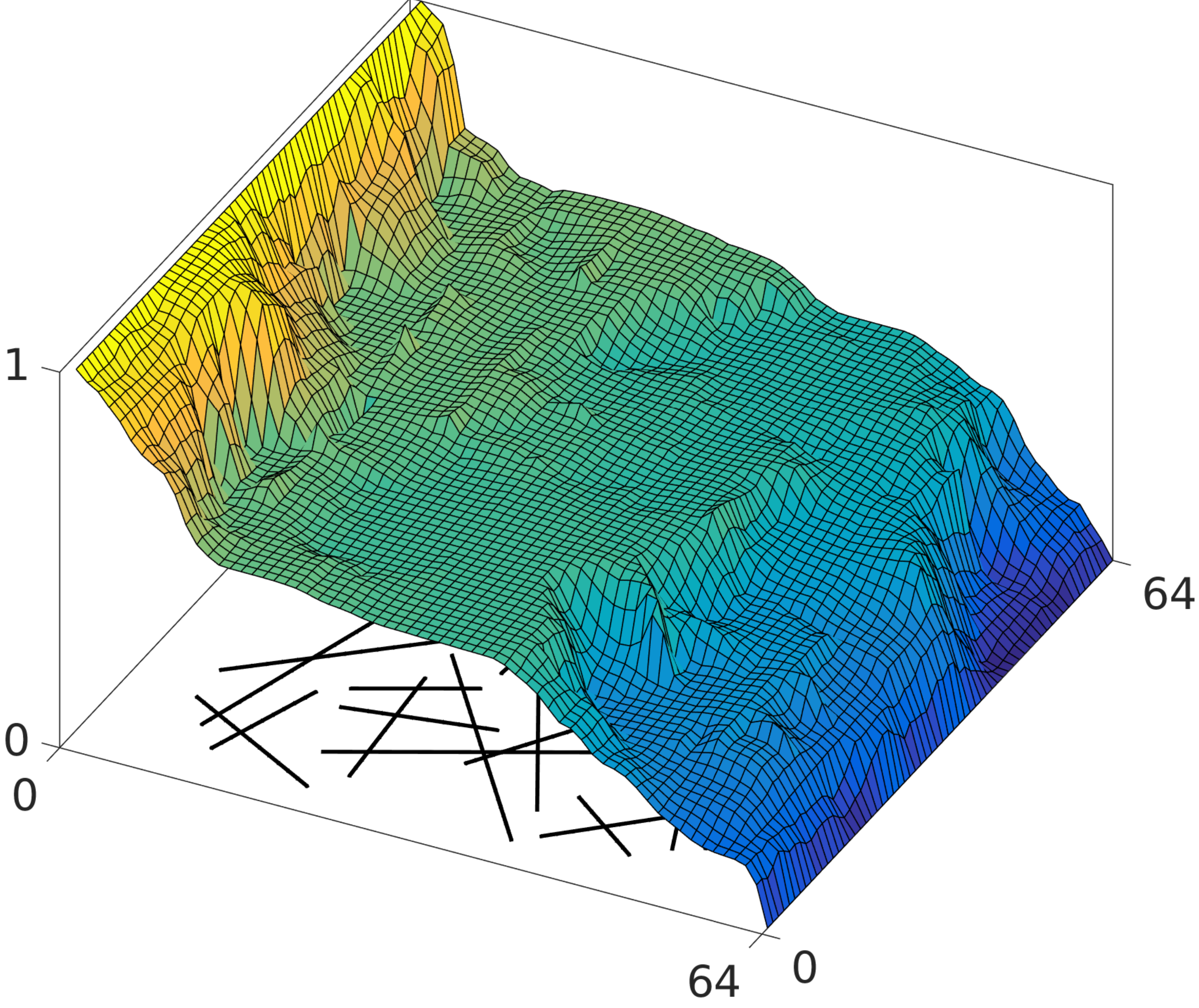}\label{fig:coupledAMS_1dof}}%
\caption{Pressure solutions after a single multiscale iteration with 1 fracture coarse cell, with 4 different prolongation coupling strategies. In all cases, the matrix coarse grid contains $8 \times 8$ blocks and the 2-norm of the error is indicated.}%
\label{fig:1MSiter_1dof}%
\end{figure}

\begin{figure}[htb!]%
\subfigure[Decoupled-AMS 15 fracture DOF,\newline $\|\epsilon\|_2 = 3.6173$]{\includegraphics[width=0.5\textwidth,height=0.4\textwidth]{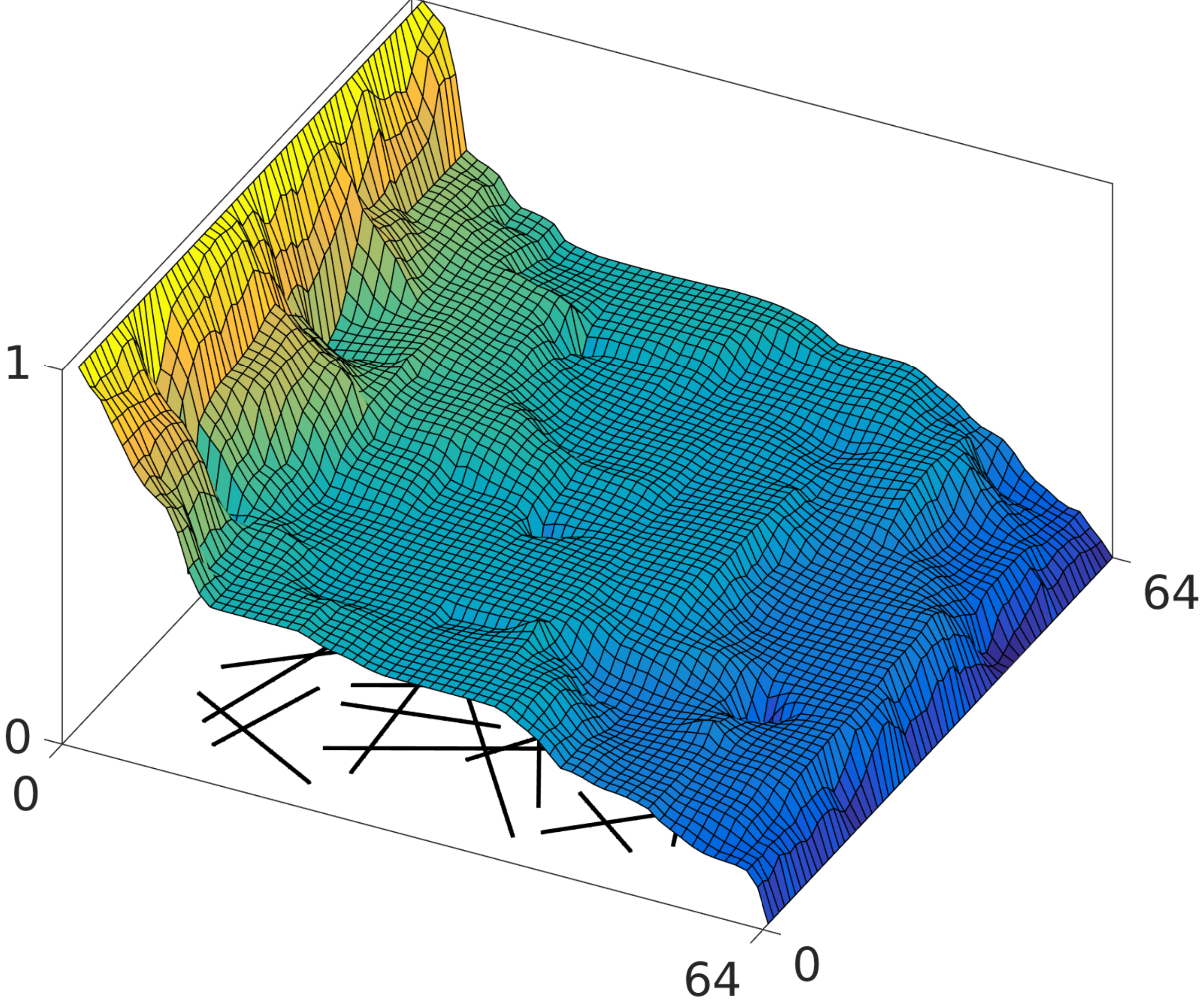}\label{fig:decoupledAMS_15dof}}\hspace{0.3cm}%
\subfigure[Frac-AMS 15 fracture DOF,\newline $\|\epsilon\|_2 = 0.9502$]{\includegraphics[width=0.5\textwidth,height=0.4\textwidth]{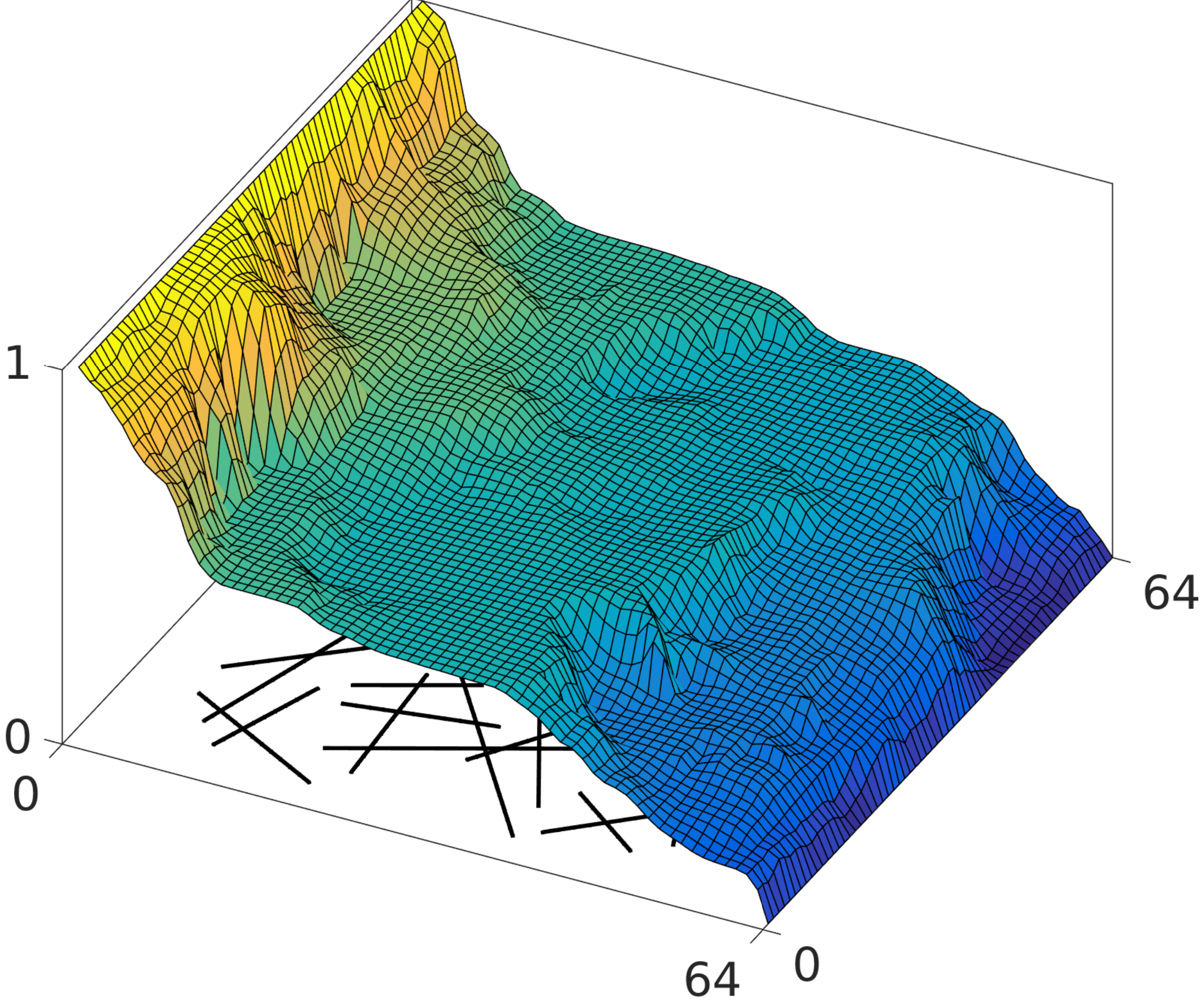}\label{fig:fracAMS_15dof}}\\%
\subfigure[Rock-AMS 15 fracture DOF,\newline $\|\epsilon\|_2 = 2.6708$]{\includegraphics[width=0.5\textwidth,height=0.4\textwidth]{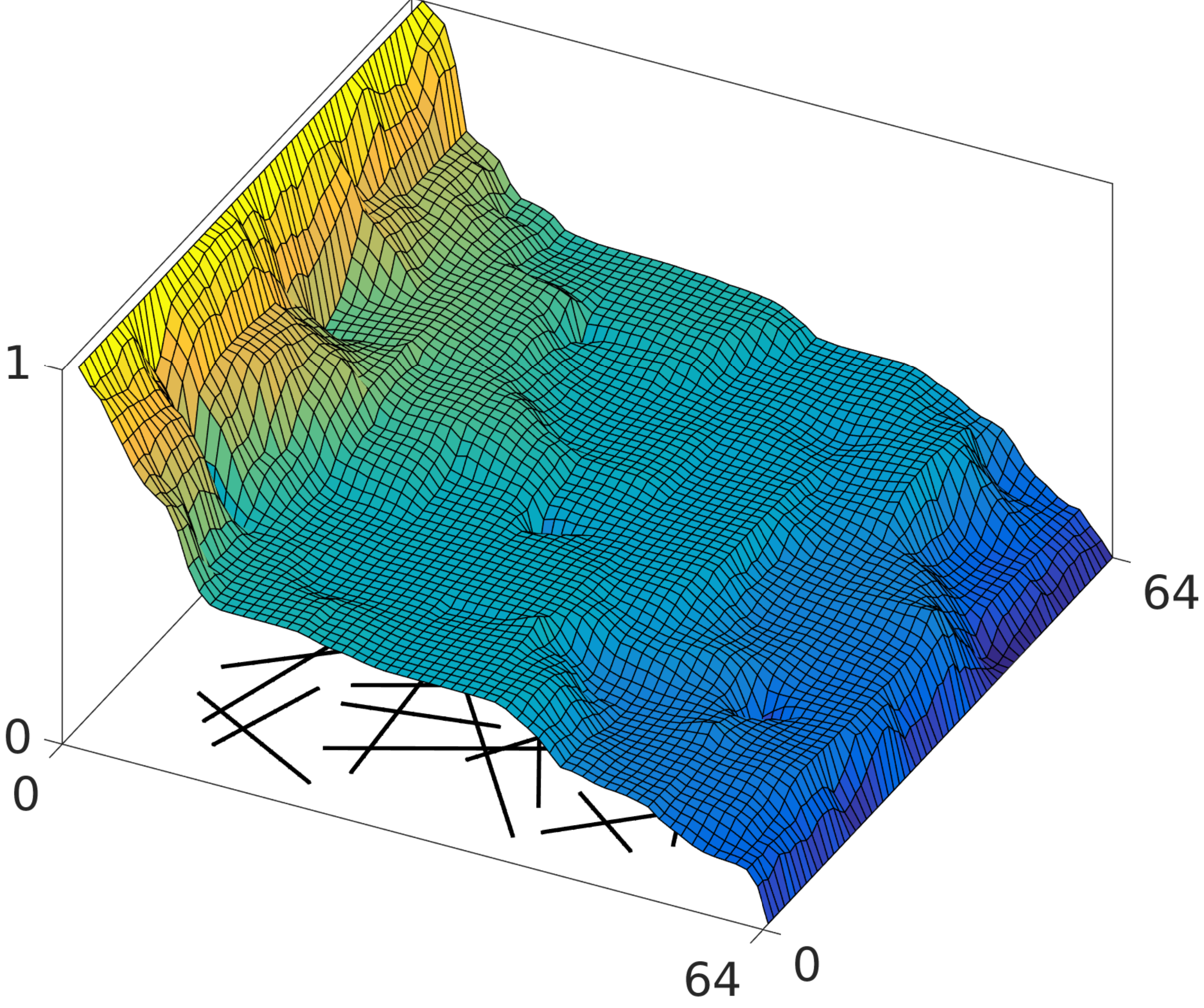}\label{fig:rockAMS_15dof}}\hspace{0.3cm}%
\subfigure[Coupled-AMS 15 fracture DOF,\newline $\|\epsilon\|_2 = 1.8185$]{\includegraphics[width=0.5\textwidth,height=0.4\textwidth]{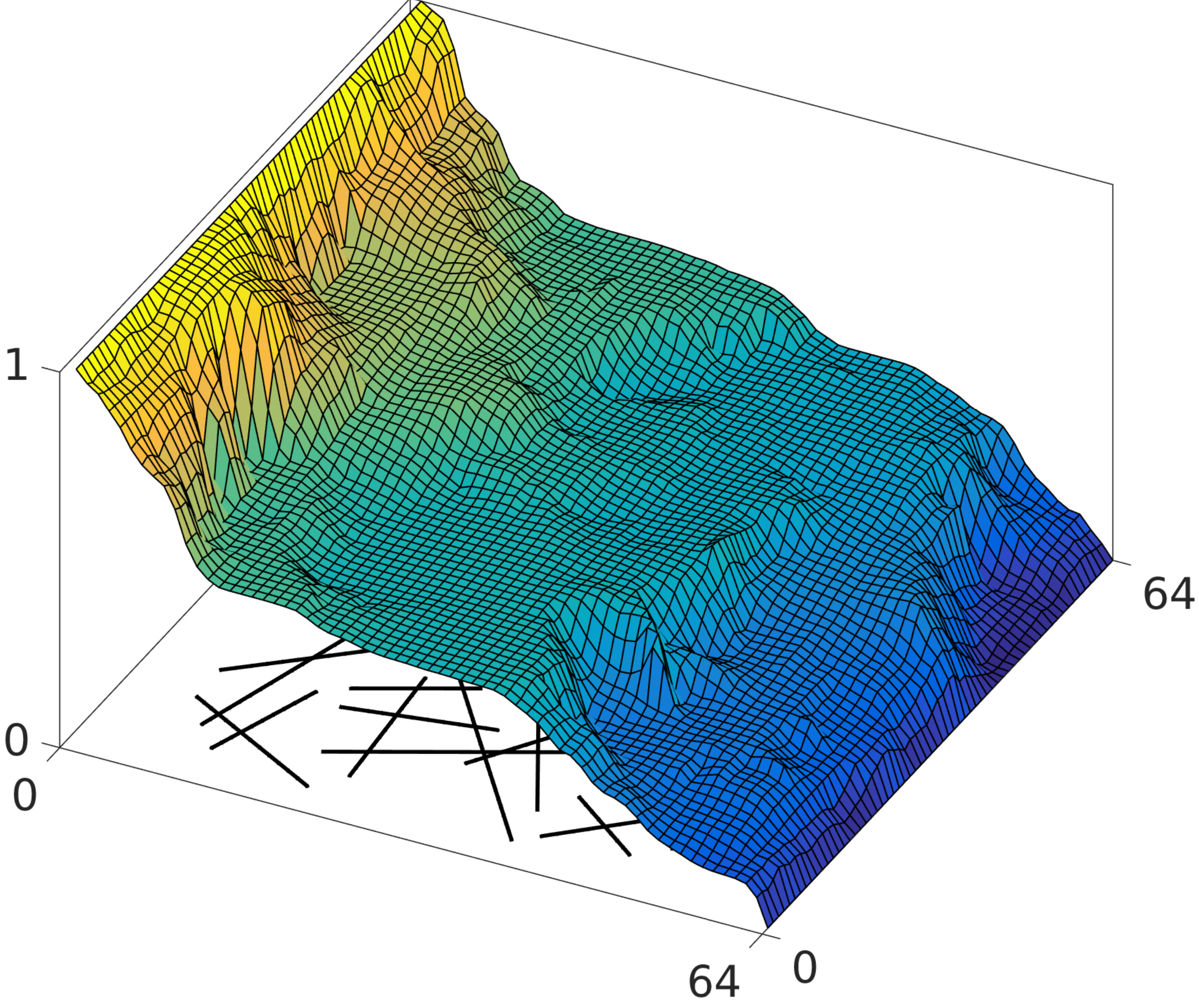}\label{fig:coupledAMS_15dof}}%
\caption{Pressure solutions after a single multiscale iteration with 15 fracture coarse cells, with 4 different prolongation coupling strategies. In all cases, the matrix coarse grid contains $8 \times 8$ blocks and the 2-norm of the error is indicated.}%
\label{fig:1MSiter_15dof}%
\end{figure}

In order to test the accuracy of the basis functions as pressure interpolators, F-AMS was stopped after Step 5 of its first iteration (see Table~\ref{fams-algorithm}). The solution is shown in Figs.~\ref{fig:1MSiter_1dof} and \ref{fig:1MSiter_15dof}. Using a single fracture coarse DOF leads to a poor approximation of the pressure distribution, especially for Decoupled- and Frac-AMS. Figure~\ref{fig:fracAMS_1dof} depicts the results of Frac-AMS with 1 fracture coarse DOF. In this setup, F-AMS treats fractures similar to \cite{hadi-frac-jcp}. Note that, due to the large length scale of the network, having a constant interpolator for the pressure along the network results in an initial solution which lacks a lot of the fine-scale features. In contrast, the Rock-AMS and Coupled-AMS (Figs.~\ref{fig:rockAMS_1dof} and \ref{fig:coupledAMS_1dof}) place a lot more emphasis on the matrix basis functions and, since, in this test case, the rock heterogeneity is the main source of approximation error, their results are more accurate. It may seem unexpected that the Coupled-AMS performs slightly worse than Rock-AMS. This can be attributed to the fact that a single fracture DOF is not sufficient to accurately capture the pressure distribution in the large fracture network, especially under the localization assumption. However, when only few additional coarse DOF are added in the fracture domain (as shown in Fig.~\ref{fig:1MSiter_15dof}), the situation improves dramatically for Decoupled-AMS, Frac-AMS and Coupled-AMS. Note that Rock-AMS, on the other hand, is insensitive to this change. 

\begin{figure}[htb!]%
\centering%
\subfigure{\includegraphics[width=0.5\textwidth,height=0.33\textwidth]{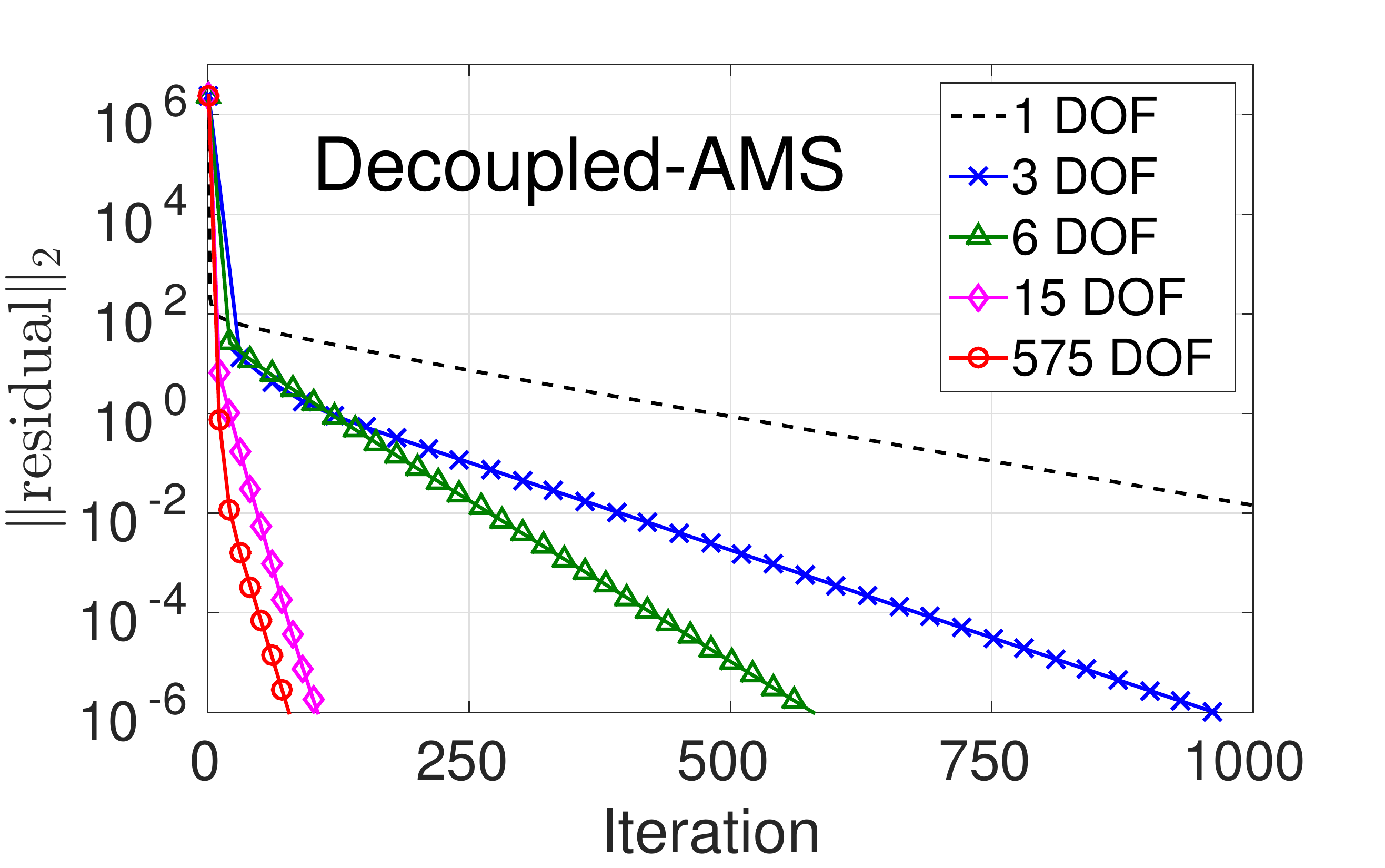}\label{fig:decoupledAMS_conv}}%
\subfigure{\includegraphics[width=0.5\textwidth,height=0.33\textwidth]{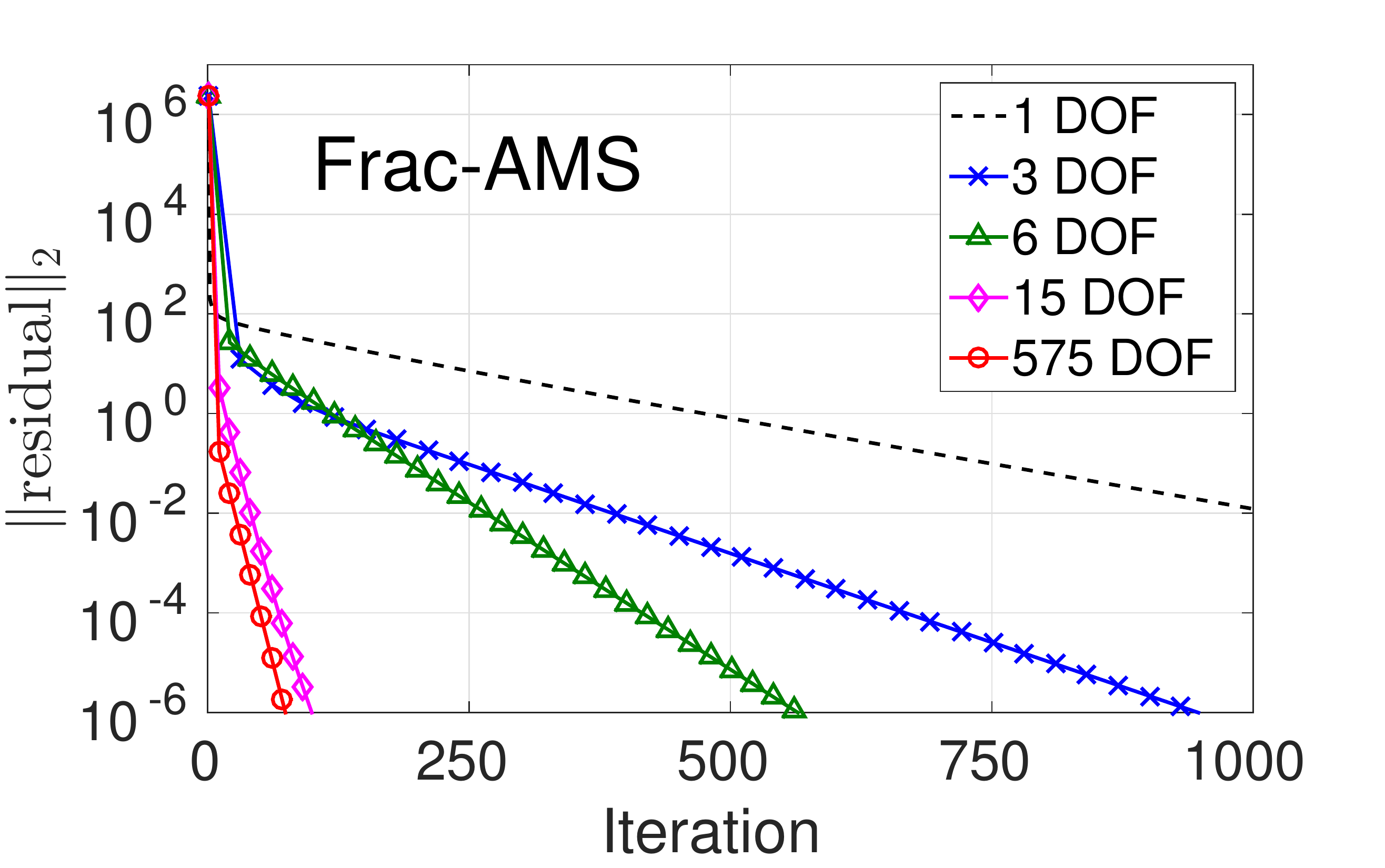}\label{fig:fracAMS_conv}}\\%
\subfigure{\includegraphics[width=0.5\textwidth,height=0.33\textwidth]{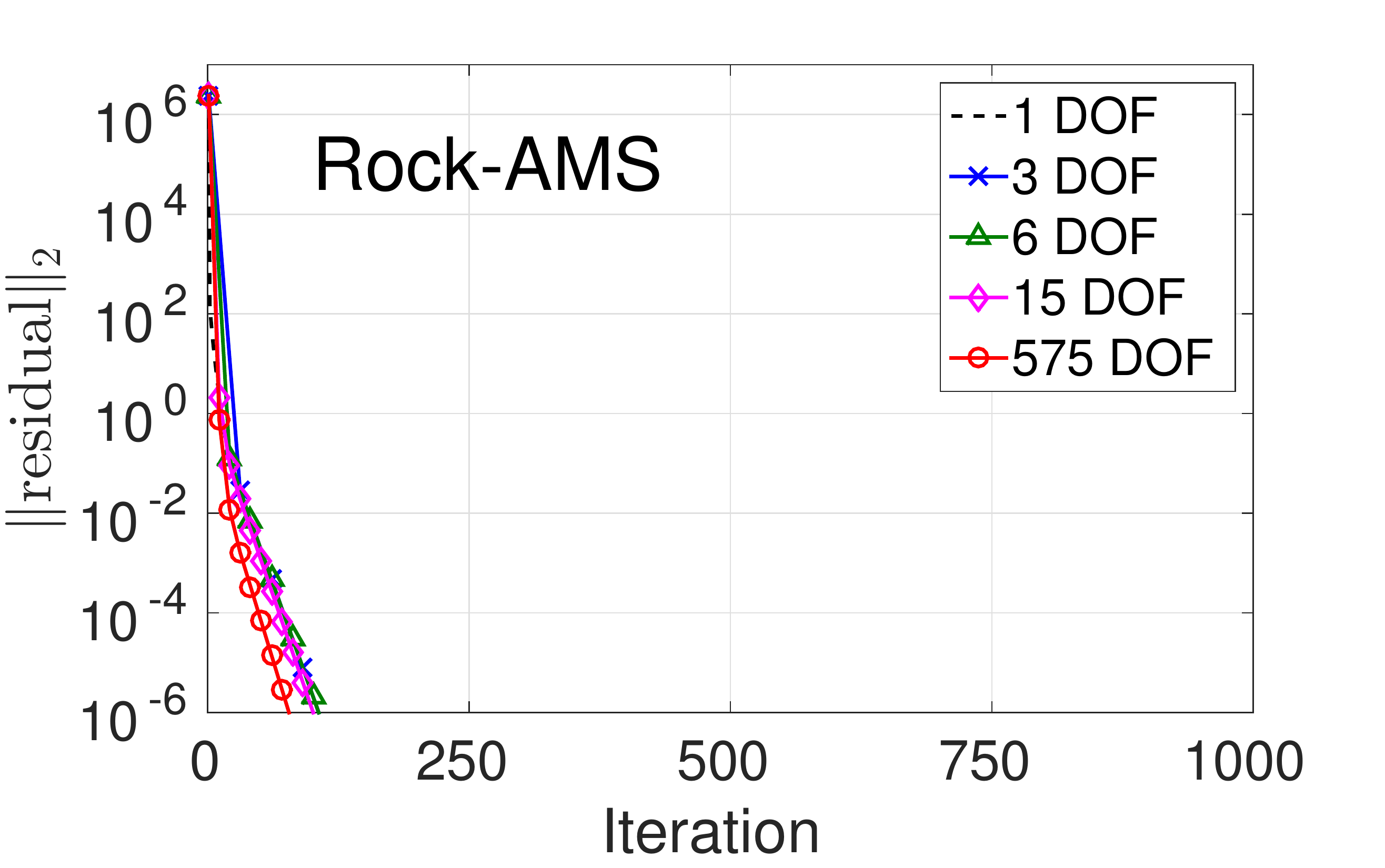}\label{fig:rockAMS_conv}}%
\subfigure{\includegraphics[width=0.5\textwidth,height=0.33\textwidth]{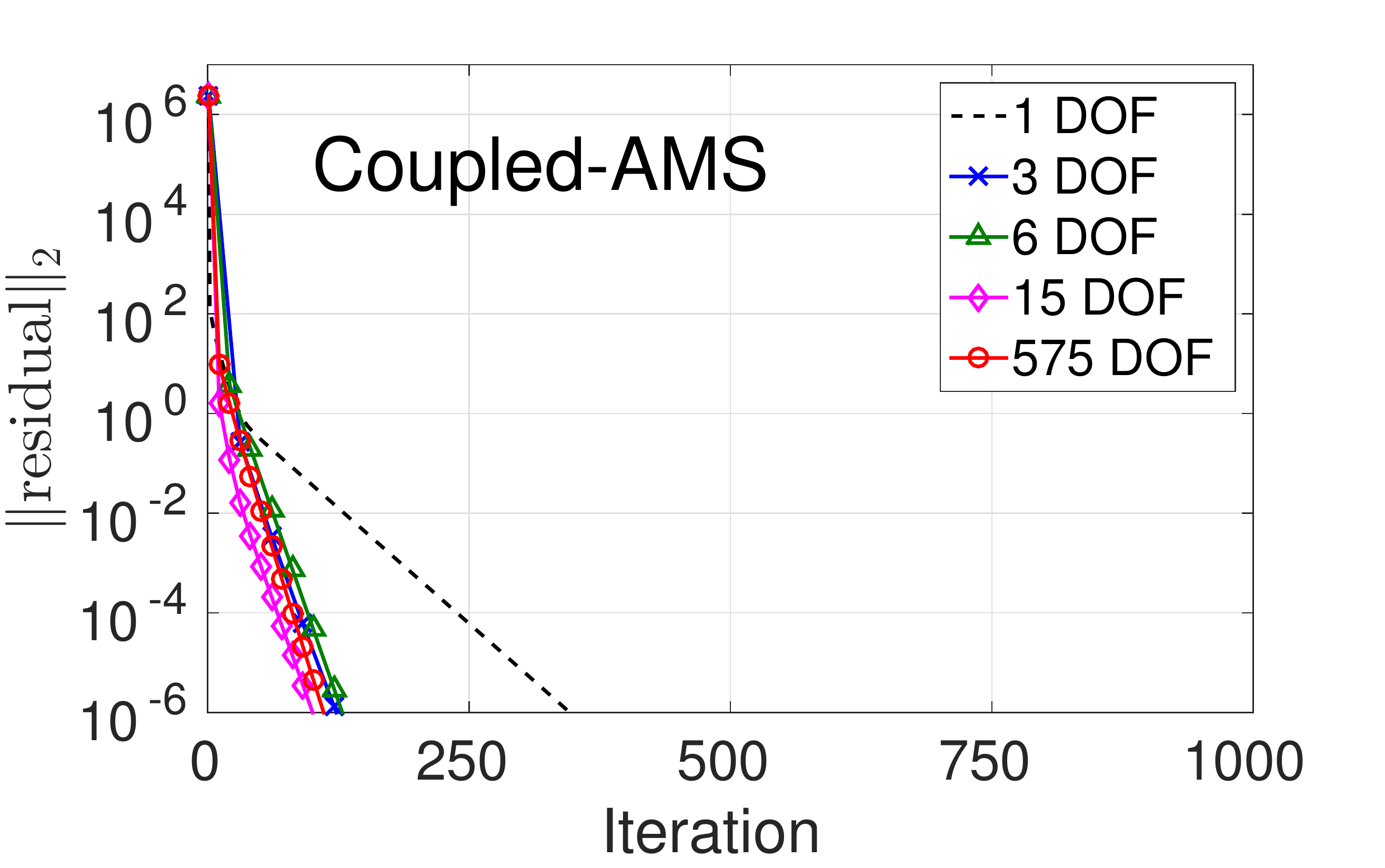}\label{fig:coupledAMS_conv}}%
\caption{Convergence history of F-AMS for the 2D test case with four basis function coupling strategies and different number of coarse DOF in the fracture network.}%
\label{fig:conv}
\end{figure}

For the results in Fig.~\ref{fig:conv}, as well as all subsequent experiments in this paper, F-AMS was iterated until converged to a residual 2-norm of $10^{-6}$. It is clear that Rock-AMS shows a good convergence rate on this particular 2D test case, regardless of the fracture coarsening factor. Also, the other strategies reach a similar behaviour when the fracture network is enriched with only few additional coarse-scale DOF.
 
\subsection{Basis function truncation}

\begin{figure}[htb!]%
\subfigure[Wells and fractures]{\includegraphics[height=0.42\textwidth,width=0.42\textwidth]{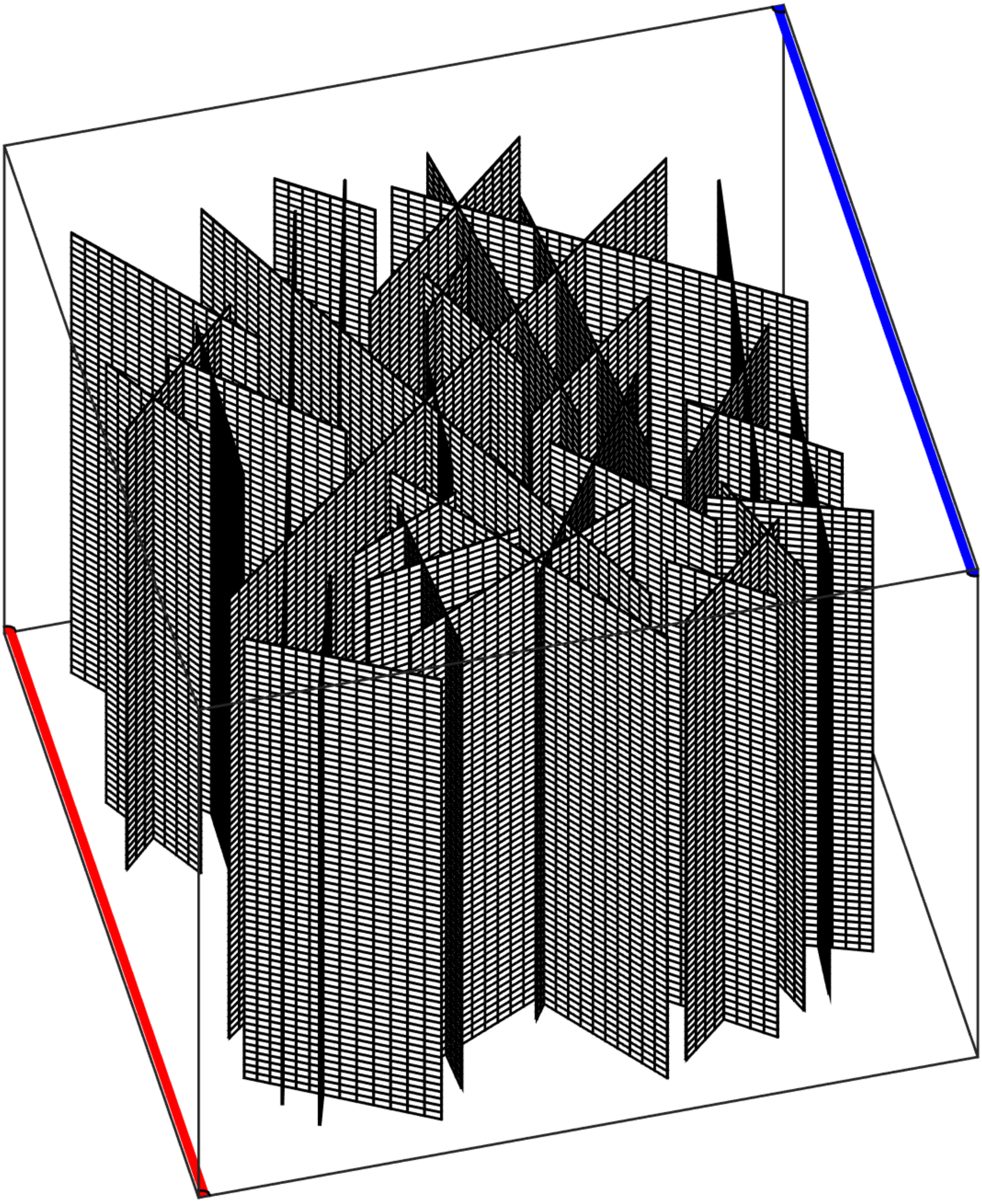}\label{fig:fracs3d}}\hspace{1.8cm}%
\subfigure[$log_{10}(k^m)$]{\includegraphics[height=0.42\textwidth,width=0.42\textwidth]{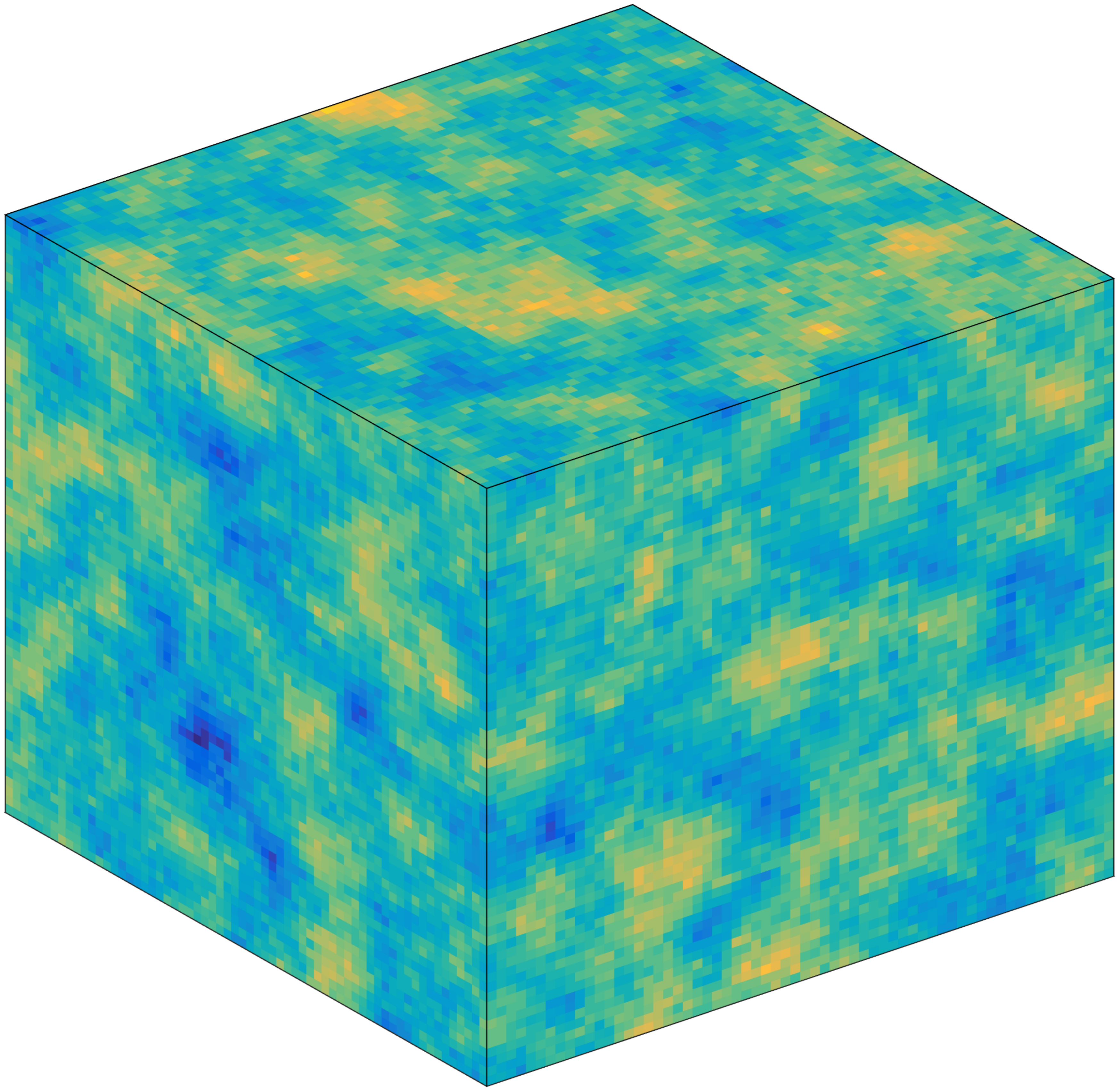}\hspace{0.2cm}\raisebox{0.7cm}{\includegraphics[height=0.3\textwidth,width=0.4cm]{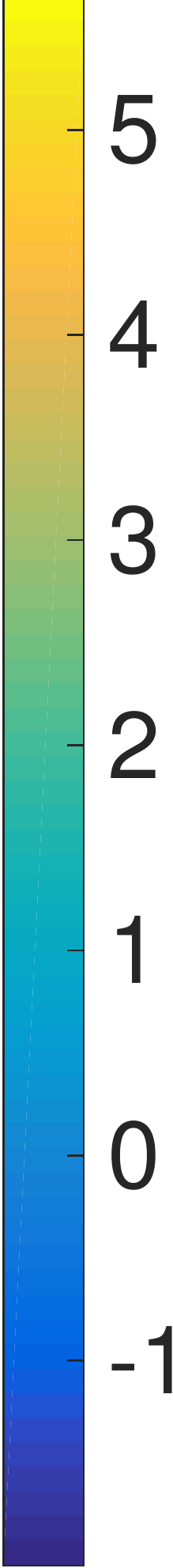}}\label{fig:perm3d}}\\%
\subfigure[Fracture pressure]{\includegraphics[height=0.42\textwidth,width=0.42\textwidth]{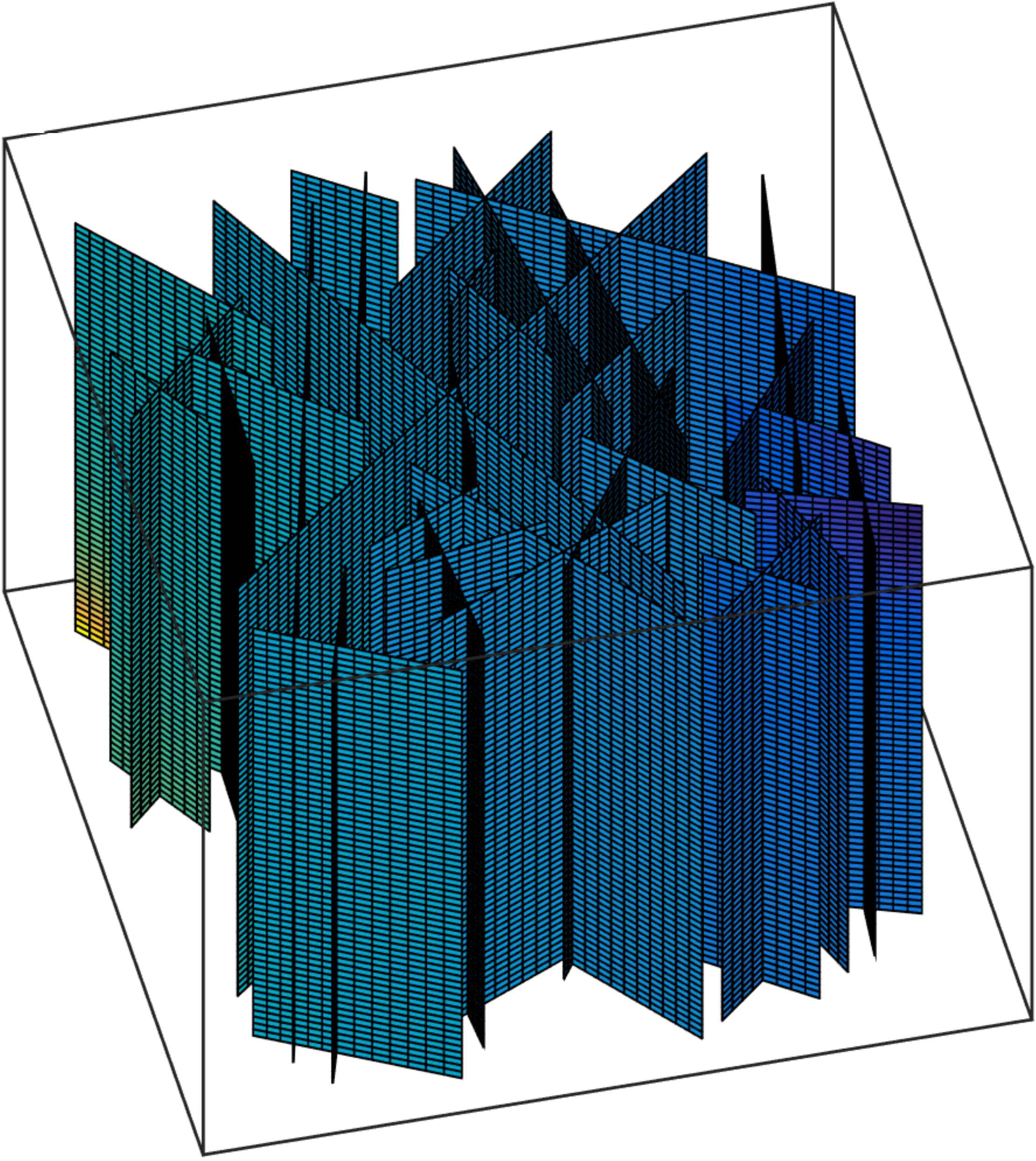}\hspace{0.2cm}\raisebox{0.7cm}{\includegraphics[height=0.3\textwidth,width=0.7cm]{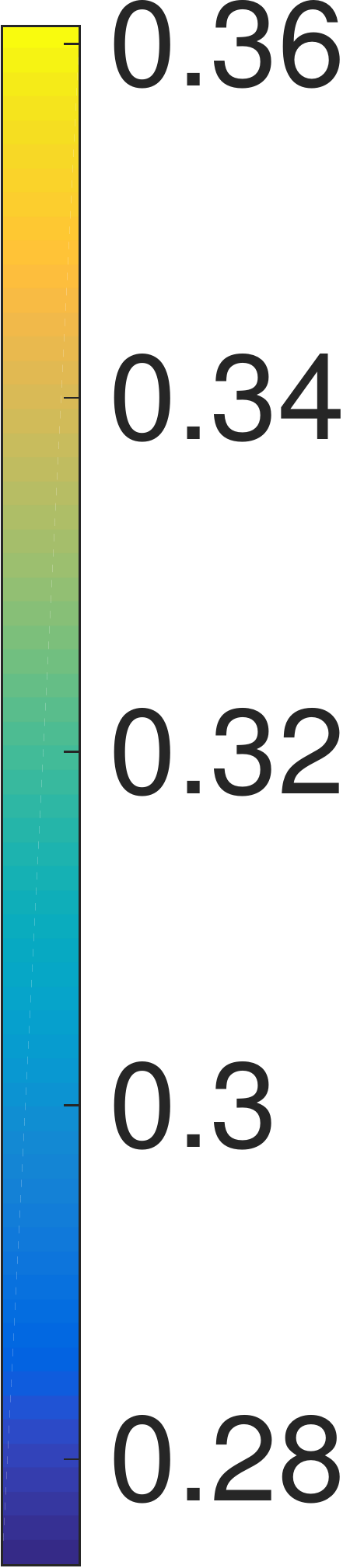}}\label{fig:fracPres3d}}\hspace{0.9cm}%
\subfigure[Matrix pressure]{\includegraphics[height=0.42\textwidth,width=0.42\textwidth]{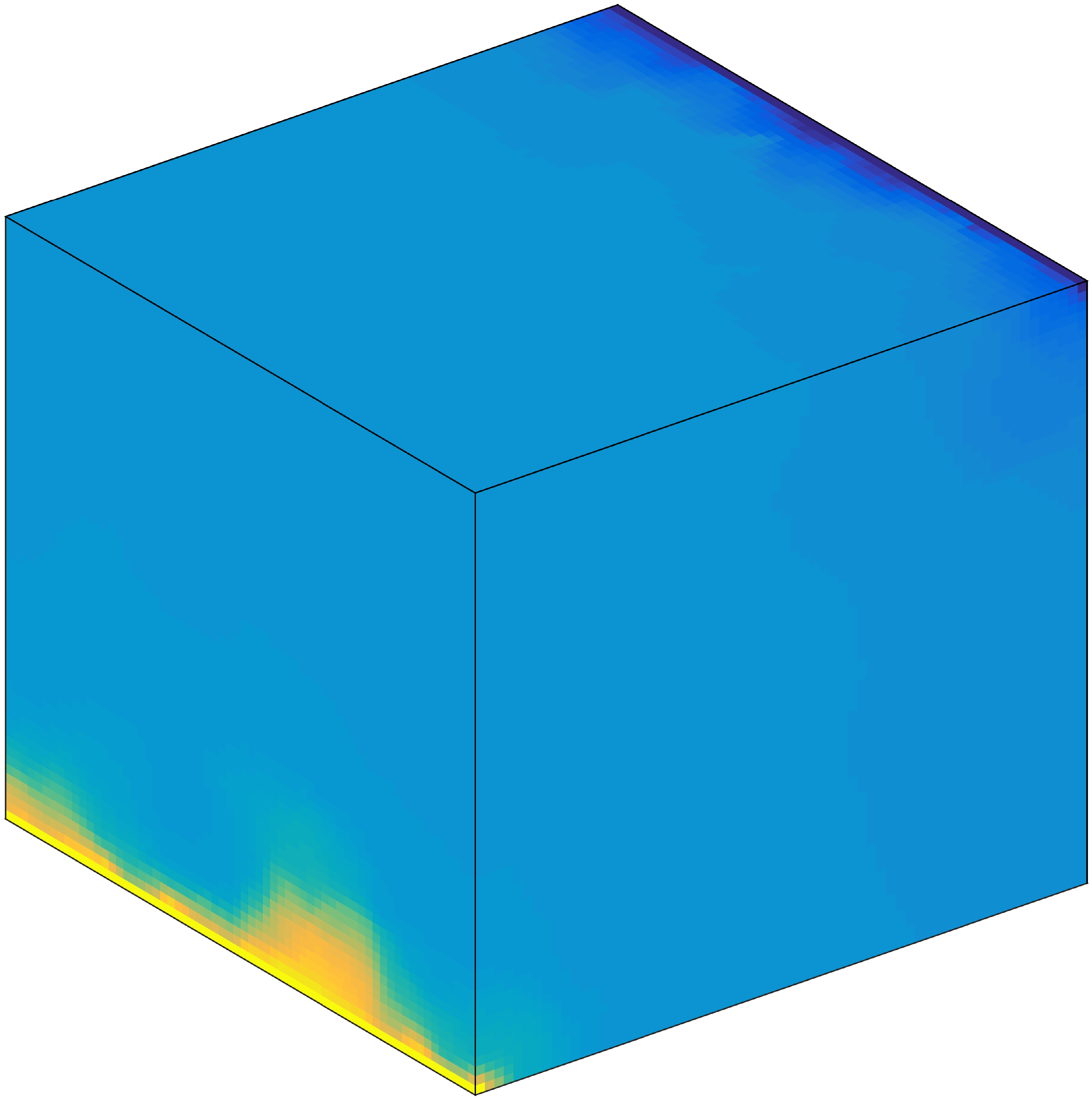}\hspace{0.2cm}\raisebox{0.7cm}{\includegraphics[height=0.3\textwidth,width=0.6cm]{Fig/method/colorbar.pdf}}\label{fig:pres3d}}%
\caption{Illustration of the 3D test case with $575 \times 64$ fracture and $64 \times 64 \times 64$ matrix grid cells. The logarithm of the heterogeneous permeability map is provided in (b). Two pressure-constrained wells are placed on opposite edges, as shown in (a). The pressure solution is shown for $T_{ratio}=100$ for fracture (c) and matrix (d). F-AMS employs the coarsening ratios of $8\times8\times8$ for matrix and $8\times8$ for fractures.}%
\label{fig:3dcase}%
\end{figure}

In order to get an idea of the performance of F-AMS on realistic fractured reservoirs, Fig.~\ref{fig:3dcase} introduces a 3D scenario, where the fracture network from Fig.~\ref{fig:2dcase} was extruded and discretized along the Z axis. Two pressure-constrained horizontal wells are placed on opposite edges of the domain boundary. Figure~\ref{fig:pres3d} shows the fine-scale pressure solution obtained on the heterogeneous (patchy) matrix permeability field shown in Fig.~\ref{fig:perm3d}. Note that, even though the matrix-fracture conductivity contrast is of only two orders of magnitude, this is enough to make the pressure distribution in the fracture network insensitive to the matrix heterogeneity (see the approximately constant pressure in Fig.~\ref{fig:fracPres3d}).

The procedure described in Table~\ref{fractureCoarsening} is followed to determine the fracture coarse nodes along the projection of the network on the X-Y plane. Then, the resulting coarse grids are extruded along the Z-axis uniformly, with the vertical distribution of the coarse nodes honouring the user-specified coarsening ratio (see Fig.~\ref{fig:coarseGrids}).

The increased number of cells in both the matrix ($64 \times 64 \times 64$), as well as the fracture network ($575 \times 64$), compared to the 2D case, can lead to a larger number of non-zeros for the basis functions which take into account the coupling between the two media. This can be particularly severe for Coupled-AMS (Fig.~\ref{fig:coupledAMS}), since the high density of fractures can cause a large number of dual blocks to be merged. The resulting basis functions have a wider support and can potentially lead to more accurate interpolations, however, the added density to the prolongation operator will also increase the computational effort necessary to construct and solve the coarse-scale linear system (i.e., Step 5 in Table~\ref{fams-algorithm}).

\begin{figure}[htb!]%
\centering%
\setlength{\fboxsep}{0pt}%
\setlength{\fboxrule}{0.2pt}%
\setlength{\unitlength}{1cm}%
\begin{picture}(15,0.4)%
  \put(2.4,0){\large Decoupled-AMS}%
  \put(9.9,0){\large Frac-AMS}%
  \end{picture}\\%
  \begin{picture}(1.2,1.5)%
	\put(0.0,0.3){\rotatebox{90}{CPU time (sec)}}%
		\put(0.9,-0.1){0}%
		\put(0.9,1.7){8}%
		\put(0.7,3.3){16}%
\end{picture}%
\subfigure{\fbox{\includegraphics[width=0.41\textwidth,height=0.25\textwidth]{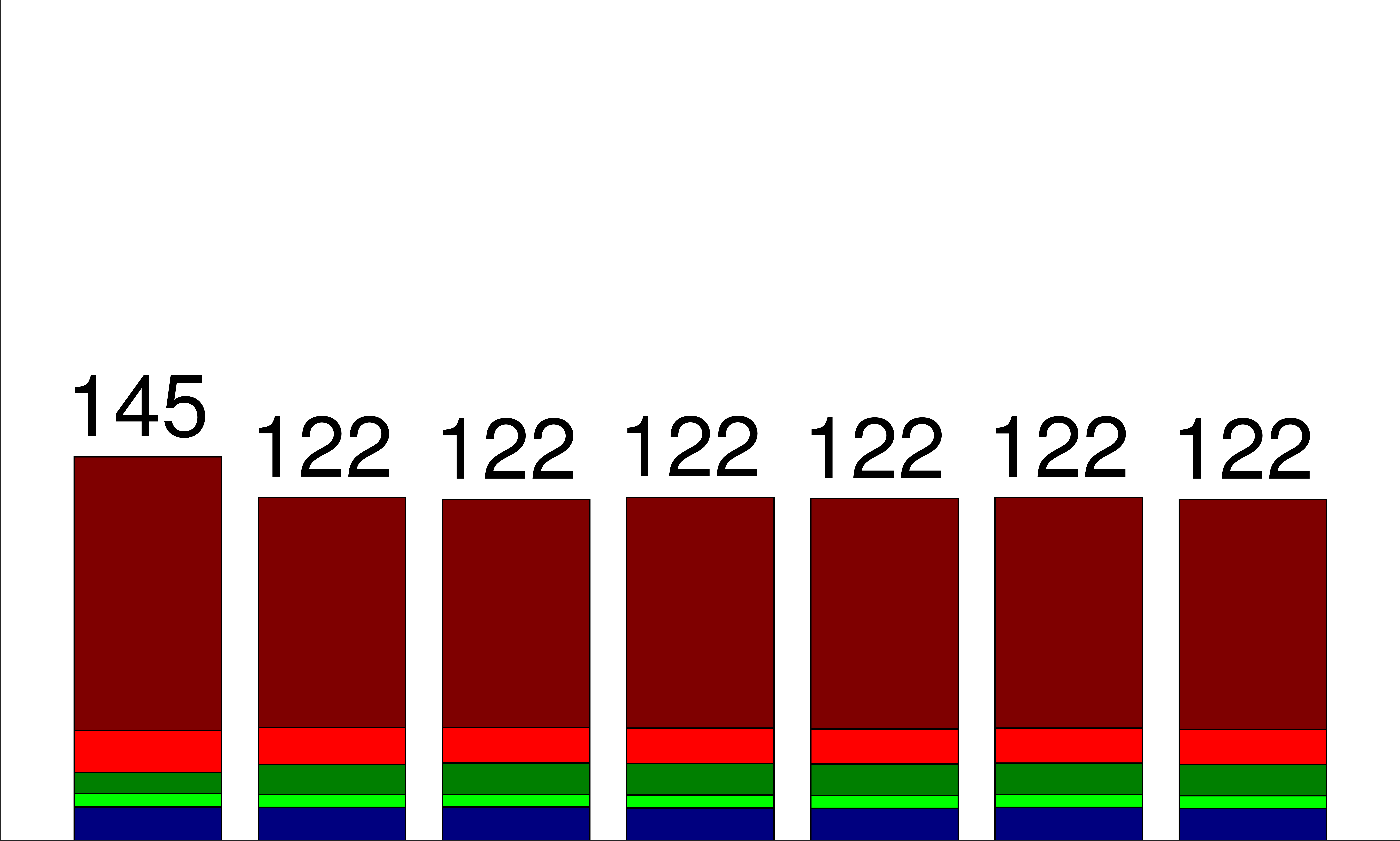}}\label{fig:cpu_trunc_decoupledAMS}}\hspace{0.7cm}%
\begin{picture}(0.5,1.5)%
		\put(0.2,-0.1){0}%
		\put(0.2,1.7){8}%
		\put(0.0,3.3){16}%
\end{picture}%
\subfigure{\fbox{\includegraphics[width=0.41\textwidth,height=0.25\textwidth]{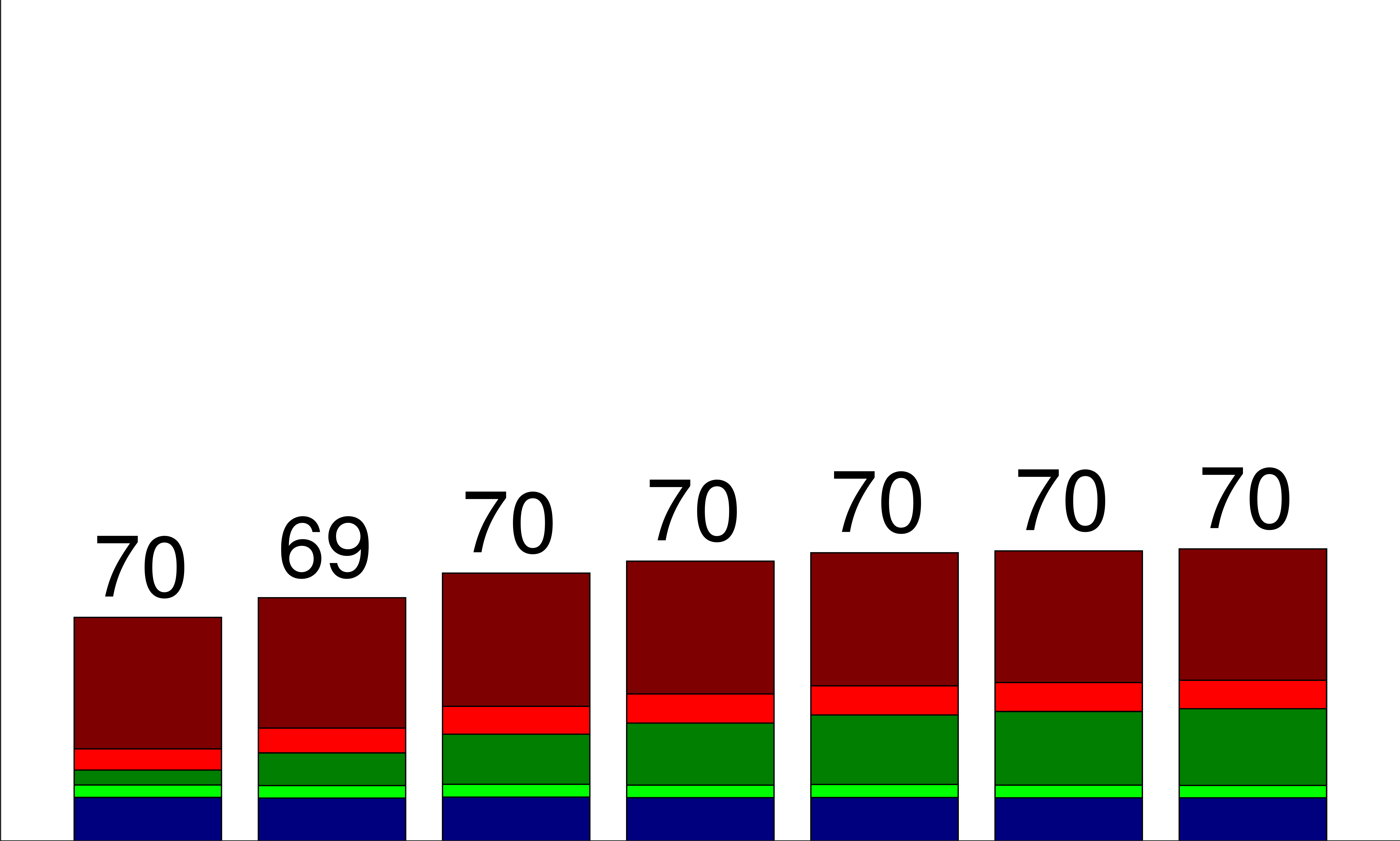}}\label{fig:cpu_trunc_fracAMS}}\\%
\begin{picture}(15,0.7)%
 \put(1.3,0.0){\rotatebox{45}{$10^{-1}$}}%
 \put(2.0,0.0){\rotatebox{45}{$10^{-2}$}}%
 \put(2.8,0.0){\rotatebox{45}{$10^{-3}$}}%
 \put(3.55,0.0){\rotatebox{45}{$10^{-4}$}}%
 \put(4.3,0.0){\rotatebox{45}{$10^{-5}$}}%
 \put(5.05,0.0){\rotatebox{45}{$10^{-6}$}}%
 \put(6.1,0.0){\rotatebox{45}{$0$}}%
  \put(8.2,0.0){\rotatebox{45}{$10^{-1}$}}%
 \put(8.9,0.0){\rotatebox{45}{$10^{-2}$}}%
 \put(9.7,0.0){\rotatebox{45}{$10^{-3}$}}%
 \put(10.45,0.0){\rotatebox{45}{$10^{-4}$}}%
 \put(11.2,0.0){\rotatebox{45}{$10^{-5}$}}%
 \put(11.95,0.0){\rotatebox{45}{$10^{-6}$}}%
 \put(13.0,0.0){\rotatebox{45}{$0$}}%
  \end{picture}\\%
  \begin{picture}(15,0.2)%
   \put(4.0,0.0){$\alpha$}%
   \put(11.0,0.0){$\alpha$}%
   \end{picture}\\%
\begin{picture}(15,0.7)%
  \put(3.0,0){\large Rock-AMS}%
  \put(9.5,0){\large Coupled-AMS}%
  \end{picture}\\%
  \begin{picture}(1.2,1.5)%
	\put(0.0,0.3){\rotatebox{90}{CPU time (sec)}}%
		\put(0.9,-0.1){0}%
		\put(0.9,1.7){8}%
		\put(0.7,3.3){16}%
\end{picture}%
\subfigure{\fbox{\includegraphics[width=0.41\textwidth,height=0.25\textwidth]{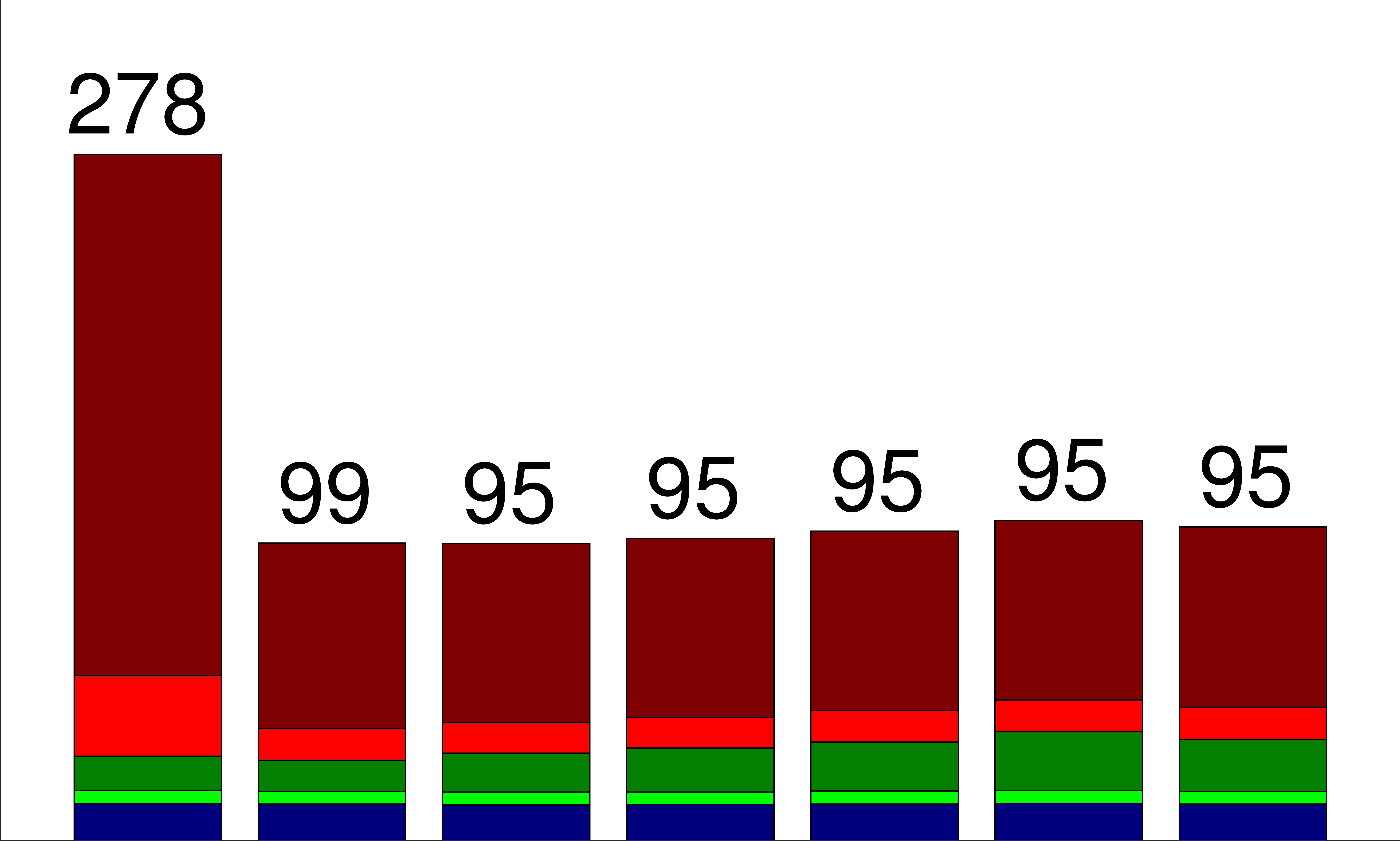}}\label{fig:cpu_trunc_rockAMS}}\hspace{0.7cm}%
  \begin{picture}(0.5,1.5)%
		\put(0.2,-0.1){0}%
		\put(0.2,1.7){8}%
		\put(0.0,3.3){16}%
\end{picture}%
\subfigure{\fbox{\includegraphics[width=0.41\textwidth,height=0.25\textwidth]{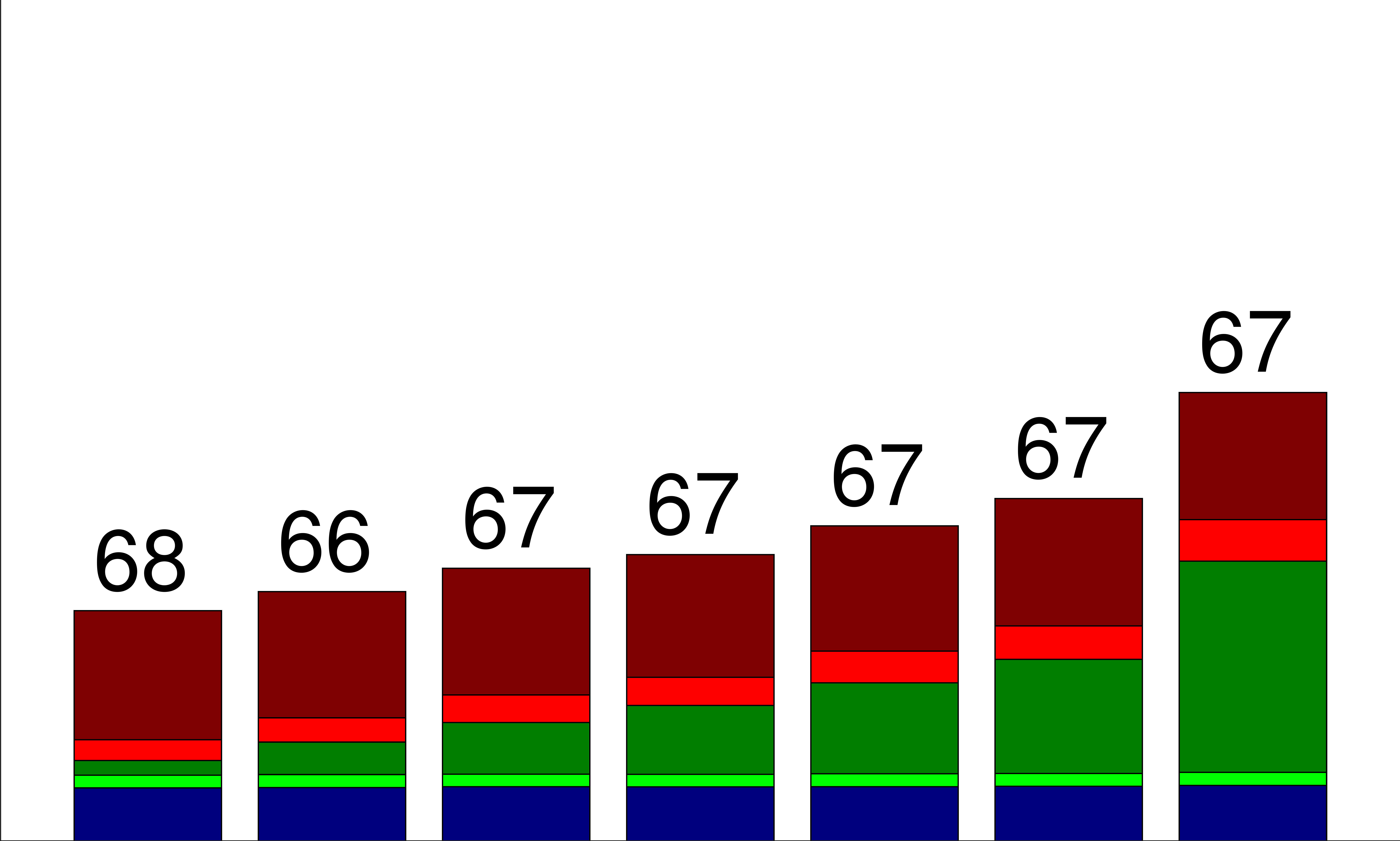}}\label{fig:cpu_trunc_coupledAMS}}\\%
\begin{picture}(15,0.7)%
 \put(1.3,0.0){\rotatebox{45}{$10^{-1}$}}%
 \put(2.0,0.0){\rotatebox{45}{$10^{-2}$}}%
 \put(2.8,0.0){\rotatebox{45}{$10^{-3}$}}%
 \put(3.55,0.0){\rotatebox{45}{$10^{-4}$}}%
 \put(4.3,0.0){\rotatebox{45}{$10^{-5}$}}%
 \put(5.05,0.0){\rotatebox{45}{$10^{-6}$}}%
 \put(6.1,0.0){\rotatebox{45}{$0$}}%
  \put(8.2,0.0){\rotatebox{45}{$10^{-1}$}}%
 \put(8.9,0.0){\rotatebox{45}{$10^{-2}$}}%
 \put(9.7,0.0){\rotatebox{45}{$10^{-3}$}}%
 \put(10.45,0.0){\rotatebox{45}{$10^{-4}$}}%
 \put(11.2,0.0){\rotatebox{45}{$10^{-5}$}}%
 \put(11.95,0.0){\rotatebox{45}{$10^{-6}$}}%
 \put(13.0,0.0){\rotatebox{45}{$0$}}%
  \end{picture}\\%
  \begin{picture}(15,0.2)%
   \put(4.0,0.0){$\alpha$}%
   \put(11.0,0.0){$\alpha$}%
   \end{picture}\\%
\subfigure{\fbox{\includegraphics[width=\textwidth]{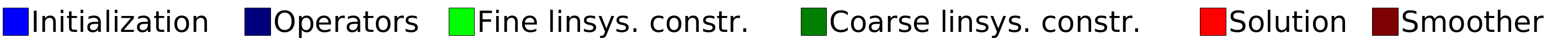}}}%
\caption{F-AMS performance for different values of the basis function truncation threshold, $\alpha$. Choosing $\alpha = 0$ invalidates the truncation procedure. The number of performed iterations to reach $10^{-6}$ residual 2-norm is given on top of each bar.}%
\label{fig:trunc}%
\end{figure}

\begin{figure}[htb!]%
\centering%
\subfigure{\includegraphics[width=0.48\textwidth]{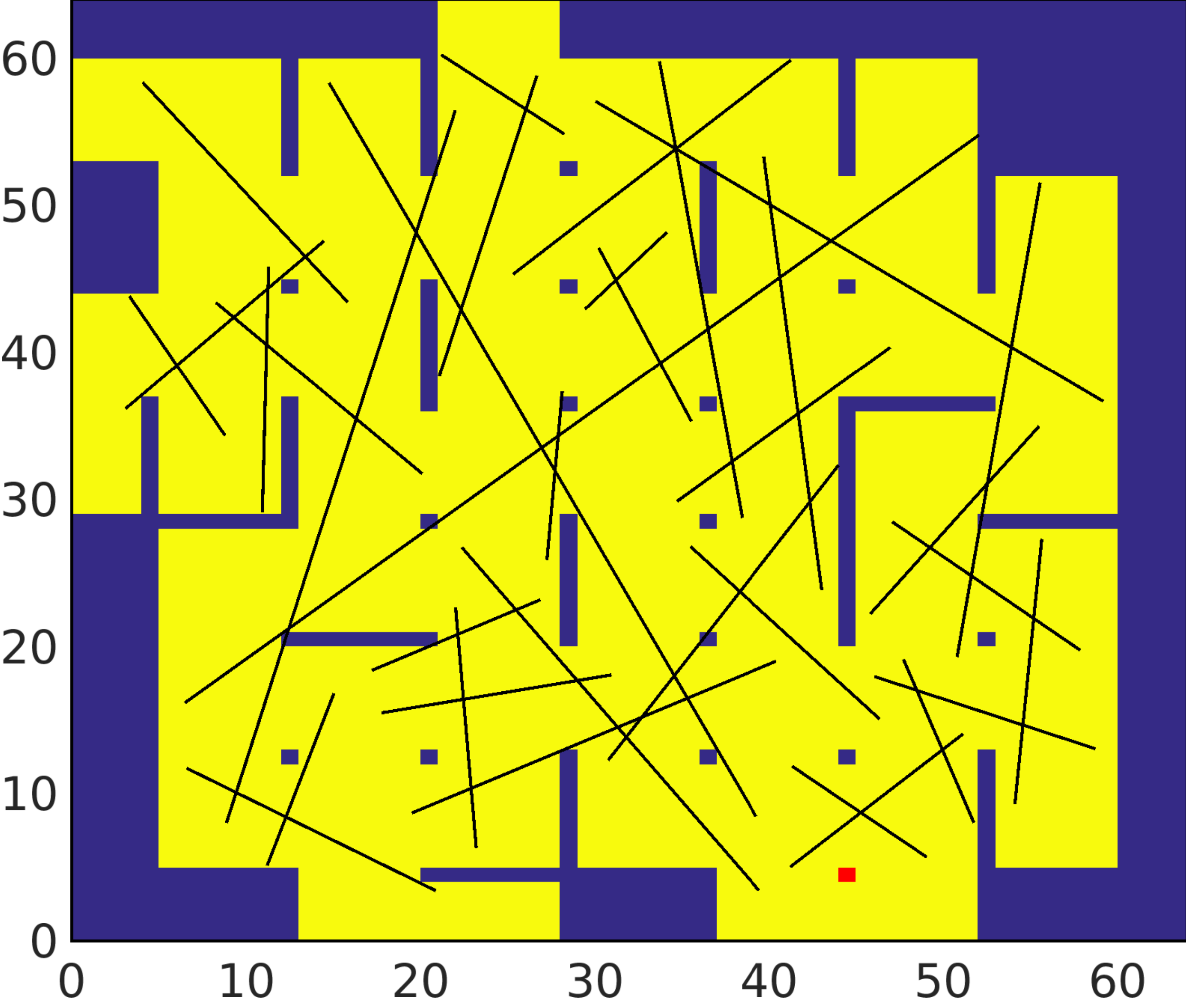}\label{fig:basis_support_untrunc}}\hspace{0.5cm}%
\subfigure{\includegraphics[width=0.48\textwidth]{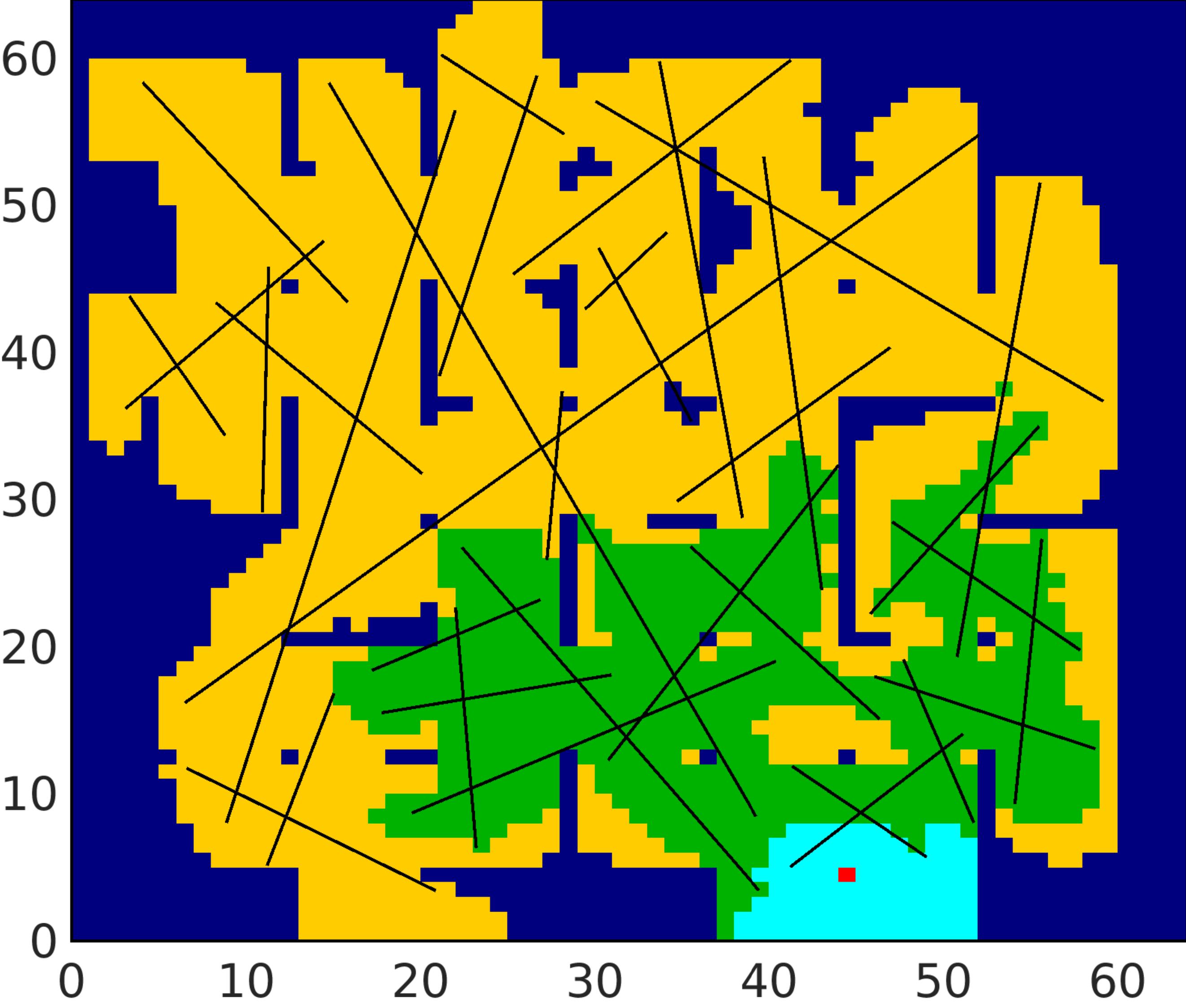}\label{fig:basis_support_trunc}}%
\caption{Coupled-AMS basis function support before (left) and after truncation (right). The coarse node (red) is located in a low permeable region of the 2D test case. The fracture network contains only 1 coarse-scale DOF. The colors on the right plot correspond to $\alpha = 10^{-5}$ (orange), $10^{-3}$ (green) and $10^{-2}$ (light blue).}%
\label{fig:basis_support}%
\end{figure}

One can limit the density of $\bm{\mathcal{P}}$ by truncating basis function values below a specified threshold, $\alpha \in [0,1)$. However, in order to preserve partition of unity, the affected rows in $\bm{\mathcal{P}}$ need to be rescaled by dividing the remaining values by the row sum. Figure~\ref{fig:trunc} shows the CPU time spent by F-AMS on the 3D test case, while varying the value of $\alpha$. Notice that the very restrictive value of $\alpha = 10^{-1}$ leads to an increase in the number of overall iterations, because the smoother needs to compensate for the induced inaccuracy of the basis functions. However, starting with $\alpha=10^{-2}$, the convergence is no longer much affected and the algorithm gains efficiency from the reduced number of FLOPS necessary to perform $\bm{\mathcal{R}} \bm A^\nu \bm{\mathcal{P}}$ and invert the result. The truncation has the biggest impact on the Coupled-AMS strategy, which experiences a speed-up factor of 2, compared to the un-truncated case (last bar in Fig.~\ref{fig:cpu_trunc_coupledAMS}, where $\alpha=0$). Figure~\ref{fig:basis_support} shows that, for this coupling strategy, when only a single DOF is considered for fracture network, the support of basis functions can be as big as the span of the fracture network. Also, this figure shows that after the truncation stage the locality of the basis function support can be maintained. As conclusion to this study, the subsequent experiments will use a value of $\alpha = 10^{-2}$, regardless of the choice of basis function coupling strategy.

\subsection{Sensitivity to the coarsening factor}
\label{sec:coarsening_ratio}

\begin{figure}[htb!]%
\centering%
\setlength{\fboxsep}{0pt}%
\setlength{\fboxrule}{0.2pt}%
\setlength{\unitlength}{1cm}%
\begin{picture}(15,1.5)%
  \put(5.4,1.3){\large Rock coarsening}%
	\put(1.9,0.65){$16 \times 16 \times 16$ \hspace{1.9cm} $8 \times 8 \times 8$ \hspace{2.2cm} $4 \times 4 \times 4$}%
  \put(2.15,0.1){(64 DOF) \hspace{2.0cm} (512 DOF) \hspace{1.8cm} (4096 DOF)}%
\end{picture}\\%
\begin{picture}(1.3,1.5)%
	\put(0,1.0){\rotatebox{90}{CPU time (sec)}}%
	\put(0.9,-0.1){0}%
	\put(0.75,1.1){25}%
	\put(0.75,2.3){50}%
	\put(0.75,3.5){75}%
	\put(0.6,4.7){100}%
\end{picture}%
\subfigure{\fbox{\includegraphics[width=0.28\textwidth,height=0.35\textwidth]{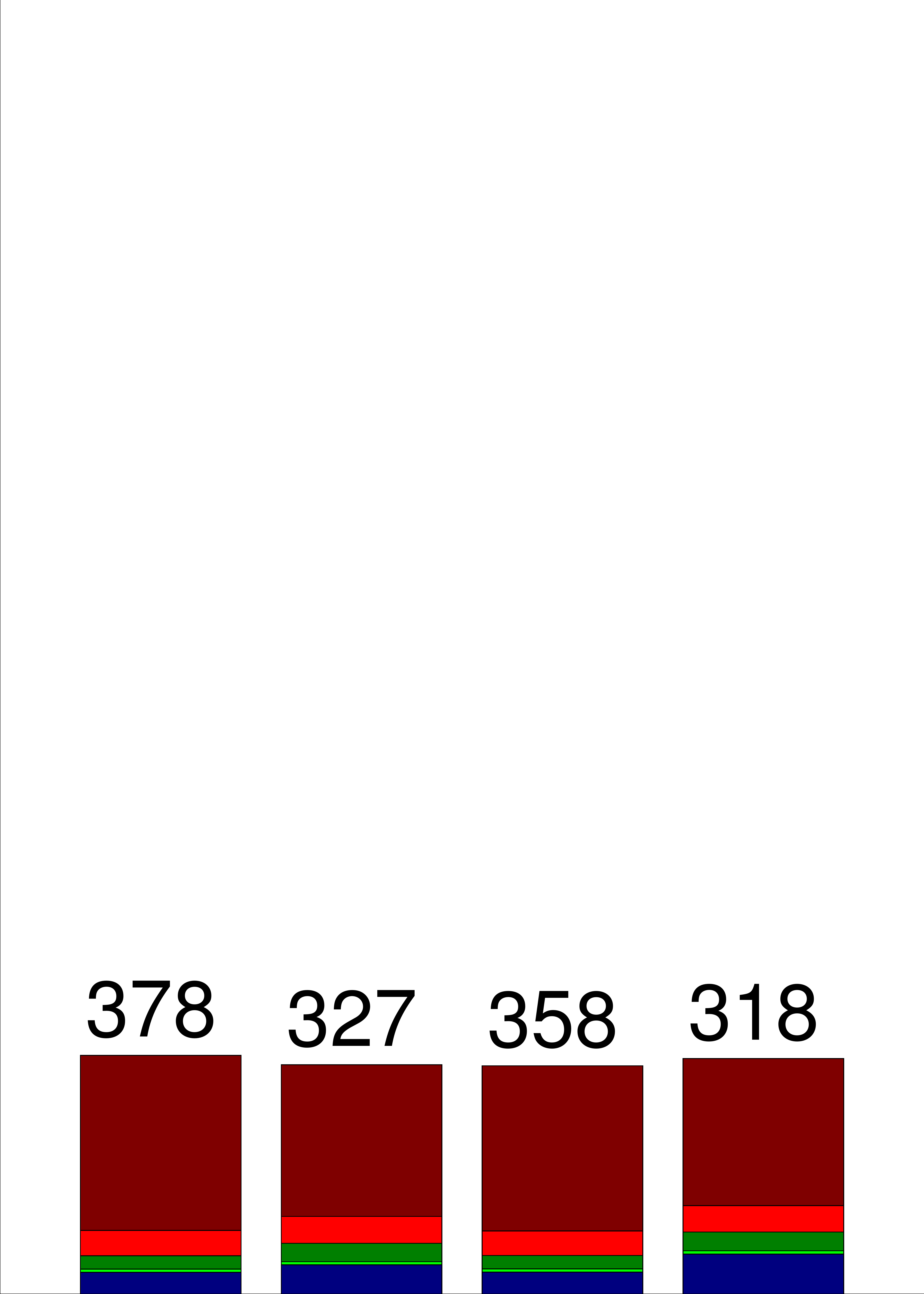}}%
\fbox{\includegraphics[width=0.28\textwidth,height=0.35\textwidth]{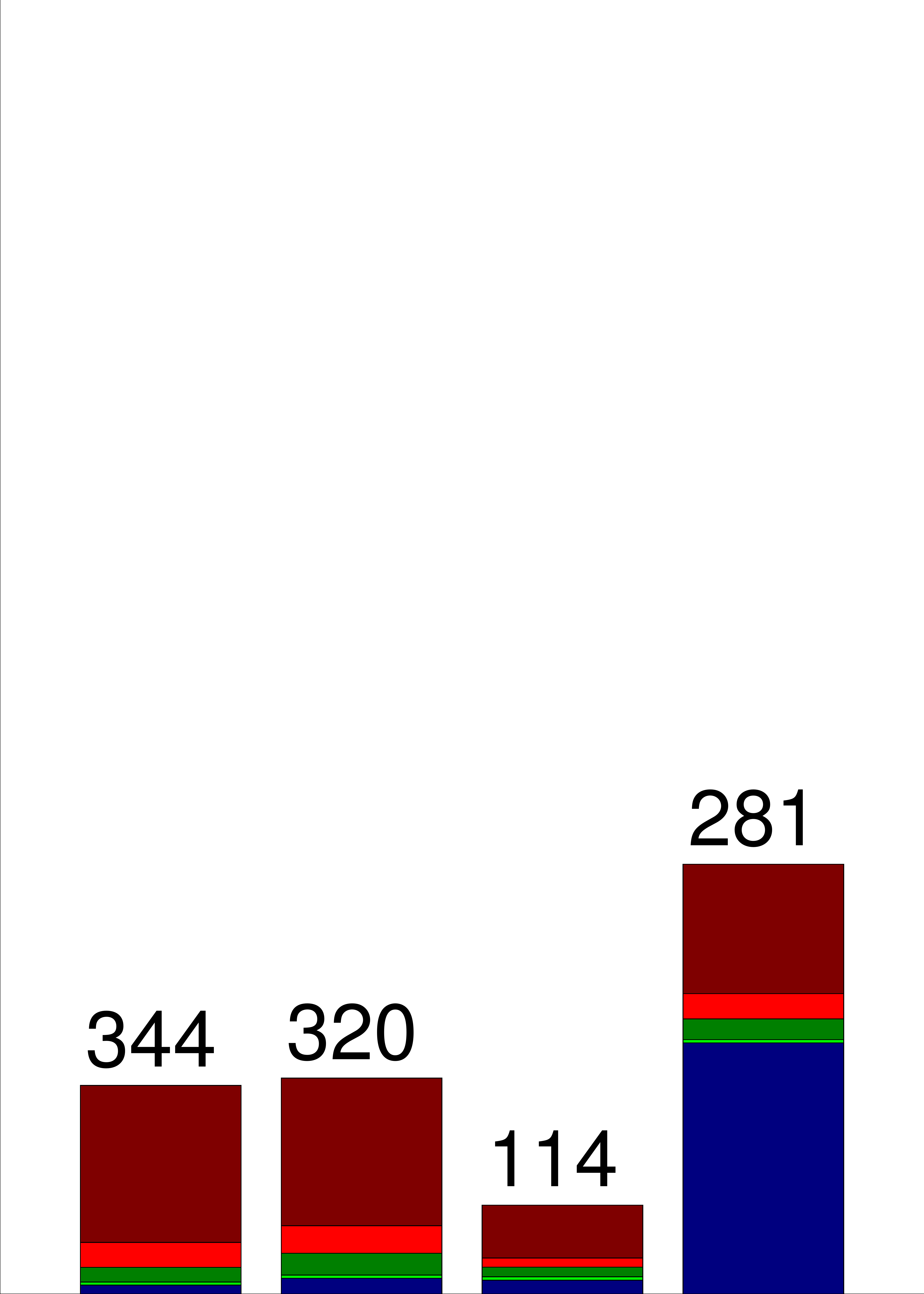}}%
\fbox{\includegraphics[width=0.28\textwidth,height=0.35\textwidth]{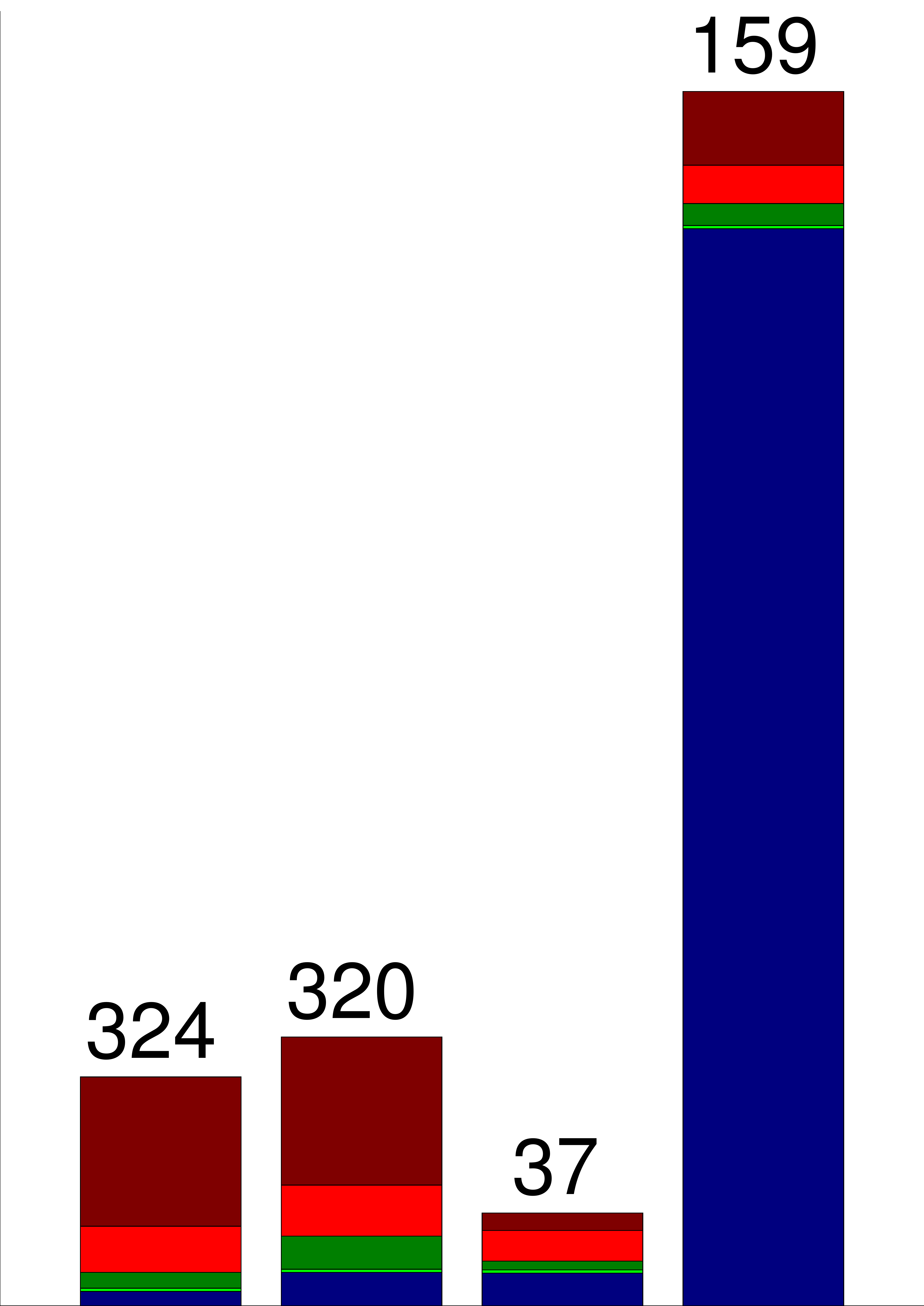}}}%
\begin{picture}(1.3,1.5)%
	\put(0.2,1.8){\rotatebox{90}{$16 \times 16$}}%
   \put(0.65,1.6){\rotatebox{90}{(76 DOF)}}%
\end{picture}\\%
\begin{picture}(1.3,1.5)%
	\put(0,1.0){\rotatebox{90}{CPU time (sec)}}%
	\put(0.9,-0.1){0}%
	\put(0.75,1.1){25}%
	\put(0.75,2.3){50}%
	\put(0.75,3.5){75}%
	\put(0.6,4.7){100}%
\end{picture}%
\subfigure{\fbox{\includegraphics[width=0.28\textwidth,height=0.35\textwidth]{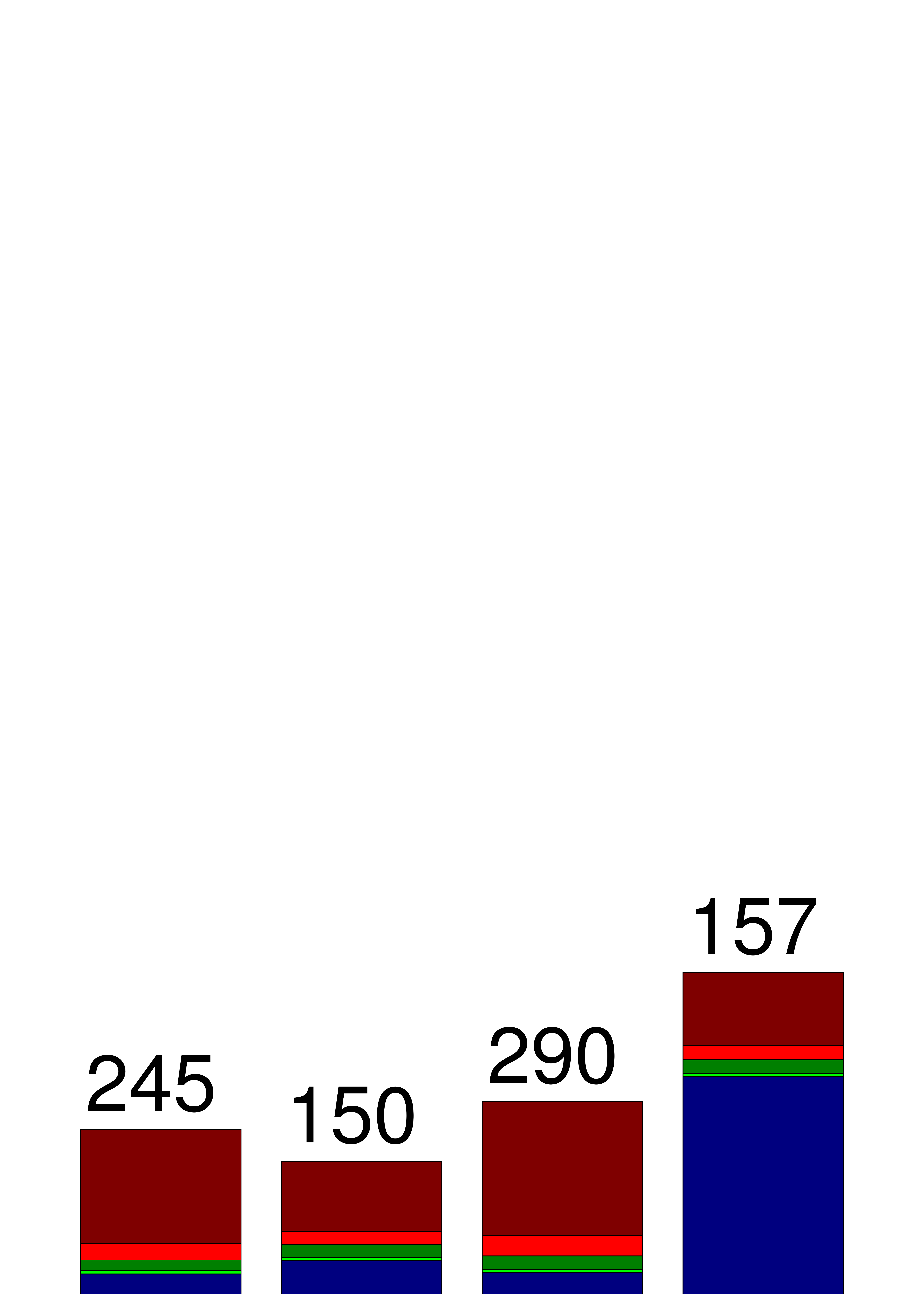}}%
\fbox{\includegraphics[width=0.28\textwidth,height=0.35\textwidth]{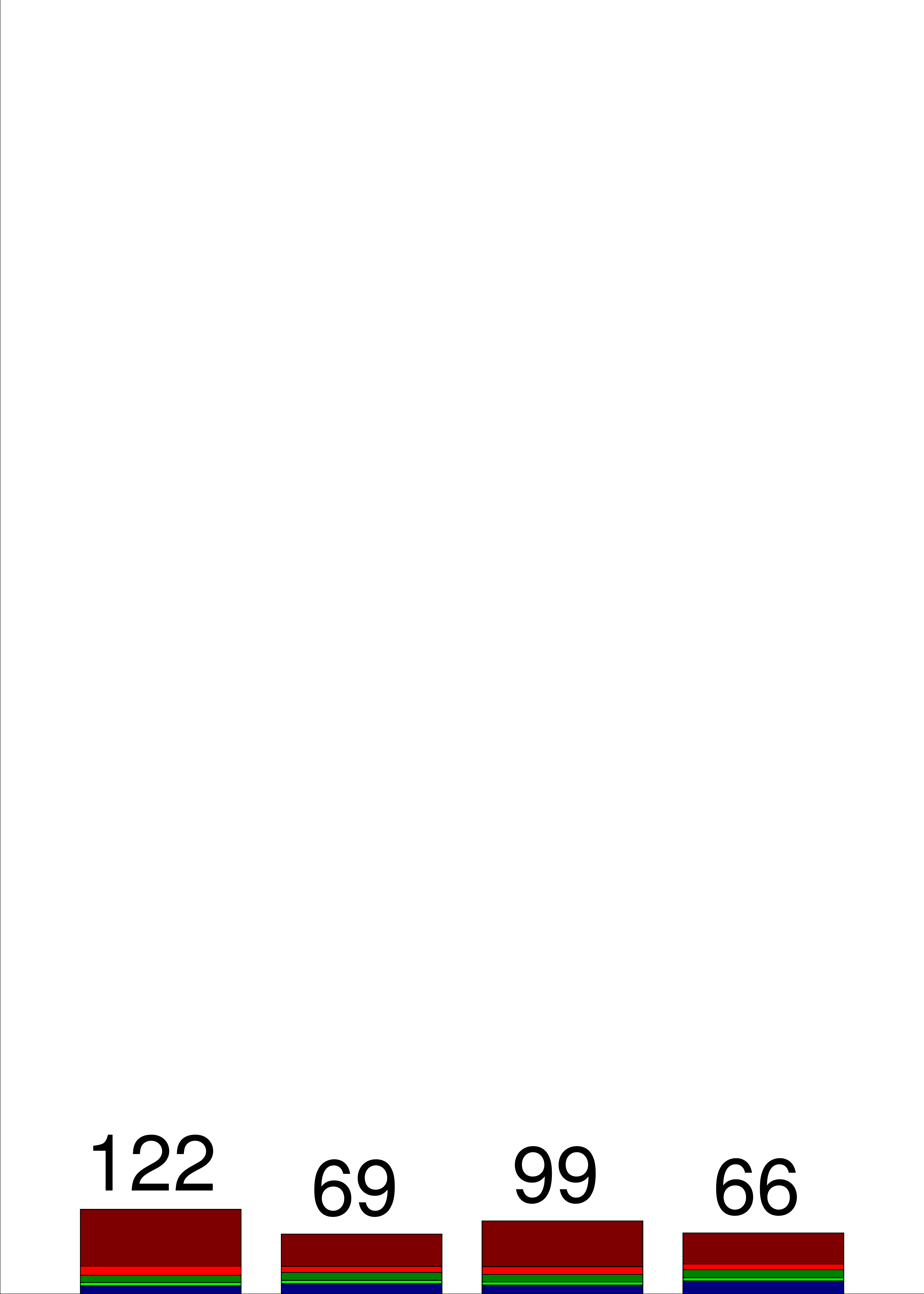}}%
\fbox{\includegraphics[width=0.28\textwidth,height=0.35\textwidth]{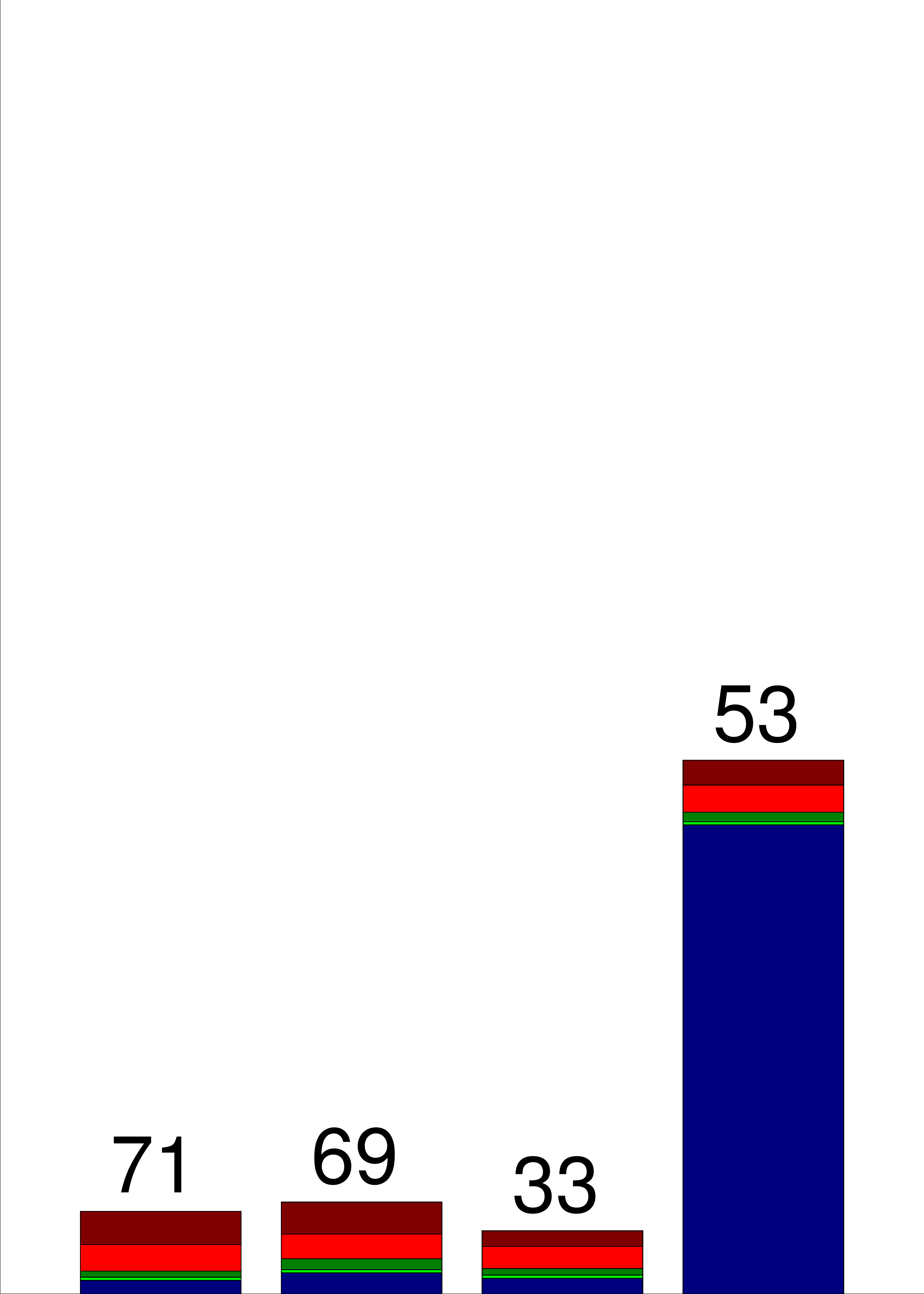}}}%
\begin{picture}(1.3,3.0)%
	\put(0.2,2.0){\rotatebox{90}{$8 \times 8$}}%
   \put(0.65,1.5){\rotatebox{90}{(368 DOF)}}%
	\put(1.3,0.2){\rotatebox{90}{\large Fracture coarsening}}%
\end{picture}\\%
\begin{picture}(1.3,1.5)%
	\put(0,1.0){\rotatebox{90}{CPU time (sec)}}%
	\put(0.9,-0.1){0}%
	\put(0.75,1.1){25}%
	\put(0.75,2.3){50}%
	\put(0.75,3.5){75}%
	\put(0.6,4.7){100}%
\end{picture}%
\subfigure{\fbox{\includegraphics[width=0.28\textwidth,height=0.35\textwidth]{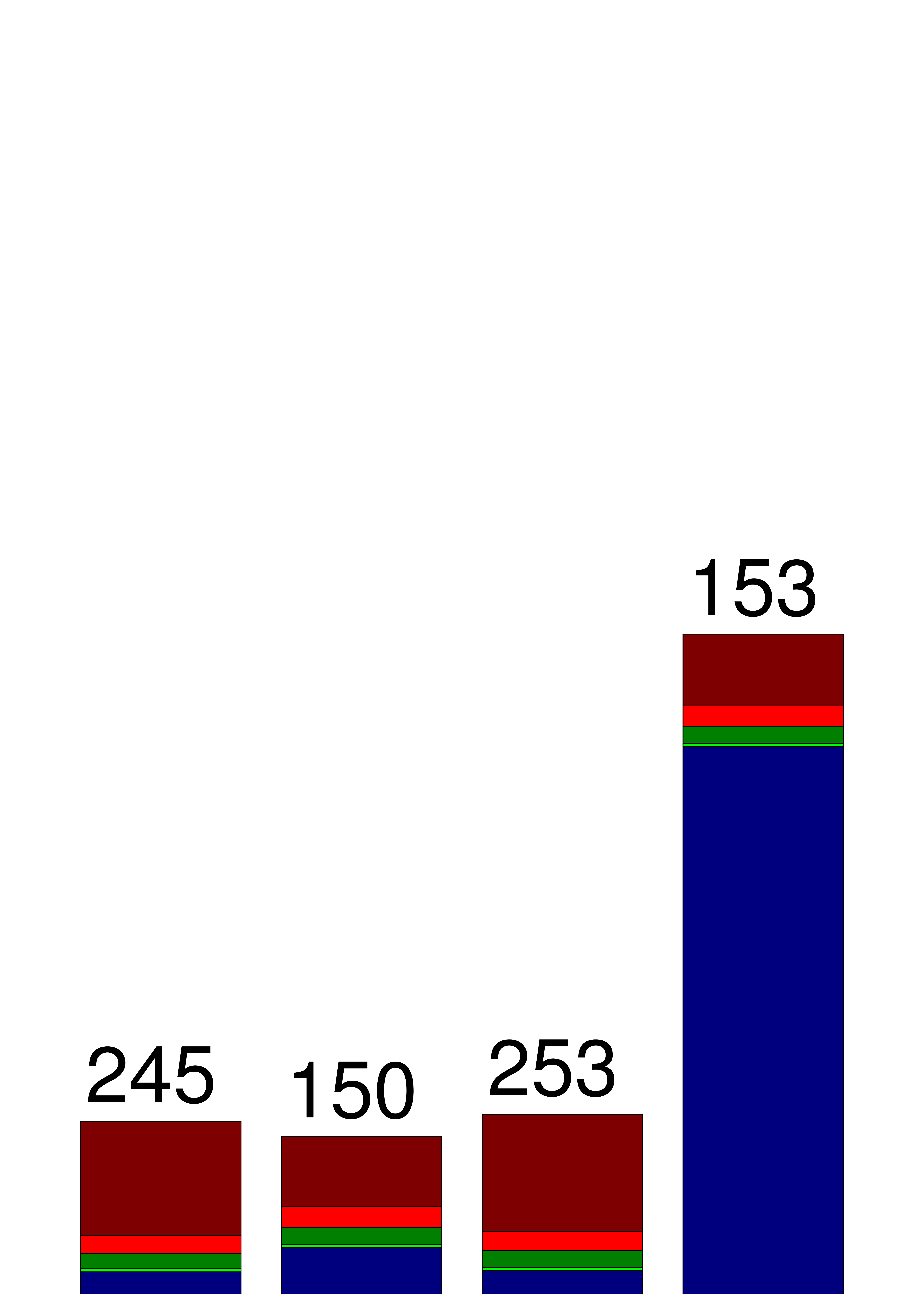}}%
\fbox{\includegraphics[width=0.28\textwidth,height=0.35\textwidth]{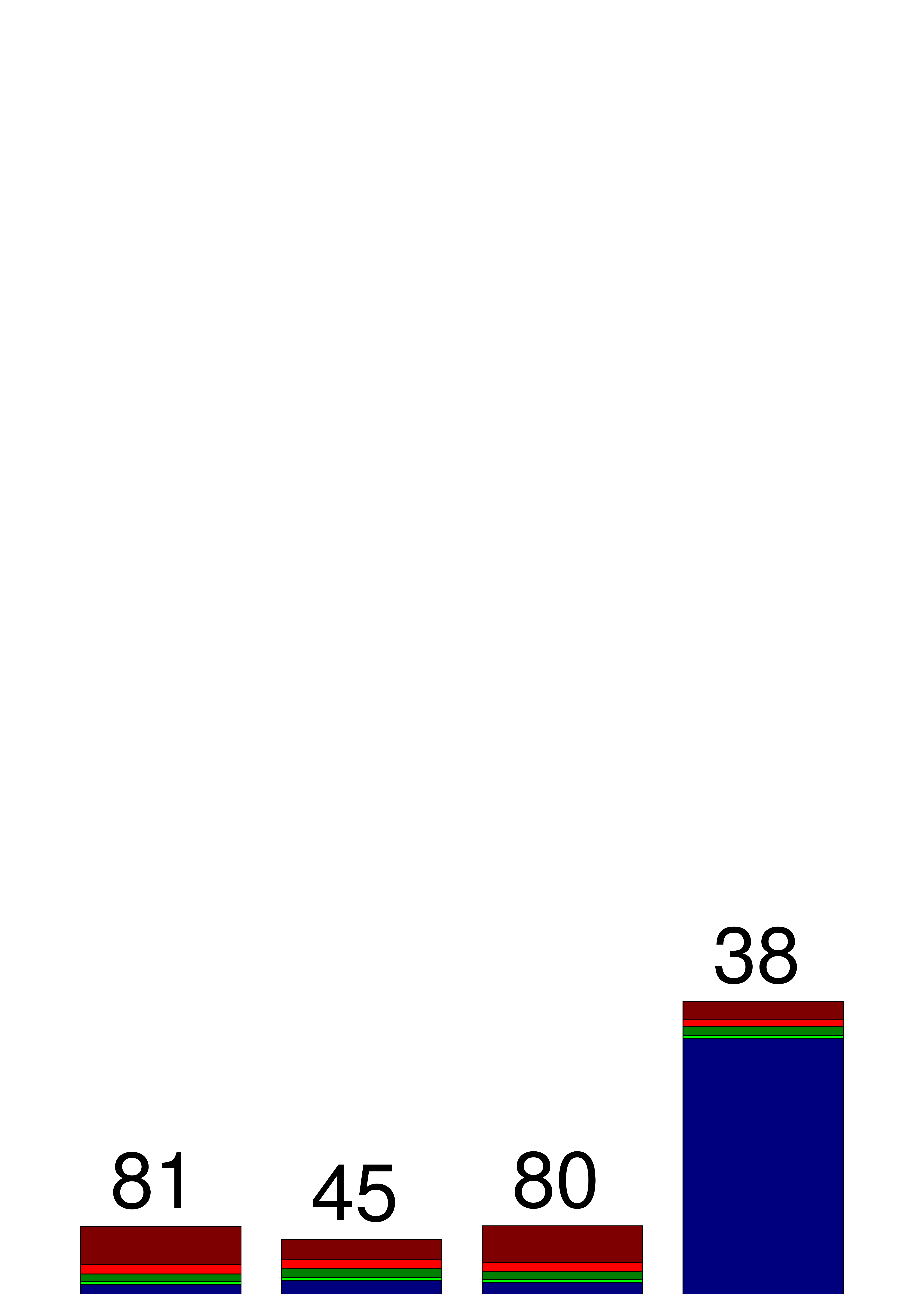}}%
\fbox{\includegraphics[width=0.28\textwidth,height=0.35\textwidth]{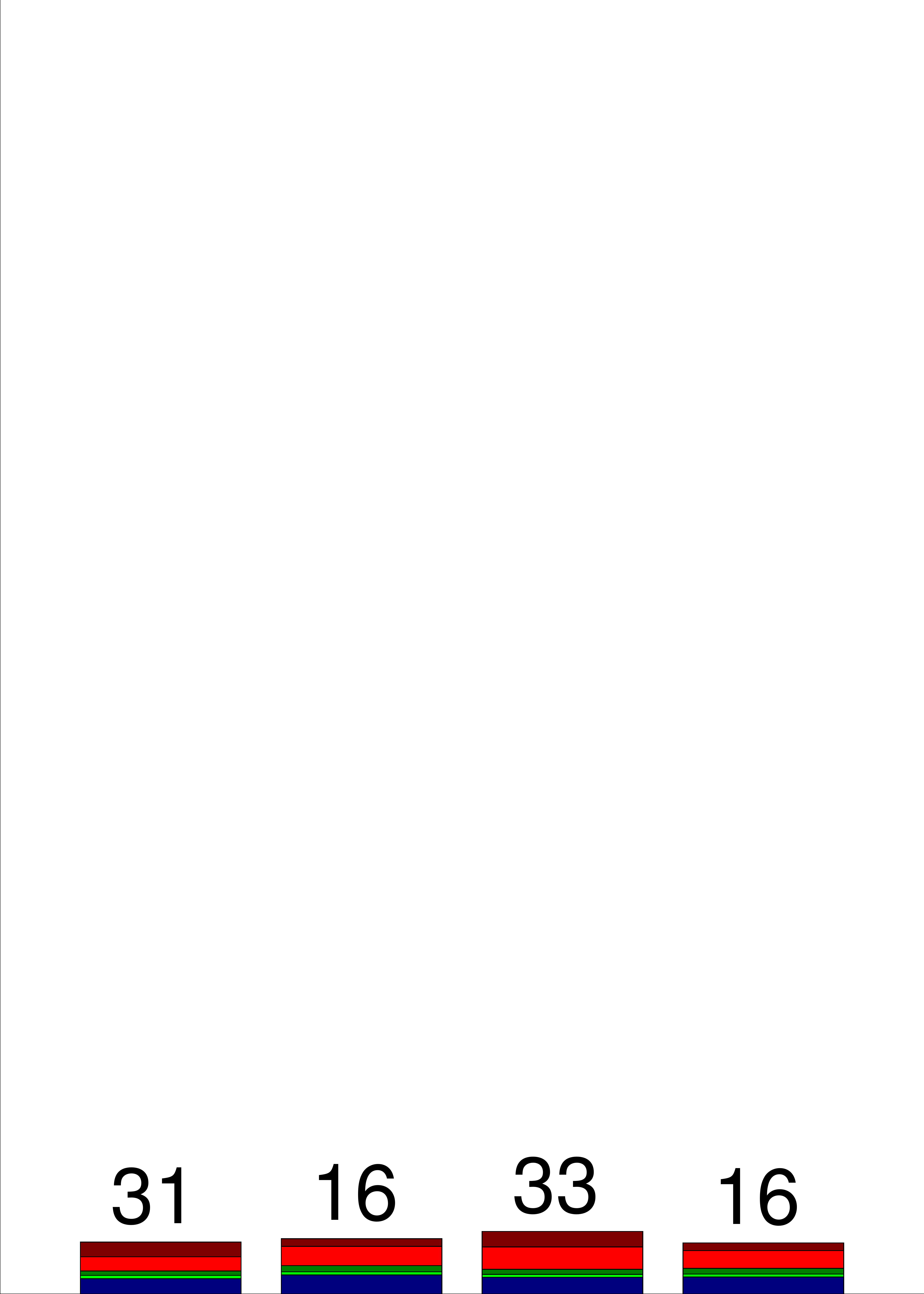}}}%
\begin{picture}(1.3,1.5)%
	\put(0.2,2.0){\rotatebox{90}{$4 \times 4$}}%
   \put(0.65,1.4){\rotatebox{90}{(1968 DOF)}}%
\end{picture}\\%
\begin{picture}(15,2.3)%
	\put(0,0.0){\rotatebox{45}{Decoupled-AMS}}%
	\put(0.8,0.0){\rotatebox{45}{\hspace{1.1cm}Frac-AMS}}%
	\put(1.7,0.0){\rotatebox{45}{\hspace{1.0cm}Rock-AMS}}%
	\put(2.5,0.0){\rotatebox{45}{\hspace{0.4cm}Coupled-AMS}}%
	\put(3.9,0.0){\rotatebox{45}{Decoupled-AMS}}%
	\put(4.7,0.0){\rotatebox{45}{\hspace{1.1cm}Frac-AMS}}%
	\put(5.6,0.0){\rotatebox{45}{\hspace{1.0cm}Rock-AMS}}%
	\put(6.4,0.0){\rotatebox{45}{\hspace{0.4cm}Coupled-AMS}}%
	\put(7.8,0.0){\rotatebox{45}{Decoupled-AMS}}%
	\put(8.6,0.0){\rotatebox{45}{\hspace{1.1cm}Frac-AMS}}%
	\put(9.4,0.0){\rotatebox{45}{\hspace{1.0cm}Rock-AMS}}%
	\put(10.3,0.0){\rotatebox{45}{\hspace{0.4cm}Coupled-AMS}}%
\end{picture}\\%
\subfigure{\fbox{\includegraphics[width=\textwidth]{Fig/method/legend.png}}}%
\caption{F-AMS performance for different matrix (decreasing from left to right) and fracture coarsening ratios (decreasing from top to bottom). The number of performed iterations to reach $10^{-6}$ residual 2-norm is given on top of each bar.}%
\label{fig:dof}%
\end{figure}%

The sensitivity of F-AMS to the number of coarse DOF in the fracture network, as well as the matrix coarsening ratio is studied for the 3D test case shown in Fig.~\ref{fig:3dcase}. The coarsening factor is defined as the average number of fine cells contained in one (matrix or fracture) primal-coarse block, along each axis. Recall that, along the fracture length, this is given by the $d_{min}$ parameter in Table~\ref{fractureCoarsening}. The experiment designed for this purpose is focused on ``isotropic'' coarsening factors, mainly due to the point-wise nature of ILU(0), which was chosen as global smoothing stage for the implementation of F-AMS used here (Table~\ref{fams-algorithm}).

Figure~\ref{fig:dof} shows the F-AMS CPU times obtained with three different coarsening factors for the matrix, as well as the fracture network. From this figure, the Coupled-AMS  is found to automatically adapt itself to the coarsening ratio. Its convergence rate is surpassed only in cases where there is a large discrepancy between the rock and fracture coarsening ratios. However, this comes with the additional computational cost of having basis functions with wider support. Still, due to the truncation factor $\alpha = 10^{-2}$, the construction and solution of coarse system for Coupled-AMS is not significantly higher than that of the alternative strategies. In addition, in each row, the optimum simulation results are obtained when fracture and matrix coarsening ratios are similar.

Based on this study and unless otherwise stated, the coarsening ratio of $8$ in each direction for both matrix and fracture media is employed in the experiments presented in the following subsections. Note that this option leads to more efficient coarse-scale systems than the alternative option of using coarsening ratios of $4$. In addition, the number of linear iterations can be significantly reduced when F-AMS is employed as preconditioner for GMRES \cite{Saad}.

\subsection{Sensitivity to the transmissibility ratio}

\begin{figure}[htb!]%
\centering%
\setlength{\fboxsep}{0pt}%
\setlength{\fboxrule}{0.2pt}%
\setlength{\unitlength}{1cm}%
\begin{picture}(15,0.2)%
  \put(2.6,0){\large Decoupled-AMS}%
  \put(9.9,0){\large Frac-AMS}%
  \end{picture}\\%
  \begin{picture}(1.2,1.5)%
	\put(0.0,0.3){\rotatebox{90}{CPU time (sec)}}%
		\put(0.9,-0.1){0}%
		\put(0.9,1.7){6}%
		\put(0.7,3.3){12}%
\end{picture}%
\subfigure{\fbox{\includegraphics[width=0.41\textwidth,height=0.25\textwidth]{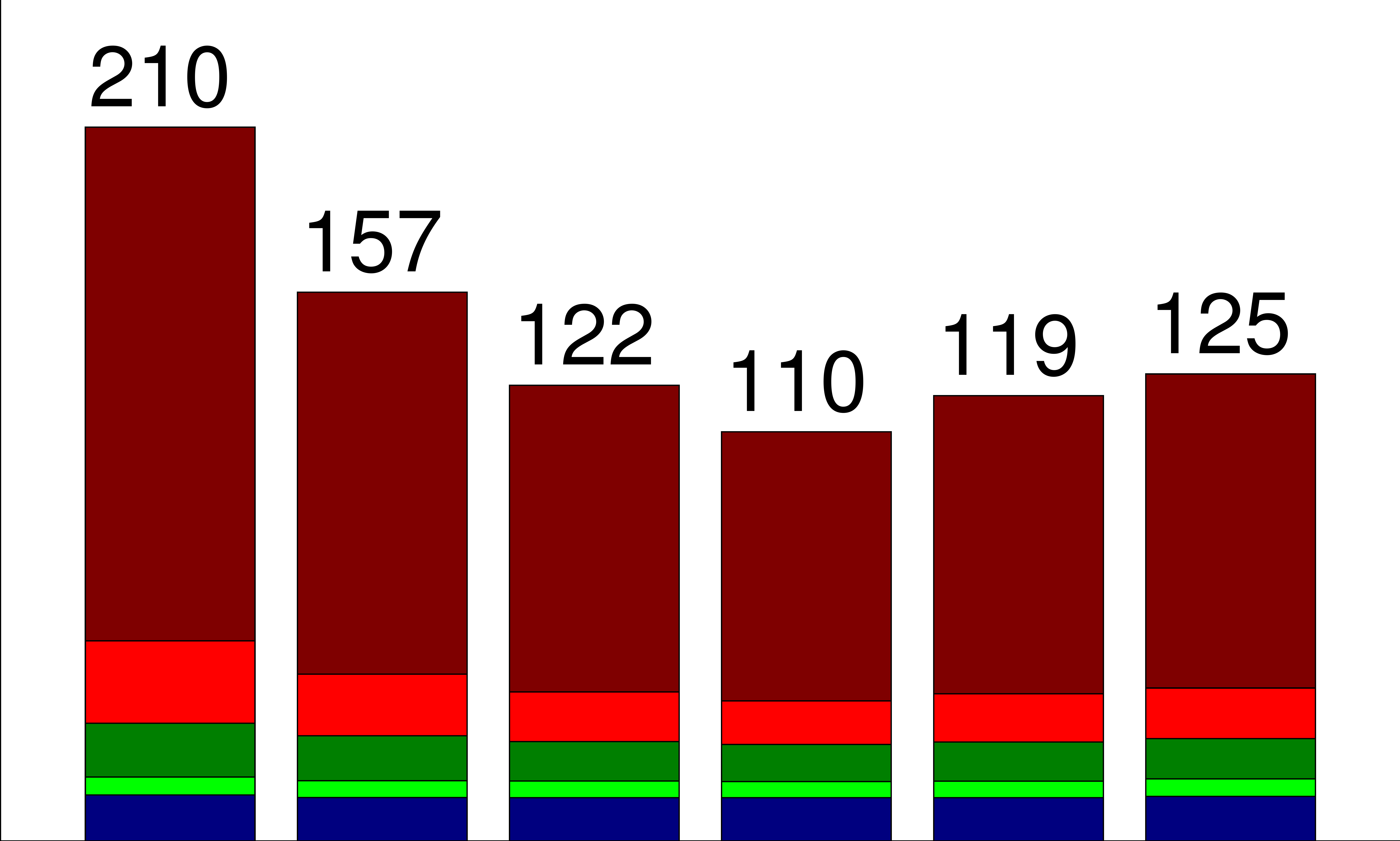}}\label{fig:cpu_transm_decoupledAMS}}\hspace{0.7cm}%
\begin{picture}(0.5,1.5)%
		\put(0.2,-0.1){0}%
		\put(0.2,1.7){6}%
		\put(0.0,3.3){12}%
\end{picture}%
\subfigure{\fbox{\includegraphics[width=0.41\textwidth,height=0.25\textwidth]{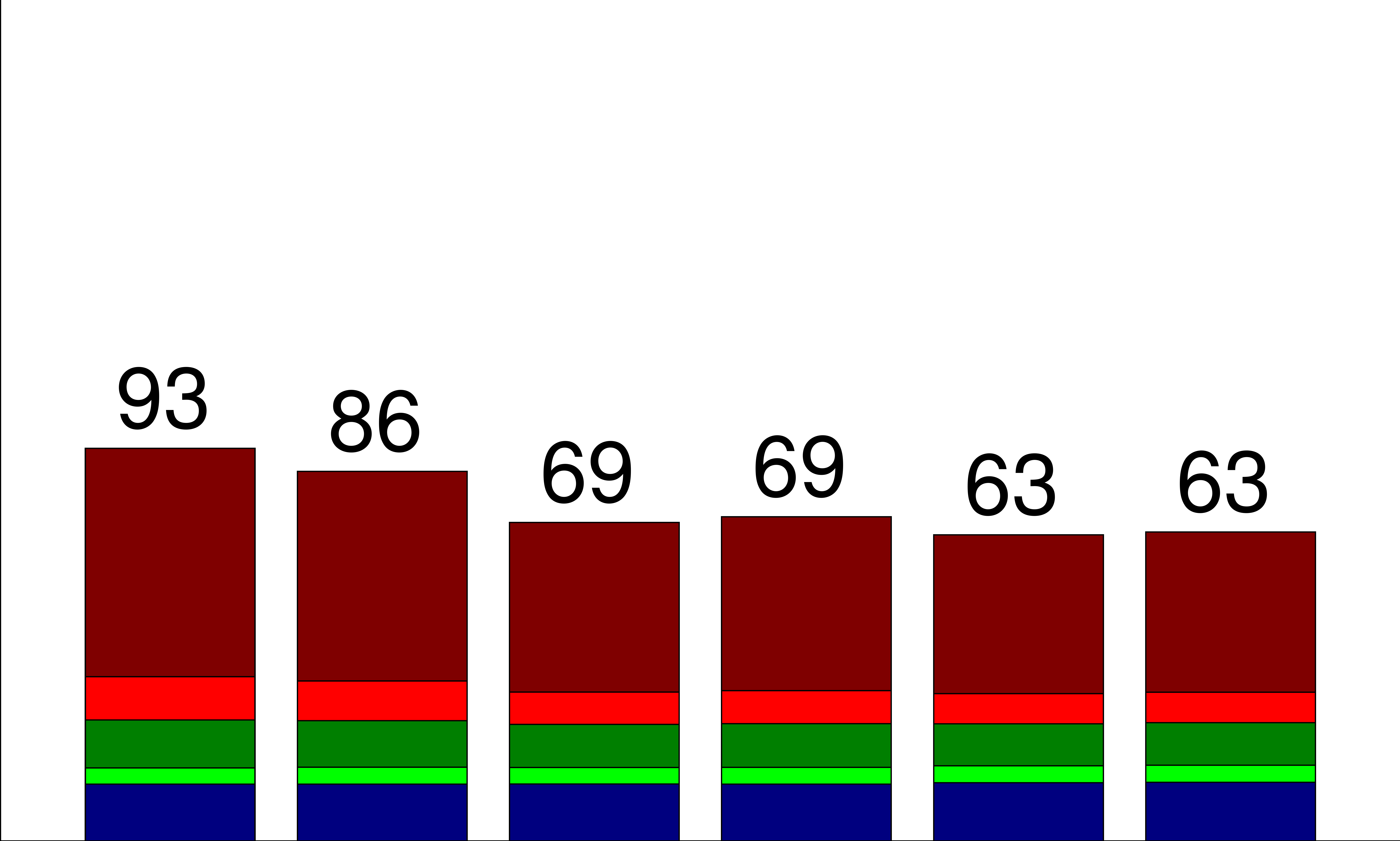}}\label{fig:cpu_transm_fracAMS}}\\%
\begin{picture}(15,0.2)%
 \put(1.7,0.0){$10^0$}%
 \put(2.6,0.0){$10^1$}%
 \put(3.45,0.0){$10^2$}%
 \put(4.3,0.0){$10^3$}%
 \put(5.15,0.0){$10^4$}%
 \put(6.05,0.0){$10^5$}%
  \put(8.6,0.0){$10^0$}%
 \put(9.5,0.0){$10^1$}%
 \put(10.35,0.0){$10^2$}%
 \put(11.2,0.0){$10^3$}%
 \put(12.05,0.0){$10^4$}%
 \put(12.95,0.0){$10^5$}%
  \end{picture}\\%
  \begin{picture}(15,0.2)%
   \put(3.8,.0){$T_{ratio}$}%
   \put(10.7,.0){$T_{ratio}$}%
   \end{picture}\\%
\begin{picture}(15,0.9)%
  \put(3.2,0){\large Rock-AMS}%
  \put(9.4,0){\large Coupled-AMS}%
  \end{picture}\\%
  \begin{picture}(1.2,1.5)%
	\put(0.0,0.3){\rotatebox{90}{CPU time (sec)}}%
		\put(0.9,-0.1){0}%
		\put(0.9,1.7){6}%
		\put(0.7,3.3){12}%
\end{picture}%
\subfigure{\fbox{\includegraphics[width=0.41\textwidth,height=0.25\textwidth]{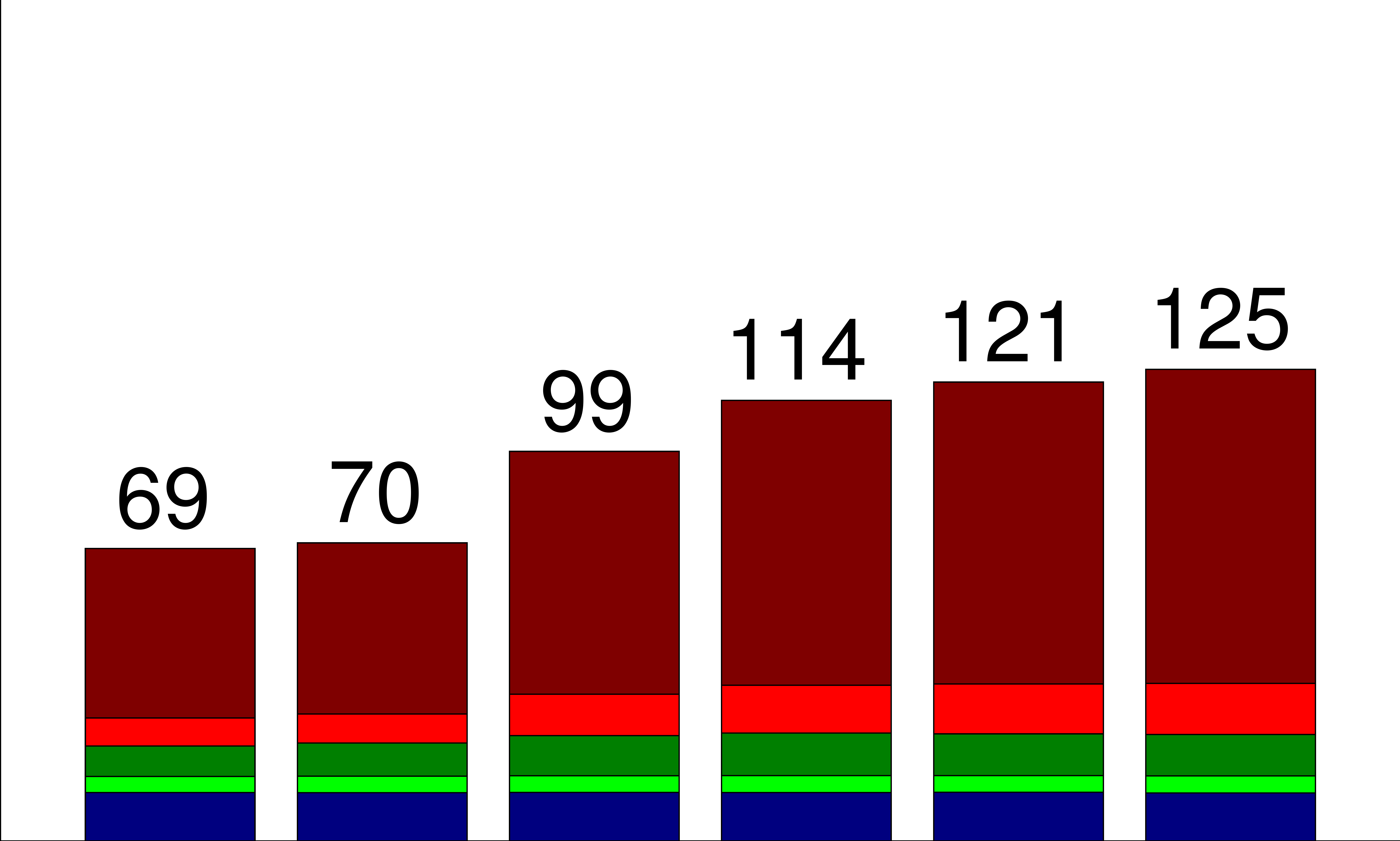}}\label{fig:cpu_transm_rockAMS}}\hspace{0.7cm}%
\begin{picture}(0.5,1.5)%
		\put(0.2,-0.1){0}%
		\put(0.2,1.7){6}%
		\put(0.0,3.3){12}%
\end{picture}%
\subfigure{\fbox{\includegraphics[width=0.41\textwidth,height=0.25\textwidth]{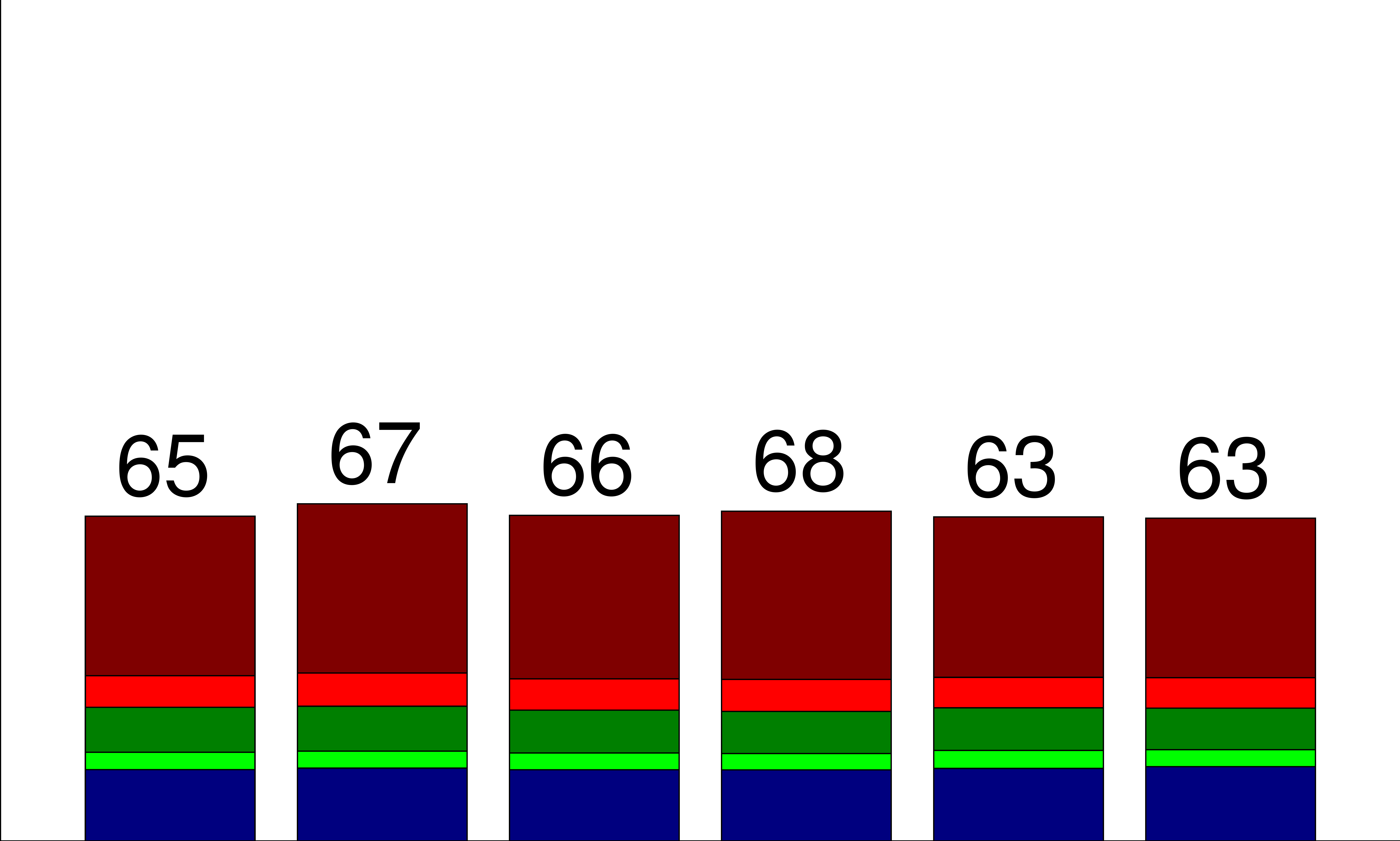}}\label{fig:cpu_transm_coupledAMS}}\\%
\begin{picture}(15,0.2)%
 \put(1.7,0.0){$10^0$}%
 \put(2.6,0.0){$10^1$}%
 \put(3.45,0.0){$10^2$}%
 \put(4.3,0.0){$10^3$}%
 \put(5.15,0.0){$10^4$}%
 \put(6.05,0.0){$10^5$}%
  \put(8.6,0.0){$10^0$}%
 \put(9.5,0.0){$10^1$}%
 \put(10.35,0.0){$10^2$}%
 \put(11.2,0.0){$10^3$}%
 \put(12.05,0.0){$10^4$}%
 \put(12.95,0.0){$10^5$}%
  \end{picture}\\%
  \begin{picture}(15,0.2)%
   \put(3.8,.0){$T_{ratio}$}%
   \put(10.7,.0){$T_{ratio}$}%
   \end{picture}\\%
\subfigure{\fbox{\includegraphics[width=\textwidth]{Fig/method/legend.png}}}%
\caption{F-AMS performance for different matrix-fracture transmissibility ratios. The number of performed iterations to reach $10^{-6}$ residual 2-norm is given on top of each bar.}%
\label{fig:transm_ratio}%
\end{figure}%

\begin{figure}[htb!]%
\centering%
\includegraphics[width=0.8\textwidth,height=0.4\textwidth]{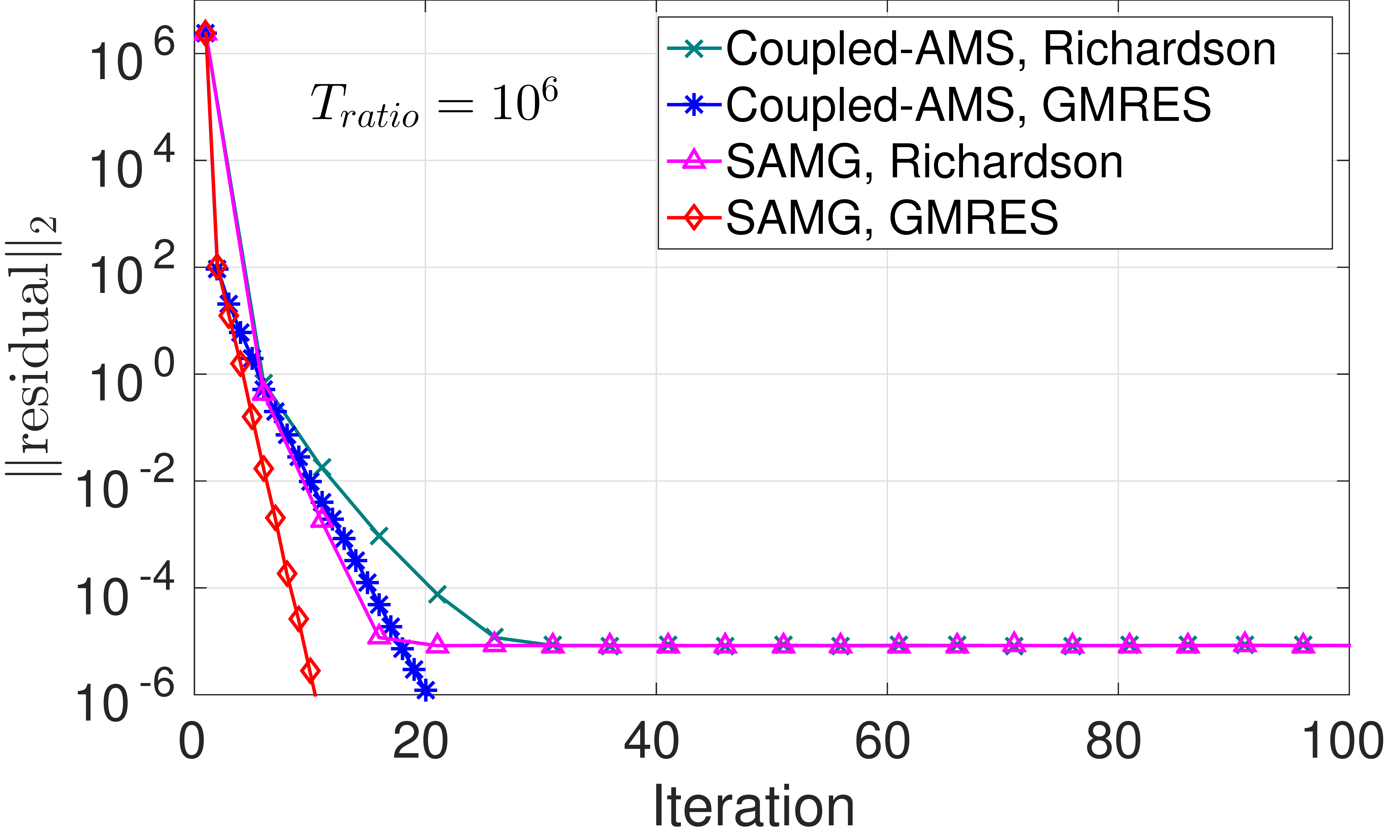}%
\caption{Convergence history of F-AMS and SAMG on the 3D test case with a fracture/matrix transmissibility ratio of $T_{ratio}=10^6$. Notice that neither method can converge when iterated in a Richardson's loop. Instead,  as preconditioners to GMRES, both methods converge after a few iterations.}%
\label{fig:conv_perm1e6}%
\end{figure}%

The next set of experiments aims to investigate the sensitivity of F-AMS to the conductivity contrast between the matrix and fractures. The transmissibility ratio $T_{ratio}$, as defined in Eq.~\eqref{tratio}, is varied over several orders of magnitude (Fig.~\ref{fig:transm_ratio}) while measuring CPU times and number of linear iterations performed by the Richardson loop. A coarsening factor of $8$ was chosen for both media, based on the results from the previous subsection. As the network becomes more conductive, the influence of the matrix heterogeneity on the fracture pressure decreases (see Figs.~\ref{fig:fracPres3d} and \ref{fig:pres3d}). As such, Rock-AMS exhibits a degradation in performance for higher $T_{ratio}$ (Fig.~\ref{fig:cpu_transm_rockAMS}), while the reverse is true for Frac-AMS (Fig.~\ref{fig:cpu_transm_fracAMS}). On the other hand, by automatically adapting to the change, the Coupled-AMS  strategy remains relatively insensitive to $T_{ratio}$, as shown in Fig.~\ref{fig:cpu_transm_coupledAMS}. Finally, Decoupled-AMS requires the most number of iterations when the fracture and matrix transmissibility values are close (Fig.~\ref{fig:cpu_transm_decoupledAMS}), since in this case the two-way coupling between the media is the most pronounced.

Note that the solver could not converge to the chosen tolerance, of $10^{-6}$ residual norm, when F-AMS was iterated in a Richardson's loop for $T_{ratio} \ge 10^6$ and the same holds for SAMG \cite{SAMG}. However, as shown in Fig.~\ref{fig:conv_perm1e6}, both methods converge successfully when employed as preconditioners for GMRES \cite{Saad}.

\subsection{CPU benchmark study}

This final subsection presents the results of a benchmark study between F-AMS and SAMG \cite{SAMG} on 3D heterogeneous fractured reservoirs. Both methods are employed as preconditioners to GMRES \cite{Saad} and iterated until converged with a residual 2-norm below $10^{-6}$. 

Unlike the Richardson loop, similar performance was observed for all experiments when F-AMS is used as preconditioner to GMRES, regardless of which coupling strategy was chosen. Therefore, the presentation of the results is restricted to Decoupled-AMS, for conciseness. Note that the Coupled-AMS strategy provides a more general framework, however, with a more complex implementation. 

At each GMRES iteration, SAMG employs a single V-cycle. It is important to note that SAMG is a commercial black-box package. Thus, it is not possible to measure its CPU breakdown as accurately as for F-AMS. Instead, the time spent on its first iteration is considered as ``Initialization'', while subsequent iterations were labelled as ``Solution''. Finally, for both SAMG and F-AMS, the setup and construction of the operators are performed only once, at the beginning of the iteration procedure.

This study is aimed only to demonstrate the scalability of the F-AMS method. Note that a unique advantage of F-AMS over SAMG is that a fine-scale mass conservative flux field can be reconstructed after any iteration stage, once the coarse-scale system with $\bm{\mathcal{R}}^{FV}$ restriction operator is solved.

\subsubsection{Transmissibility contrast}

\begin{figure}[htb!]%
\centering%
\setlength{\fboxsep}{0pt}%
\setlength{\fboxrule}{0.2pt}%
\setlength{\unitlength}{1cm}%
\begin{picture}(15,0.8)%
\put(0.2,0.2){$T_{ratio}$}%
\put(2.4,0.2){$10^1$}%
\put(5.7,0.2){$10^2$}%
\put(9.0,0.2){$10^4$}%
\put(12.2,0.2){$10^8$}%
\end{picture}\\%
\begin{picture}(1.0,3.5)%
	\put(0.0,0.3){\rotatebox{90}{CPU time (sec)}}%
		\put(0.75,-0.1){0}%
		\put(0.7,1.7){2}%
		\put(0.7,3.3){4}%
\end{picture}%
\subfigure{\fbox{\includegraphics[width=0.19\textwidth,height=0.25\textwidth]{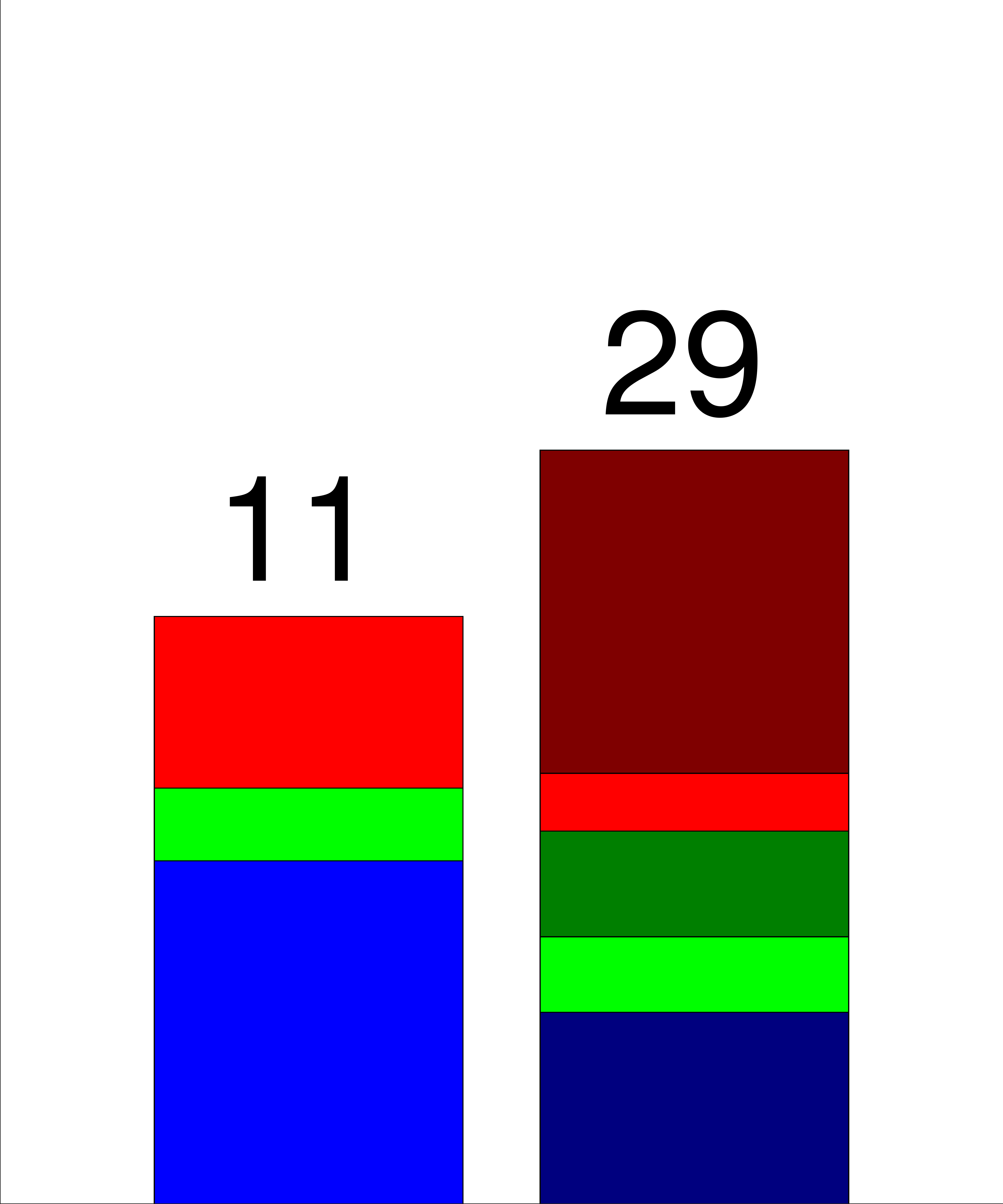}}\label{fig:cpu_bench_perm1e1}}%
\begin{picture}(0.6,3.5)%
		\put(0.3,-0.1){0}%
		\put(0.3,1.7){2}%
		\put(0.3,3.3){4}%
\end{picture}%
\subfigure{\fbox{\includegraphics[width=0.19\textwidth,height=0.25\textwidth]{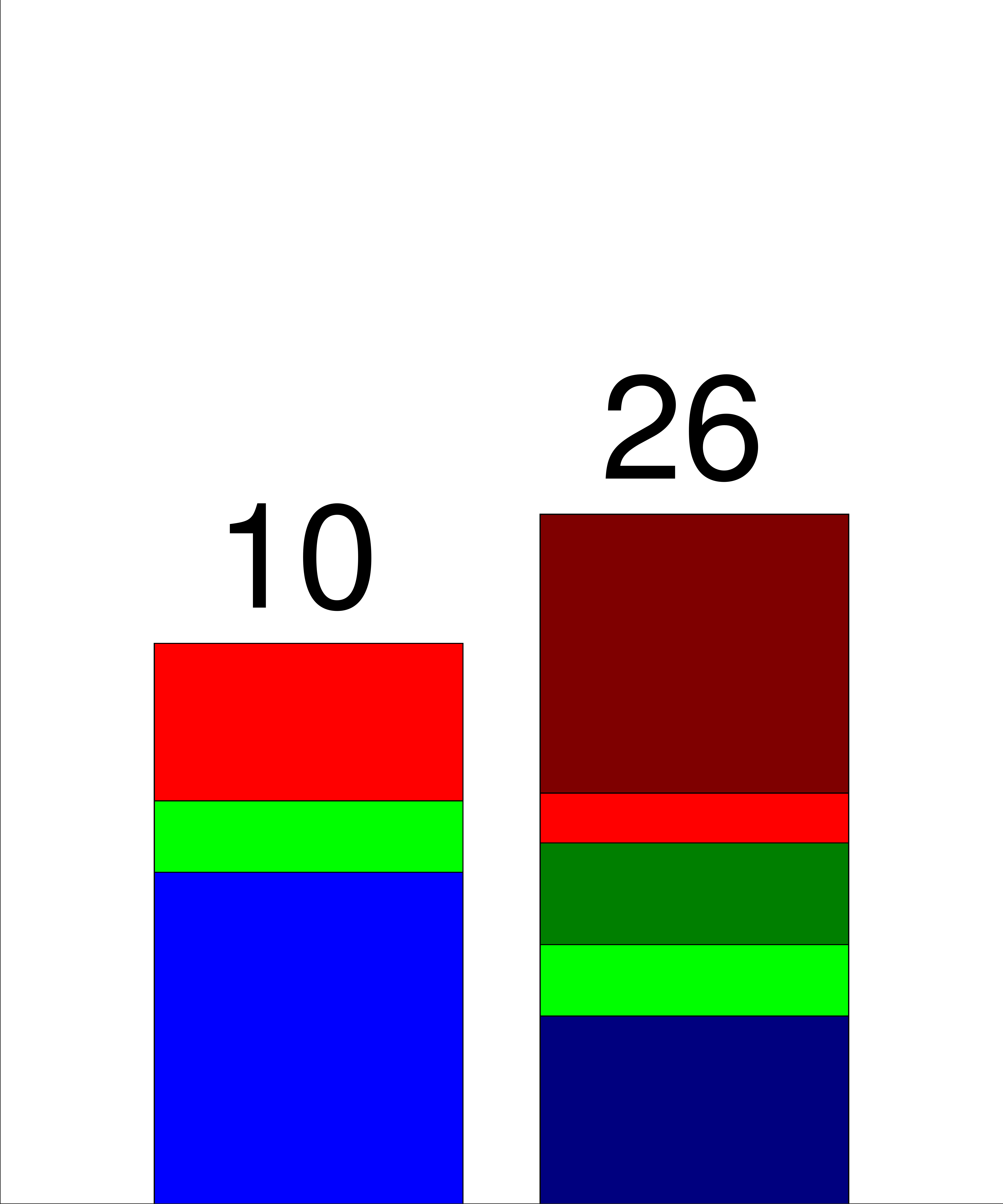}}\label{fig:cpu_bench_perm1e2}}%
\begin{picture}(0.6,3.5)%
		\put(0.3,-0.1){0}%
		\put(0.3,1.7){2}%
		\put(0.3,3.3){4}%
\end{picture}%
\subfigure{\fbox{\includegraphics[width=0.19\textwidth,height=0.25\textwidth]{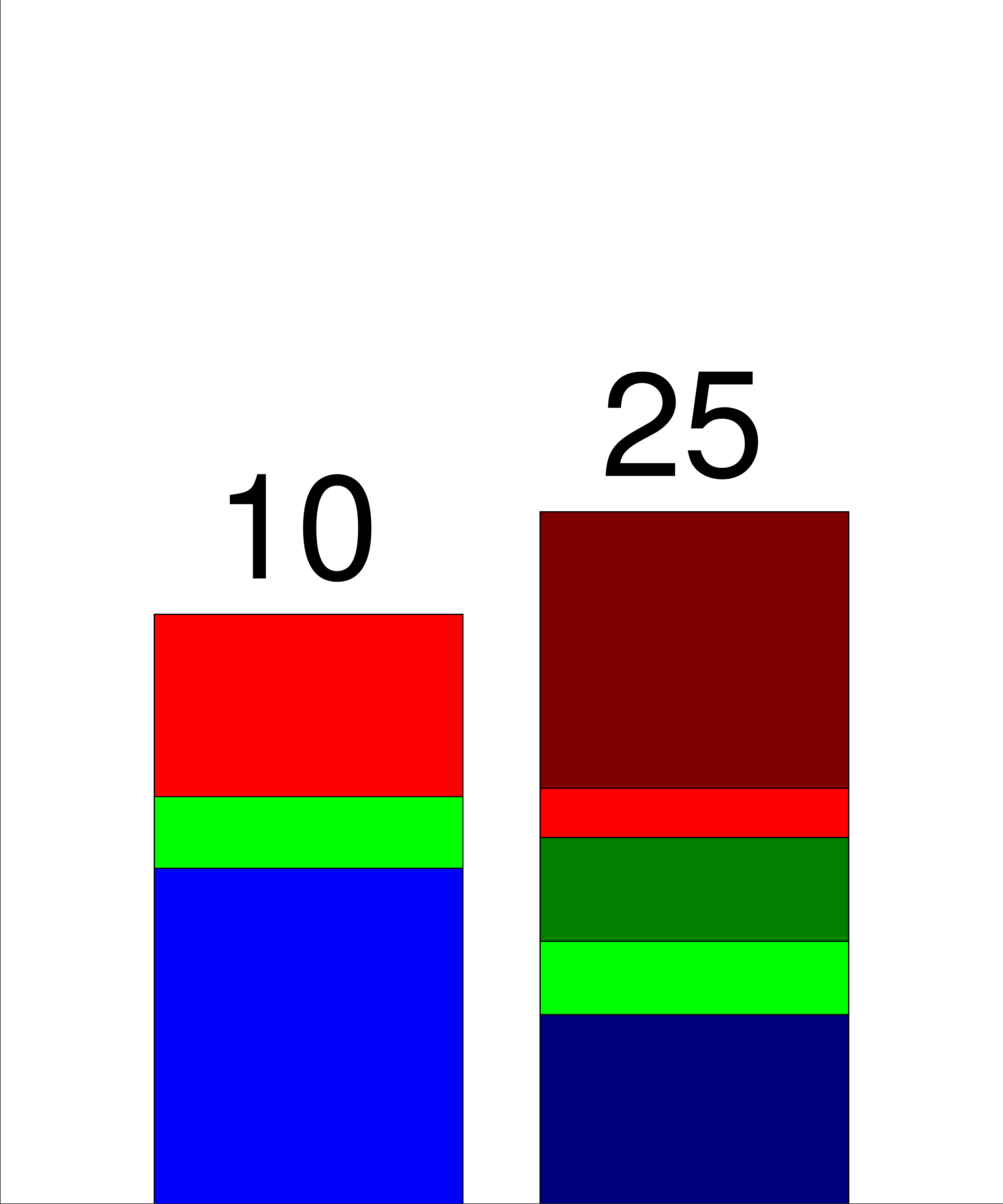}}\label{fig:cpu_bench_perm1e4}}%
\begin{picture}(0.6,3.5)%
		\put(0.3,-0.1){0}%
		\put(0.3,1.7){2}%
		\put(0.3,3.3){4}%
\end{picture}%
\subfigure{\fbox{\includegraphics[width=0.19\textwidth,height=0.25\textwidth]{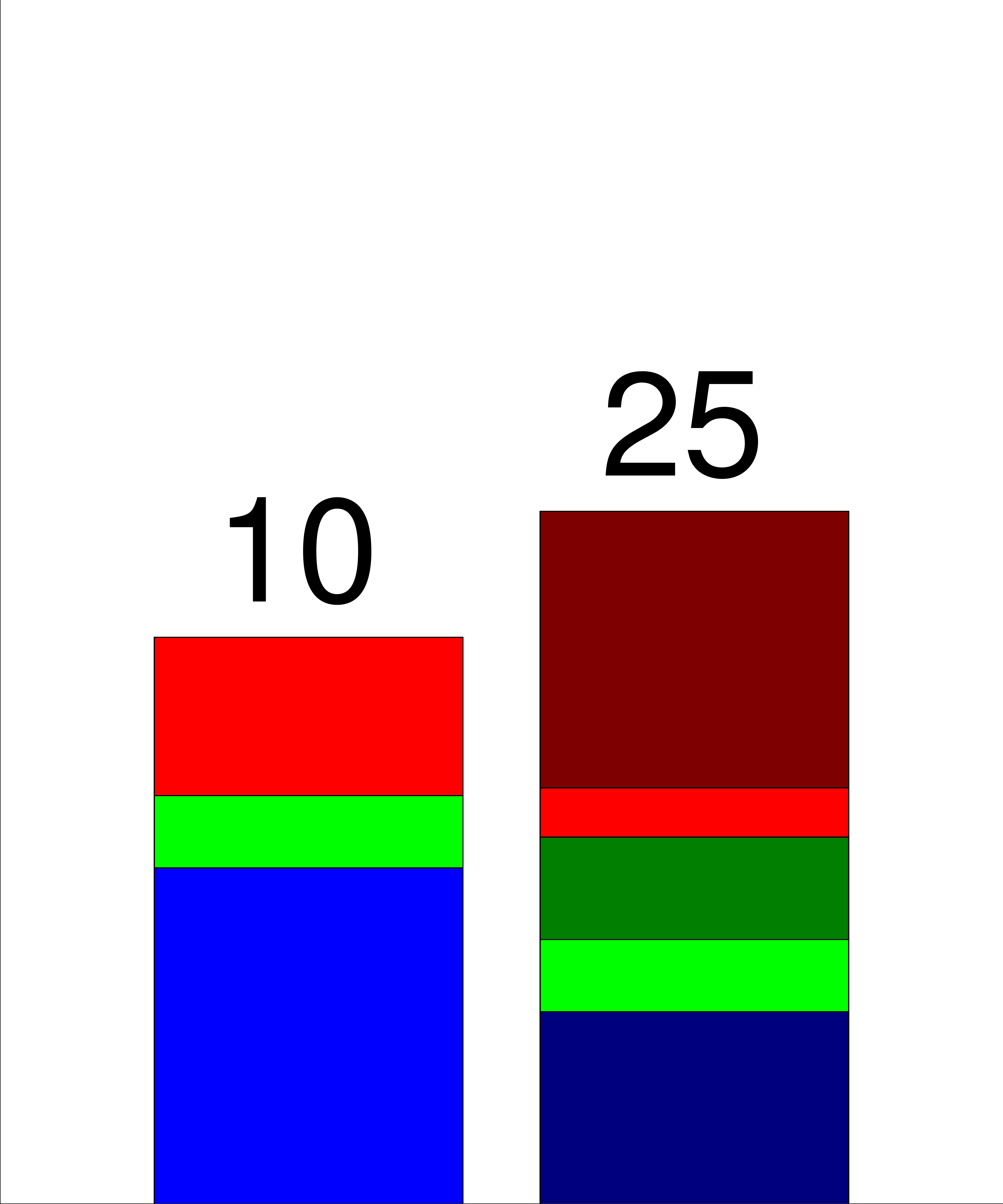}}\label{fig:cpu_bench_perm1e8}}\\%
\begin{picture}(15,1.0)%
\put(1.2,0.0){\rotatebox{45}{SAMG}}%
\put(2.2,0.0){\rotatebox{45}{F-AMS}}%
\put(4.4,0.0){\rotatebox{45}{SAMG}}%
\put(5.4,0.0){\rotatebox{45}{F-AMS}}%
\put(7.75,0.0){\rotatebox{45}{SAMG}}%
\put(8.75,0.0){\rotatebox{45}{F-AMS}}%
\put(11.00,0.0){\rotatebox{45}{SAMG}}%
\put(12.00,0.0){\rotatebox{45}{F-AMS}}%
\end{picture}\\%
\subfigure{\fbox{\includegraphics[width=\textwidth]{Fig/method/legend.png}}}%
\caption{F-AMS (Decoupled-AMS) performance compared with SAMG for different matrix-fracture transmissibility ratios. The number of performed iterations to reach $10^{-6}$ residual 2-norm is given on top of each bar. For these experiments, both methods were used as preconditioners for GMRES. Similar performance was observed for other F-AMS coupling strategies.}%
\label{fig:bench}%
\end{figure}

The test case from Fig.~\ref{fig:3dcase} is used, with different values of fracture-matrix transmissibility contrasts, i.e., $T_{ratio}$ in Eq.~\eqref{tratio}. As can be seen in Fig.~\ref{fig:bench}, the number of iterations and the CPU time for both F-AMS and SAMG is insensitive to the contrast. This is a significant achievement for F-AMS, compared to \cite{hadi-frac-jcp}. 

\subsubsection{Fracture density}

\begin{figure}[htb!]%
\subfigure{\includegraphics[height=0.25\textwidth,width=0.25\textwidth]{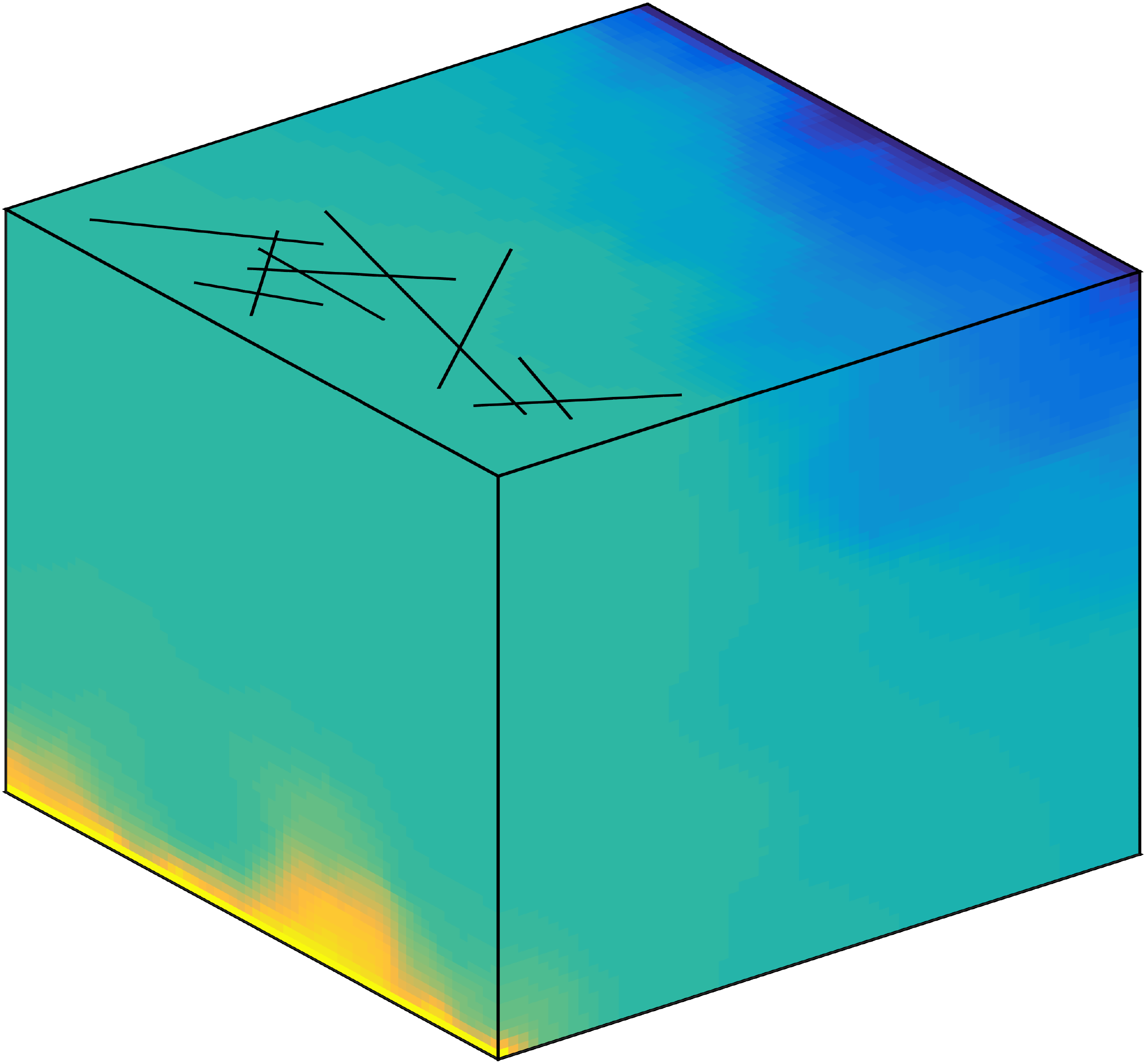}\label{fig:phase1}\hspace{0.1cm}}%
\subfigure{\includegraphics[height=0.25\textwidth,width=0.25\textwidth]{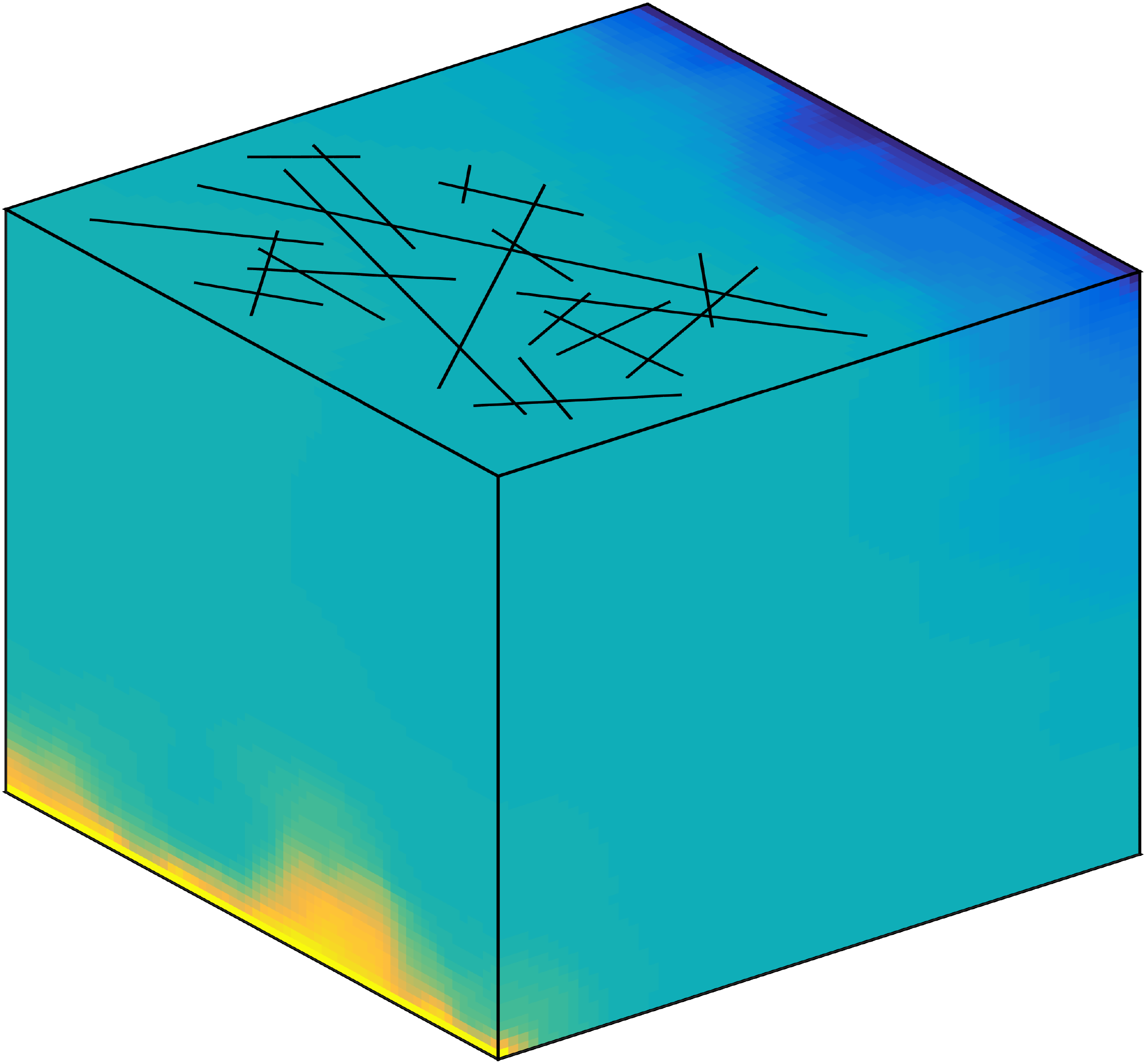}\label{fig:phase2}\hspace{0.1cm}}%
\subfigure{\includegraphics[height=0.25\textwidth,width=0.25\textwidth]{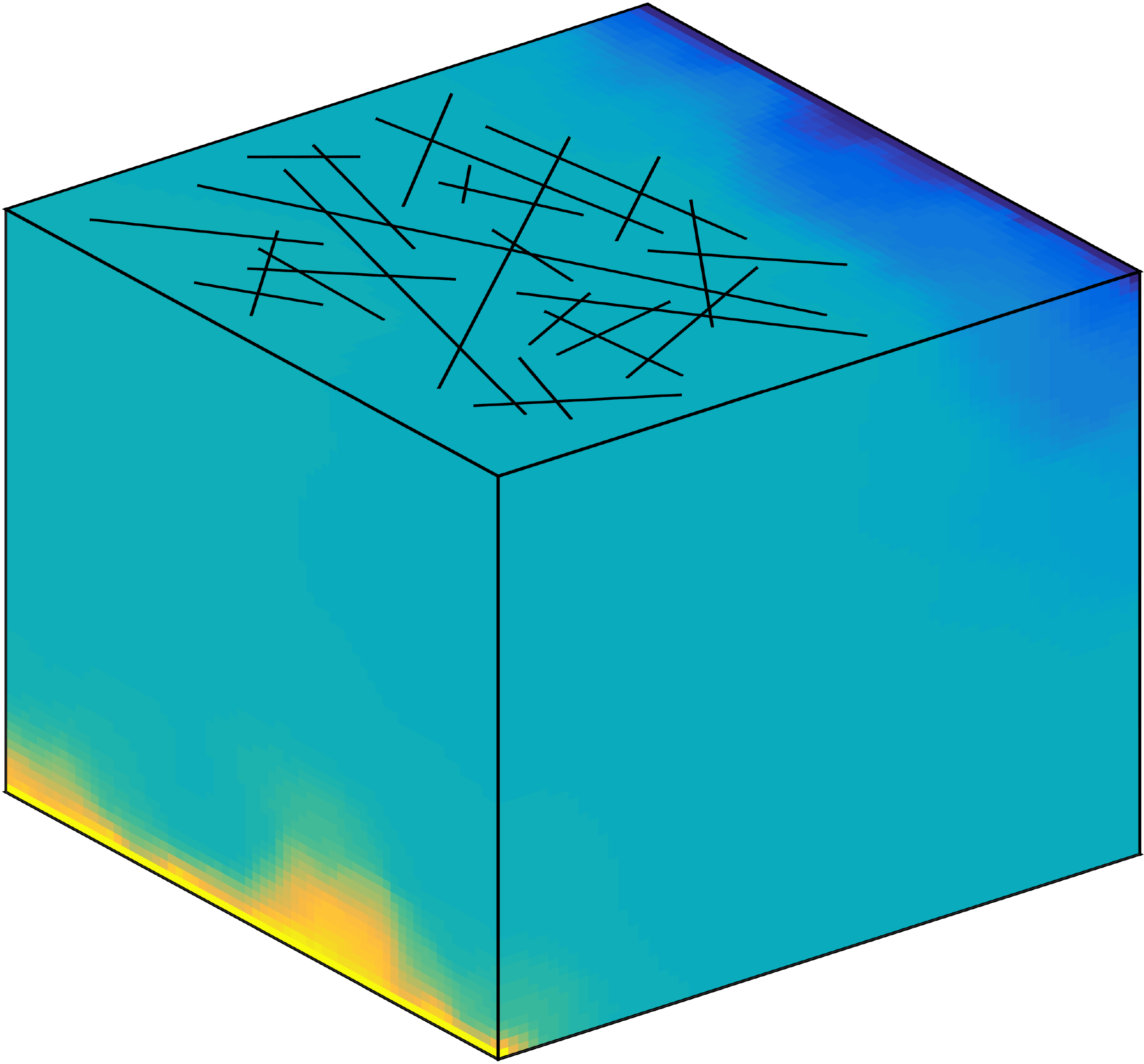}\label{fig:phase3}\hspace{0.1cm}}%
\subfigure{\includegraphics[height=0.25\textwidth,width=0.25\textwidth]{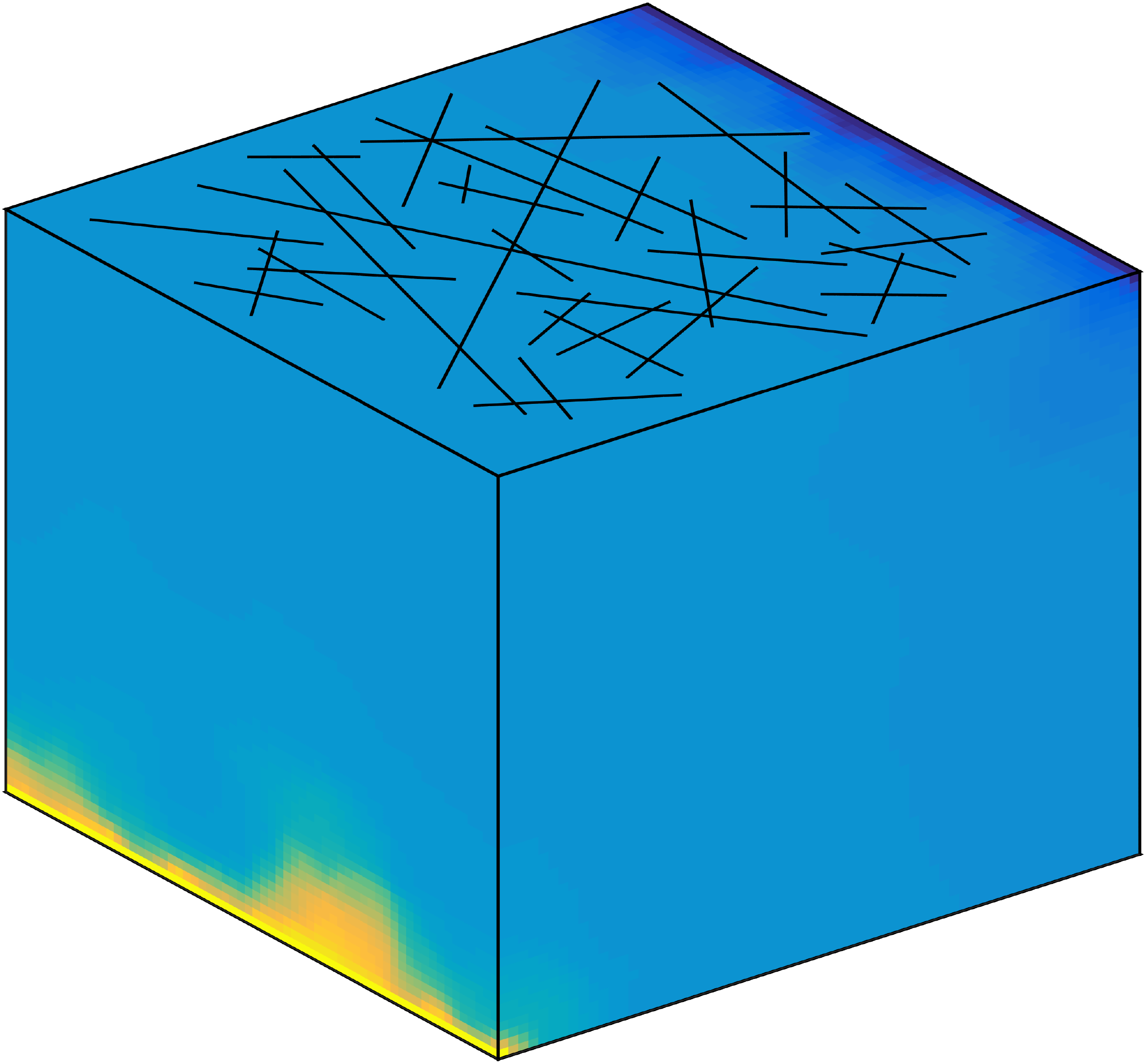}\hspace{0.1cm}\includegraphics[height=0.25\textwidth,width=0.5cm]{Fig/method/colorbar.pdf}\label{fig:phase4}}%
\caption{Pressure solution obtained for different fracture densities. The left-most contains 27 fracture plates ($136 \times 64$ cells), followed by one with 76 fracture plates ($334 \times 64$ cells), the next has 96 fracture plates ($438 \times 64$ cells), and, finally, the right-most is perforated by 127 fracture plates ($575 \times 64$ cells).}%
\label{fig:slicedFracs}%
\end{figure}

\begin{figure}[htb!]%
\centering%
\setlength{\fboxsep}{0pt}%
\setlength{\fboxrule}{0.2pt}%
\setlength{\unitlength}{1cm}%
\begin{picture}(15,0.8)%
\put(0.2,0.2){\#fracs}%
\put(2.4,0.2){$27$}%
\put(5.7,0.2){$76$}%
\put(8.9,0.2){$96$}%
\put(12.1,0.2){$127$}%
\end{picture}\\%
\begin{picture}(1.0,3.5)%
	\put(0.0,0.3){\rotatebox{90}{CPU time (sec)}}%
		\put(0.75,-0.1){0}%
		\put(0.7,1.7){2}%
		\put(0.7,3.3){4}%
\end{picture}%
\subfigure{\fbox{\includegraphics[width=0.19\textwidth,height=0.25\textwidth]{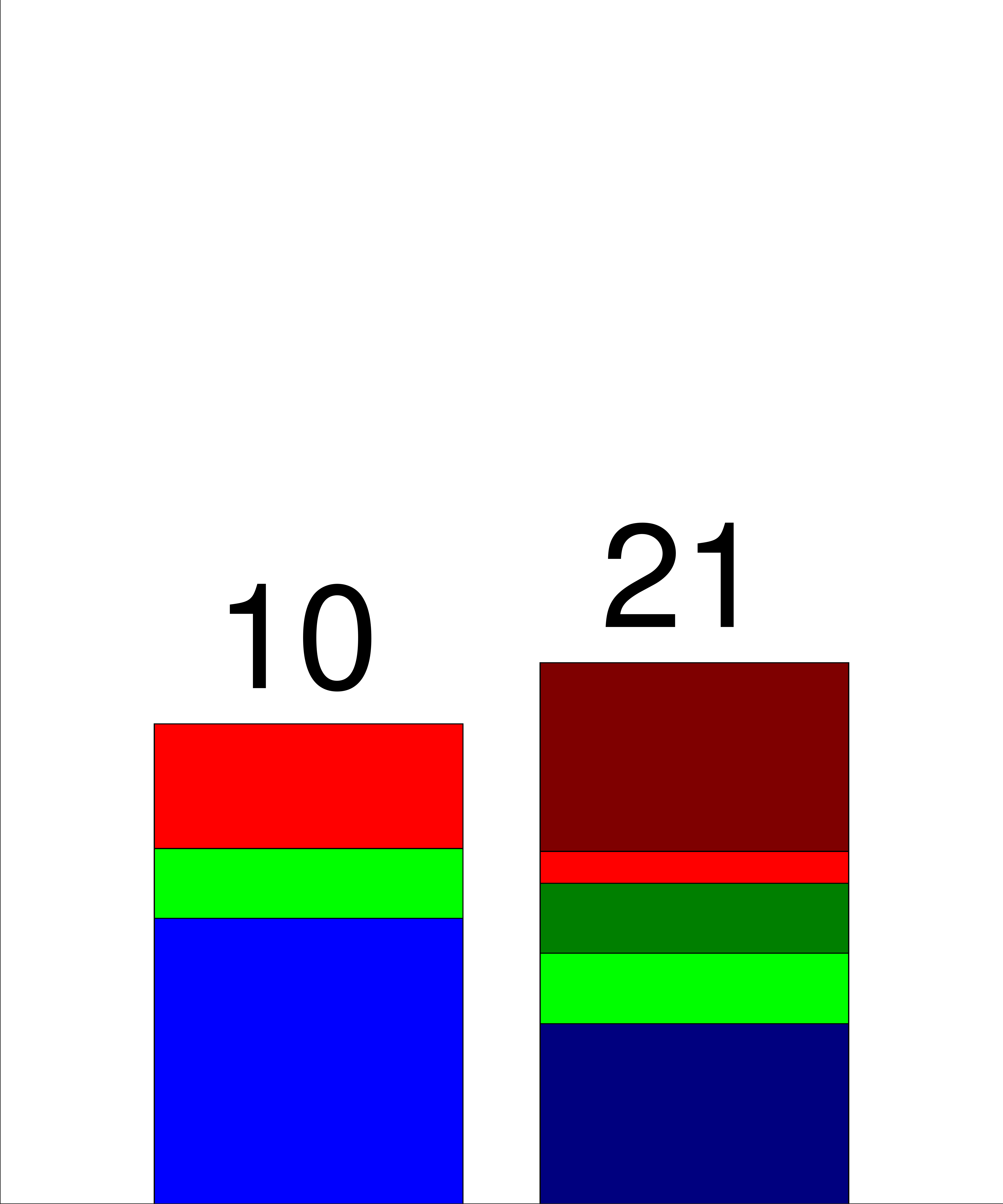}}\label{fig:cpu_nfracs1}}%
\begin{picture}(0.6,3.5)%
		\put(0.3,-0.1){0}%
		\put(0.3,1.7){2}%
		\put(0.3,3.3){4}%
\end{picture}%
\subfigure{\fbox{\includegraphics[width=0.19\textwidth,height=0.25\textwidth]{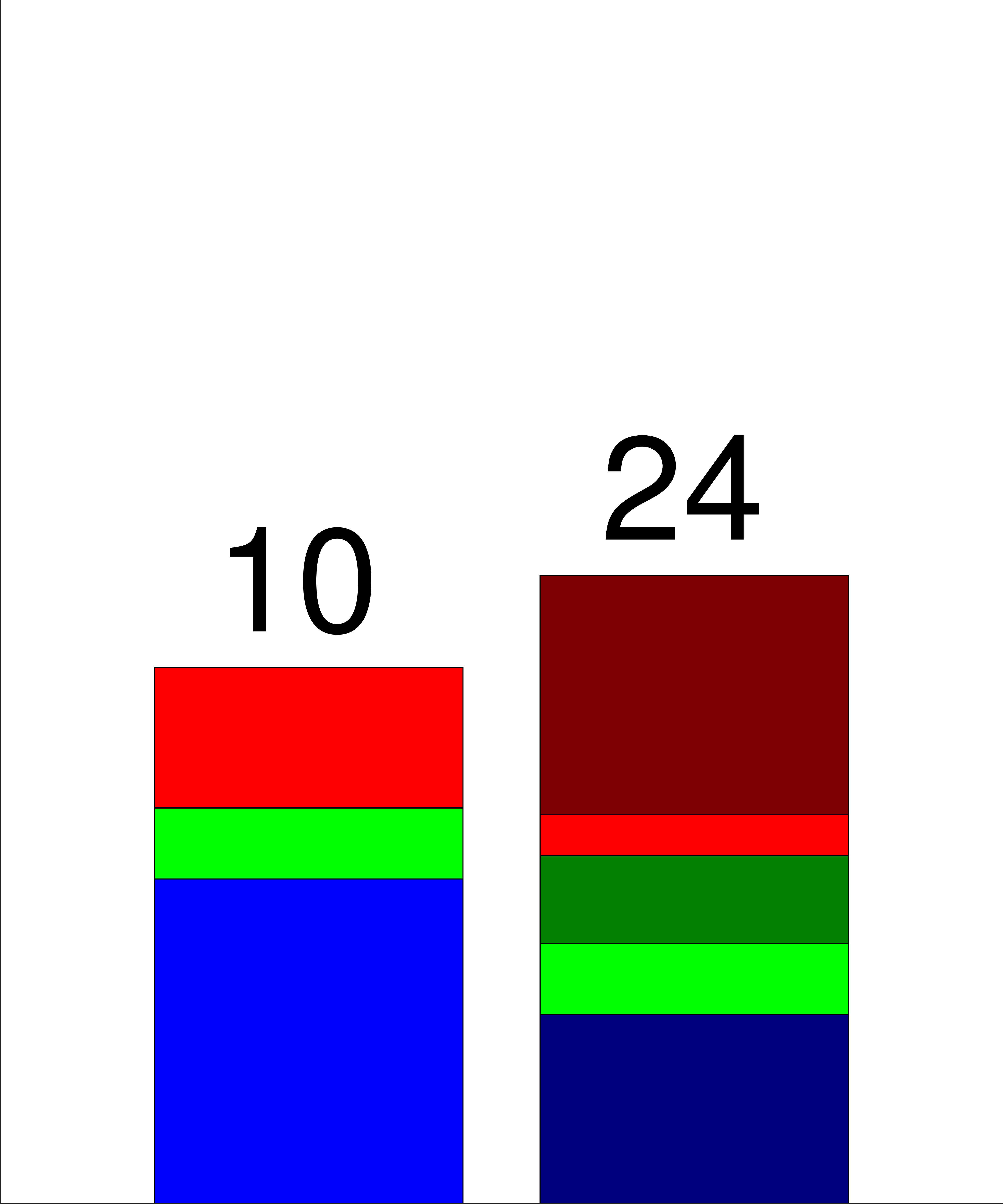}}\label{fig:cpu_nfracs2}}%
\begin{picture}(0.6,3.5)%
		\put(0.3,-0.1){0}%
		\put(0.3,1.7){2}%
		\put(0.3,3.3){4}%
\end{picture}%
\subfigure{\fbox{\includegraphics[width=0.19\textwidth,height=0.25\textwidth]{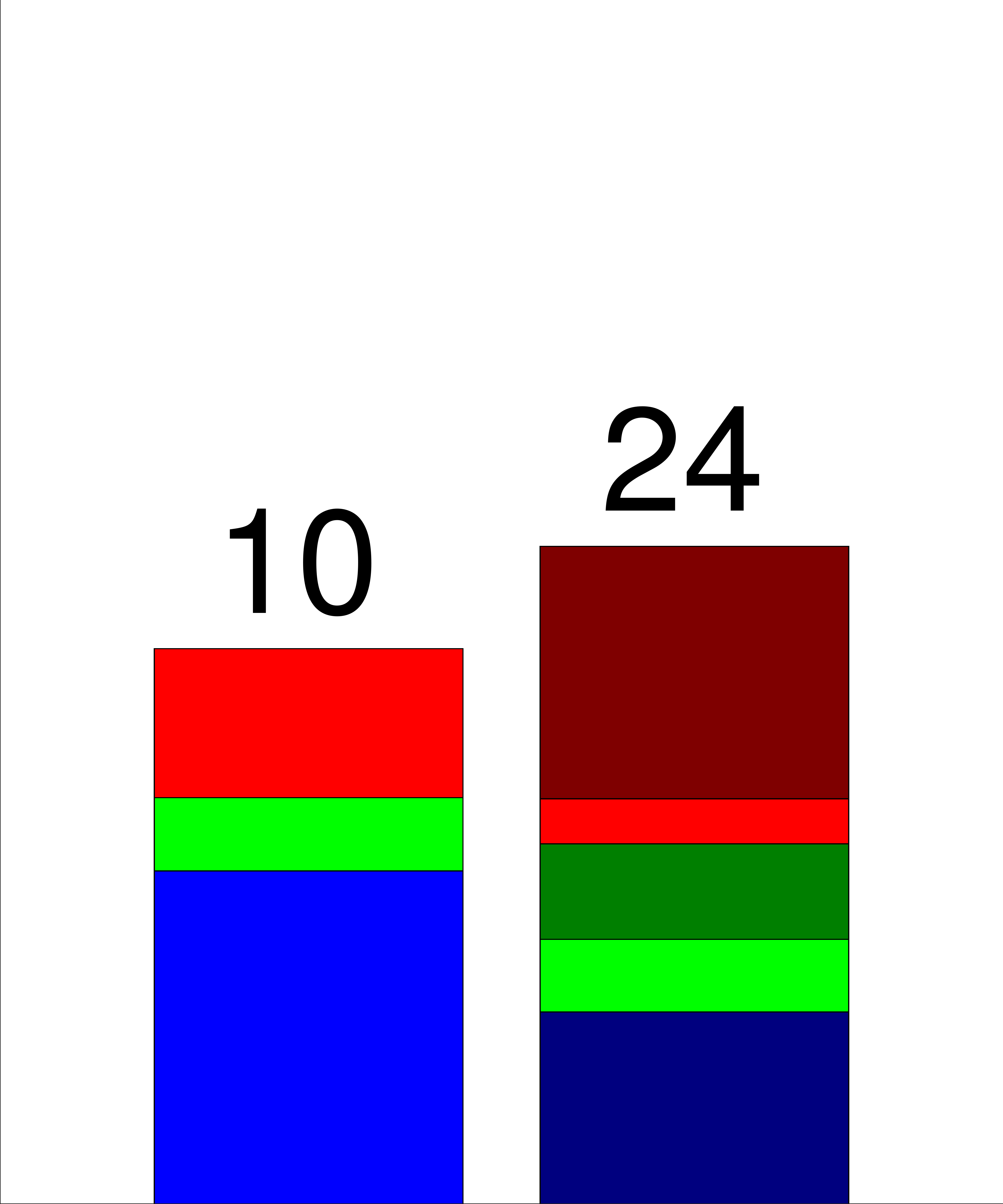}}\label{fig:cpu_nfracs3}}%
\begin{picture}(0.6,3.5)%
		\put(0.3,-0.1){0}%
		\put(0.3,1.7){2}%
		\put(0.3,3.3){4}%
\end{picture}%
\subfigure{\fbox{\includegraphics[width=0.19\textwidth,height=0.25\textwidth]{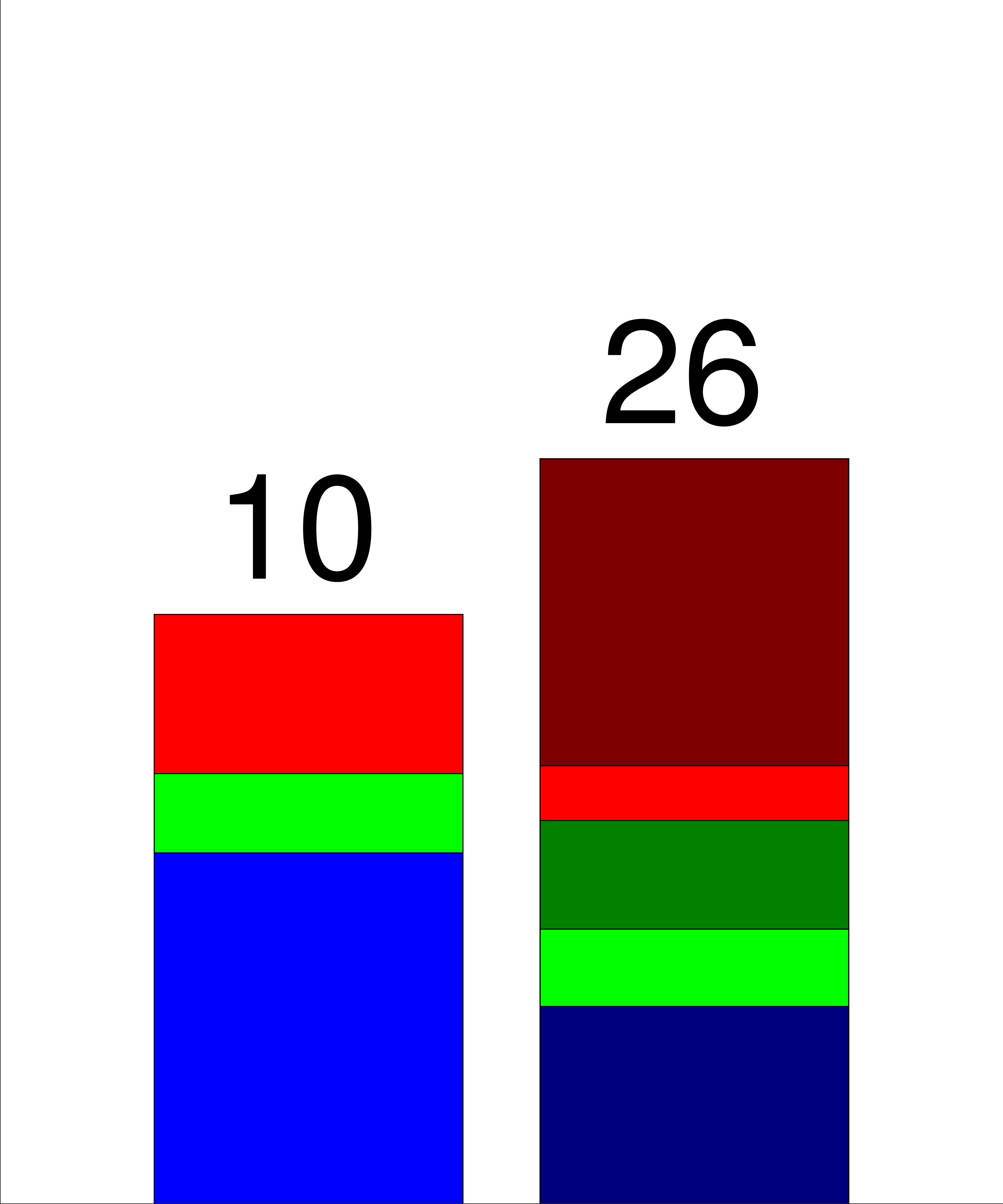}}\label{fig:cpu_nfracs4}}\\%
\begin{picture}(15,1.0)%
\put(1.2,0.0){\rotatebox{45}{SAMG}}%
\put(2.2,0.0){\rotatebox{45}{F-AMS}}%
\put(4.4,0.0){\rotatebox{45}{SAMG}}%
\put(5.4,0.0){\rotatebox{45}{F-AMS}}%
\put(7.75,0.0){\rotatebox{45}{SAMG}}%
\put(8.75,0.0){\rotatebox{45}{F-AMS}}%
\put(11.00,0.0){\rotatebox{45}{SAMG}}%
\put(12.00,0.0){\rotatebox{45}{F-AMS}}%
\end{picture}\\%
\subfigure{\fbox{\includegraphics[width=\textwidth]{Fig/method/legend.png}}}%
\caption{F-AMS (Decoupled-AMS) performance for reservoirs with different number of fractures. The number of performed iterations to reach $10^{-6}$ residual 2-norm is given on top of each bar. For these experiments, both methods were employed as preconditioners for GMRES. Similar performance was observed for other F-AMS coupling strategies.}%
\label{fig:nfracs}%
\end{figure}

What follows is a study of the scalability of F-AMS when faced with a dynamic fracture network, where the number of fracture plates is increased step by step. The 3D fracture map shown in Fig.~\ref{fig:fracs3d} is considered, where the network is now created through 4 phases, as presented in Fig.~\ref{fig:slicedFracs}. Note that, as new fractures are added, not only the number of DOF increases, but also the pressure variation along the network can increase. The detailed description of the CPU times obtained using F-AMS in these four cases are depicted in Fig.~\ref{fig:nfracs}. It is clear that, by maintaining the prescribed fracture coarsening factor of $8 \times 8$, F-AMS maintains virtually the same convergence rate. The slight increase in CPU time is mainly due to computation of extra fracture basis functions, as well as the construction and solution of a slightly larger coarse-scale linear system. 

In consequence, by having multiple coarse-scale DOF in each fracture network, F-AMS can automatically scale with fracture length and density. This is in contrast to \cite{hadi-frac-jcp}, where the use of a single fracture basis function would lead to a drastic deterioration of the multiscale convergence for test cases containing fracture networks with large length scales.

\subsubsection{Domain scale}
\label{sec:domain_scale}

\begin{figure}[htb!]%
\subfigure{\includegraphics[height=0.17\textwidth,width=0.17\textwidth]{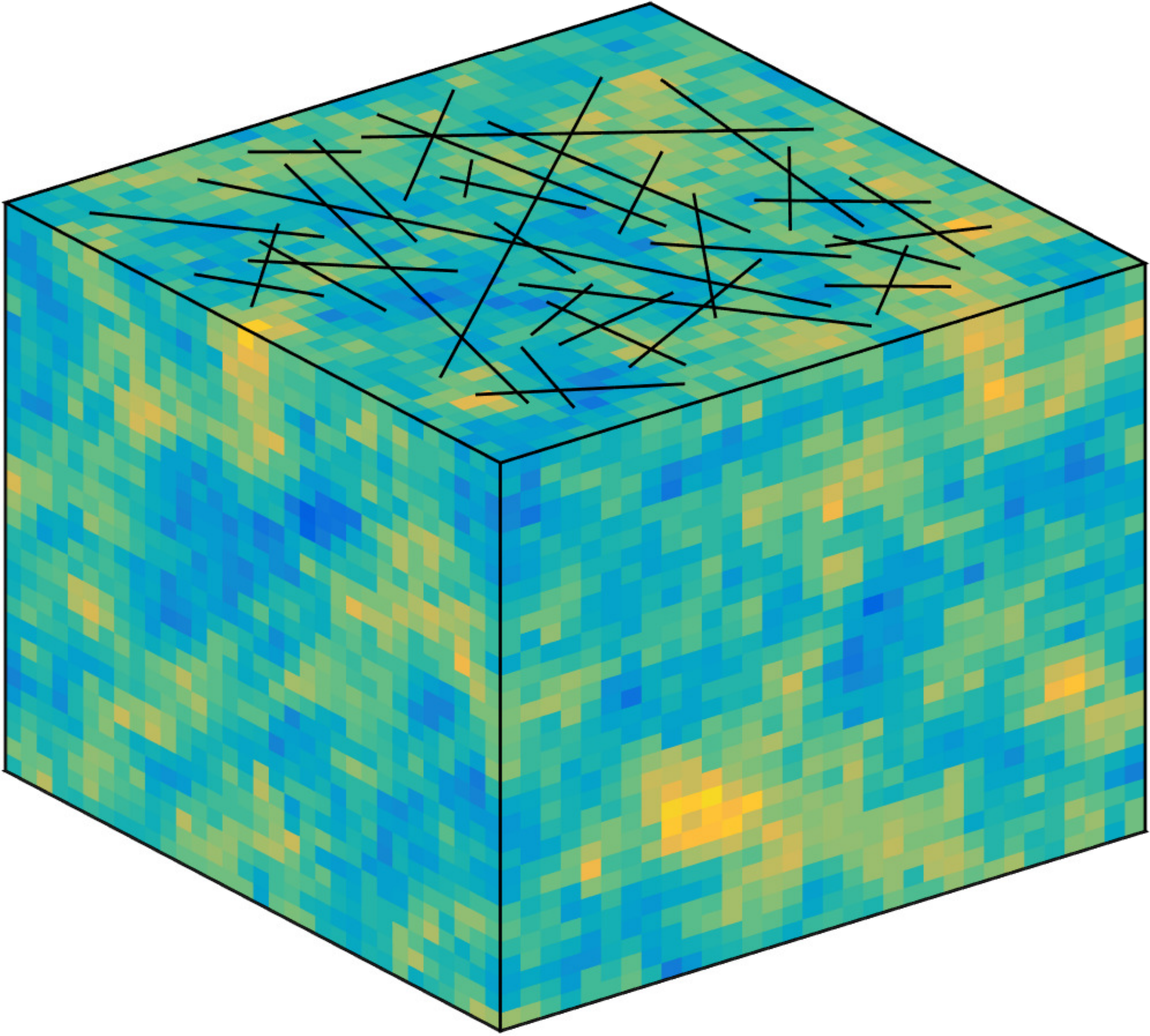}\label{fig:perm32}\hspace{0.1cm}}%
\subfigure{\includegraphics[height=0.22\textwidth,width=0.22\textwidth]{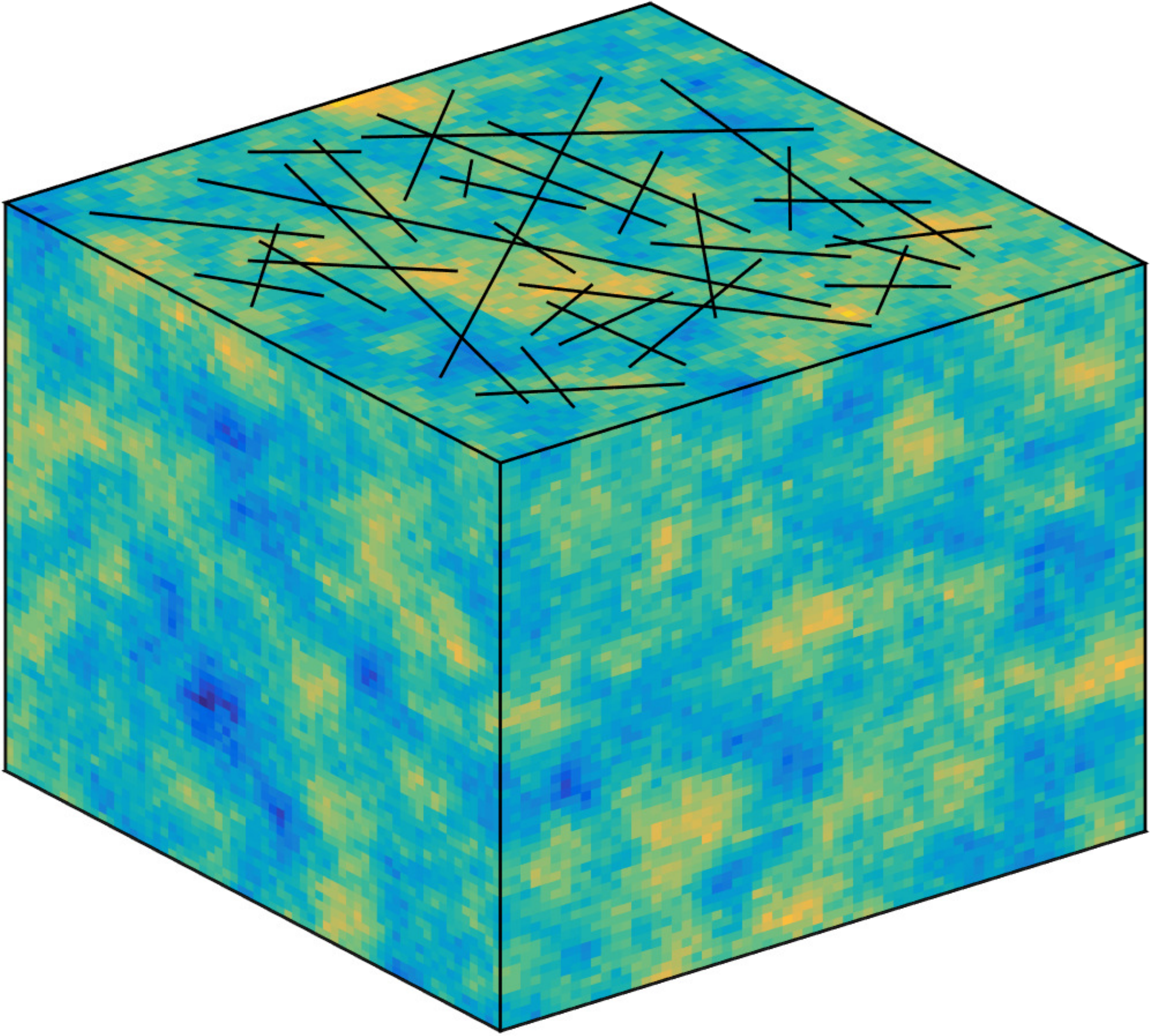}\label{fig:perm64}\hspace{0.1cm}}%
\subfigure{\includegraphics[height=0.27\textwidth,width=0.27\textwidth]{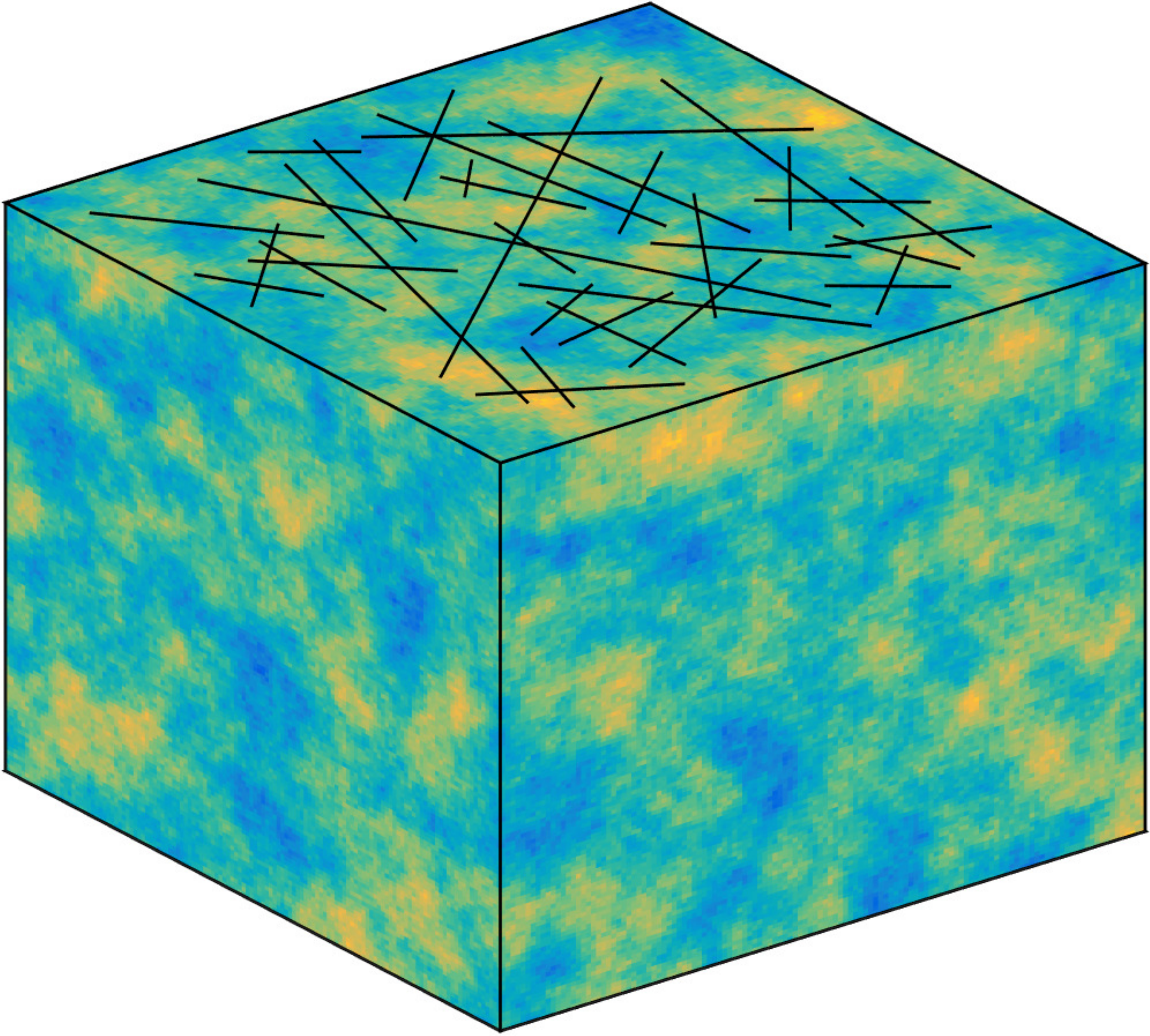}\label{fig:perm128}\hspace{0.1cm}}%
\subfigure{\includegraphics[height=0.32\textwidth,width=0.32\textwidth]{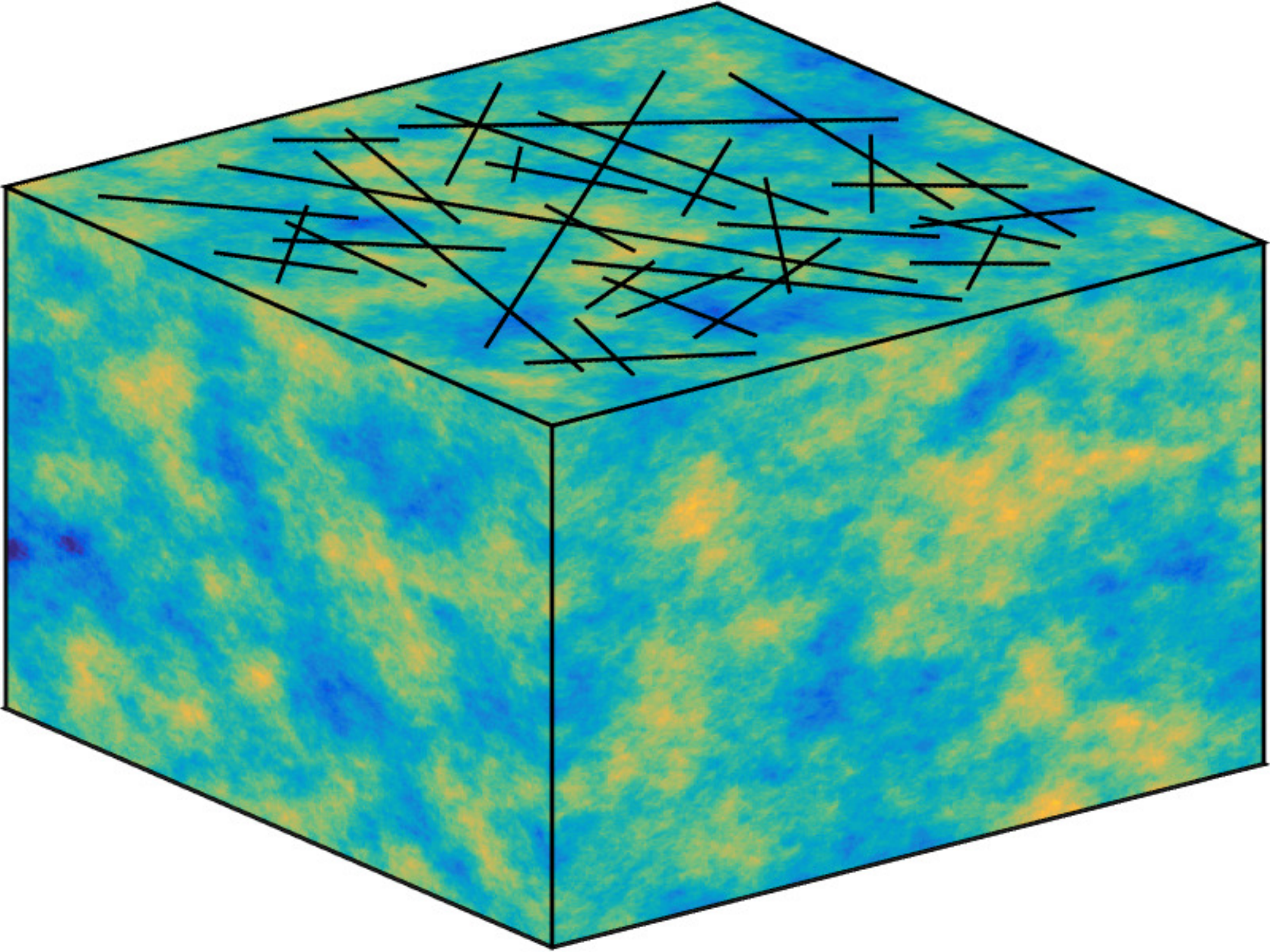}\hspace{0.1cm}\raisebox{0.1cm}{\includegraphics[height=0.3\textwidth,width=0.4cm]{Fig/3D/colorbar_perm3d.pdf}}\label{fig:perm256}}%
\caption{Logarithm of the permeability ($log_{10}(k^m)$) and fracture distribution in four reservoirs of increasing size. The left-most has $32^3$ matrix and $320 \times 32$ fracture cells, followed by $64^3$ matrix and $575 \times 64$ fracture cells, then $128^3$ matrix and $1087 \times 128$ fracture cells, and, finally, $256^3$ matrix and $2117 \times 256$ fracture cells.}%
\label{fig:scaledFracs}%
\end{figure}

\begin{figure}[htb!]%
\centering%
\setlength{\fboxsep}{0pt}%
\setlength{\fboxrule}{0.2pt}%
\setlength{\unitlength}{1cm}%
\begin{picture}(15,0.8)%
\put(0.0,0.2){Scale}%
\put(2.4,0.2){$32^3$}%
\put(5.6,0.2){$64^3$}%
\put(8.9,0.2){$128^3$}%
\put(12.4,0.2){$256^3$}%
\end{picture}\\%
  \begin{picture}(1.2,1.5)%
	\put(0.0,0.3){\rotatebox{90}{CPU time (sec)}}%
		\put(0.95,-0.1){0}%
		\put(0.55,1.7){0.5}%
		\put(0.9,3.3){1}%
\end{picture}%
\subfigure{\fbox{\includegraphics[width=0.19\textwidth,height=0.25\textwidth]{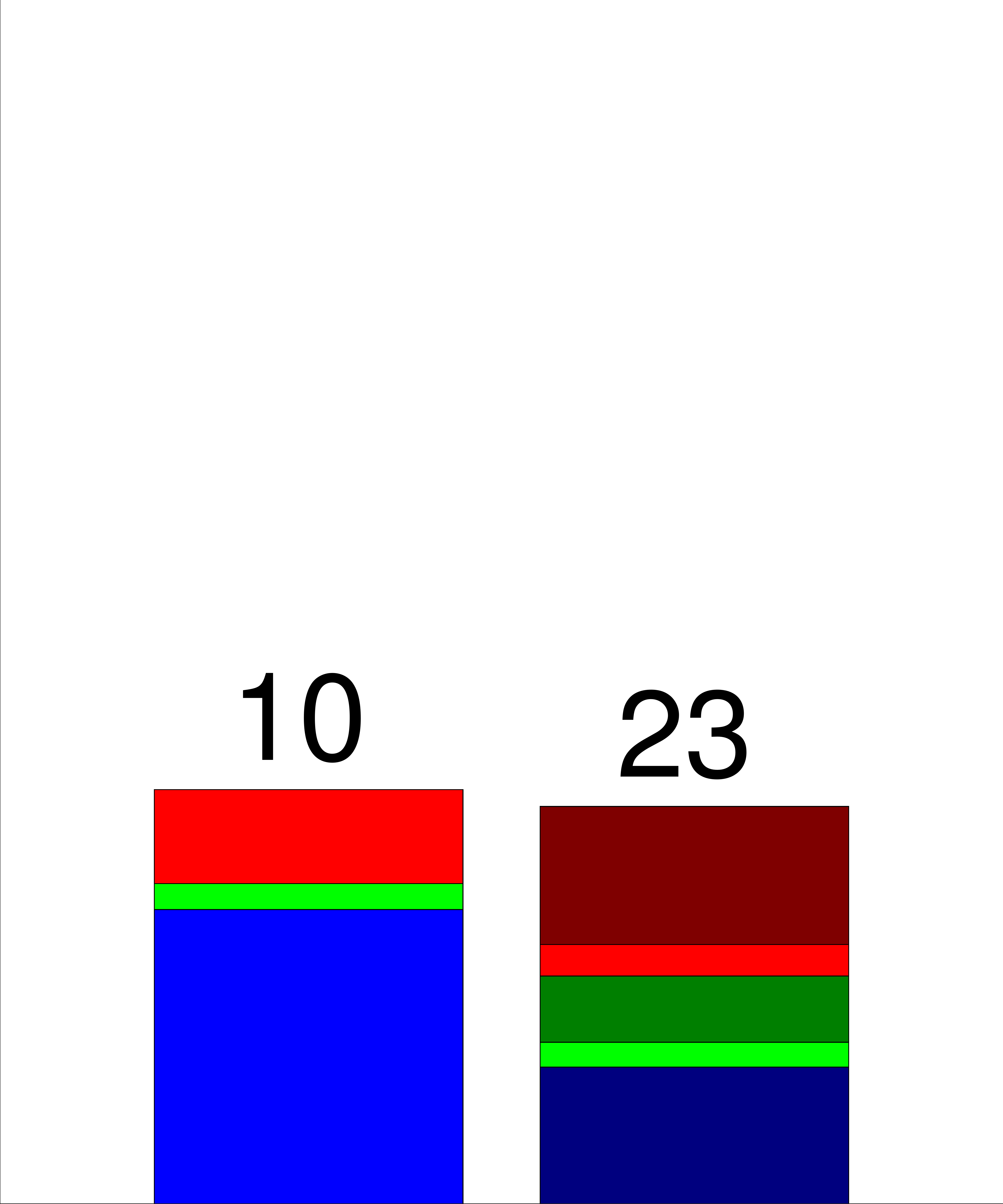}}\label{fig:cpu_scale32}}%
\begin{picture}(0.6,1.5)%
		\put(0.3,-0.1){0}%
		\put(0.3,1.7){2}%
		\put(0.3,3.3){4}%
\end{picture}%
\subfigure{\fbox{\includegraphics[width=0.19\textwidth,height=0.25\textwidth]{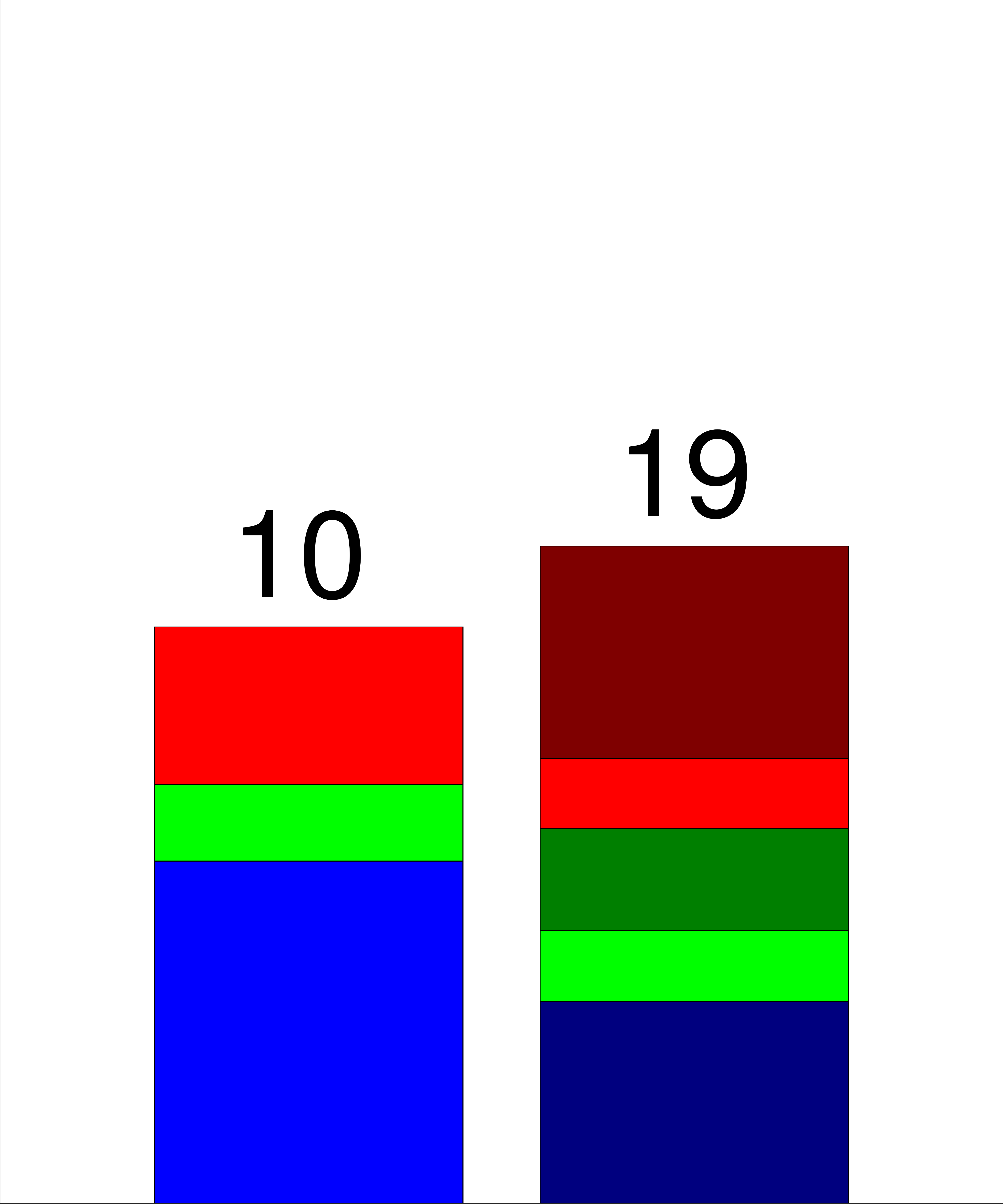}}\label{fig:cpu_scale64}}%
  \begin{picture}(0.7,1.5)%
		\put(0.4,-0.1){0}%
		\put(0.2,1.7){12}%
		\put(0.2,3.3){24}%
\end{picture}%
\subfigure{\fbox{\includegraphics[width=0.19\textwidth,height=0.25\textwidth]{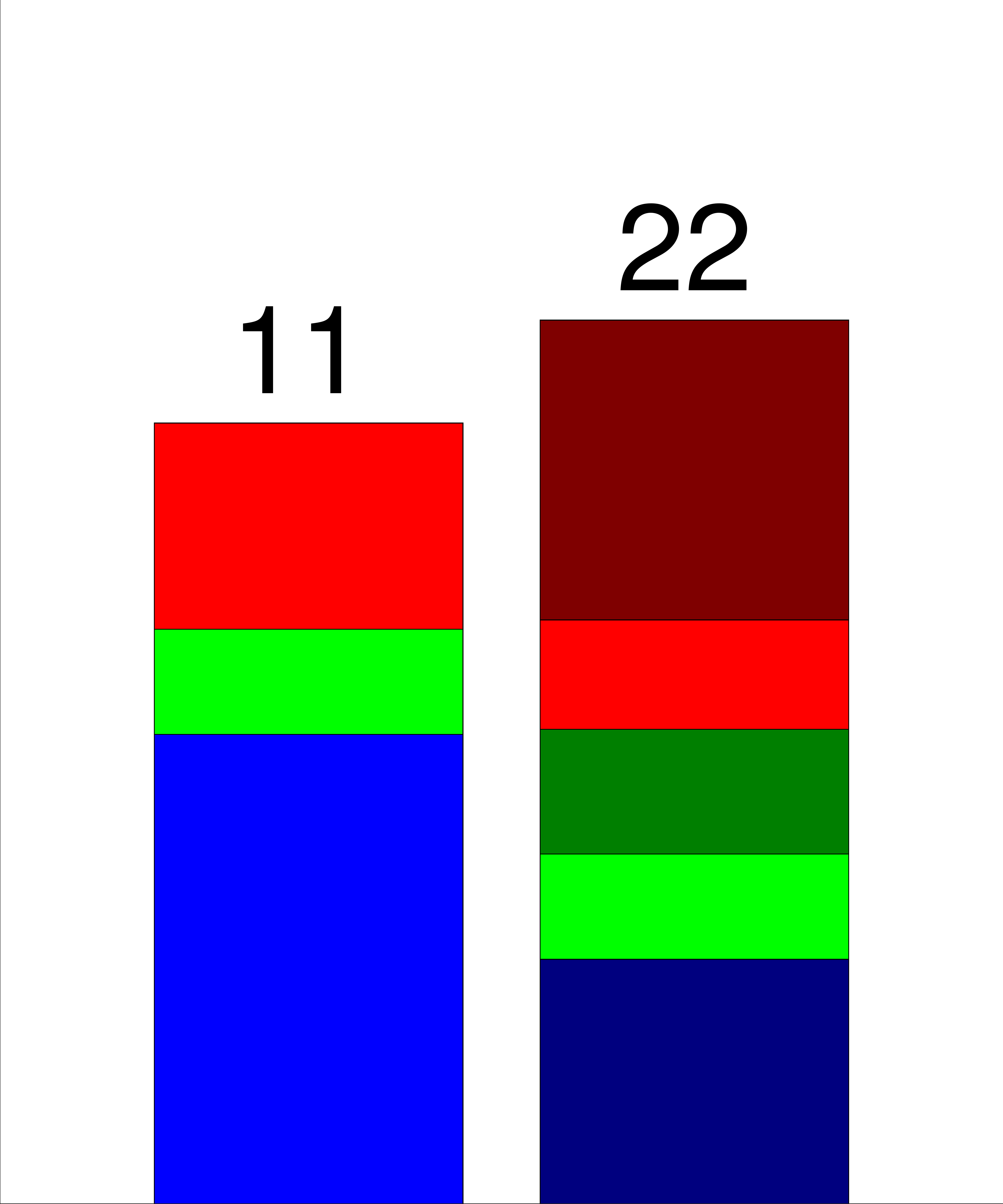}}\label{fig:cpu_scale128}}%
\begin{picture}(0.8,1.5)%
		\put(0.5,-0.1){0}%
		\put(0.3,1.7){90}%
		\put(0.1,3.3){180}%
\end{picture}%
\subfigure{\fbox{\includegraphics[width=0.19\textwidth,height=0.25\textwidth]{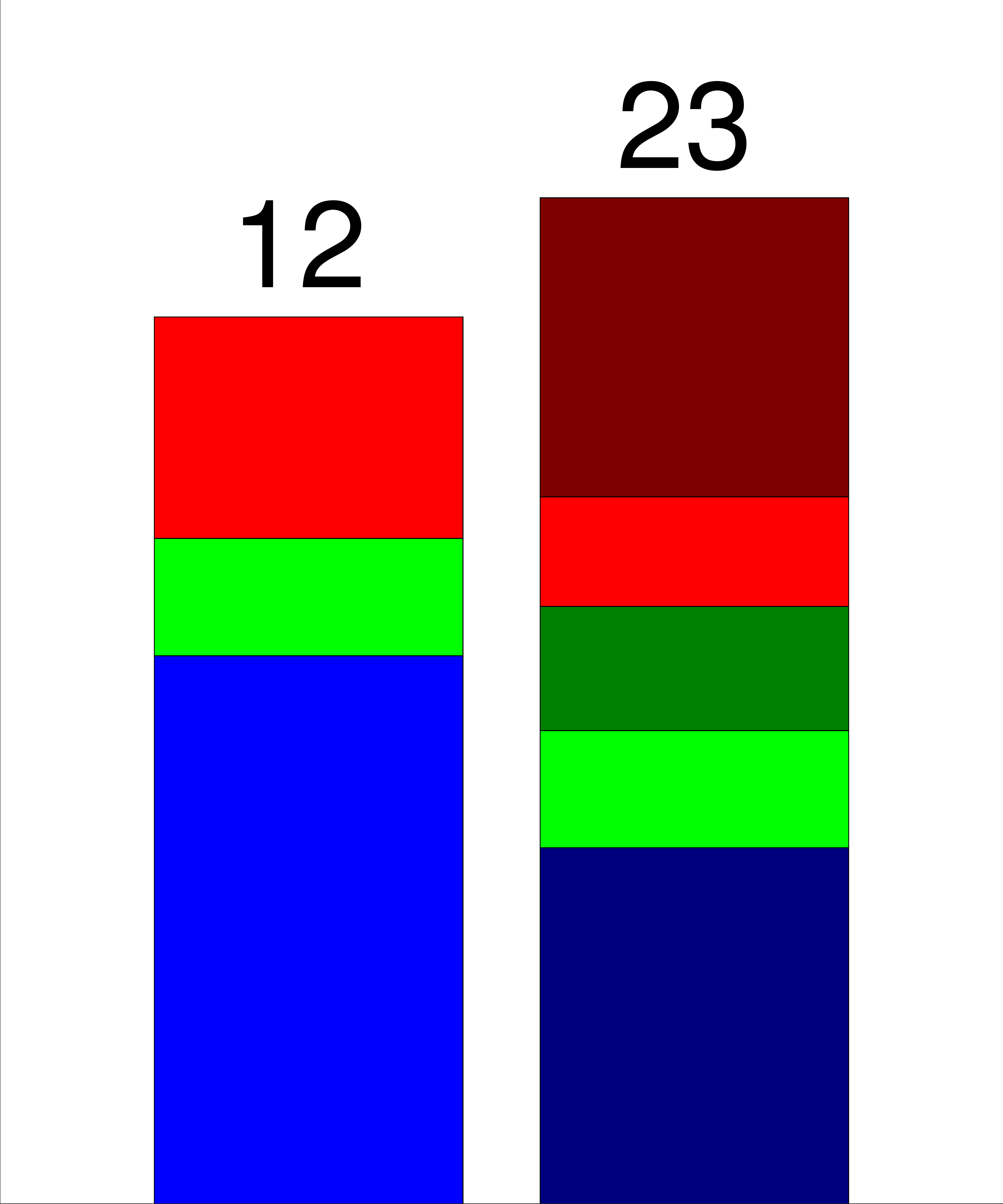}}\label{fig:cpu_scale256}}\\%
 \begin{picture}(15,1.0)%
 \put(1.2,0.0){\rotatebox{45}{SAMG}}%
 \put(2.2,0.0){\rotatebox{45}{F-AMS}}%
 \put(4.4,0.0){\rotatebox{45}{SAMG}}%
 \put(5.4,0.0){\rotatebox{45}{F-AMS}}%
 \put(7.75,0.0){\rotatebox{45}{SAMG}}%
 \put(8.75,0.0){\rotatebox{45}{F-AMS}}%
 \put(11.25,0.0){\rotatebox{45}{SAMG}}%
 \put(12.25,0.0){\rotatebox{45}{F-AMS}}%
  \end{picture}\\%
\subfigure{\fbox{\includegraphics[width=\textwidth]{Fig/method/legend.png}}}%
\caption{F-AMS (Decoupled-AMS) performance compared with SAMG for different domain sizes. The coarsening factors are $6$, $8$, $8$ and $10$ respectively. Also, a transmissibility ratio of  $T_{ratio} = 10^2$ was considered. The number of iterations to reach $10^{-6}$ residual 2-norm is given on top of each bar. For these experiments, both methods were employed as preconditioners for GMRES. Similar performance was observed for other F-AMS coupling strategies.}%
\label{fig:scale}%
\end{figure}

The scalability of F-AMS, benchmarked with SAMG, is investigated for heterogeneous (patchy) reservoir of increasing size. To this end, both the matrix and fracture fine-scale grid resolution is varied from $32^3$ matrix and $320\times32$ fractures (smallest) up to $256^3$ matrix and $2117\times256$ fracture cells (see Fig.~\ref{fig:scaledFracs}). The transmissibility ratio between the two media is $T_{ratio} = 10^2$. Figure~\ref{fig:scale} shows the obtained CPU times. F-AMS and SAMG both maintain their convergence rates and experience a similar level of scalability, in terms of CPU time, i.e., they grow linearly with the problem size. During these experiments, and as previously studied in Subsection~\ref{sec:coarsening_ratio}, F-AMS was found very sensitive to the coarsening strategy used. The reported results use the optimum coarsening ratio found during repeated experiments (some of which are detailed in Table~\ref{tab:coarsening_ratio}). Note that SAMG uses adaptive coarsening at each coarse level in its V-cycles.

\subsubsection{Heterogeneous fractures}

\begin{figure}[htb!]%
\centering%
\setlength{\fboxsep}{0pt}%
\setlength{\fboxrule}{0.2pt}%
\setlength{\unitlength}{1cm}%
\begin{minipage}{0.55\linewidth}%
\includegraphics[height=0.8\textwidth,width=0.8\textwidth]{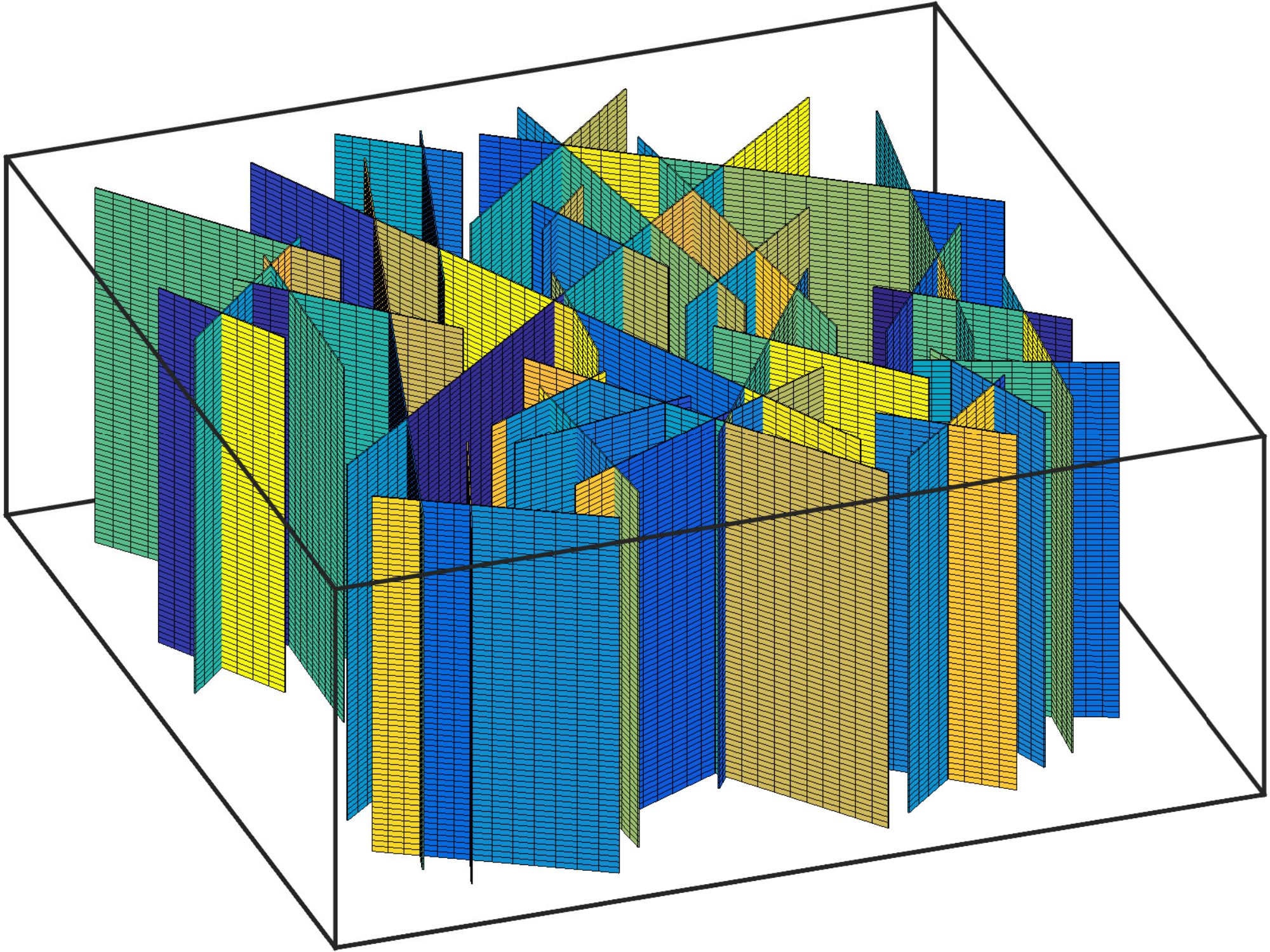}\hspace{0.2cm}\raisebox{0.7cm}{\includegraphics[height=0.5\textwidth,width=0.6cm]{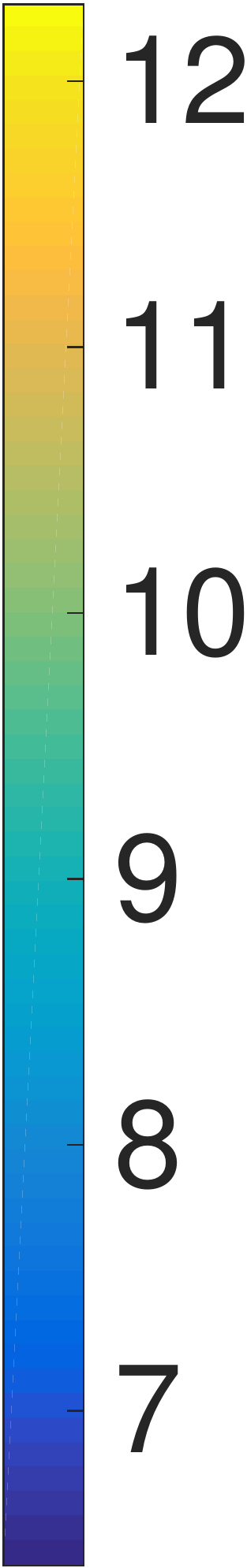}}\label{fig:perm_hetero}\\%
 \begin{picture}(5,0.5)%
 \put(3.2,0.2){$log_{10}(k^f)$}%
 \end{picture}%
\end{minipage} \hspace{0.2cm}%
\begin{minipage}{.25\linewidth}%
\vspace{0.8cm}%
  \begin{picture}(1.0,3.5)%
	\put(0.0,0.3){\rotatebox{90}{CPU time (sec)}}%
		\put(0.75,-0.1){0}%
		\put(0.7,1.7){2}%
		\put(0.7,3.3){4}%
\end{picture}%
\fbox{\includegraphics[width=0.8\textwidth,height=\textwidth]{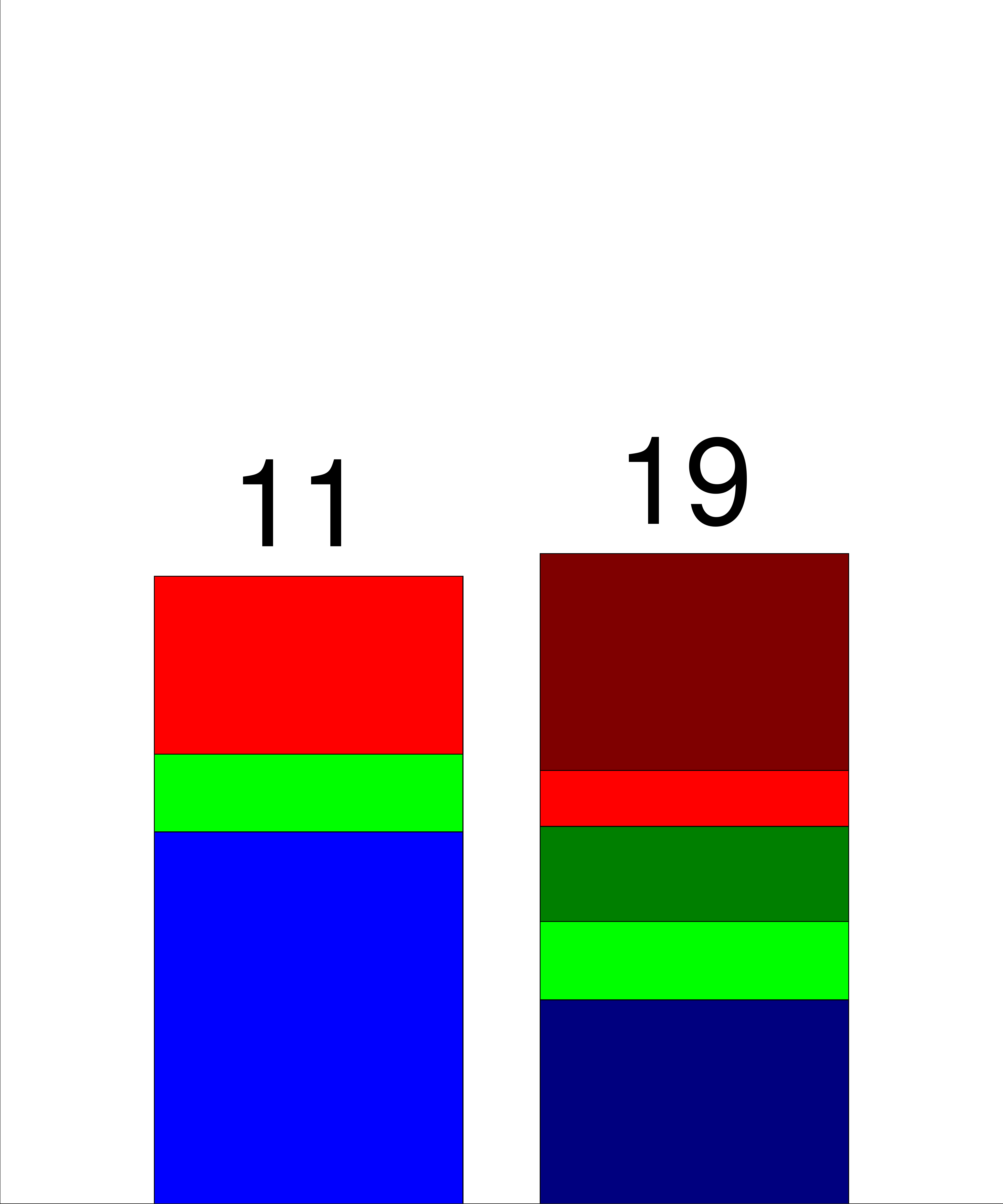}}\label{fig:cpu_hetero}\\%
 \begin{picture}(5,1.2)%
 \put(1.1,0.0){\rotatebox{45}{SAMG}}%
 \put(2.1,0.0){\rotatebox{45}{F-AMS}}%
  \end{picture}%
  \end{minipage}\\%
\fbox{\includegraphics[width=\textwidth]{Fig/method/legend.png}}%
\caption{F-AMS (Decoupled-AMS) performance compared with SAMG on a test case with heterogeneous fracture permeability. The logarithm of the fracture permeability is shown on the left. The number of iterations to reach $10^{-6}$ residual 2-norm is given on top of each bar. For these experiments, both methods were employed as preconditioners for GMRES.}%
\label{fig:hetero}%
\end{figure}

Finally, the sensitivity of the F-AMS method (benchmarked with SAMG) to strongly heterogeneous fracture properties is investigated. The permeability of of the $127$ fracture plates from Fig.~\ref{fig:fracs3d} is randomly perturbed to a span $6$ orders of magnitude. Figure~\ref{fig:hetero} shows that, if an appropriate coarsening ratio is chosen --in this case $6 \times 6 \times 8$ in the matrix and $4 \times 8$ in the fractures--, then the F-AMS and SAMG performances are comparable. In addition to the coarsening ratio, multiscale methods are also sensitive to the heterogeneity contrasts (here, in both fractures and matrix). Improvements can be achieved by adapting the coarse grid geometry to follow the fracture and matrix conductivity distribution, or by enriching the prolongation operator with additional basis functions \cite{Yalchin-enriched1,Yalchin-enriched2,Davide14}. These are subjects of future studies.

The results of the experiments presented in this section show that the performance and scalability of F-AMS is comparable to that of SAMG. As such, even in its current proof of concept implementation, F-AMS is found a promising multiscale method for naturally fractured porous media. Note that, for realistic multiphase test cases, simulations can be further accelerated by employing only few iterations of F-AMS, followed by a flux reconstruction stage, leading to efficient approximate solutions \cite{hadi-erest-spej}.

\section{Conclusions}
\label{sec:conclusions}

In this paper, a novel general multiscale framework, F-AMS, was devised for efficient and accurate simulation of flow through heterogeneous porous media with embedded fractures of various length scales. For the first time in the multiscale community, the possibility to prescribe an arbitrary coarse grid in each fracture network was presented. Then, for each coarse node (from both matrix and fractures), a locally-supported basis function was defined, by considering one of the four cross-media coupling strategies (Decoupled-AMS, Frac-AMS, Rock-AMS or Coupled-AMS). All of these flexibilities allow the user to tweak the trade-off between the computational budget of the setup stage and the convergence rate.  

Aligned with the EDFM approach of having independent grids for fracture and matrix \cite{SeongFrac2}, this work also introduced a distance-based automatic coarsening algorithm for the fracture domain. This allows the user to specify the desired (uniform) coarsening factor for the fracture domain, in similar fashion to the matrix. In addition,  the effect of truncating small non-zeros from the prolongation operator was studied, in order to maintain efficiency, especially for the Coupled-AMS strategy. For all test cases considered, the truncation value of $\alpha = 10^{-2}$ was found optimum.

The numerical results illustrate that F-AMS (similar as SAMG) is insensitive to the fracture-matrix conductivity contrast, and - importantly - shares the same scalability with respect to the fracture density, domain scale and heterogeneous properties. However, the performance of F-AMS was found to dramatically degrade if a sub-optimal coarsening strategy is chosen for challenging test cases. The method can be further extended to address this, e.g. by considering enhanced coarsening strategies, different choices for the smoother or employment of enriched prolongation operators \cite{Yalchin-enriched1,Yalchin-enriched2}.

Another important finding of this paper was that all basis function coupling strategies perform similar when F-AMS is used as preconditioner to GMRES. This recommends the Decoupled-AMS approach for commercial reservoir simulation, due to its convenient implementation, an attractive feature for real-field applications.

In summary, it is concluded that F-AMS is an important multiscale development for flow in heterogeneous media with embedded fractures. It was shown that only few fracture coarse nodes are required to deliver good approximate pressure solutions, at the original fine-scale resolution. Future developments of F-AMS will include consideration of complex physics such as stress-dependent fracture networks and also capillarity and gravity effects for multiphase flow simulations. 

\section*{Acknowledgements}
\label{sec:acknowledgements}
Financial support of PI/ADNOC is acknowledged. Also, many thanks are due towards the members of the DARSim research group of TU Delft, for the useful discussions during the development of F-AMS method.

\appendix
\section{Algebraic Formulation of the F-AMS prolongation operators}

Consider the fine-scale system for matrix and fractures, $\bm{A} p = q$, i.e.,  
\begin{align}\label{A-dis-mf}
\underbrace{\left[\begin{array}{cc} \bm A^{mm} & \bm A^{mf} \\ \bm A^{fm} & \bm A^{ff} \end{array}\right]}_{\bm A}
\underbrace{\left[\begin{array}{c}p^m \\p^f\end{array}\right]}_{p}
=
\underbrace{\left[\begin{array}{c}q^m \\q^f \end{array}\right]}_{q}.
\end{align}
The permutation operator $\bm\wp$ containing matrix and fracture permutation block operators $\bm\wp^{m}$ and $\bm\wp^{f}$, respectively,
\begin{align}
{\bm \wp} \equiv {\left[\begin{array}{cc}\bm\wp^{m} & 0 \\0 & \bm\wp^{f} \end{array}\right]},
\end{align}
is defined such that it reorders the linear system \eqref{A-dis-mf} based on the wirebasket ordering \cite{NordbottenBjorstad07,wirebasket,yixuan-ams} of Internal (I), Face (F), Edge (E) and Vertex (V) for both matrix (superscript $^m$) and fracture (superscript $^f$) unknowns, i.e.,
\begin{align}\label{A-reordered}
\underbrace{\left[\begin{array}{cccc|ccc}
\bm A^{{I^m}{I^m}} & \bm A^{{I^m}{F^m}} & 0 & 0 &
\bm A^{{I^m}{F^f}} & \bm A^{{I^m}{E^f}} & \bm A^{{I^m}{V^f}} \\
\bm A^{{F^m}{I^m}} & \bm A^{{F^m}{F^m}} & \bm A^{{F^m}{E^m}} & 0&
\bm A^{{F^m}{F^f}} & \bm A^{{F^m}{E^f}} & \bm A^{{F^m}{V^f}} \\
0 & \bm A^{{E^m}{F^m}} & \bm A^{{E^m}{E^m}} & \bm A^{{E^m}{V^m}} &
\bm A^{{E^m}{F^f}} & \bm A^{{E^m}{E^f}} & \bm A^{{E^m}{V^f}} \\
0 & 0 & \bm A^{{V^m}{E^m}} & \bm A^{{V^m}{V^m}} &
\bm A^{{V^m}{F^f}} & \bm A^{{V^m}{E^f}} & \bm A^{{V^m}{V^f}} \\
\hline
\bm A^{{F^f}{I^m}} & \bm A^{{F^f}{F^m}} & \bm A^{{F^f}{E^m}} & \bm A^{{F^f}{V^m}} &
\bm A^{{F^f}{F^f}} & \bm A^{{F^f}{E^f}} & 0 \\
\bm A^{{E^f}{I^m}} & \bm A^{{E^f}{F^m}} & \bm A^{{E^f}{E^m}} & \bm A^{{E^f}{V^m}} &
\bm A^{{E^f}{F^f}} & \bm A^{{E^f}{E^f}} & \bm A^{{E^f}{V^f}} \\
\bm A^{{V^f}{I^m}} & \bm A^{{V^f}{F^m}} & \bm A^{{V^f}{E^m}} & \bm A^{{V^f}{V^m}} &
0 & \bm A^{{V^f}{E^f}} & \bm A^{{V^f}{V^f}}
\end{array}\right]}_{\bm{\wp  A  \wp^T}}
\underbrace{\left[\begin{array}{c}
p^{I^m}\\ 
p^{F^m}\\ 
p^{E^m}\\ 
p^{V^m}\\ 
\hline
p^{F^f}\\ 
p^{E^f}\\ 
p^{V^f}\\ 
\end{array}\right]}_{\bm\wp p}
=
\underbrace{\left[\begin{array}{c}
q^{I^m}\\ 
q^{F^m}\\ 
q^{E^m}\\ 
q^{V^m}\\ 
\hline
q^{F^f}\\ 
q^{E^f}\\ 
q^{V^f}\\ 
\end{array}\right]}_{\bm\wp q}.
\end{align}
Note that fractures have only Face, Edge, and Vertex cells, since they are represented in a lower-dimensional space than the matrix. Also, according to the two-point flux approximation (TPFA) stencil for structured grids, $\bm A^{{I^m}{E^m}}$, $\bm A^{{E^m}{I^m}}$, $\bm A^{{I^m}{V^m}}$, $\bm A^{{V^m}{I^m}}$, $\bm A^{{V^m}{F^m}}$, $\bm A^{{F^m}{V^m}}$, $\bm A^{{F^f}{V^f}}$, $\bm A^{{V^f}{F^f}}$ are zero. More importantly, for media with embedded fractures, the coupling off-diagonal blocks $\bm A^{mf}$ and $\bm A^{fm}$ are full, i.e., each matrix cell may overlap with fracture cells of any type (F, E, or V). This is the main reason behind the consideration of the four types of basis functions, each with a different level of matrix-fracture coupling, as previously discussed in this paper. 

The algebraic construction of the prolongation operator for each strategy is described next.

\subsection{Decoupled-AMS}
\label{sec:app-decoupledAMS}
In the Decoupled-AMS prolongation operator, the matrix-fracture coupling terms are completely neglected. To this end, all off-diagonal block matrix entries (belonging to $\bm A^{mf}$ and $\bm A^{fm}$) are set to zero. In addition, similar to the AMS \cite{yixuan-ams} and C-AMS \cite{Tene-cams} methods, the linear system is further simplified to account for the localization boundary condition within each medium (by neglecting connectivity between each cell and its lower-ranked neighbours in the wirebasket hierarchy). This leads to the following approximate linear system:
\begin{align}\label{DecoupledAMS-approxSys}
\underbrace{\left[\begin{array}{cccc|ccc}
\underline{\bm A}^{{I^m}{I^m}} & \bm A^{{I^m}{F^m}} & 0 & 0 &
0 & 0 & 0 \\
0 & \overline{\underline{\bm{A}}}^{{F^m}{F^m}} & \bm A^{{F^m}{E^m}} & 0&
0 & 0 & 0 \\
0 & 0 & \overline{\underline{\bm{A}}}^{{E^m}{E^m}} & \bm A^{{E^m}{V^m}} &
0 & 0 & 0 \\
0 & 0 & 0 & \breve{\bm{A}}^{mm} &
0 & 0 & \breve{\bm{A}}^{mf} \\
\hline
0 & 0 & 0 & 0 &
\underline{\bm A}^{{F^f}{F^f}} & \bm A^{{F^f}{E^f}} & 0 \\
0 & 0 & 0 & 0 &
0 & \overline{\underline{\bm{A}}}^{{E^f}{E^f}} & \bm A^{{E^f}{V^f}} \\
0 & 0 & 0 & \breve{\bm{A}}^{fm} &
0 & 0 & \breve{\bm{A}}^{ff}
\end{array}\right]}_{\bm{\wp  A'  \wp^T}} %
\underbrace{\left[\begin{array}{c}
p'^{I^m}\\ 
p'^{F^m}\\ 
p'^{E^m}\\ 
p'^{V^m}\\ 
\hline
p'^{F^f}\\ 
p'^{E^f}\\ 
p'^{V^f}\\ 
\end{array}\right]}_{\bm\wp p'}
=
\underbrace{\left[\begin{array}{c}
0\\ 
0\\ 
0\\ 
\breve{q}^m\\ 
\hline
0\\ 
0\\ 
\breve{q}^f\\ 
\end{array}\right]}_{\bm\wp q'}.
\end{align}
Here, the diagonal blocks marked as $\overline{\bm{A}}$ indicate that the matrix-matrix and fracture-fracture transmissibilities, neglected due to this localization assumption ($\bm A^{{F^m}{I^m}}$, $\bm A^{{E^m}{F^m}}$, $\bm A^{{V^m}{E^m}}$  and  $\bm A^{{E^f}{F^f}}$, $\bm A^{{V^f}{E^f}}$, respectively) have also been removed from the diagonal term. At the same time, the notation $\underline{\bm A}$ indicates diagonal blocks where the neglected matrix-fracture transmissibilities have been removed from the diagonal term. Finally,
\begin{align}
\breve{\bm{A}} \equiv \left[\begin{array}{cc}
\breve{\bm{A}}^{mm} & \breve{\bm{A}}^{mf} \\
\breve{\bm{A}}^{fm} & \breve{\bm{A}}^{ff} %
\end{array}\right] = %
\bm{(\mathcal{R} A \mathcal{P})}, & &%
\breve{p} \equiv \left[\begin{array}{c} %
p'^{V^m} \\
p'^{V^f} %
\end{array}\right], & &%
\breve q \equiv \left[\begin{array}{c}
\breve{q}^m \\
\breve{q}^f
\end{array}\right] = %
\bm{\mathcal{R}} q, %
\end{align}
are the components of the coarse-scale system.

After solving for the coarse-scale pressures, $\breve p = {\breve{\bm{A}}}^{-1} \breve q$, the approximate system can be inverted algebraically, due to its upper-triangular structure. Consequently, the prolongation operator, which satisfies $p' = \bm{\mathcal{P}} \breve p$, reads
\begin{align}\label{DecoupledAMS-prlg}
\bm{\mathcal{P}} = {\bm \wp^T}%
{\left[\begin{array}{c|c}
- {\underline{\bm{A}}^{{I^m}{I^m}}}^{-1} \bm{A}^{{I^m}{F^m}}  \bm{\mathcal{P}}^{F^m V^m} & 0\\ 
-{\overline{\underline{\bm{A}}}^{{F^m}{F^m}}}^{-1} \bm{A}^{{F^m}{E^m}} \bm{\mathcal{P}}^{E^m V^m} & 0\\ 
-{\overline{\underline{\bm{A}}}^{{E^m}{E^m}}}^{-1} \bm{A}^{{E^m}{V^m}} & 0\\ 
\bm{I}^{V^m V^m}& 0\\ 
\hline
0 & 
-{\underline{\bm{A}}^{{F^f}{F^f}}}^{-1} \bm{A}^{{F^f}{E^f}}  \bm{\mathcal{P}}^{E^f V^f}\\ 
0 & 
-{\overline{\underline{\bm{A}}}^{{E^f}{E^f}}}^{-1} \bm{A}^{{E^f}{V^f}} \\ 
0 & 
\bm{I}^{V^f V^f}
\end{array}\right]},
\end{align}
where $\bm{I}$ is the identity matrix and the transpose operator ${\bm \wp^T}$ back-transforms the wirebasket ordering into the natural ordering. Also, $\bm{\mathcal{P}}^{E^m V^m}$, $\bm{\mathcal{P}}^{F^m V^m}$ and $\bm{\mathcal{P}}^{E^f V^f}$ are sub-blocks of the prolongation with the corresponding rows and columns given in the superscripts. For example, $\bm{\mathcal{P}}^{E^m V^m} = -{\overline{\underline{\bm{A}}}^{{E^m}{E^m}}}^{-1} \bm{A}^{{E^m}{V^m}}$ and $ \bm{\mathcal{P}}^{F^m V^m} = -{\overline{\underline{\bm{A}}}^{{F^m}{F^m}}}^{-1} \bm{A}^{{F^m}{E^m}} \bm{\mathcal{P}}^{E^m V^m} $. Note that, once computed, the higher-rank sub-blocks of $\bm{\mathcal{P}}$ become boundary conditions for the values of basis functions in lower-rank cells, in accordance to the localization assumption (e.g. the values obtained for matrix edges, $\bm{\mathcal{P}}^{E_m V_m}$, are used to compute the prolongation in adjacent faces, $\bm{\mathcal{P}}^{F^m V^m}$).

\subsection{Frac-AMS}
\label{sec:app-fracAMS}
The Frac-AMS approach considers the effect of the $\bm A^{mf}$ transmissibilities when computing basis functions. This leads to the following approximate system operator
\begin{align}\label{FracAMS-approxSys}
\bm{\wp A' \wp^T} = {\left[\begin{array}{cccc|ccc}
\bm{A}^{{I^m}{I^m}} & \bm{A}^{{I^m}{F^m}} & 0 & 0 &
\bm{A}^{{I^m}{F^f}} & \bm{A}^{{I^m}{E^f}} & \bm{A}^{{I^m}{V^f}} \\
0 & \overline{\bm{A}}^{{F^m}{F^m}} & \bm{A}^{{F^m}{E^m}} & 0&
\bm{A}^{{F^m}{F^f}} & \bm{A}^{{F^m}{E^f}} & \bm{A}^{{F^m}{V^f}} \\
0 & 0  & \overline{\bm{A}}^{{E^m}{E^m}} & \bm{A}^{{E^m}{V^m}} &
\bm{A}^{{E^m}{F^f}} & \bm{A}^{{E^m}{E^f}} & \bm{A}^{{E^m}{V^f}} \\
0 & 0 & 0 & \breve{\bm{A}}^{mm} &
0 & 0 & \breve{\bm{A}}^{mf} \\
\hline
0 & 0 & 0 & 0 &
\underline{\bm{A}}^{{F^f}{F^f}} & \bm{A}^{{F^f}{E^f}} & 0 \\
0 & 0 & 0 &  0 &
0 & \overline{\underline{\bm{A}}}^{{E^f}{E^f}} & \bm{A}^{{E^f}{V^f}} \\
0 & 0 & 0 & \breve{\bm{A}}^{fm} &
0 & 0 & \breve{\bm{A}}^{ff}
\end{array}\right]},
\end{align}
where $\bm A^{{F^m}{I^m}}$, $\bm A^{{E^m}{F^m}}$, $\bm A^{{V^m}{E^m}}$, $\bm A^{{V^m}{F^f}}$ and $\bm A^{{V^m}{E^f}}$ are set to zero due to localization boundary condition corresponding to Frac-AMS coupling for matrix, while, at the same time, the $\bm A^{{E^f}{F^f}}$, $\bm A^{{V^f}{E^f}}$ are zero in the fracture equations.

The Frac-AMS prolongation operator reads
\begin{align}\label{FracAMS-prlg}
\bm{\mathcal{P}} = {\bm \wp^T}%
{\left[\begin{array}{c|c}
- {\bm{A}^{{I^m}{I^m}}}^{-1} \bm{A}^{{I^m}{F^m}}  \bm{\mathcal{P}}^{F^m V^m} &
- {{\bm{A}}^{{I^m}{I^m}}}^{-1} (\bm{A}^{{I^m}{F^m}} \bm{\mathcal{P}}^{F^m V^f} + \bm{A}^{{I^m}{f}} \bm{\mathcal{P}}^{f V^f}) \\
-{\overline{\bm{A}}^{{F^m}{F^m}}}^{-1} \bm{A}^{{F^m}{E^m}} \bm{\mathcal{P}}^{E^m V^m} &
- {\overline{\bm{A}}^{{F^m}{F^m}}}^{-1} (\bm{A}^{{F^m}{E^m}} \bm{\mathcal{P}}^{E^m V^f} + \bm{A}^{{F^m}{f}} \bm{\mathcal{P}}^{f V^f}) \\
-{\overline{\bm{A}}^{{E^m}{E^m}}}^{-1} \bm{A}^{{E^m}{V^m}} & 
- {\overline{\bm{A}}^{{E^m}{E^m}}}^{-1} \bm{A}^{{E^m}{f}} \bm{\mathcal{P}}^{f V^f} \\
\bm{I}^{V^m V^m}& 
0\\ 
\hline
0 & 
- {\underline{\bm{A}}^{{F^f}{F^f}}}^{-1} \bm{A}^{{F^f}{E^f}}  \bm{\mathcal{P}}^{E^f V^f}\\ 
0 & 
-{\overline{\underline{\bm{A}}}^{{E^f}{E^f}}}^{-1} \bm{A}^{{E^f}{V^f}}\\ 
0 & 
\bm{I}^{V^f V^f}\\ 
\end{array}\right]},
\end{align}
where the superscript $^f$ (e.g. from $\bm{\mathcal{P}}^{fV^f}$) corresponds to all the fracture cells, regardless of their containing dual block. Similar as in the previous case, the sub-blocks of prolongation operator $\bm{\mathcal{P}}^{F^m V^m}$, $\bm{\mathcal{P}}^{F^m V^f}$, $\bm{\mathcal{P}}^{f V^f}$, $\bm{\mathcal{P}}^{E^m V^m}$,  $\bm{\mathcal{P}}^{E^m V^f}$, $\bm{\mathcal{P}}^{f V^f}$ and $\bm{\mathcal{P}}^{E^f V^f}$ represent  the corresponding rows and columns given in their superscripts. For example, $\bm{\mathcal{P}}^{E^m V^f} = - {\overline{\bm{A}}^{{E^m}{E^m}}}^{-1} \bm{A}^{{E^m}{f}} \bm{\mathcal{P}}^{f V^f}$ and, specially,
\begin{align}\label{pfVf}
\bm{\mathcal{P}}^{f V^f} = {\left[ \begin{array}{c}
- {\underline{\bm{A}}^{{F^f}{F^f}}}^{-1} \bm{A}^{{F^f}{E^f}}  \bm{\mathcal{P}}^{E^f V^f}\\ 
-{\overline{\underline{\bm{A}}}^{{E^f}{E^f}}}^{-1} \bm{A}^{{E^f}{V^f}}\\ 
\bm{I}^{V^f V^f}\\ 
\end{array} \right]}.
\end{align}
\subsection{Rock-AMS}%
\label{sec:app-rockAMS}%
For Rock-AMS, the $\bm A^{mf}$ transmissibilities are set to zero,%
\begin{align}%
\label{RockAMS-approxSys}%
\bm{\wp A' \wp^T} = {\left[\begin{array}{cccc|ccc}%
\underline{\bm{A}}^{{I^m}{I^m}} & \bm{A}^{{I^m}{F^m}} & 0 & 0 &%
0 & 0 &0 \\
0 & \overline{\underline{\bm{A}}}^{{F^m}{F^m}} & \bm{A}^{{F^m}{E^m}} & 0&%
0 & 0 &0 \\
0 & 0  & \overline{\underline{\bm{A}}}^{{E^m}{E^m}} & \bm{A}^{{E^m}{V^m}} &%
0 & 0 &0 \\
0 & 0 & 0 & \breve{\bm{A}}^{mm} &%
0 & 0 & \breve{\bm{A}}^{mf} \\
\hline%
\bm{A}^{{F^f}{I^m}} & \bm{A}^{{F^f}{F^m}} & \bm{A}^{{F^f}{E^m}} & \bm{A}^{{F^f}{V^m}} &%
\bm{A}^{{F^f}{F^f}} & \bm{A}^{{F^f}{E^f}} & 0 \\
\bm{A}^{{E^f}{I^m}} & \bm{A}^{{E^f}{F^m}} & \bm{A}^{{E^f}{E^m}} & \bm{A}^{{E^f}{V^m}} &%
0 & \overline{\bm{A}}^{{E^f}{E^f}} & \bm{A}^{{E^f}{V^f}} \\
0 & 0 & 0 & \breve{\bm{A}}^{fm} &%
0 & 0 & \breve{\bm{A}}^{ff}%
\end{array}\right]},%
\end{align}%
where, the localization boundary condition was also applied.\\
\indent The Rock-AMS prolongation operator reads%
\begin{align}\label{RockAMS-prlg}%
\bm{\mathcal{P}} = {\bm \wp^T}%
{\left[\begin{array}{c|c}%
- {\underline{\bm{A}}^{{I^m}{I^m}}}^{-1} \bm{A}^{{I^m}{F^m}}  \bm{\mathcal{P}}^{F^m V^m} &%
0 \\
-{\overline{\underline{\bm{A}}}^{{F^m}{F^m}}}^{-1} \bm{A}^{{F^m}{E^m}} \bm{\mathcal{P}}^{E^m V^m} &%
0 \\
-{\overline{\underline{\bm{A}}}^{{E^m}{E^m}}}^{-1} \bm{A}^{{E^m}{V^m}} & %
0 \\
\bm{I}^{V^m V^m}& %
0\\ 
\hline %
-{\bm{A}^{{F^f}{F^f}}}^{-1} (\bm{A}^{{F^f}{E^f}} \bm{\mathcal{P}}^{E^f V^m} + \bm{A}^{{F^f}{m}} \bm{\mathcal{P}}^{m V^m}) & %
- {\bm{A}^{{F^f}{F^f}}}^{-1} \bm{A}^{{F^f}{E^f}}  \bm{\mathcal{P}}^{E^f V^f}\\ 
- {\overline{\bm{A}}^{{E^f}{E^f}}}^{-1} \bm{A}^{{E^f}{m}} \bm{\mathcal{P}}^{m V^m} & %
-{\overline{\bm{A}}^{{E^f}{E^f}}}^{-1} \bm{A}^{{E^f}{V^f}}\\ 
0 & %
\bm{I}^{V^f V^f}\\ 
\end{array}\right]},%
\end{align}%
where the superscript $^m$ (e.g. from $\bm{\mathcal{P}}^{mV^m}$) corresponds to all the matrix cells, regardless of their containing dual block.
The sub-blocks of the prolongation operator, i.e.,$\bm{\mathcal{P}}^{F^m V^m}$, $\bm{\mathcal{P}}^{E^m V^m}$, $\bm{\mathcal{P}}^{E^f V^m} $, $\bm{\mathcal{P}}^{m V^m} $ and $\bm{\mathcal{P}}^{E^f V^f}$, are defined similar to the previous cases in the sense that they represent  the corresponding rows and columns given in their superscripts. For example, $ \bm{\mathcal{P}}^{E^m V^m} = -{\overline{\underline{\bm{A}}}^{{E^m}{E^m}}}^{-1} \bm{A}^{{E^m}{V^m}}$ and, specially, 
\begin{align}\label{PmVm}
\bm{\mathcal{P}}^{m V^m} = 
\left[
\begin{array}{c}%
- {\underline{\bm{A}}^{{I^m}{I^m}}}^{-1} \bm{A}^{{I^m}{F^m}}  \bm{\mathcal{P}}^{F^m V^m}\\
-{\overline{\underline{\bm{A}}}^{{F^m}{F^m}}}^{-1} \bm{A}^{{F^m}{E^m}} \bm{\mathcal{P}}^{E^m V^m}\\
-{\overline{\underline{\bm{A}}}^{{E^m}{E^m}}}^{-1} \bm{A}^{{E^m}{V^m}} \\
\bm{I}^{V^m V^m}\\ 
\end{array}\right].
\end{align}

\subsection{Coupled-AMS}%
\label{sec:app-coupledAMS}%
In the Coupled-AMS approach, all adjacent Face and Edge blocks are merged between the media (Fig.~\ref{fig:mergedDuals}), i.e. $F = F^m \bigcup F^f$ and $E = E^m \bigcup E^f$. Also, let $V$ denote the set of coarse nodes, irrespective of their location. In this new setting, the $\bm\wp$-reordered approximate linear system is defined as%
\begin{align}\label{CoupledAMS-approxSys}%
\bm{\wp A' \wp^T} = {\left[\begin{array}{cccc}%
\bm{A}^{I I} & \bm{A}^{I F} & \bm{A}^{I E} & \bm{A}^{I V} \\
0 & \overline{\bm{A}}^{F F} & \bm{A}^{F E} & \bm{A}^{F V} \\
0 & 0  & \overline{\bm{A}}^{E E} & \bm{A}^{E V} \\
0 & 0 & 0 & \breve{\bm{A}}%
\end{array}\right]},%
\end{align}%
where the localization boundary condition was appropriately employed.\\
\indent Then, the Coupled-AMS prolongation operator reads
\begin{align}\label{Coupled-AMS_prolong_1}
\bm{\mathcal{P}} = {\bm \wp^T}%
{\left[\begin{array}{c}
-{\bm{A}^{I I}}^{-1} (A^{I F} \bm{\mathcal{P}}^{F V} + \bm{A}^{I E} \bm{\mathcal{P}}^{E V} + \bm{A}^{I V})\\ 
- {\overline{\bm{A}}^{F F}}^{-1} (\bm{A}^{F E} \bm{\mathcal{P}}^{E V} + \bm{A}^{F V})\\ 
-{\overline{\bm{A}}^{E E}}^{-1} \bm{A}^{E V}\\ 
\bm{I}^{V V}
\end{array}\right]}.
\end{align}

\section{Coarsening ratios used during scale sensitivity test}

The scale sensitivity experiment revealed that the performance of F-AMS is highly dependent on the coarsening ratios used. Only the optimum configuration was featured in the plots presented in the body of the manuscript. Table~\ref{tab:coarsening_ratio} lists experimental results obtained when using primal grids with more refined and more coarse blocks, respectively, for comparison purposes.

\begin{table}[htb!] 
\begin{center}
\begin{tabular}{|c|c|c|c|}
\hline
\hline
Scale & \descrcell{Coarsening ratio}{(matrix, fracs)} & \descrcell{Total CPU time}{(sec)} & \# iterations \\
\hline 
\hline
\multirow{3}{*}{$32^3$} %
& $2 \times 2 \times 2$, $2 \times 2$ & $1.904$ & $11$ \\ \cline{2-4}
& $\bm{6 \times 6 \times 6}$, $\bm{6 \times 6}$ & $\bm{0.330}$ & $\bm{23}$ \\ \cline{2-4}
& $8 \times 8 \times 8$, $8 \times 8$ & $0.351$ & $29$ \\ 
\hline
\hline
\multirow{3}{*}{$64^3$} %
& $4 \times 4 \times 4$, $4 \times 4$ & $3.270$ & $13$ \\ \cline{2-4}
& $\bm{6 \times 6 \times 6}$, $\bm{6 \times 6}$ & $\bm{2.186}$ & $\bm{19}$ \\ \cline{2-4}
& $9 \times 9 \times 9$, $9 \times 9$ & $2.371$ & $27$ \\ 
\hline
\hline
\multirow{3}{*}{$128^3$} %
& $6 \times 6 \times 6$, $6 \times 6$ & $21.790$ & $17$ \\ \cline{2-4}
& $\bm{8 \times 8 \times 8}$, $\bm{8 \times 8}$ & $\bm{17.620}$ &$\bm{22}$ \\ \cline{2-4}
& $11 \times 11 \times 11$, $11 \times 11$ & $21.350$ & $31$ \\ 
\hline
\hline
\multirow{3}{*}{$256^3$} %
& $8 \times 8 \times 8$, $8 \times 8$ & $164.600$ & $18$ \\ \cline{2-4}
& $\bm{10 \times 10 \times 10}$, $\bm{10 \times 10}$ & $\bm{150.400}$ & $\bm{23}$ \\ \cline{2-4}
& $17 \times 17 \times 17$, $17 \times 17$ & $252.100$ & $40$ \\
\hline
\hline
\end{tabular}
\vspace{0.5cm}%
\caption{Performance of F-AMS during the scale sensitivity test cases, when using different coarsening factors. The middle row for each test (shown in bold) is the optimum configuration, whose results were presented in the body of the manuscript.}
\label{tab:coarsening_ratio}
\end{center}
\end{table}

\end{document}